%% file: main.tex
\documentclass[12pt]{report}

\usepackage{tamuconfig}

\usepackage[T1]{fontenc}

\usepackage{graphicx}

\DeclareGraphicsExtensions{.png}
\DeclareGraphicsExtensions{.pdf}

\graphicspath{ {./graphic/} }

\usepackage[hidelinks]{hyperref}

\usepackage{graphicx, bm, amsmath, siunitx, longtable, tabularx, float, mathtools, booktabs, eqnarray, color, xcolor, soul, amsfonts, amssymb, appendix, tensor, subcaption, enumitem, makecell, listings, dsfont, lscape, diagbox, attachfile}
\usepackage[version=4]{mhchem}
\usepackage[export]{adjustbox}
\allowdisplaybreaks

\newif\ifDarkTheme
\DarkThemefalse

\newif\ifOGAPS
\OGAPSfalse

\newcommand{\B}[1]{{\bm #1}}
\newcommand{\T}{^{\mbox{\tiny T}}}
\newcommand{\ds}{\displaystyle}
\newcommand{\dd}{\; \text{d}}
\DeclarePairedDelimiter{\nint}\lfloor\rceil

\newcommand{\ces}{constrained expressions}
\renewcommand{\ce}{constrained expression}
\newcommand{\p}[2]{\prescript{(#1)}{}{#2}}
\newcommand{\pp}[3]{\prescript{#1}{#2}{#3}}
\newcommand{\C}[1]{\tensor*[^{}_{}]{\mathfrak{C}}{^{}_{#1}}}
\newcommand{\pC}[2]{\tensor*[^{(#1)}_{}]{\mathfrak{C}}{^{}_{#2}}}
\newcommand{\ppC}[3]{\tensor*[^{#1}_{#2}]{\mathfrak{C}}{^{}_{#3}}}
\newcommand{\R}{\mathbb{R}}
\DeclarePairedDelimiter{\ceil}\lceil\rceil
\newcommand{\andd}{\quad \text{and} \quad}
\newcommand{\Part}[2]{\frac{\partial #1}{\partial #2}}
\newcommand{\PartS}[2]{\frac{\partial^2 #1}{\partial #2^2}}

\ifOGAPS
\definecolor{titleColor}{RGB}{255,255,255}
\definecolor{backColorAlt}{RGB}{255,255,255}
\definecolor{textColor}{RGB}{0,0,0}
\definecolor{backColor}{RGB}{255,255,255}
\definecolor{titleColorText}{RGB}{0,0,0}
\else
    \definecolor{titleColor}{RGB}{75,0,130}
    \definecolor{titleColorText}{RGB}{255,255,255}
    \ifDarkTheme
        \definecolor{backColorAlt}{RGB}{0,0,0}
        \definecolor{textColor}{RGB}{230,230,230}
        \definecolor{backColor}{RGB}{110,104,128}
    \else
        \definecolor{backColorAlt}{RGB}{255,255,255}
        \definecolor{textColor}{RGB}{0,0,0}
        \definecolor{backColor}{RGB}{220,208,255}
    \fi
\fi

\usepackage{placeins, tocloft}

\newcommand{\listExamplesname}{List of Examples}
\newlistof[chapter]{Examples}{ex}{\listExamplesname}
\renewcommand{\theExamples}{\thechapter.\arabic{Examples}}
\newcommand{\Examples}[1]{%
    \par\textbf{Example \theExamples : #1}
    \addcontentsline{ex}{Examples}{\protect\numberline{\theExamples}#1}\par}

\usepackage{tcolorbox}
\tcbuselibrary{many}

\newlength\defaultparindent
\AtBeginDocument{\setlength\defaultparindent{\parindent}}
\ifOGAPS
    \newtcolorbox{exBox}[1][]{
        colback=backColor,
        colframe=titleColor,
        fonttitle=\bfseries,
        colbacktitle=titleColor,
        coltitle=titleColorText,
        coltext=textColor,
        enhanced,
        breakable,
        arc=5pt,
        attach boxed title to top left ={yshift=-\tcboxedtitleheight/2,yshifttext=-\tcboxedtitleheight/2,xshift=0.5cm},
        title=#1,
        before upper={\parindent\defaultparindent},
    }
\else
    \newtcolorbox{exBox}[1][]{
        colback=backColor,
        colframe=titleColor,
        fonttitle=\bfseries,
        colbacktitle=titleColor,
        coltitle=titleColorText,
        coltext=textColor,
        enhanced,
        breakable,
        arc=5pt,
        attach boxed title to top right ={yshift=-\tcboxedtitleheight/2,yshifttext=-\tcboxedtitleheight/2,xshift=-0.5cm},
        title=#1,
        before upper={\parindent\defaultparindent},
    }
\fi

\newenvironment{example}[1]
    {%
    \refstepcounter{Examples}
    \begingroup \allowdisplaybreaks\vspace{0.4cm}\begin{exBox}[\Examples{#1}]\relax}
    {\end{exBox}\endgroup}
    
\renewcommand{\proof}{\noindent\hrulefill\newline\noindent \textbf{Proof:} }
\newcounter{Theorems}
\newcommand{\Theorems}{%
    \par\textbf{Theorem \theTheorems}\par
}
    
\newtcolorbox{thrmBox}[1][]{
    colback=backColor,
    colframe=titleColor,
    fonttitle=\bfseries,
    colbacktitle=titleColor,
    coltitle=titleColorText,
    coltext=textColor,
    enhanced,
    breakable,
    arc=5pt,
    attach boxed title to top left ={yshift=-\tcboxedtitleheight/2,yshifttext=-\tcboxedtitleheight/2,xshift=0.5cm},
    title=#1,
    before upper={\noindent\parindent\defaultparindent},
}

\newenvironment{theorem}
    {%
    \refstepcounter{Theorems}
    \begingroup \allowdisplaybreaks\vspace{0.4cm}\begin{thrmBox}[\Theorems]\relax}
    {$\blacksquare$\end{thrmBox}\endgroup}
    
\newcounter{Properties}
\newcommand{\Properties}{%
    \par\textbf{Property \theProperties}\par
}
    
\newtcolorbox{propDefBox}[1][]{
    colback=backColorAlt,
    colframe=titleColor,
    fonttitle=\bfseries,
    colbacktitle=titleColor,
    coltitle=titleColorText,
    coltext=textColor,
    enhanced,
    breakable,
    arc=5pt,
    attach boxed title to top left ={yshift=-\tcboxedtitleheight/2,yshifttext=-\tcboxedtitleheight/2,xshift=0.5cm},
    title=#1,
}

\newenvironment{property}
    {%
    \refstepcounter{Properties}
    \begingroup \allowdisplaybreaks\vspace{0.4cm}\begin{propDefBox}[\Properties]\relax}
    {\end{propDefBox}\endgroup}

\newcounter{Definitions}
\newcommand{\Definitions}{%
    \par\textbf{Definition \theDefinitions}\par
}

\newenvironment{definition}
    {%
    \refstepcounter{Definitions}
    \begingroup \allowdisplaybreaks\vspace{0.4cm}\begin{propDefBox}[\Definitions]\relax}
    {\end{propDefBox}\endgroup}

\ifDarkTheme
    \color{textColor}
    \pagecolor{black}
\fi

\begin{document}

\renewcommand{\tamumanuscripttitle}{The Multivariate Theory of Functional Connections: An \texorpdfstring{$n$}{n}-Dimensional Constraint Embedding Technique Applied to Partial Differential Equations}

\renewcommand{\tamupapertype}{Dissertation Proposal}

\renewcommand{\tamufullname}{Carl Dakota Leake}

\renewcommand{\tamudegree}{Doctor of Philosophy}
\renewcommand{\tamuchairone}{Daniele Mortari}

\renewcommand{\tamumemberone}{John Hurtado}
\newcommand{\tamumembertwo}{Junuthula Reddy}
\newcommand{\tamumemberthree}{Moble Benedict}
\renewcommand{\tamudepthead}{Srinivas Rao Vadali}

\renewcommand{\tamugradmonth}{August}
\renewcommand{\tamugradyear}{2021}
\renewcommand{\tamudepartment}{Aerospace Engineering}

\include{Data/titlepage} 
\include{Data/abstract}
\include{Data/DedicationAndAcknowledgments}
\include{Data/contributors}
\include{Data/nomenclature}

\include{Data/lists}  

\include{Data/Introduction}
\include{Data/TfcTheory}

\include{Data/DeApplications}

\include{Data/FlexBodyApplications}
\include{Data/conclusions}


\let\oldbibitem\bibitem
\renewcommand{\bibitem}{\setlength{\itemsep}{0pt}\oldbibitem}
\bibliographystyle{ieeetr}

\phantomsection
\addcontentsline{toc}{chapter}{REFERENCES}

\renewcommand{\bibname}{{\normalsize\rm REFERENCES}}

\bibliography{Data/Refs}

\counterwithin*{Definitions}{chapter}
\renewcommand{\theDefinitions}{\thechapter.\arabic{Definitions}}

\include{Data/appendices}

\end{document}

%% file: Data/titlepage.tex
%
%
%
%


\providecommand{\tabularnewline}{\\}

\begin{titlepage}
\begin{center}
\MakeUppercase{\tamumanuscripttitle}
\vspace{4em}

A \tamupapertype

by

\MakeUppercase{\tamufullname}

\vspace{4em}

\begin{singlespace}

Submitted to the Office of Graduate and Professional Studies of \\
Texas A\&M University \\

in partial fulfillment of the requirements for the degree of \\
\end{singlespace}

\MakeUppercase{\tamudegree}
\par\end{center}
\vspace{2em}
\begin{singlespace}
\begin{tabular}{ll}
 & \tabularnewline
& \cr
Chair of Committee, & \tamuchairone\tabularnewline
Committee Members, & \tamumemberone\tabularnewline
 & \tamumembertwo\tabularnewline
 & \tamumemberthree\tabularnewline
Head of Department, & \tamudepthead\tabularnewline

\end{tabular}
\end{singlespace}
\vspace{3em}

\begin{center}
\tamugradmonth \hspace{2pt} \tamugradyear

\vspace{3em}

Major Subject: \tamudepartment \par
\vspace{3em}
Copyright \tamugradyear \hspace{.5em}\tamufullname 
\par\end{center}
\end{titlepage}
\pagebreak{}

%% file: Data/abstract.tex
%

\chapter*{ABSTRACT}
\addcontentsline{toc}{chapter}{ABSTRACT} 

\pagestyle{plain} 
\pagenumbering{roman} 
\setcounter{page}{2}

The Theory of Functional Connections (TFC) is a functional interpolation framework founded upon the so-called \ce: a functional that expresses the family of all possible functions that satisfy some user-specified, linear constraints. These \ces\ can be utilized to transform constrained problems into unconstrained ones. The benefits of doing so include faster solution times, more accurate solutions, and more robust convergence. This dissertation contains a comprehensive, self-contained presentation of the TFC theory beginning with simple univariate point constraints and ending with general linear constraints in $n$-dimensions; relevant mathematical theorems and clarifying examples are included throughout the presentation to expand and solidify the reader's understanding. Furthermore, this dissertation describes how TFC can be applied to estimate differential equations' solutions, its primary application to date. In addition, comparisons with other state-of-the-art algorithms that estimate differential equations' solutions are included to showcase the advantages and disadvantages of the TFC approach. Lastly, the aforementioned concepts are leveraged to estimate solutions of differential equations from the field of flexible body dynamics.

\pagebreak{}

%% file: Data/DedicationAndAcknowledgments.tex
%
%
%
%

\chapter*{DEDICATION}
\addcontentsline{toc}{chapter}{DEDICATION}  

\begin{center}
This dissertation is dedicated to the ocean of human knowledge: \\ what follows is a molecule of $H_2O$.
\end{center}
\vspace{-1.0cm}

{\let\clearpage\relax
\chapter*{ACKNOWLEDGMENTS}}
\addcontentsline{toc}{chapter}{ACKNOWLEDGMENTS}  

It is a nearly impossible task to list all of those who have positively nudged one's trajectory through life. Even the smallest push can have a massive effect, as I conjecture life is a chaotic system. Moreover, I am fortunate enough to have been incredibly blessed with supportive and encouraging family, friends, teachers, and mentors throughout my education. I would love to list them all, but I fear I will miss one or more, and the list would be much too long to expect any reader to trudge through. Hence, I have elected to name just six people who I felt played the largest part. To everyone else who has encouraged me, shaped me, helped me, pushed me, loved me, and challenged me throughout this exciting journey, please know that even if you are not called out by name I am still eternally grateful and I love and appreciate you: this would not be possible without you.

I believe a good parent is hard to find, and an excellent parent is exceedingly rare. Yet, I find myself with two excellent parents: my mother, Donna Leake, and father, Robert Leake. Without your guidance, love, and support I would not be where I am today. To Nina Rogerson, the woman I love, thank you for your companionship and support. It is has been an extraordinary six years, and I look forward to many, many more together. To Hunter Johnson and Daniele Mortari, without your friendship and intellectual support it would surely have been a dull and markedly less productive four years. To Aaron Schutte, thank you for your mentorship and support throughout two internships and an NSTRF fellowship.
\pagebreak{}

%% file: Data/contributors.tex
%
%
%
%

\chapter*{CONTRIBUTORS AND FUNDING SOURCES}
\addcontentsline{toc}{chapter}{CONTRIBUTORS AND FUNDING SOURCES}  

\subsection*{Contributors}
This work was supported by a dissertation committee consisting of Professor Daniele Mortari (advisor) and Professors John Hurtado and Moble Benedict of the Department of Aerospace Engineering and Professor Junuthula Reddy of the Department of Mechanical Engineering.

The Theory of Functional Connections was collaboratively developed by Daniele Mortari (advisor), Hunter Johnston (PhD candidate), and Carl Leake (author/PhD candidate). To clarify the major contributions of each, the following figure is included.

\begin{center}
    \includegraphics[width=\linewidth]{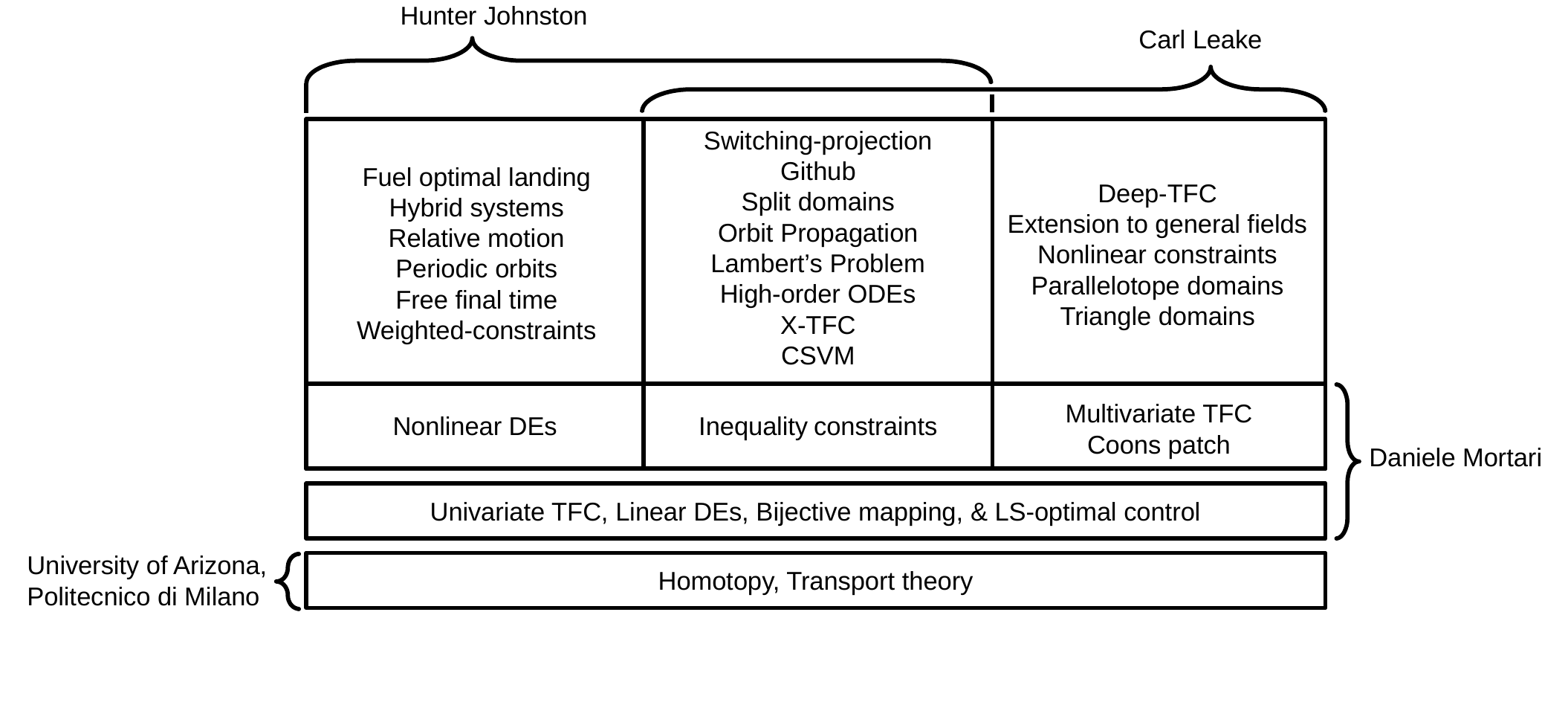}
\end{center}

In addition, Jonathan Cameron and Kevin Webb provided valuable insight and guidance on the natural balloon shape boundary-value problem, and those on the Venus Variable Altitude Aerobot project shared the Venus planetary data used to create Table \ref{tab:BalloonAtmData} and the balloon data used to create Table \ref{tab:BalloonData}.  All other work conducted for the dissertation was completed by the author independently.

\subsection*{Funding Sources}
Graduate study was supported by Texas A\&M University teaching and research assistantships from  August 2017 - August 2019 and by the NASA Space Technology Research Fellowship (NSTRF), grant number 80NSSC19K1152, from August 2019 - August 2021.
\pagebreak{}

%% file: Data/nomenclature.tex
%
%
%
%


\chapter*{NOMENCLATURE}
\addcontentsline{toc}{chapter}{NOMENCLATURE}  


\hspace*{-1.25in}
\vspace{12pt}
\begin{spacing}{1.0}
	\begin{longtable}[htbp]{@{}p{0.25\textwidth} p{0.72\textwidth}@{}}
	   API & Application Programming Interface \\ [2ex]
	   BFGS & Broyden-Fletcher-Goldfarb-Shanno \\ [2ex]
	   CSVM & Constrained Support Vector Machine \\ [2ex]
	   DE & Differential equation\\	[2ex]
	   Deep-TFC & Deep Theory of Functional Connections \\ [2ex]
	   ELM & Extreme Learning Machine \\ [2ex]
	   FEM & Finite Element Method \\ [2ex]
	   i.i.d. & Independently and identically distributed \\ [2ex]
	   JIT & Just-in-time (compiler) \\ [2ex]
	   JPL & Jet Propulsion Laboratory \\ [2ex]
	   LS & Least-squares\\ [2ex]
	   LS-SVM & Least-squares support vector machine \\ [2ex]
	   NN & Neural network \\ [2ex]
	   NSTRF  & NASA Space Technology Research Fellowship\\ [2ex]
	   ODE  & Ordinary differential equation\\ [2ex]
	   PDE  & Partial differential equation\\ [2ex]
	   SVM & Support vector machine \\ [2ex]
	   TFC  & Theory of Function Connections\\ [2ex]
	   XLA & Accelerated Linear Algebra \\ [2ex]
	   X-TFC & Extreme Theory of Functional Connections \\ [2ex]
	   $c_k$  & Slope in the linear map for the $k$-th independent variable that maps the basis function domain to the problem domain \\ [2ex]
	   $\pC{k}{i}$ & Constraint operator for the $i$-th constraint of the $k$-th independent variable \\ [2ex]
	   $g(\B{x})$ & Free function $\mathbb{R}^n\mapsto\mathbb{R}$. Note that a superscript may be used to denote the free function for a specific dependent variable, e.g., $g^u(\B{x})$ is the free function for the dependent variable $u$. \\ [2ex]
	   $\mathcal{J}$ & Jacobian matrix\\ [2ex]
	   $\mathbb{L}$ & Loss function $\mathbb{R}^m\mapsto\mathbb{R}^n$ \\ [2ex]
	   $L^1$ & Space of Lebesgue-integrable functions, i.e., $\int_\Omega |f| \dd{\mu} < \infty$\hspace*{\fill}\\ [2ex]
	   $L^2$ & Space of square-Lebesgue-integrable functions, i.e., $\int_\Omega |f|^2 \dd{\mu} < \infty$\hspace*{\fill}\\ [2ex]
	   $\mathbb{R}$ & Set of real numbers \\ [2ex]
	   $\mathbb{S}_{ij}$ & Support matrix \\ [2ex]
	   $\B{x}$ & A vector of the independent variables, i.e., $\B{x} = \{x_1,x_2,\cdots,x_n\}$, where $n$ is the number of independent variables. \\ [2ex]
	   $x_k$ & The $k$-th independent variable \\  [2ex]
	   $\mathbb{Z}$ & Set of integers \\ [2ex]
	   $\mathbb{Z}^+$ & Set of positive integers \\ [2ex]
	   $\mathbb{Z}/n\mathbb{Z}$ & Set of integers modulo $n$ \\ [2ex]
	   $z_k$ & Basis function domain variable for the $k$-th independent variable \\ [2ex]
	   $\delta_{ij}$ & Kronecker delta \\ [2ex]
	   $\p{k}{\kappa_i}(\B{x})$ & Portion of the $i$-th constraint on the $k$-th independent variable that does not contain the dependent variable. Note the pre-superscript is dropped in the univariate formulation  as there is only one independent variable. \\ [2ex]
	   $\mu(z)$ & Measure function $\mathbb{R}\mapsto\mathbb{R}$\\ [2ex]
	   $\p{k}{\rho}_i(\B{x},g(\B{x}))$ & Projection functional for the $i$-th constraint of the $k$-th independent variable. Note that in the univariate formulation, the pre-superscript is dropped as there is only one independent variable. \\ [2ex]
	   $\p{k}{\phi}_i(x_k)$ & Switching function for the $i$-th constraint on the $k$-th independent variable. Note that in the univariate formulation, the pre-superscript is dropped as there is only one independent variable. \\ [2ex]
	   $\Omega$ & Domain \\ [2ex]
	   $\mathds{1}(x,x_1)$ & Heaviside function, $\mathbb{R}\mapsto\mathbb{R}$ \\ [2ex]
	   $\mathds{1}_0(x)$ & Heaviside function where $x_1 = 0$ \\ [2ex]
	   $\ceil{x}$ & Rounds $x$ to the next largest integer \\ [2ex]
	   $\nint{x}$ & Rounds $x$ to the nearest integer \\ [2ex]
	\end{longtable}
\end{spacing}

\pagebreak{}

%% file: Data/lists.tex
%
%
%
%

\phantomsection
\addcontentsline{toc}{chapter}{TABLE OF CONTENTS}  

\begin{singlespace}
\renewcommand\contentsname{\normalfont} {\centerline{TABLE OF CONTENTS}}

\setcounter{tocdepth}{4} 

\setlength{\cftaftertoctitleskip}{1em}
\renewcommand{\cftaftertoctitle}{%
\hfill{\normalfont {Page}\par}}

\tableofcontents

\end{singlespace}

\pagebreak{}


\phantomsection
\addcontentsline{toc}{chapter}{LIST OF FIGURES}  

\renewcommand{\cftloftitlefont}{\center\normalfont\MakeUppercase}

\setlength{\cftbeforeloftitleskip}{-12pt} 
\renewcommand{\cftafterloftitleskip}{12pt}

\renewcommand{\cftafterloftitle}{%
\\[4em]\mbox{}\hspace{2pt}Figure\hfill{\normalfont Page}\vskip\baselineskip}

\begingroup

\begin{center}
\begin{singlespace}
\setlength{\cftbeforechapskip}{0.4cm}
\setlength{\cftbeforesecskip}{0.30cm}
\setlength{\cftbeforesubsecskip}{0.30cm}
\setlength{\cftbeforefigskip}{0.4cm}
\setlength{\cftbeforetabskip}{0.4cm}



\listoffigures

\end{singlespace}
\end{center}

\pagebreak{}

%
\phantomsection
\addcontentsline{toc}{chapter}{LIST OF TABLES}  

\renewcommand{\cftlottitlefont}{\center\normalfont\MakeUppercase}

\setlength{\cftbeforelottitleskip}{-12pt} 

\renewcommand{\cftafterlottitleskip}{1pt}

\renewcommand{\cftafterlottitle}{%
\\[4em]\mbox{}\hspace{2pt}Table\hfill{\normalfont Page}\vskip\baselineskip}

\begin{center}
\begin{singlespace}

\setlength{\cftbeforechapskip}{0.4cm}
\setlength{\cftbeforesecskip}{0.30cm}
\setlength{\cftbeforesubsecskip}{0.30cm}
\setlength{\cftbeforefigskip}{0.4cm}
\setlength{\cftbeforetabskip}{0.4cm}

\listoftables 

\end{singlespace}
\end{center}

\pagebreak{}

%
\phantomsection
\addcontentsline{toc}{chapter}{LIST OF EXAMPLES}  

\renewcommand{\cftextitlefont}{\center\normalfont\MakeUppercase}

\setlength{\cftbeforeextitleskip}{-12pt} 

\renewcommand{\cftafterextitleskip}{1pt}

\renewcommand{\cftafterextitle}{%
\\[4em]\mbox{}\hspace{2pt}Example\hfill{\normalfont Page}\vskip\baselineskip}

\begin{center}
\begin{singlespace}

\setlength{\cftbeforechapskip}{0.4cm}
\setlength{\cftbeforesecskip}{0.30cm}
\setlength{\cftbeforesubsecskip}{0.30cm}
\setlength{\cftbeforefigskip}{0.4cm}
\setlength{\cftbeforetabskip}{0.4cm}
\setlength{\cftbeforeExamplesskip}{0.4cm}

\setlength{\cftExamplesindent}{1.5em}
\setlength{\cftExamplesnumwidth}{2.3em}

\listofExamples

\end{singlespace}
\end{center}
\endgroup
\pagebreak{}  

%% file: Data/Introduction.tex

\pagestyle{plain} 
\pagenumbering{arabic} 
\setcounter{page}{1}

\chapter{INTRODUCTION AND MOTIVATION}

Interpolation is a concept that many use every day---for example, to extract an estimated value between experimental data points or between data points of tabulated values for a computationally expensive function---but pay little attention to. This is not surprising given that most imagine a discrete set of points when considering interpolation, a fairly simplistic problem, and the age of the algorithms used to solve them dates as far back as Waring polynomial interpolation, which was published in 1779 \cite{Waring}. Yet, if the interpolation complexity is increased, it quickly becomes difficult or impossible to point to an algorithm that can easily perform the interpolation.

Consider the advantages of satisfying properties at points other than their value, for example, their derivatives or a linear combination of derivatives and values between points. Rather than thinking of properties at points, it may be easier to imagine them as constraints, for example, $u_x(x_0) + \pi u(x_1) = v(x_2)$, where $x$ is an independent variable, $x_0$, $x_1$, and $x_2$ are some specific values in the domain, $u$ and $v$ are dependent variables, and $u_x$ denotes a derivative of $u$ with respect to $x$. Dream bigger. What if one could do this in $n$-dimensional domains or write \emph{all} possible functions that satisfy the constraints rather than just one function that satisfies the constraints? At this point, one is describing something much more complicated than simple point-wise interpolation; rather, they are describing a sort of function-based interpolation.

A rich framework for function-based interpolation could transform problems with linear constraints into unconstrained problems. In terms of optimization-type problems, this would mean one could use simpler optimizers, as an optimizer that handles constraints would no longer be needed, and/or the function to be minimized would not need to be augmented to include the constraints. Indeed, such a framework would enhance one's ability to solve such problems and is the driving motivation behind the Theory of Functional Connections (TFC): a general framework for function-based interpolation. 

The concept of function-based interpolation itself is not new, and numerous methods exist \cite{interp_1,interp_2,interp_3,interp_x,interp_4,interp_5}; however, these previous techniques only work for a class or sub-class of functions and cannot be used to describe all functions that satisfy a set of constraints. Therefore, their scope of applications is limited; TFC does not have this restriction.

\section{Original Idea}
The idea that sparked the Theory of Functional Connections (TFC) was conceived by Daniele Mortari while teaching the Waring, better known as Lagrange, polynomial interpolation method \cite{Waring}, which is used to generate an interpolating function that passes through a set of points. For example, the Lagrange polynomial, $y(x)$, for a set of $n$ points, $(x_1,y_1),\dots,(x_k,y_k),\dots,(x_n,y_n)$, can be written as,
\begin{equation*}
    y (x) = \ds\sum_{k=1}^n  y_k \prod_{i\ne k} \dfrac{x - x_i}{x_k - x_i}.
\end{equation*}
The Lagrange polynomial represents one function that passes through the $n$ points. In other words, the Lagrange polynomial is an interpolating function for these points: a function that satisfies the constraints $y(x_k) = y_k$.

Mortari's original insight was that by replacing $x$ with an arbitrary function $g(x)$, $x_k$ with $g(x_k)$, and so on, one could write the the family of \emph{all} possible functions that passes through the set of points. That is,
\begin{equation*}
    y (x) = \ds\sum_{k=1}^n  y_k \prod_{i\ne k} \dfrac{x - x_i}{x_k - x_i} \qquad \to \qquad y (x, g(x)) = \ds\sum_{k=1}^n y_k \prod_{i\ne k} \dfrac{g(x) - g(x_i)}{g(x_k) - g(x_i)}.
\end{equation*}
From this seed of an idea sprouted Mortari's 2016 seminal article \cite{U-TFC} that demonstrated how to embed univariate value and derivative constraints, and constraints consisting of linear combinations of values and derivatives at points, into the TFC framework. Since then, the TFC framework has grown to encompass a larger variety of increasingly exotic constraints, including integral \cite{Integral-TFC}, component \cite{Component-TFC}, and inequality \cite{Inequality-TFC} constraints. Moreover, the original univariate framework has been extended to multiple variables \cite{M-TFC2,M-TFC} and some non-rectangular domains \cite{TFC-Selected,BijectiveMapping-TFC}. 

\section{Overview of the Remaining Chapters}
In lieu of a large literature review concentrated at the beginning of the dissertation, most chapters contain their own smaller literature review whose contents pertain specifically to that chapter. The rest of this dissertation is structured as follows. 

\subsection*{Chapter \ref{chap:tfcTheory}. Theory of Functional Connections}
This chapter describes the theory behind the TFC functional interpolation framework. It is split into two major sections: the univariate theory and the multivariate theory. The univariate theory is introduced first and describes how to construct \ces\ for value, derivative, integral, and component constraints, and linear combinations thereof. In addition, it includes mathematical theorems that pertain to the univariate \ce. The multivariate theory section generalizes this to $n$-dimensions. Examples are included throughout to help solidify the reader's understanding. 

\subsection*{Chapter \ref{chap:deApplications}. Applications in Differential Equations}
This chapter utilizes the TFC framework introduced in the previous chapter to solve differential equations by embedding the differential equation constraints into the constrained expression and using the free function to minimize the differential equation's residual at a discrete set of points. In addition, the chapter discusses useful free function choices and optimization methods and includes a summary of the TFC numerical implementation. To strengthen the reader's understanding, a simple PDE is numerically estimated using each of the common free function choices; additional examples are provided that highlight the strengths and weaknesses of each free function choice. 

\subsection*{Chapter \ref{chap:flexBodyApplications}. Applications in Flexible Body Problems}
Building on the foundations of the previous chapters, this chapter utilizes the TFC framework introduced in Chapter \ref{chap:tfcTheory} and its application to differential equations introduced in Chapter \ref{chap:deApplications} to apply the method to differential equations that appear in flexible body problems. In other words, this chapter contains flexible-body-related ODEs and PDEs that are solved via TFC.

\subsection*{Chapter \ref{chap:Conclusions}. Summary and Conclusions}
This chapter summarizes the major ideas covered in the dissertation and draws conclusions based on the content discussed throughout. In addition, this chapter presents ideas for future study. 

\subsection*{Appendices}
The appendices include more detailed explanations of some of the topics covered in the main body of the text. In addition, they also include the following extensions of the TFC framework: nonlinear constraints, inequality constraints, parallelotope domains, lower-dimensional constraints in $n$-dimensions, and an extension to general fields, i.e., beyond the field of real numbers. 

%% file: Data/TfcTheory.tex

\chapter{THEORY OF FUNCTIONAL CONNECTIONS}\label{chap:tfcTheory}
The seminal article on the Theory of Functional Connections\footnote{This theory was originally published under the name ``Theory of Connections.'' However, this name conflicted with a specific theory in differential geometry and was not the most accurate description of the functional interpolation method. Therefore, in 2019, this name was changed to the ``Theory of Functional Connections'' to highlight the tie to functional interpolation and the fact that it provides \emph{all} functions satisfying a set of linear constraints in $n$-dimensional space.} (TFC) introduced the notion of a \emph{\ce} \cite{U-TFC}: a mathematical expression that utilizes a function that can be chosen by the user, the so-called free function, that can describe \emph{all} possible functions satisfying a given set of linear constraints. At the time the article was written, this statement was merely a conjecture but has since been proven mathematically. Since their conception, the process for deriving \ces\ and the language used to discuss them has changed, but what they are in mathematical terms has remained constant: \ces\ are functionals.\footnote{In other literature, functionals are also referred to as ``functions of functions'' or ``higher-order functions.''} Therefore, it is useful to define a functional and investigate some of its properties before delving further into TFC \ces.

\section{Functionals}
To begin, consider the following definition of a functional.

\begin{definition}
A functional, e.g., $f(x,g(x))$, has independent variable(s) and function(s) as inputs and produces a function as an output.
\end{definition}

\noindent Note that a functional as defined here coincides with the computer science definition of a functional. One can think of a functional as a map for functions. That is, the functional takes a variable or variables and a function or functions as inputs and produces a function as its output, e.g., $f^*(x) = f(x,g(x))$. This dissertation is focused on constraint embedding, or in other words, functional interpolation; hence, for now, there is no need to concern oneself with the domains and ranges of the input and output functions. Rather, functionals can be discussed in the context of their potential input functions, hereon referred to as the domain of the functional, and potential output functions, hereon referred to as the codomain of the functional. 

\begin{figure}[!ht]
    \centering\includegraphics[width=.65\linewidth]{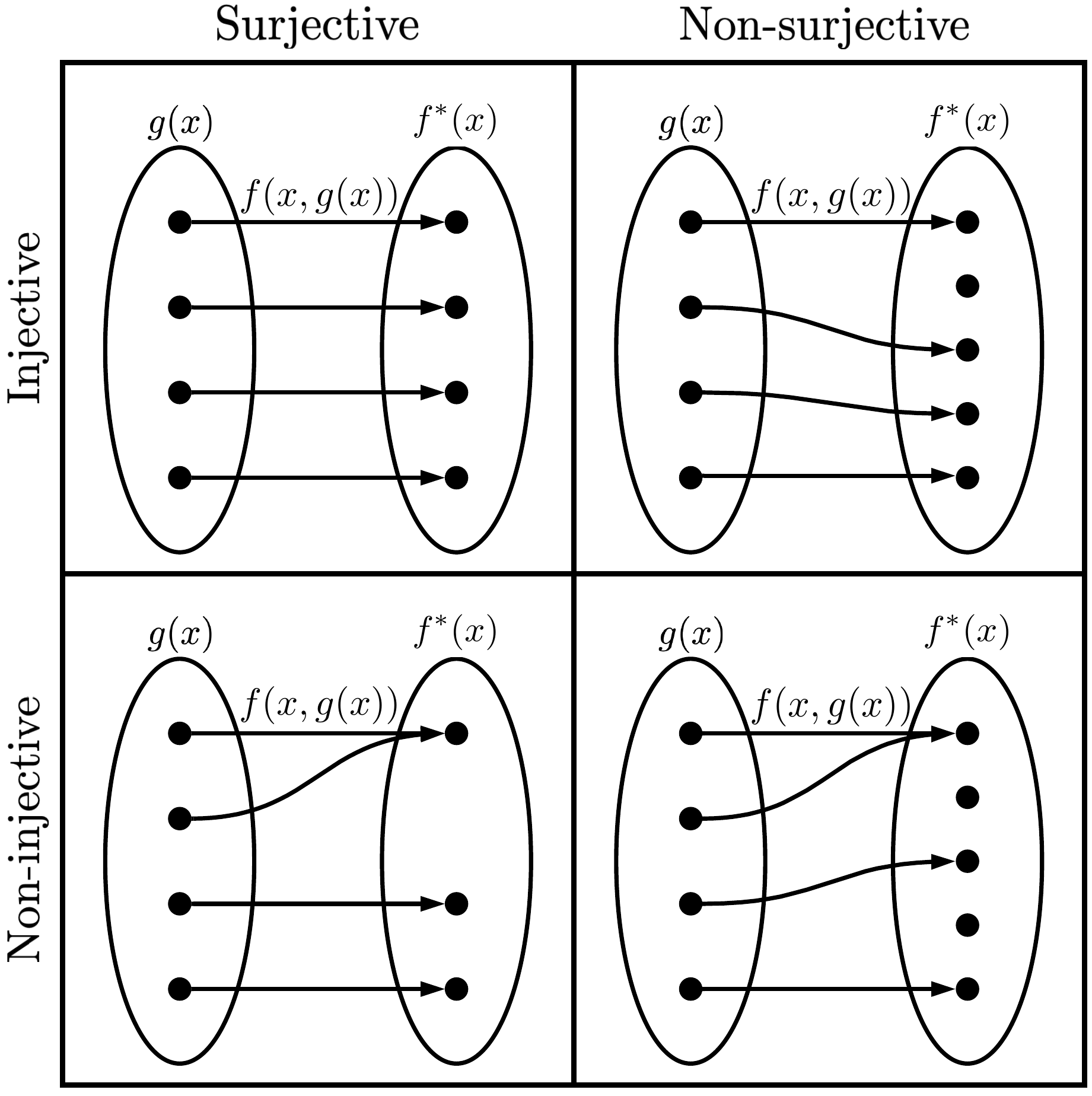}
    \caption{Graphical representation of injective and surjective functionals.}
    \label{fig:Injective_Surjective}
\end{figure}

Next, the definitions of injective, surjective, and bijective are extended from functions to functionals. 
\begin{definition}
A functional is injective if every function in its codomain is the image of at most one function in its domain. 
\end{definition}
\begin{definition}
A functional, $f(x,g(x))$, is surjective if for every function in the codomain, $f^*(x)$, there exists at least one function, $g(x)$, in the domain such that $f^*(x) = f(x,g(x))$.
\end{definition}
\begin{definition}
A functional is bijective if it is both injective and surjective. 
\end{definition}
To elaborate, Figure \ref{fig:Injective_Surjective} gives a graphical representation of each of these functionals, and examples of each of these functionals follow. Note that the phrase ``smooth functions'' is used here to denote continuous, infinitely differentiable, real-valued functions. 

Consider the functional $f(x,g(x)) = e^{-g(x)}$ whose domain is all smooth functions and whose codomain is all smooth functions. The functional is injective because for every $f^*(x)$ in the codomain there is at most one $g(x)$ that maps $f(x,g(x))$ to $f^*(x)$. However, the functional is not surjective, because the functional does not span the codomain. For example, consider the desired output function $f^*(x) = -2$: there is no $g(x)$ that produces this output. 

Next, consider the functional $f(x,g(x)) = g(x)-g(0)$ whose domain is all smooth functions and whose codomain is all smooth functions $f^*(x)$ such that $f^*(0) = 0$. This functional is surjective, because it spans the set of all smooth functions that are $0$ when $x=0$, but it is not injective. For example, the functions $g(x) = x$ and $g(x) = x+3$ produce the same result, i.e., $f(x,x) = f(x,x+3) = x$.
%

Finally, consider the functional $f(x,g(x)) = g(x)$ whose domain is all smooth functions and whose codomain is all smooth functions. This functional is bijective, because it is both injective and surjective. 

In addition, the notion of projection is extended to functionals. Consider an analogy to vector projection wherein a projection matrix, i.e., an idempotent matrix $P^n = P\ \forall n\in\mathbb{Z}^+$, projects a vector from one vector space to another. In other words, the properties of $P$ are (1) it transports vectors from one vector space to another, and (2) when it operates on itself (the operator being matrix multiplication), it produces itself ($P^n = P$). A projection property for functionals can be defined similarly. Functionals already have the first property: they transport functions from one set, their domain, to another set, their codomain, e.g., the \ce\ transports functions from the set of all real-valued functions defined at the constraints to the set of functions that satisfy the constraints. Following the analogy, if a functional produces itself when operating on itself, where the operator is using the functional's output as its input function, then that functional is said to be a projection functional.

\begin{definition}
A functional is said to be a projection functional if it produces itself when operating on itself.
\end{definition}
For example, consider a functional operating on itself, $f(x,f(x,g(x)))$. If\\$f(x,f(x,g(x))) = f(x,g(x))$, then the functional is a projection functional. Note that proving $f(x,f(x,g(x))) = f(x,g(x))$ automatically extends to a functional operating on itself $n$ times where $n\in\mathbb{Z}^+$: for example, $f(x,f(x,f(x,g(x))) = f(x,f(x,g(x))) = f(x,g(x))$, and so on.

\section{Univariate Theory}

The majority of this dissertation focuses on multivariate TFC; however, the multivariate TFC framework is built by recursively applying univariate TFC. Hence, it is paramount the reader understands univariate TFC before moving to the multivariate case. First, the original form of the univariate \ce\ from Reference \cite{U-TFC} will be presented via Example \ref{ex:UniPointConstraints}. Then, the constrained expression will be manipulated to expose an underlying structure made up of so-called projection functionals and switching functions \cite{M-TFC2}. Throughout the remainder of this section, that structure will be utilized to create \ces\ for various types of linear constraints and prove mathematical theorems related to univariate \ces.  

\begin{example}{Constraints at a point}\mbox{}\label{ex:UniPointConstraints}
Constraints at a point consist of constraints on the value and derivatives at the point. Consider the follow constraints,
\begin{equation*}
    y(0) = 1, \quad y_x(1) = 2, \andd y(2) = 3.
\end{equation*}

Given a set of $k$ point constraints, the univariate \ce\ takes the following form \cite{U-TFC},
\begin{equation}\label{eq:uniCE}
    y (x, g (x)) = g (x) + \sum_{j = 1}^k s_j (x) \, \eta_j(x,g(x)),
\end{equation}
where $g(x)$ is a free function, $s_j(x)$ are $k$ linearly independent functions called support functions, and $\eta_j(x,g(x))$ are $k$ coefficient functionals that are solved by imposing the constraints. The free function $g(x)$ can be chosen to be any function provided that it is defined at the constraints' locations. 

For this example, the support functions are chosen to be $s_1(x) = 1$, $s_2(x) = x^2$, and $s_3(x) = x^3$. Following Equation \eqref{eq:uniCE} and imposing the three constraints leads to the simultaneous set of equations
\begin{align*}
    y(0) &= 1 = g(0) + \eta_1(x,g(x))\\
    y_x(1) &= 2 = g_x(1) + 2\eta_2(x,g(x)) + 3\eta_3(x,g(x))\\
    y(2) &= 3 = g(2) + \eta_1(x,g(x)) + 4 \eta_2(x,g(x)) + 8 \eta_3(x,g(x)).
\end{align*}

\noindent Solving this set of equations for the unknowns $\eta_j(x,g(x))$ leads to the solution,
\begin{align*}
    \eta_1(x,g(x)) &= 1 - g(0) \\
    \eta_2(x,g(x)) &= \frac{10 - 3 g (0) + 3 g (2) - 8 g_x (1)}{4} \\
    \eta_3(x,g(x)) &= \frac{g (0) - g (2) + 2 g_x (1)}{2}.
\end{align*}

\noindent Substituting the coefficient functionals back into Equation \eqref{eq:uniCE} and simplifying yields,
\begin{equation}\label{eq:uniEx1Soln}
\begin{aligned}
    y(x,g(x)) = g(x) &+ \frac{-2 x^3 + 3 x^2 + 4}{4} \Big(1 - g (0)\Big) + \Big(-x^3 + 2 x^2\Big) \Big(2 - g_x (1)\Big)\\
    &+\frac{2 x^3 - 3 x^2}{4}\Big(3 - g (2)\Big).
\end{aligned}
\end{equation}

\noindent It is simple to verify that regardless of how $g (x)$ is chosen, provided $g (x)$ exists at the constraint points, Equation \eqref{eq:uniEx1Soln} always satisfies the given constraints. 

The support functions in the previous example were selected as $s_1(x) = 1$, $s_2(x) = x^2$, and $s_3(x) = x^3$. However, these support functions could have been any linearly independent set of functions that permits a solution for the coefficient functionals $\eta_j(x,g(x))$; to clarify the latter of these requirements, consider the same constraints with support functions $s_1(x) = 1$, $s_2(x) = x$, and $s_3(x) = x^2$. Then, the set of equations with unknowns $\eta_j(x,g(x))$ is,
\begin{equation*}
    \begin{bmatrix} 1 & 0 & 0 \\ 0 & 1 & 2 \\ 1 & 2 & 4\end{bmatrix} \begin{Bmatrix} \eta_1(x,g(x)) \\ \eta_2(x,g(x)) \\ \eta_3(x,g(x))\end{Bmatrix} = \begin{Bmatrix} 1 - g (0) \\ 2 - g_x (1) \\ 3 - g (2)\end{Bmatrix}.
\end{equation*}
Notice that when using these support functions, the matrix that multiplies the coefficient functionals is singular. Thus, no solution exists, and therefore, the support functions $s_1(x) = 1$, $s_2(x) = x$, and $s_3(x) = x^2$ are an invalid set for these constraints. Note that the matrix singularity does not depend on the free function. This means that the singularity arises when a linear combination of the selected support functions cannot be used to interpolate the constraints. Therefore, the support function matrix's singularity depends on both the support functions chosen and the specific constraints to be embedded. This raises another important restriction on the expression of the support functions: not only must they be linearly independent, but they must constitute an interpolation model that is consistent with the specified constraints.

Notice that each term, except the term containing only the free function, in the constrained expression is associated with a specific constraint and has a particular structure. To illustrate, examine the first constraint term from Equation \eqref{eq:uniEx1Soln}, 
\begin{equation*}
    \underbrace{\frac{-2 x^3 + 3 x^2 + 4}{4}}_{\phi_1 (x)}\underbrace{(1 - g (0))}_{\rho_1(x,g(x))}.
\end{equation*}
The first term in the product, $\phi_1 (x)$, is called a {\it switching function}\footnote{Reference \cite{U-TFC} introduced these switching functions as ``coefficient'' functions, $\B{\beta}_k$, but they were not used in the same way the switching-projection form uses them.} and is a function that is equal to $1$ when evaluated at the constraint it is referencing and equal to $0$ when evaluated at all the other constraints. For example, when evaluating the switching function $\phi_1 (x)$ at the constraint it is referencing it is equal to 1, i.e., $\phi_1 (0) = 1$, and when it is evaluated at the other constraints it is equal to $0$, i.e., $\frac{\partial \phi_1}{\partial x} (1) = 0$ and $\phi_1 (2) = 0$. The second term of the product, $\rho_1 (x, g (x))$, is called a {\it projection functional}, and is derived by setting the constraint function equal to zero and replacing $y(x)$ with $g(x)$. In the case of constraints at a point, this is simply the difference between the constraint value and the free function evaluated at the constraint point. It is called the projection functional because it projects the free function to the set of functions that vanish at the constraint. 

The switching-projection structure is important because it shows up in other constraint types too. Based on this structure, an alternate way to define the constrained expression can be derived,
\begin{equation}\label{eq:uniCeAlt}
    y (x, g (x)) = g (x) + \sum_{j = 1}^k \phi_j (x) \, \rho_j(x,g(x)).
\end{equation}

For this case, the projection functionals are simple to derive, but the switching functions require some attention. From their definition, these functions must go to $1$ at their associated constraint and $0$ at all other constraints. Hence, the following algorithm for deriving the switching functions is proposed:
\begin{enumerate}
    \item Choose $k$ support functions, $s_k(x)$.
    \item Write each switching function as a linear combination of the support functions with unknown coefficients.
    \item Based on the switching function definition, write a system of equations to solve for the unknown coefficients. 
\end{enumerate}

To validate that this algorithm works, consider the same constraints and support functions and rederive the \ce\ shown in Equation \eqref{eq:uniEx1Soln}. Hence, $\phi_1 (x) = s_i (x) \, \alpha_{i1}$, $\phi_2 (x) = s_i (x) \, \alpha_{i2}$, and $\phi_3 (x) = s_i (x) \, \alpha_{i3}$, for some as yet unknown coefficients $\alpha_{ij}$. Note that in the previous mathematical expressions and throughout the remainder of the dissertation, the Einstein summation convention is used to improve readability. Now, the definition of the switching function is used to come up with a set of equations. For example, the first switching function has the three equations,
\begin{equation*}
    \phi_1(0) = 1, \quad \frac{\partial \phi_1}{\partial x}(1) = 0, \quad \text{and} \quad \phi_1(2) = 0.
\end{equation*}

\noindent These equations are expanded in terms of the support functions,
\begin{align*}
    \phi_1(0) &= (1) \cdot \alpha_{11} + (0) \cdot\alpha_{21} + (0) \cdot\alpha_{31} = 1\\
    \frac{\partial \phi_1}{\partial x}(1) &= (0)\cdot \alpha_{11} + (2) \cdot\alpha_{21} + (3) \cdot\alpha_{31} = 0\\
    \phi_1(2) &= (1) \cdot\alpha_{11} + (4) \cdot\alpha_{21} + (8) \cdot\alpha_{31} = 0,
\end{align*}
which can be compactly written as,
\begin{equation*}
    \begin{bmatrix} 1 & 0 & 0 \\ 0 & 2 & 3 \\ 1 & 4 & 8\end{bmatrix} \begin{Bmatrix} \alpha_{11} \\ \alpha_{21} \\ \alpha_{31} \end{Bmatrix} = \begin{Bmatrix} 1 \\ 0 \\ 0 \end{Bmatrix}.
\end{equation*}

\noindent The same is done for the other two switching functions to produce a set of equations that can be solved by matrix inversion.
\begin{align*}
    \begin{bmatrix} 1 & 0 & 0 \\ 0 & 2 & 3 \\ 1 & 4 & 8\end{bmatrix} \begin{bmatrix} \alpha_{11} & \alpha_{12} & \alpha_{13} \\ \alpha_{21} & \alpha_{22} & \alpha_{23} \\ \alpha_{31} & \alpha_{32} & \alpha_{33} \end{bmatrix} &= \begin{bmatrix} 1 & 0 & 0 \\ 0 & 1 & 0 \\ 0 & 0 & 1\end{bmatrix} \\
     \begin{bmatrix} \alpha_{11} & \alpha_{12} & \alpha_{13} \\ \alpha_{21} & \alpha_{22} & \alpha_{23} \\ \alpha_{31} & \alpha_{32} & \alpha_{33} \end{bmatrix} &= \begin{bmatrix} 1 & 0 & 0 \\ 0 & 2 & 3 \\ 1 & 4 & 8\end{bmatrix}^{-1} = \begin{bmatrix} 1 & 0 & 0 \\ \frac{3}{4} & 2 & -\frac{3}{4} \\ -\frac{1}{2} & -1 & \frac{1}{2}\end{bmatrix}.
\end{align*}

\noindent Substituting the constants back into the switching functions and simplifying yields,
\begin{equation*}
    \phi_1 (x) = \frac{-2 x^3 + 3 x^2 + 4}{4}, \quad \phi_2 (x) =-x^3 + 2 x^2, \quad \text{and} \quad \phi_3 (x) = \frac{2 x^3 - 3 x^2}{4}.
\end{equation*}

\noindent Substituting the projection functionals and switching functions back into the \ce\ shown in Equation \eqref{eq:uniCeAlt} yields,
\begin{align*}
    y(x,g(x)) =\ &g(x) + \frac{-2 x^3 + 3 x^2 + 4}{4} \Big(1 - g (0)\Big) + \Big(-x^3 + 2 x^2\Big) \Big(2 - g_x (1)\Big)\\
    &+\frac{2 x^3 - 3 x^2}{4}\Big(3 - g (2)\Big),
\end{align*}
which is identical to Equation \eqref{eq:uniEx1Soln}. 
\end{example}

As demonstrated in Example \ref{ex:UniPointConstraints}, the switching-projection approach, Equation \eqref{eq:uniCeAlt}, is a valid method for deriving \ces; although it was only demonstrated for one set of constraints here, this \ce\ derivation technique will be proven mathematically in Section \ref{subsec:UniProofs}. Similar to the original approach, Equation \eqref{eq:uniCE}, there is a risk of obtaining a singular matrix when solving for $\alpha_{ij}$ if the support functions selected are not able to interpolate the constraints. However, as will be demonstrated in the sections that follow, the switching-projection approach can be used for many constraint types, easily extended to multivariate domains via recursive applications of the univariate theory, and lends itself nicely to mathematical proofs. Before moving to these other topics, it is useful to first examine the anatomy of a constraint and define the so-called constraint operator. This analysis will prove invaluable as it provides a method to unify the way linear constraints are written. Consequently, mathematical analyses can be done on this unified form and thereby applied to all linear constraints: without this method, one would need to conduct the same mathematical analysis for each constraint type separately.

\subsection{Anatomy of a Linear Constraint}
Linear constraints can be conveniently dissected into two portions: (1) an operator that operates on a dependent variable and (2) the remaining constants and functions of the constraint. Let the former be called the constraint operator and denoted by the symbol $\C{}$ and the latter denoted by the symbol $\kappa$. Using this nomenclature, a constraint on the dependent variable $y$ would typically be written in the form,
\begin{equation*}
    \kappa = \C{}[y].
\end{equation*}
For example, the constraint $3 = 2 y(2) - \pi y_{xx}(0)$ consists of $\kappa = 3$ and $\C{}[y] = 2 y(x) - \pi y_{xx}(0)$. Definition \ref{def:constraint_operator} defines the constraint operator more rigorously.

\begin{definition}\label{def:constraint_operator}
The constraint operator, $\C{i}$, is a linear operator that operates on a function and returns the function evaluated at the $i$-th specified constraint.
\end{definition}
\noindent The word evaluation in the previous definition requires some elaboration; evaluation means to evaluate the operand function in the same way as the dependent variable in the constraint. Notice that this means the constraint operator is not affected by terms in the constraint that do not contain the dependent variable. As an example, again consider the constraint $3 = 2 y(2) + \pi y_{xx}(0)$, and suppose it is the first constraint in the set  ($i = 1$). For this constraint, the constraint operator operates as follows,
\begin{equation*}
    \C{1} [f(x)] =  2 f(2) + \pi f_{xx}(0).
\end{equation*}

In addition, notice that the constraint operator satisfies the two properties of a linear operator:
\begin{enumerate}
    \item $\C{i} [f(x) + g(x)] = \C{i}[f(x)] + \C{i}[g(x)]$
    \item $\C{i}[a g(x)] = a\C{i}[g(x)]$
\end{enumerate}
For example, again consider the linear constraint $3 = 2 y(2) + \pi y_{xx}(0)$,
\begin{align*}
    \C{1} [f(x)+g(x)] &= \C{1} [f(x)] + \C{1}[g(x)] = 2 f(2) + \pi f_{xx}(0) + 2 g(2) + \pi g_{xx}(0) \\
    \C{1}[a f(x)] &= a \C{1} [f(x)] = a \Big( 2 f(2) + \pi f_{xx}(0)\Big).
\end{align*}

Naturally, the constraint operator has specific properties when operating on the support functions, switching functions, and projection functionals.
\begin{property}
The constraint operator acting on the support functions $s_j (x)$ produces the support matrix \begin{equation*}
    \mathbb{S}_{ij} = \C{i}[s_j(x)].
\end{equation*}
\end{property}
Consider the example given in Example \ref{ex:UniPointConstraints} where the support functions were $s_1(x) = 1$, $s_2(x) = x^2$, and $s_3(x) = x^3$. By applying the constraint operator,
\begin{align*}
    \mathbb{S}_{ij} &= \C{i}[s_j(x)] = \begin{bmatrix}\C{1}[s_1(x)] & \C{1}[s_2(x)] & \C{1}[s_3(x)] \\ \C{2}[s_1(x)] & \C{2}[s_2(x)] & \C{2}[s_3(x)] \\ \C{3}[s_1(x)] & \C{3}[s_2(x)] & \C{3}[s_3(x)] \end{bmatrix} \\
    &= \begin{bmatrix}s_1(0) & s_2(0) & s_3(0) \\ \frac{\partial s_1}{\partial x}(1) & \frac{\partial s_2}{\partial x}(1) & \frac{\partial s_3}{\partial x}(1) \\ s_1(2) & s_2(2) & s_3(2) \end{bmatrix} = \begin{bmatrix} 1 & 0 & 0 \\ 0 & 2 & 3 \\ 1 & 4 & 8 \end{bmatrix},
\end{align*}
which is identical to the support matrix from Example \ref{ex:UniPointConstraints}. It follows that $\mathbb{S}_{ij} \, \alpha_{jk} =  \alpha_{ij} \, \mathbb{S}_{jk} = \delta_{ik}$, where $\delta_{ik}$ is the Kroneker delta, and the solution of the $\alpha_{ij}$ coefficients can be determined by simply inverting the support matrix. 

\begin{property}\label{prop:SwitchDefinition}
The constraint operator acting on the switching functions $\phi_j(x)$ produces the Kronecker delta.  
\begin{equation*}
    \C{i}[\phi_j(x)] = \C{i}[s_k(x) \alpha_{kj}] = \C{i}[s_k(x)] \alpha_{kj} = \mathbb{S}_{ik}\alpha_{kj} = \delta_{ij}
\end{equation*}
\end{property}

\noindent This property is just a mathematical restatement of the linguistic definition of the switching function given earlier. One can intuit this property from the switching function definition, since they evaluate to $1$ at their specified constraint condition, i.e., $i=j$, and to $0$ at all other constraint conditions, i.e., $i \neq j$.

Using the constraint operator definition, one can define the projection functional in a compact and precise manner.
\begin{definition}\label{def:projection_function}
The projection functional is the difference between the numerical portion of the constraint and the constraint operator acting on the free function. Mathematically,
\begin{equation*}
    \rho_i(x,g(x)) = \kappa_i - \C{i}[g(x)].
\end{equation*}
\end{definition}
\begin{definition}
The univariate free function is any function $g(x)\colon\mathbb{R}\mapsto\mathbb{R}$ such that $\C{i}[g]$ is defined. 
\end{definition}

\noindent Again, consider the constraint $3 = 2 y(2) + \pi y_{xx}(0)$,
\begin{equation*}
    \rho_1(x,g(x)) = \kappa_1 - \C{1}[g(x)] = 3 - 2 g(2) - \pi g_{xx}(0).
\end{equation*}

\noindent Note that in the univariate case, $\kappa_i$ is a scalar value, i.e., $\kappa_i\in\mathbb{R}$, but in the multivariate case, $\kappa_i$ can be a function. In addition, notice what happens if $g(x)$ is a function that already satisfies the constraints.
\begin{property}\label{prop:projZero}
If $g(x)$ is a function that satisfies the constraints, then the projection functional is equal to zero. 
\end{property}
Property \ref{prop:projZero} follows from the definition of the projection functional; if $g(x)$ satisfies the constraints, then,
\begin{align*}
    \rho_i(x,g(x)) &= \kappa_i-\C{i}[g(x)] \\
    &= \kappa_i - \kappa_i \\
    &= 0.
\end{align*}

Now that the constraint operator has been defined, and consequently, rigorous definitions for the projection functionals and switching functions have been provided, other constraint types become easy to embed into univariate \ces: For example, integral constraints.
\begin{example}{Integral constraints}
Consider the following set of constraints,
\begin{equation*}
    \int_{-2}^3 y(x) \dd{x} = 5 \andd \int_0^2 3 y(x) \dd{x} = 2.
\end{equation*}
Based on Definition \ref{def:projection_function}, the projection functionals for these constraints can be written as,
\begin{align*}
    \rho_1(x,g(x)) &= 5-\int_{-2}^3 g(\tau) \dd{\tau}\\
    \rho_2(x,g(x)) &= 2-\int_0^2 3g(\tau) \dd{\tau}.
\end{align*}
Notice that the integrals in the projection functionals use a dummy variable, $\tau$, rather than $x$. Furthermore, based on Property \ref{prop:SwitchDefinition}, the switching function equations can be written as,
\begin{align*}
    \int_{-2}^3 \phi_1(x) \dd{x} &= 1, &&\int_0^2 3\phi_1(x) \dd{x} = 0,\\
    \int_{-2}^3 \phi_2(x) \dd{x} &= 0, &&\int_0^2 3\phi_2(x) \dd{x} = 1.
\end{align*}
Setting $\phi_1(x)$ and $\phi_2(x)$ to be a linear combination of the support functions $s_1(x) = 1$  and $s_2(x) = x$ with unknown coefficients $\alpha_{ij}$ yields,
\begin{align*}
    \begin{bmatrix} 5 & \frac{5}{2} \\ 6 & 6\end{bmatrix} \begin{bmatrix} \alpha_{11} & \alpha_{12} \\ \alpha_{21} & \alpha_{22} \end{bmatrix} &= \begin{bmatrix} 1 & 0 \\ 0 & 1 \end{bmatrix}\\
    \begin{bmatrix} \alpha_{11} & \alpha_{12} \\ \alpha_{21} & \alpha_{22} \end{bmatrix} &= \begin{bmatrix} \frac{2}{5} & -\frac{1}{2} \\ -\frac{1}{6} & \frac{1}{3} \end{bmatrix}.
\end{align*}
Hence,
\begin{equation*}
    \phi_1(x) = \frac{2-2x}{5} \andd \phi_2(x) = \frac{2x-1}{6}.
\end{equation*}
Thus, following Equation \eqref{eq:uniCeAlt}, the \ce\ for these constraints is,
\begin{equation*}
    y(x,g(x)) = g(x) + \frac{2-2x}{5}\Big(5-\int_{-2}^3 g(\tau) \dd{\tau}\Big) + \frac{2x-1}{6}\Big( 2-\int_0^2 3g(\tau) \dd{\tau} \Big).
\end{equation*}
The previous \ce\ will always satisfy the constraints regardless of how the free function, $g(x)$, is chosen.
\end{example}

\subsection{Component Constraints}
When handling component constraints, one must decide which dependent variable's \ce\ the component constraint will be embedded into. This dependent variable will define the constraint operator, and all other dependent variables will become part of the constraint's $\kappa$ term. Regardless of which dependent variable is chosen, a valid \ce\ will be produced.

\begin{example}{Component constraints}\label{ex:UniComponentConstraints}
Consider the following set of constraints,
\begin{equation*}
    u(0)+v(0) = 5 \andd u_x(2) + v(3) = 4.
\end{equation*}
Two different sets of \ces\ will be produced: one where the component constraints are embedded into the \ce\ for $u$, and the second where the component constraints are embedded into the \ce\ for $v$. If the constraints are embedded into $u$, then the projection functionals are,
\begin{align*}
    \rho_1(x,g^u(x),g^v(x)) &= 5-g^u(0)-v(0,g^v(x))\\ \rho_2(x,g^u(x),g^v(x)) &= 4-g^u_x(2)-v(3,g^v(x)),
\end{align*}
where $g^u(x)$ is the free function used in the $u$ \ce; similarly, $g^v(x)$ will be the free function used in the $v$ \ce. The equations for the switching function are,
\begin{align*}
    \phi_1(0) = 1,\quad & \quad \frac{\partial \phi_1}{\partial x}(2) = 0\\
    \phi_2(0) = 0,\quad & \quad \frac{\partial \phi_2}{\partial x}(2) = 1.
\end{align*}
Let the support functions be $s_1(x) = 1$ and $s_2(x) = x$, then
\begin{align*}
    \begin{bmatrix} 1 & 0 \\ 0 & 1\end{bmatrix} \begin{bmatrix} \alpha_{11} & \alpha_{12} \\ \alpha_{21} & \alpha_{22} \end{bmatrix} &= \begin{bmatrix} 1 & 0 \\ 0 & 1 \end{bmatrix}\\
    \begin{bmatrix} \alpha_{11} & \alpha_{12} \\ \alpha_{21} & \alpha_{22} \end{bmatrix} &= \begin{bmatrix} 1 & 0 \\ 0 & 1 \end{bmatrix}.
\end{align*}
Thus, the switching functions are,
\begin{equation*}
    \phi_1(x) = 1 \andd \phi_2(x) = x,
\end{equation*}
and the first set of \ces, where the component constrains are embedded into $u$ is,
\begin{equation}\label{eq:ComponentExSet1}
\begin{aligned}
    u(x,g^u(x),g^v(x)) &= g^u(x) + 5-g^u(0)-v(0,g^v(x)) \\
    &\quad +  x\Big(4-g^u_x(2)-v(3,g^v(x))\Big)\\
    v(x,g^v(x)) &= g^v(x).
\end{aligned}
\end{equation}
A similar derivation yields the second set of constrained expressions, where the component constraints are embedded into $v$,
\begin{equation}\label{eq:ComponentExSet2}
\begin{aligned}
    u(x,g^u(x)) &= g^u(x)\\
    v(x,g^v(x),g^u(x)) &= g^v(x) + \frac{3-x}{3}\Big( 5-u(0,g^u(x))-g^v(0)\Big) \\
    &\quad + \frac{x}{3}\Big(4-u_x(2,g^u(x))-g^v(3)\Big).
\end{aligned}
\end{equation}
Notice that regardless of how $g^u(x)$ and $g^v(x)$ are chosen, Equations \eqref{eq:ComponentExSet1} and \eqref{eq:ComponentExSet2} will always satisfy the constraints.
\end{example}

Example \ref{ex:UniComponentConstraints} shows that component constraints can be placed on either dependent variable. However, notice that in the previous example, one could not put one component constraint on one dependent variable and the other component constraint on the other; doing so would result in an infinite recursion whenever trying to evaluate either \ce\ because each \ce\ would require an evaluation of the other. For example, suppose one tried to embed the first component constraint in $u$ and the second in $v$, then, the \ces\ would be,
\begin{align*}
    u(x,g^u(x)) &= g^u(x)+5-g^u(0)-v(0,g^v(x))\\
    v(x,g^v(x),g^u(x)) &= g^v(x)+4-g^v(3)-u_x(2,g^u(x)).
\end{align*}
Notice that evaluating either \ce\ requires an evaluation of the other; hence, an infinite recursion is encountered. 

The aforementioned infinite recursions can be avoided, in general, by choosing to embed as many component constraints as possible into one dependent variable, then embed as many component constraints that remain as possible into the second dependent variable, and so on, until all component constraints are accounted for. However, there may be instances when one is interested in all ways in which a set of component constraints can be embedded. Fortunately, graph theory provides a succinct method to do just that. For readers unfamiliar with the basics of graph theory, see Appendix \ref{app:GraphTheory}.

For a given set of constraints, consider a directed graph whose nodes are composed of all dependent variables that contain component constraints and whose edges connect nodes if there is a constraint between them. The direction of the edges will denote dependency in the processing order, i.e., for every edge, the target must be processed before its source is processed. Thus, to determine the order in which to create the \ces, one need only trace the graph backwards, starting at the leaf node(s) and working towards the root node(s). Infinite recursions can be avoided by checking that the resultant graph is acyclic. As mentioned in Appendix \ref{app:GraphTheory}, if a directed graph's adjacency matrix is nilpotent, then the graph is acyclic \cite{GraphTheory}. Hence, one can create all possible graphs for a given set of constraints by considering all permutations of all source/target pairs---$2^n$ possibilities where $n$ denotes the number of source/target pairs---and then reduce the set to those that do not contain infinite recursions by using the adjacency matrix. 

\begin{example}{Component constraint graphs}\label{ex:UniComponentConstraintGraphs}
Consider the follow set of component constraints,
\begin{gather*}
    u(0)+v(0)+w(0) = 5, \quad u_x(1) + v(2) = \pi \\
    u_x(3)+v_x(4) = e, \andd v(1)+w(2) = 1.
\end{gather*}
\begin{figure}[H]
    \centering
    \begin{subfigure}{0.3\linewidth}
        \centering
        \includegraphics[width=\linewidth]{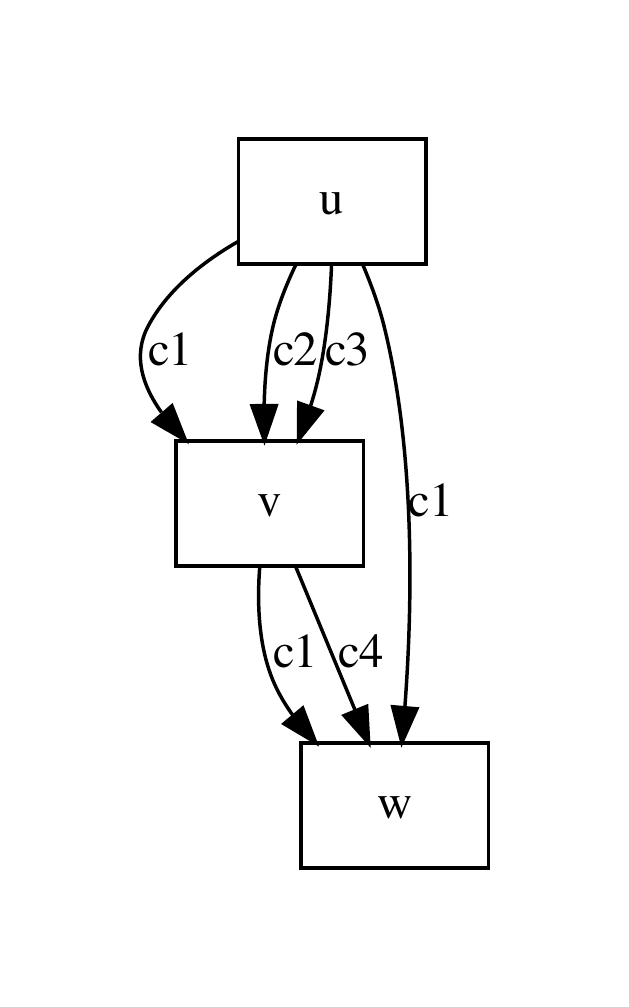}
    \end{subfigure}
    \begin{subfigure}{0.3\linewidth}
        \centering
        \includegraphics[width=\linewidth]{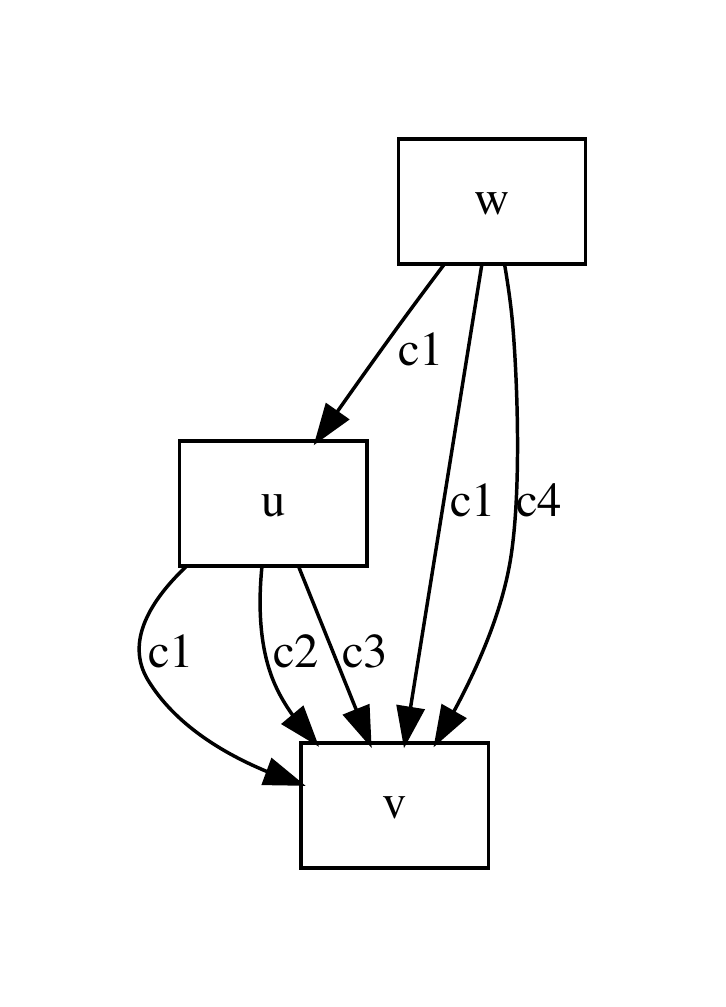}
    \end{subfigure}
    \begin{subfigure}{0.3\linewidth}
        \centering
        \includegraphics[width=\linewidth]{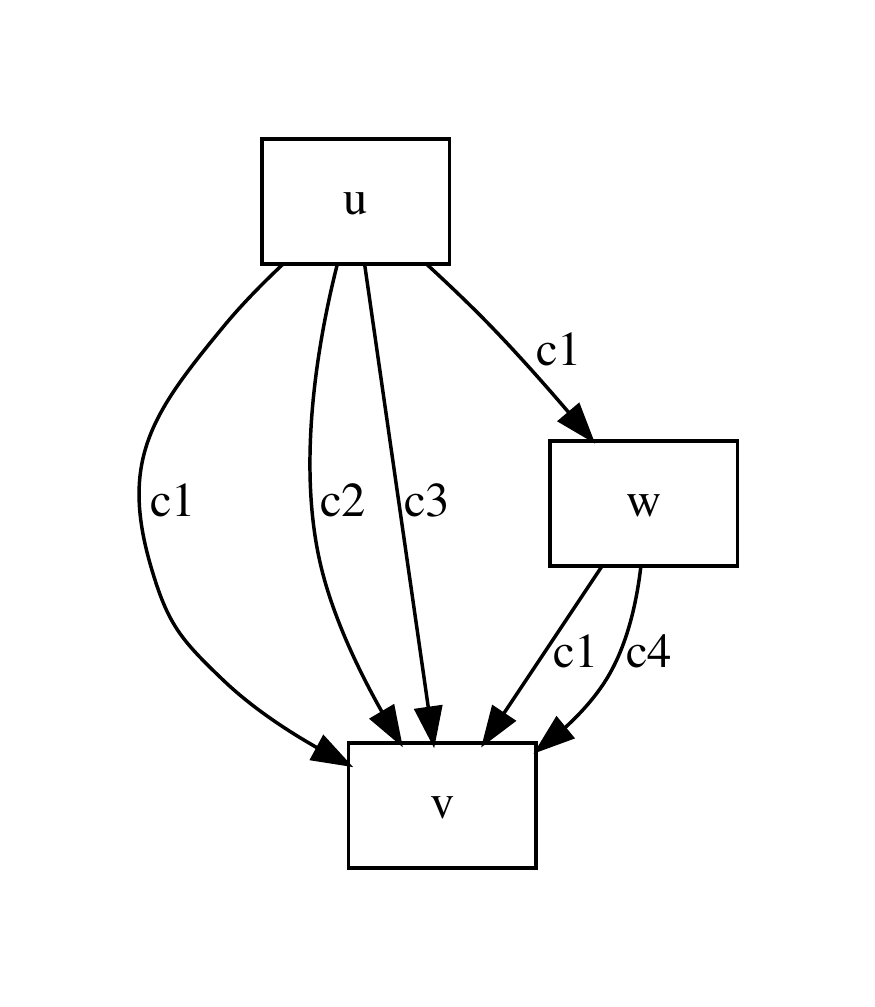}
    \end{subfigure}
\end{figure}
\begin{figure}[H]\ContinuedFloat
    \centering
    \begin{subfigure}{0.3\linewidth}
        \centering
        \includegraphics[width=\linewidth]{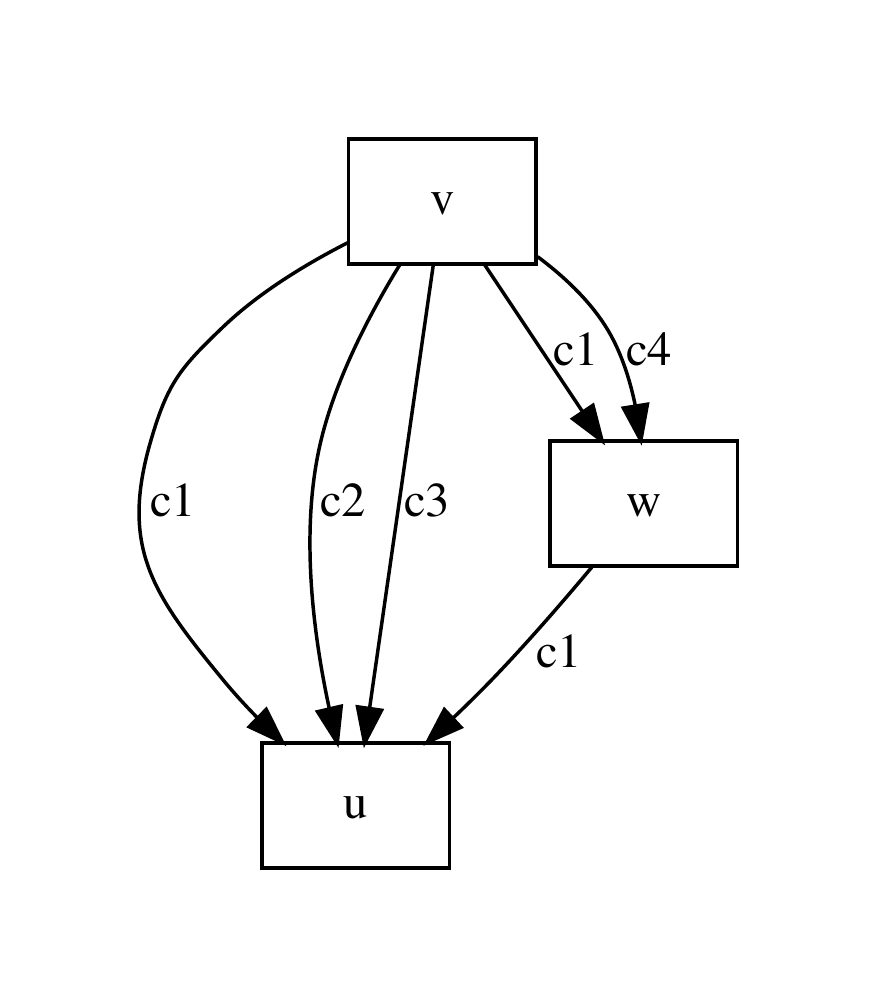}
    \end{subfigure}
    \begin{subfigure}{0.3\linewidth}
        \centering
        \includegraphics[width=\linewidth]{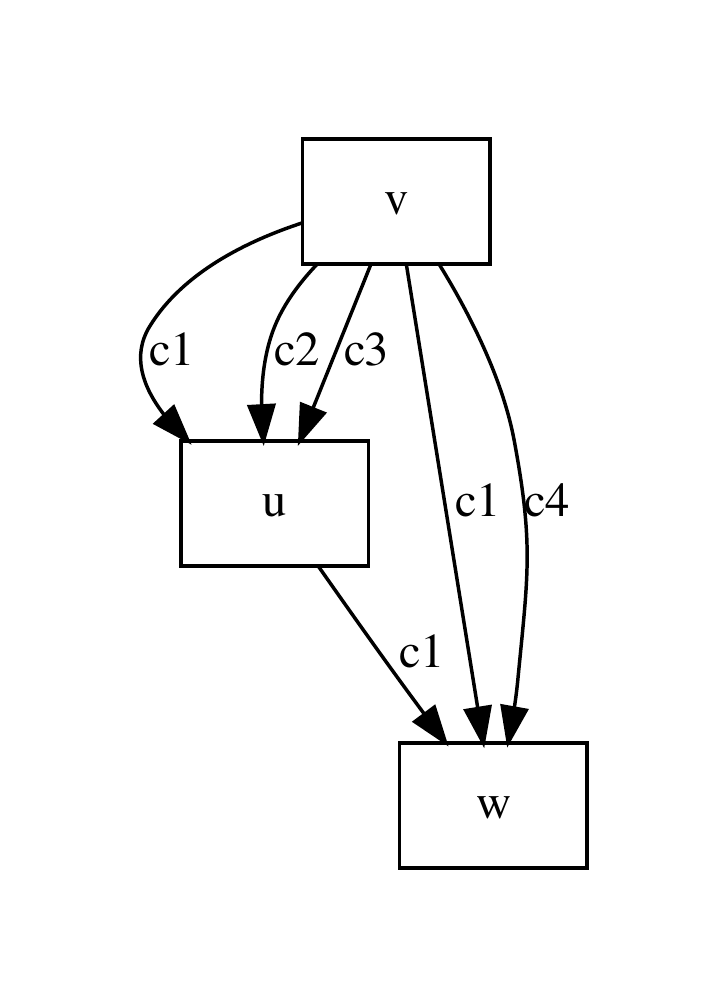}
    \end{subfigure}
    \begin{subfigure}{0.3\linewidth}
        \centering
        \includegraphics[width=\linewidth]{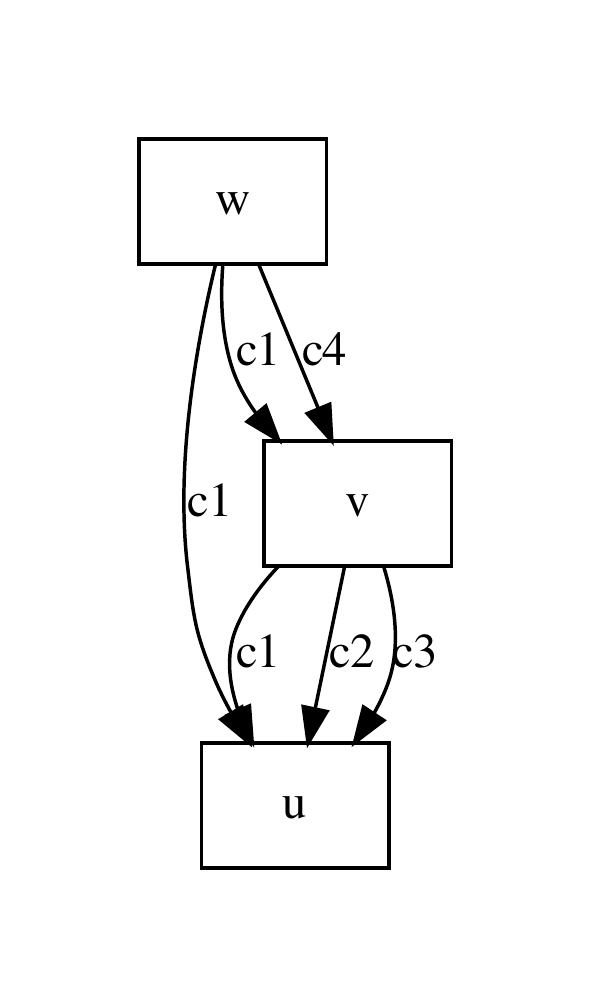}
    \end{subfigure}
    \caption{Valid component constraint graphs.}
    \label{fig:ExUniComponentGraph}
\end{figure}
Let $c_1$ denote the first component constraint, $c_2$ the second component constraint, and so on. Using the graph theory method just described, a set of directed, acyclic graphs can be created that show all possible ways in which the component constraints can be embedded; these graphs are shown in Figure \ref{fig:ExUniComponentGraph}. For example, the first graph in Figure \ref{fig:ExUniComponentGraph} is associated with embedding the first, second, and third component constraints into the \ce\ for $u$, and the fourth component constraint into the constrained expression for $v$. Moreover, based on the same graph, when constructing the constrained expressions, one must build the \ce\ for $w$ first, the \ce\ for $v$ second, and the \ce\ for $u$ last: this order was determined by traversing the graph backwards from leaf node to root node.
\end{example}

\subsection{Linear Constraints}
The term linear constraints refers to constraints that consist of linear combinations of the constraint types shown earlier. These constraints can be embedded by simply applying the techniques introduced previously.
\begin{example}{Linear constraints}
Consider the following set of constraints,
\begin{gather*}
    u(0)+u_x(0) = \pi, \quad u(1)+v(1) = 5, \\
    \int_{-1}^1 v(x) \dd{x} + v(1) = 6, \andd  v(2) = e,
\end{gather*}
and suppose the choice is made to embed the component constraint inside the $u$ \ce: although in this case, it would be equally valid to choose to embed it inside the \ce\ for $v$. Thus, for $u$, the projection functionals are,
\begin{equation*}
    \rho_1(x,g^u(x)) = \pi-g^u(0)-g^u_x(0) \andd \rho_2(x,g^u(x),g^v(x)) = 5-g^u(1)-v(1,g^v(x)),
\end{equation*}
and for $v$ they are,
\begin{equation*}
    \rho_1(x,g^v(x)) = 6 - \int_{-1}^1 g^v(\tau) \dd{\tau} - g^v(1) \andd \rho_2(x,g^v(x)) = e-g^v(2).
\end{equation*}

For $u$, the switching function equations are,
\begin{align*}
    \phi_1(0)+\frac{\partial \phi_1}{\partial x}(0) = 1, \quad & \quad \phi_1(1) = 0,\\
    \phi_2(0)+\frac{\partial \phi_2}{\partial x}(0) = 0, \quad & \quad \phi_2(1) = 1,
\end{align*}
and for $v$ they are,
\begin{align*}
    \int_{-1}^1 \phi_1(x) \dd{x}+ \phi_1(1) = 1, \quad & \quad \phi_1(2) = 0,\\
    \int_{-1}^1 \phi_2(x) \dd{x}+ \phi_2(1) = 0, \quad & \quad \phi_2(2) = 1.
\end{align*}
For $u$, let the support functions be $s_1(x) = x$ and $s_2(x) = x^2$. Then,
\begin{align*}
    \begin{bmatrix} 1 & 0 \\ 1 & 1\end{bmatrix} \begin{bmatrix} \alpha_{11} & \alpha_{12} \\ \alpha_{21} & \alpha_{22} \end{bmatrix} &= \begin{bmatrix} 1 & 0 \\ 0 & 1 \end{bmatrix}\\
    \begin{bmatrix} \alpha_{11} & \alpha_{12} \\ \alpha_{21} & \alpha_{22} \end{bmatrix} &= \begin{bmatrix} 1 & 0 \\ -1 & 1 \end{bmatrix},
\end{align*}
and the switching functions are,
\begin{equation*}
    \phi_1(x) = x-x^2 \andd \phi_2(x) = x^2.
\end{equation*}
For $v$, let the support functions be $s_1(x) = 1$ and $s_2(x) = x$. Then,
\begin{align*}
    \begin{bmatrix} 3 & 1 \\ 1 & 2\end{bmatrix} \begin{bmatrix} \alpha_{11} & \alpha_{12} \\ \alpha_{21} & \alpha_{22} \end{bmatrix} &= \begin{bmatrix} 1 & 0 \\ 0 & 1 \end{bmatrix}\\
    \begin{bmatrix} \alpha_{11} & \alpha_{12} \\ \alpha_{21} & \alpha_{22} \end{bmatrix} &= \begin{bmatrix} \frac{2}{5} &  -\frac{1}{5} \\ -\frac{1}{5} & \frac{3}{5}\end{bmatrix},
\end{align*}
and the switching functions are,
\begin{equation*}
    \phi_1(x) = \frac{2-x}{5} \andd \phi_2(x) = \frac{3x-1}{5}.
\end{equation*}

Putting the projection functionals and switching functions together yields the \ces,
\begin{align*}
    u(x,g^u(x),g^v(x)) &= g^u(x) +(x-x^2)\Big(\pi-g^u(0)-g^u_x(0)\Big) \\
    &\quad + x^2\Big(5-g^u(1)-v(1,g^v(x))\Big)\\
    v(x,g^v(x)) &= g^v(x) + \frac{2-x}{5}\Big(6 - \int_{-1}^1 g^v(\tau) \dd{\tau} - g^v(1)\Big) + \frac{3x-1}{5}\Big(e-g^v(2)\Big).
\end{align*}
As before, regardless of how $g^u(x)$ and $g^v(x)$ are chosen, these two \ces\ will always satisfy the constraints.
\end{example}

In addition to actual linear constraints, this technique can be applied to constraints that can be rewritten as linear constraints. For example, consider the nonlinear constraints shown in Appendix \ref{app:SimpleNLConstraints}; these nonlinear constraints can be rewritten as a set of linear constraints, which ultimately means they can be embedded into \ces.

\subsection{Univariate Constrained Expression Theorems}\label{subsec:UniProofs}
This section presents important theorems related to univariate TFC \ces. Theorem \ref{thrm:UniCe} shows that the \ce\ form given in the previous section satisfies the constraints regardless of how the free function is chosen. This theorem is critical, as this is the objective of \ces. 
\begin{theorem}\label{thrm:UniCe}%
The switching-projection form of the univariate \ce,
\begin{equation*}
    y(x,g(x)) = g(x) + \phi_j(x) \rho_j(x,g(x)),
\end{equation*}
satisfies the user-specified constraints for any free function.

\proof
One must show that $\C{i}[y(x,g(x))] = \kappa_i$. Apply $\C{i}$ to $y(x,g(x))$ and drop the $x$ and $g(x)$ arguments for clarity.
\begin{equation*}
    \C{i}[y] = \C{i}[g] + \C{i}[\phi_j\rho_j]
\end{equation*}
Expand $\rho_j$ and simplify,
\begin{align*}
    \C{i}[y] &= \C{i}[g] + \C{i}\Big[\phi_j(\kappa_j - \C{j}[g])\Big]\\
    \C{i}[y] &= \C{i}[g] + \C{i}[\phi_j](\kappa_j - \C{j}[g])\\
    \C{i}[y] &= \C{i}[g] + \delta_{ij}(\kappa_j - \C{j}[g])\\
    \C{i}[y] &= \C{i}[g] + \kappa_i - \C{i}[g]\\
    \C{i}[y] &= \kappa_i.
\end{align*}
Therefore, $\C{i}[y(x,g(x))] = \kappa_i$ for an any free function $g(x)$.
\end{theorem}

The natural question that arises after learning that the \ce\ satisfies the constraints for any free function is, can the \ce\ represent any function that satisfies the constraints? In other words, does the \ce\ represent the family of all possible functions that satisfy the constraints? Theorem \ref{thrm:UniGExists} shows that indeed it does.

\begin{theorem}\label{thrm:UniGExists}%
For any function satisfying the constraints, $f (x) \colon \mathbb{R}\mapsto\mathbb{R}$, there exists at least one free function, $g (x)$, such that the \ce\ $y(x,g(x)) = f(x)$. In other words, \ces\ are surjective functionals whose domain is all free functions and whose codomain is all functions that satisfy the constraints. 

\proof
As highlighted in Property \ref{prop:projZero}, the projection functionals are equal to zero whenever $g(x)$ satisfies the constraints. Thus, if $g(x)$ is a function that satisfies the constraints, then the \ce\ becomes,
\begin{align*}
    y (x, g (x)) &= g (x) + \rho_i (x, g(x)) \phi_i (x) \\
    &= g(x) + 0 \\
    &= g (x).
\end{align*} Hence, by choosing $g (x) = f (x)$, the \ce\ becomes $y (x, f (x)) = f (x)$. Therefore, for any function satisfying the constraints, $f(x)$, there exists at least one free function, $g (x) = f (x)$, such that the constrained expression is equal to the function satisfying the constraints, i.e., $y (x, f (x)) = f (x)$.
\end{theorem}

Given that the codomain of the \ce\ functional is the set of all functions satisfying the constraints, but the domain is the set of all functions, one might hypothesize that there may be multiple free function choices that produce the same output: Theorem \ref{thrm:NonUniG} shows that there are.

\begin{theorem}\label{thrm:NonUniG}%
For a given function satisfying the constraints, $f (x) \colon \mathbb{R}\mapsto\mathbb{R}$, the free function, $g (x)$, such that the \ce\ $y(x,g(x)) = f(x)$ is not unique. In other words, \ces\ are not injective functionals over the domain of all free functions and codomain of all functions that satisfy the constraints. 

\proof
Consider the free function $g (x) = f (x) + \beta_j \,  s_j (x)$ where $\beta_j$ are scalar values on $\R$ and $s_j (x)$ are the support functions used to construct the switching functions. Substituting this free function into the \ce,
\begin{equation*}
    y (x,g(x)) = g (x) + \phi_i (x) \, \rho_i (x, g (x)),
\end{equation*}
and dropping the $x$ and $g(x)$ arguments for clarity yields,
\begin{equation*}
    y = f + \beta_j s_j  + \phi_i \rho_i.
\end{equation*}

\noindent Now, expand the projection functionals and simplify,
\begin{align*}
    y &= f + \beta_j \, s_j + \phi_i \Big(\kappa_i - \C{i} [f + \beta_j \, s_j]\Big)\\
    y &= f + \beta_j s_j +  \phi_i\Big(\kappa_i - \C{i}[f] -  \C{i}[s_j]\beta_j\Big)\\
    y &= f + \beta_j s_j -  \phi_i\C{i}[s_j]\beta_j.
\end{align*}
Next, decompose the switching functions and simplify,
\begin{align*}
    y &= f + \beta_j  s_j -  \alpha_{ki}  s_k \mathbb{S}_{ij} \beta_j\\
    y &= f + \beta_j \Big(\delta_{jk} - \alpha_{ki} \mathbb{S}_{ij}\Big) s_k\\
    y &= f + \beta_j \Big(\delta_{jk} - \delta_{jk} \Big) s_k\\
    y &= f.
\end{align*}
The result obtained is independent of the $\beta_js_j(x)$ terms in the free function. Therefore, for any function, $f (x)$, satisfying the constraints, the free function, $g(x)$, that produces $f(x)$ via the \ce, i.e., $y(x,g(x)) = f(x)$, is not unique.
\end{theorem}

Notice that the non-uniqueness of $g(x)$ depends on the support functions used in the \ce, which has an immediate consequence when using \ces\ in optimization. If any terms in $g (x)$ are linearly dependent to the support functions used to construct the \ce, their contribution is negated and thus arbitrary. For some optimization techniques, it is critical that the linearly dependent terms that do not contribute to the final solution be removed; else, the optimization technique becomes impaired. For example, when solving differential equations using a linear combination of basis functions as the free function and least-squares as the optimization process \cite{M-TFC2,LDE,NDE}, the basis functions that are linearly dependent to the support functions have to be omitted from the free function to maintain full rank matrices in the least-squares.

Based on the previous results, one convenient way to think of the \ce\ is a functional that projects the free function to the set of functions that satisfy the constraints. As Theorem \ref{thrm:UniProj} shows, thinking of the \ce\ as a projection functional is a valid perspective. 

\begin{theorem}\label{thrm:UniProj}%
The \ce\ is a projection functional.

\proof
One must show that $y (x, y (x, g (x))) = y (x, g (x))$. Theorem \ref{thrm:UniCe} states that the constrained expression returns a function that satisfies the constraints. In other words, for any $g (x)$ that is defined at the constraints, $y (x, g (x))$ is a function that satisfies the constraints. From Theorem \ref{thrm:UniGExists}, if the free function used in the \ce\ satisfies the constraints, then the \ce\ returns that free function exactly. Hence, if the \ce\ functional is given itself as the free function, it will simply return itself.
\end{theorem}

The previous proofs coupled with the functional-related definitions given earlier provide a more rigorous definition for the univariate \ce: the univariate \ce\ is a surjective, projection functional whose domain is the set of all free functions and whose codomain is the set of all functions that satisfy the constraints. It is surjective because it spans the set of all functions that satisfy the constraints, its codomain, based on Theorem \ref{thrm:UniGExists}, but it is not injective because Theorem \ref{thrm:NonUniG} shows that functions in the codomain are the image of more than one function in the domain; \ces\ are thus not bijective either because they are not injective. Moreover, the \ce\ is a projection functional as shown in Theorem \ref{thrm:UniProj}. 

\section{Multivariate Theory}
This section utilizes the univariate theory introduced in the previous section to extend TFC to the multivariate case. As such, one should ensure they have a firm grasp of the concepts introduced in the univariate section before moving on. The section begins by introducing the recursive method: a method for generating multivariate \ces\ by using the univariate \ce\ for one independent variable as the free function in the univariate \ce\ for a different independent variable. Afterward, the mathematical theorems presented for univariate \ces\ are also extended to the multivariate case. Finally, a compact tensor form of the multivariate \ce\ is presented.

\subsection{Recursive Method}
Oftentimes, the constraints of a problem do not include integral constraints. In these cases, one independent variable's constraints will not interfere with another independent variable's constraints. Consequently, using the univariate \ce\ for one independent variable as the free function in the univariate \ce\ of another independent variable produces a function that satisfies both independent variables' constraints. To prove this, one must first understand how the constraint operator of one independent variable affects the \ce\ of another independent variable; the following discussion and properties will help achieve this understanding. A pre-superscript will be used to distinguish the operators, functions, and functionals of one independent variable from another. For example, $\pC{k}{j}$ denotes the constraint operator for the $j$-th constraint of the $k$-th independent variable.

\begin{property}\label{prop:ConstraintOperatorOnOtherVariables}
    For non-integral constraints, the constraint operator for the $k$-th independent variable operating on a product of functions wherein one function is not a function of the $k$-th independent variable and the other is leads to,
    \begin{align*}
        \pC{k}{j}&[f(x_1,\dots,x_{k-1},x_{k+1},\dots,x_n) h(x_1,\dots,x_k,\dots,x_n)] = \\
        &f(x_1,\dots,x_{k-1},x_{k+1},\dots,x_n) \pC{k}{j}[h(x_1,\dots,x_k,\dots,x_n)],
    \end{align*}
    where $f$ is not a function of the $k$-th independent variable, $x_k$, but $h$ is.
\end{property}

\noindent Property \ref{prop:ConstraintOperatorOnOtherVariables} holds for non-integral constraints because $\pC{k}{j}$ operates on the $k$-th independent variable only, and $f$ is not a function of the $k$-th independent variable, i.e., it is effectively a constant. In particular, this property is useful in multivariate expressions, which oftentimes contain such products. For example, $\pC{k}{i}[\p{k}{\phi}_j\p{k}{\kappa}_j] = \pC{k}{i}[\p{k}{\phi}_j]\p{k}{\kappa}_i$.

\begin{property}\label{prop:ConstraintOperatorOnKappas}
A set of non-integral constraints is consistent if and only if $\pC{k}{j}[\p{l}{\kappa}_i] = \pC{l}{i}[\p{k}{\kappa}_j]$.
\end{property}
\noindent Property \ref{prop:ConstraintOperatorOnKappas} is easiest to understand via an example of inconsistent constraints:
\begin{equation*}
    z(x,0) = 5 \andd z(0,y) = 4.
\end{equation*}
Clearly, these constraints cannot simultaneously be satisfied at the intersection point $z(0,0)$. 

In addition, for multivariate constraints, the free function must be locally $C^m$ in the neighborhood of the geometric intersection of constraints, where $m$ is the sum of the orders of derivatives of the intersecting constraints. This restriction on the free function is necessary for generating \ces\ using recursive applications of univariate expressions, as it ensures that Clairaut's theorem holds for the free function, and thus, $\pC{l}{i}\Big[\pC{k}{j}[g]\Big] = \pC{k}{j}\Big[\pC{l}{i}[g]\Big]$. 
\begin{definition}\label{def:multiFree}
The multivariate free function is any function $g (x)\colon \mathbb{R}^n\mapsto\mathbb{R}$ such that $\pC{i_j}{k}[g]$ is defined and $\pC{i_j}{m} \Big[ \cdots \big[\pC{i_k}{n}[g]\big]\cdots\Big]$ is defined, where the latter consists of at most one constraint operator from each dimension; the latter must be freely permutable, e.g., $\pC{i_j}{m} \Big[ \cdots \big[\pC{i_k}{n}[g]\big]\cdots\Big] = \pC{i_k}{n} \Big[ \cdots \big[\pC{i_j}{m}[g]\big]\cdots\Big]$, for any non-integral constraints.
\end{definition}
\noindent These properties and restriction on the free function are utilized in Theorem \ref{thrm:ConsructingNdCes} to show that the recursive method produces a valid multivariate \ce\ for non-integral constraints.

\begin{theorem}\label{thrm:ConsructingNdCes}%
For non-integral constraints, a valid multivariate \ce\ can be constructed by recursively applying the univariate \ce\ from one independent variable as the free function in the \ce\ for another independent variable. In this recursion, all univariate \ces\ must be used once and only once, and the first univariate \ce\ is built using a regular free function.

\proof
First, show that $\p{k}{u}(\B{x},\p{l}{u}(\B{x},g(\B{x})))$ is a valid bivariate \ce\ that satisfies both sets of constraints, where $\B{x}$ represents a vector of the independent variables, i.e., $\B{x} = \{ x_1, x_2, \cdots, x_n \}$. Then, apply it $n$ times recursively to produce an $n$-dimensional, multivariate \ce\ that satisfies the constraints on all $n$ dimensions. Consider two univariate \ces:
\begin{align*}
    \p{k}{u}(\B{x},g(\B{x})) &= g(\B{x}) + \p{k}{\phi}_j(x_k)\p{k}{\rho}_j(\B{x},g(\B{x})),\\
    \p{l}{u}(\B{x},g(\B{x})) &= g(\B{x}) + \p{l}{\phi}_i(x_l)\p{l}{\rho}_i(\B{x},g(\B{x})).
\end{align*}
Substitute the univariate \ce\ for the $l$-th independent variable as the free function in the univariate \ce\ for the $k$-th independent variable,
\begin{equation*}
    \p{k}{u}(\B{x},\p{l}{u}(\B{x},g(\B{x}))) = \p{l}{u}(\B{x},g(\B{x})) + \p{k}{\phi}_j(x_k)\p{k}{\rho}_j(\B{x},\p{l}{u}(\B{x},g(\B{x}))).
\end{equation*}
Clearly, from Theorem \ref{thrm:UniCe}, which shows that a univariate \ce\ satisfies the constraints for any free function that is defined at the constraints, the constraints of the $k$-th independent variable must be satisfied, as $\p{l}{u}(\B{x},g(\B{x}))$ is a valid free function. Next, expand $\p{l}{u}(\B{x},g(\B{x}))$ and $\p{k}{\rho}_j(\B{x},g(\B{x}))$ and drop the $\B{x}$ and $g(\B{x})$ arguments for clarity.
\begin{align*}
    \p{k}{u} &= g + \p{l}{\phi}_i\p{l}{\rho}_i +  \p{k}{\phi}_j\Big(\p{k}{\kappa}_j-\pC{k}{j}[g] - \pC{k}{j}\Big[\p{l}{\phi}_i\p{l}{\rho}_i\Big]\Big) \\
    &= g + \p{l}{\phi}_i\Big(\p{l}{\kappa}_i -\pC{l}{i}[g]\Big) +  \p{k}{\phi}_j\bigg(\p{k}{\kappa}_j-\pC{k}{j}[g] \\
    &\quad\quad - \p{l}{\phi}_i\Big(\pC{k}{j}[\p{l}{\kappa}_i] - \pC{k}{j}\Big[\pC{l}{i}[g] \Big]\Big)\bigg),
\end{align*}
where Property \ref{prop:ConstraintOperatorOnOtherVariables} has been used to simplify the expression. Now, evaluate $\p{k}{u}$ at the $m$-th constraint for the $l$-th independent variable,
\begin{align*}
    \pC{l}{m}[\p{k}{u}] &= \pC{l}{m}[g] + \delta_{mi}\Big(\p{l}{\kappa}_i -\pC{l}{i}[g]\Big)+\p{k}{\phi}_j\bigg(\pC{l}{m}[\p{k}{\kappa}_j]-\pC{l}{m}\Big[\pC{k}{j}[g]\Big] \\
    &\quad-\delta_{mi}\Big(\pC{k}{j}[\p{l}{\kappa}_i] - \pC{k}{j}\Big[\pC{l}{i}[g] \Big]\Big)\bigg)\\
    &= \p{l}{\kappa}_m \\
    &\quad + \p{k}{\phi}_j\Big(\pC{l}{m}[\p{k}{\kappa}_j] - \pC{k}{j}[\p{l}{\kappa}_m] - \pC{l}{m}\Big[\pC{k}{j}[g]\Big]
    + \pC{k}{j}\Big[\pC{l}{m}[g] \Big]\Big)\\
    &= \p{l}{\kappa}_m.
\end{align*}
Therefore, $\pC{l}{m}[\p{k}{u}] = \p{l}{\kappa}_m$ as required, and the expression $\p{k}{u}$ satisfies both sets of original univariate constraints. 
\end{theorem}

Example \ref{ex:MultiExSimple} demonstrates Theorem \ref{thrm:ConsructingNdCes}.

\begin{example}{Multivariate non-integral constraints}\label{ex:MultiExSimple}
Consider the following set of constraints,
\begin{gather*}
   u(0,y) = y^2\sin(\pi y), \quad u(1,y)+u(2,y) = y\sin(\pi y),\\
   u_y(x,0) = 0, \andd u(x,0) = u(x,1).
\end{gather*}
The univariate \ces\ for the constraints on $x$ and $y$ are,
\begin{align*}
    \p{1}{u}(x,y,g(x,y)) &= g(x,y) + \frac{3-2x}{3}\Big(y^2\sin(\pi y)-g(0,y)\Big)\\
    &\quad\quad+\frac{x}{3}\Big(y\sin(\pi y)-g(2,y)-g(1,y)\Big)\\
    \p{2}{u}(x,y,g(x,y)) &= g(x,y)-(y-y^2)g_y(x,0)-y^2\Big(g(x,1)-g(x,0)\Big).
\end{align*}

\noindent Then, $\p{1}{u}$ is used as the free function in $\p{2}{u}$,
\begin{align*}
    \p{2}{u}(x,y,\p{1}{u}(x,y,g(x,y))) &= \p{1}{u}(x,y,g(x,y))-(y-y^2)\p{1}{u}_y(x,0,g(x,y))\\
    &\quad -y^2\Big(\p{1}{u}(x,1,g(x,y))-\p{1}{u}(x,0,g(x,y))\Big).
\end{align*}

\noindent Substituting in $\p{1}{u}$ and simplifying yields,
\begin{equation}\label{eq:ceMultEx2}
\begin{aligned}
    u(x,y,&g(x,y)) = g(x,y)+\left(y-y^2\right) \Big(\frac{3-2 x}{3} g_y(0,0)-\frac{x}{3} \left(-g_y(1,0)-g_y(2,0)\right)\\
    &-g_y(x,0)\Big)-y^2 \Big(\frac{3-2 x}{3} g(0,0)-\frac{3-2 x}{3} g(0,1)-\frac{x}{3}  (-g(1,0)-g(2,0))\\
    &+\frac{x}{3}  (-g(1,1)-g(2,1)) -g(x,0)+g(x,1)\Big)+\frac{3-2 x}{3} \left(y^2 \sin (\pi  y)-g(0,y)\right)\\
    &+\frac{x}{3}  \Big(-g(1,y) -g(2,y)+y \sin (\pi  y)\Big).
\end{aligned}
\end{equation}

Note that substituting $\p{2}{u}$ as the free function in $\p{1}{u}$, after simplifying, yields the same result given in Equation \eqref{eq:ceMultEx2}. Equation \eqref{eq:ceMultEx2} satisfies the constraints for any $g(x,y)$ satisfying Definition \ref{def:multiFree}. Figure \ref{fig:MultiEx2} shows the constrained expression when $g (x, y) = x^2 \cos y + \sin (2 x)$, where the blue line signifies the constraint on $u (0, y)$, the black lines signify the derivative constraint on $u_y (x, 0)$, and the magenta lines signify the relative constraint $u (x, 0) = u (x, 1)$. The linear constraint $u (1, y) + u (2, y) = y \sin(\pi y)$ is not easily visualized but is nonetheless satisfied by the \ce.
\begin{figure}[H]
    \centering
    \includegraphics[width=.85\linewidth]{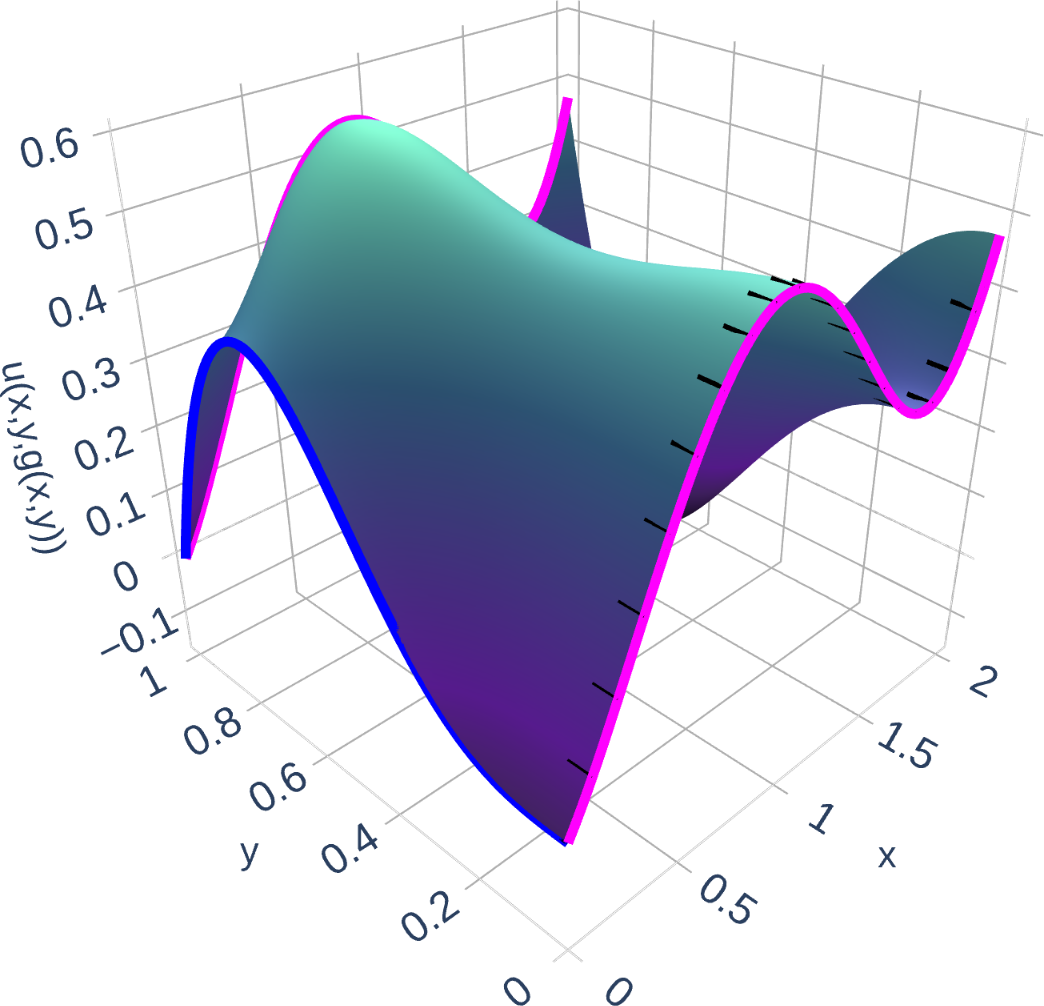}
    \caption{Constrained expression evaluated using $g (x, y) = x^2 \cos y + \sin (2 x)$. The blue line signifies the constraint on $u (0, y)$, the black lines signify the derivative constraint on $u_y (x, 0)$, and the magenta lines signify the relative constraint $u (x, 0) = u (x, 1)$. The linear constraint $u (1, y) + u (2, y) = y \sin(\pi y)$ is not easily visualized but is nonetheless satisfied by the constrained expression.}
    \label{fig:MultiEx2}
\end{figure}

\end{example}

\subsubsection{Integral Constraints}\label{sec:RecursiveIntegralConstraints}
Theorem \ref{thrm:ConsructingNdCes} is proven using Property \ref{prop:ConstraintOperatorOnOtherVariables} applied to expressions such as $\pC{l}{i}\Big[\p{k}{\phi}_j\pC{k}{j}[g]\Big] = \p{k}{\phi}_j \pC{l}{i}\Big[\pC{k}{j}[g]\Big]$. These expressions are true so long as the constraint operator of the $l$-th independent variable does not affect functions that do not contain $l$. This is true for all constraint types introduced thus far, except integral constraints. Integral constraints may have constraint operators like,
\begin{equation*}
    \pC{l}{i}[f(x_1,\dots,x_l,\dots,x_n)] = \int_a^b f(x_1,\dots,c,\dots,x_n) \dd x_k,
\end{equation*}
where $a,b,c\in\R$, which affect functions of the $k$-th independent variable, even though the constraint operator is for the $l$-th independent variable. Although this type of constraint is rare in PDEs, in the interest of introducing a general function interpolation technique, the next section presents a method to embed integral constraints into multivariate \ces. 

The interference between independent variables introduced by integral constraints can be avoided by modifying the switching functions and processing order of the univariate \ces. 

\begin{theorem}\label{thrm:IntegralConstraintModification}%
Processing the independent variables that appear as integration variables in integral constraints after the independent variables associated with the integral constraints and modifying the switching functions of all constraints of the variables of integration such that they yield zero when operated on by the constraint operators of said integral constraints is a valid method for embedding integral constraints into multivariate \ces.

\proof
Let the constraints of the $l$-th independent variable contain an integral constraint whose variable of integration is the $k$-th independent variable. Moreover, as per the theorem statement, let the switching functions of the $k$-th independent variable be defined such that $\pC{l}{m}[\p{k}{\phi}_j] = 0$ if the $m$-th constraint is the integral constraint. Now, following the recursive method, the bivariate constrained expression for the $k$-th and $l$-th independent variables is,
\begin{equation*}
    \p{k}{u}(\B{x},\p{l}{u}(\B{x},g(\B{x}))) = \p{l}{u}(\B{x},g(\B{x})) + \p{k}{\phi}_j(x_k)\p{k}{\rho}_j(\B{x},\p{l}{u}(\B{x},g(\B{x}))).
\end{equation*}
From Theorem \ref{thrm:ConsructingNdCes}, the constraints of the $k$-th independent variable and non-integral constraints of the $l$-th independent variable are satisfied. As in Theorem \ref{thrm:ConsructingNdCes}, expand the expression for $\p{k}{u}$ and drop the $\B{x}$ and $g(\B{x})$ arguments for clarity,
\begin{align*}
    \p{k}{u} =\ &g + \p{l}{\phi}_i\Big(\p{l}{\kappa}_i -\pC{l}{i}[g]\Big) \\
    &+\p{k}{\phi}_j\bigg(\p{k}{\kappa}_j-\pC{k}{j}[g] - \p{l}{\phi}_i\Big(\pC{k}{j}[\p{l}{\kappa}_i] - \pC{k}{j}\Big[\pC{l}{i}[g] \Big]\Big)\bigg).
\end{align*}
Apply the constraint operator for the integral constraint,
\begin{align*}
    \pC{l}{m}[\p{k}{u}] &= \pC{l}{m}[g] + \delta_{mi}\Big(\p{l}{\kappa}_i -\pC{l}{i}[g]\Big)\\
    &\quad+\pC{l}{m}\bigg[\p{k}{\phi}_j\underbrace{\bigg(\p{k}{\kappa}_j-\pC{k}{j}[g] - \p{l}{\phi}_i\Big(\pC{k}{j}[\p{l}{\kappa}_i] - \pC{k}{j}\Big[\pC{l}{i}[g] \Big]\Big)\bigg)}_{\text{Not a function of the $k$-th independent variable}}\bigg].
\end{align*}
As noted in the above expression, the function highlighted by the underbrace is not a function of the $k$-th independent variable; hence, the integration portion of the $\pC{l}{m}$ constraint operator only acts on $\p{k}{\phi}_j$. Moreover, recall that this method redefines the switching functions such that $\pC{l}{m}[\p{k}{\phi}_j] = 0$. Thus,
\begin{align*}
    \pC{l}{m}[\p{k}{u}] &= \pC{l}{m}[g] + \delta_{mi}\Big(\p{l}{\kappa}_i -\pC{l}{i}[g]\Big)\\
    &= \pC{l}{m}[g] + \p{l}{\kappa}_m - \pC{l}{m}[g]\\
    &= \kappa_m,
\end{align*}
as desired. Therefore, all constraints, integral and non-integral, on both the $k$-th and $l$-th independent variables are satisfied. Applying this proof recursively shows that this is a valid method for constructing multivariate \ces\ that contain integral constraints. 
\end{theorem}

Example \ref{ex:MultiInt} provides a concrete demonstration of Theorem \ref{thrm:IntegralConstraintModification}.

\begin{example}{Multivariate integral constraints}\label{ex:MultiInt}
Consider the following set of constraints,
\begin{equation*}
    u(x,0) = 2 u_y(x,1), \quad u(x,2) = \sin(x), \andd \int_{-1}^1 u(2,y) \dd y = 5.
\end{equation*}
Based on the previous discussion, the $x$ independent variable will be processed first, because it has an integral constraint with integration variable $y$, and the switching functions for the constraints in $y$ must be created such that they are equal to zero when evaluated with the constraint operator for the integral constraint. That is, the equations for the $y$ switching functions are,
\begin{align*}
    &\p{2}{\phi}_1(0) - 2 \frac{\partial \p{2}{\phi}_1}{\partial y}(1) = 1, &&\p{2}{\phi}_1(2) = 0, &&\int_{-1}^1 \p{2}{\phi}_1(y) \dd y = 0,\\
    &\p{2}{\phi}_2(0) - 2 \frac{\partial \p{2}{\phi}_2}{\partial y}(1) = 0, &&\p{2}{\phi}_2(2) = 1, &&\int_{-1}^1 \p{2}{\phi}_2(y) \dd y = 0.
\end{align*}
Even though there are only two constraints in the $y$-dimension, each switching function must satisfy three sets of equations. Therefore, each switching function should be a linear combination of three linearly independent support functions with unknown coefficients, $\p{2}{\phi}_i(y) = \alpha_{ij}s_j(y)$ where $i\in\{1,2\}$ and $j\in\{1,2,3\}$. As before, the equations can be written in a compact matrix form and solved via matrix inversion.
\begin{align*}
    \begin{bmatrix} 1 & -2 & -4 \\ 1 & 2 & 4\\ 2 & 0 & \frac{2}{3} \end{bmatrix} \begin{bmatrix} \alpha_{11} & \alpha_{12} \\ \alpha_{21} & \alpha_{22} \\ \alpha_{31} & \alpha_{32} \end{bmatrix} &= \begin{bmatrix} 1 & 0 \\ 0 & 1 \\ 0 & 0 \end{bmatrix} \\
    \begin{bmatrix} \alpha_{11} & \alpha_{12} \\ \alpha_{21} & \alpha_{22} \\ \alpha_{31} & \alpha_{32} \end{bmatrix} &= \begin{bmatrix} \frac{1}{2} & \frac{1}{2} \\ \frac{11}{4} & \frac{13}{4} \\ -\frac{3}{2} & -\frac{3}{2} \end{bmatrix}
\end{align*}
Next, the univariate constrained expressions for each of the independent variables can be written as,
\begin{align*}
    \p{1}{u}(x,y,g(x,y)) &= g(x,y) + \frac{1}{2} \bigg(5 - \int_{-1}^1 g(2,\tau) \dd \tau\bigg),\\
    \p{2}{u}(x,y,g(x,y)) &= g(x,y) +\frac{2+11y-6y^2}{4}\Big(2g_y(x,1)-g(x,0)\Big)\\
    &\quad+\frac{2+13y-6y^2}{4}\Big(\sin(x)-g(x,2)\Big),
\end{align*}
and following the method outlined earlier, the full multivariate \ce\ can be written as,
\begin{equation}\label{eq:MultiIntEx}
\begin{aligned}
    u(x,y,&g(x,y)) = \p{2}{u}(x,y,\p{1}{u}(x,y,g(x,y))) \\
    &= g(x,y)+\frac{1}{4} \Big[2(2-y)\left((6 y+1) g_y(x,1)+3 y \left(\int_{-1}^1 g(2,\tau ) \, d\tau -5\right)\right)\\
    &\quad+(y-2) (6 y+1) g(x,0)+\big(y (6 y-13)-2\big) g(x,2)\\\
    &\quad+(y (13-6 y)+2) \sin (x)\Big].
\end{aligned}
\end{equation}
As expected, Equation \eqref{eq:MultiIntEx} satisfies the constraints for any valid free function $g(x,y)$.
\end{example}

It is important to note that this method cannot embed sets of integral constraints whose independent variables refer to one another, such as,
\begin{equation*}
    \int_0^1u(x,0) \dd{x} = 1 \andd \int_0^1u(0,y) \dd{y} = 1.
\end{equation*}
The reason is that the first integral constraint requires that the $y$ independent variable be processed before $x$, but the second integral constraint requires that the $x$ independent variable be processed before $y$: obviously, these two requirements cannot be satisfied simultaneously. Therefore, this method cannot embed such constraints. 

\subsubsection{Component Constraints}
As in the univariate case, one must choose which dependent variable a component constraint will be assigned to. Again, graph theory can be used in the same manner as before to determine all possible ways in which a set of component constraints can be embedded, see Example \ref{ex:UniComponentConstraintGraphs}. However, when moving to the multivariate case, one must be cautious of the intersections between component constraints and other constraints: Example \ref{ex:MultiComponent} highlights this nuance.

\begin{example}{Multivariate component constraints}\label{ex:MultiComponent}
Consider the following set of constraints,
\begin{equation*}
    u(x,0) = 5 \andd u(0,y)+v(0,y) = 3.
\end{equation*}
If one chose to embed the component constraint into $u$ and process the constraints on $x$ first, then the \ces\ would be,
\begin{align*}
    u(x,y,g^u(x,y),g^v(x,y)) &= g^u(x,y) + 5 - g^u(x,0) - g^u(0,y)-v(0,y,g^v(x,y))\\
    &\quad+v(0,0,g^v(x,y))+g^u(0,0)\\
    v(x,y,g^v(x,y)) &= g^v(x,y).
\end{align*}
Clearly, the \ce\ for $u$ does not satisfy the two constraints for any valid free function, e.g., choosing $g^u(x,y) = 3$ and $g^v(x,y) = 2$ yields $u(0,y,3) + v(0,y,2) = 7 \neq 3$; the reason stems from the intersection between the two constraints.\footnote{Although this example only shows the intersection issue when processing the constraints on $x$ first, the same issue arises even if the constraints on $y$ are processed first.} At the intersection, $u$ must be equal to $5$, but simultaneously be equal to $3-v(0,y,g^v(0,y))$. If $g^v(x,y)$ was chosen in such a way that $g^v(0,0) = -2$, then all constraints would be satisfied, but of course, the objective of \ces\ is to provide a functional that satisfies the constraints wherein the free function can be chosen without restriction. In other words, since at the intersection of the constraints, $(x,y) = (0,0)$, $u=5$ as specified by the first constraint, the only way to simultaneously satisfy the component constraint is to change $v$; hence, the component constraint must be placed on $v$. Doing so results in the \ces,
\begin{equation*}
\begin{aligned}    
    u(x,y,g^u(x,y)) &= g^u(x,y) + 5 - g^u(x,0)\\
    v(x,y,g^v(x,y),g^u(x,y)) &= g^v(x,y) + 3 - g^v(0,y) - u(0,y,g^u(x,y)).
\end{aligned}
\end{equation*}
which satisfy the constraints for any valid $g^u(x,y)$ and $g^v(x,y)$.
\end{example}

As demonstrated in Example \ref{ex:MultiComponent}, component constraints must only be placed on dependent variables that do not have other constraints that intersect with the component constraint. Of course, if each variable in the component constraint has a constraint at the intersection point, then the component constraint can be placed on either variable. For example, if the constraints in Example \ref{ex:MultiComponent} were,
\begin{equation*}
    u(x,0) = 5, \quad v(x,0) = -2, \andd u(0,y)+v(0,y) = 3,
\end{equation*}
then the component constraint could have been embedded into either $u$ or $v$. Therefore, while the graph theory introduced in the univariate section can be used to determine component constraint embeddings that avoid infinite recursions when evaluating the \ces, it is up to the user to further reduce this set of graphs to those that avoid the intersection issues described above.

After considering the previous restrictions on component constraints, one may contrive a set of equations where it is impossible to meet the aforementioned conditions. For example, consider the following constraints,
\begin{equation*}
    u(0,y) = 5, \quad v(1,y) = 2, \andd u(x,0)+v(x,0) = 3.
\end{equation*}
The component constraint cannot be placed on $u$ because of the intersection at $(0,0)$, but it also cannot be placed on $v$ because of the intersection at $(1,0)$. The only option here is to split the domain along the $x$-axis for some $x_{\text{split}}\in(0,1)$. Then, in the left sub-domain, $x < x_{\text{split}}$, the component constraint will be embedded into $v$, and in the right sub-domain, $x > x_{\text{split}}$, the component constraint will be embedded into $u$. At the intersection of these two sub-domains, $x = x_{\text{split}}$, one can enforce $C^n$ continuity---$n$ is chosen by the user or dictated by the problem---by adding constraints at the intersection that can ultimately be embedded into the \ces. For more information and an example on splitting the domain, see Appendix \ref{app:DivideDomain}.

\subsubsection{Linear Constraints}
Multivariate linear constraints consist of linear combinations of the previously introduced constraint types. Thus, one must be conscientious of the nuances of both integral and component constraints if they appear in the linear constraints. The following step-by-step procedure can be used to construct multivariate \ces:
\begin{enumerate}
    \item Generate the directed, acyclic graphs that show all the valid ways that the component constraints can be embedded. Of these graphs, either choose one that satisfies the intersection restriction discussed in the previous section or choose one and split the domain as needed. The chosen graph will dictate the order in which the dependent variables' \ces\ are created.
    \item For each dependent variable, choose the order in which the univariate \ces\ will be processed. This order is dictated in part by the presence of integral constraints. 
    \item Build the multivariate \ces. 
\end{enumerate}
This step-by-step procedure can be used for any embeddable set of constraints. Of course, steps in this procedure may be omitted depending on the types of constraints present. For example, if the set does not contain any component constraints, then there is no need to perform step 1, and the dependent variables' \ces\ can be created in any order. Example \ref{ex:MultiLinearConstraints} demonstrates this process.

\begin{example}{Multivariate linear constraints}\mbox{}\label{ex:MultiLinearConstraints}
Consider the following set of constraints,
\begin{gather*}
    u(0,y) = \cos(\pi y), \quad \int_{-1}^2 u(1,y) \dd{y} = e, \quad u(x,1)-u(x,2) = -2,\\
    u(x,0)+v(x,0) = 5, \andd v(0,y) = 5-\cos(\pi y).
\end{gather*}
As per the procedure outlined above, step 1 is to generate the directed graphs that dictate the valid ways in which the component constraint can be embedded. In this example, those graphs are trivial: the single component constraint can be embedded into either $u$ or $v$ without producing a set of \ces\ that require an infinite recursion upon evaluation. However, $u$ contains an integral constraint along the $x=1$ line, and there is no corresponding constraint at $x=1$ in $v$. Hence, the component constraint is embedded into the $v$ \ce, and the intersection issue is avoided. It follows that the $u$ \ce\ will be created before the $v$ \ce.

Next, the processing order for the independent variables must be decided. For $v$, the order does not matter as no integral constraints are present. In contrast, for $u$, the constraints on $x$ must be processed before those on $y$, as the constraints on $x$ contain an integral constraint wherein $y$ is an integration variable. 

Now the \ces\ can be created. First the $u$ \ce\ is created. The univariate \ces\ for the constraints on $x$ and $y$ are,
\begin{align*}
    \p{1}{u}(x,y,g^u(x,y)) &= g^u(x,y)) +(1-x)\Big(\cos(\pi y)-g^u(0,y)\Big) + \frac{x}{3}\Big(e-\int_{-1}^2 g^u(1,\tau) \dd{\tau}\Big)\\
    \p{2}{u}(x,y,g^u(x,y)) &= g^u(x,y) + \frac{1-2y}{2}\Big(g^u(x,2)-g^u(x,1)-2\Big),
\end{align*}
where monomials have been used as the support functions. These two \ces\ are used in the processing order defined above to produce the multivariate \ce\ for $u$,
\begin{align*}
    u(x,y,&g^u(x,y)) = g^u(x,y)+\frac{1}{3} x \left(e-\int_{-1}^2 g^u(1,\tau ) \dd{\tau} \right)+\frac{1-2 y}{2}\big((1-x) (1+g^u(0,1))\\
    &+(1-x) (1-g^u(0,2))-g^u(x,1)+g^u(x,2)-2\Big)+(1-x) (\cos (\pi  y)-g^u(0,y)).
\end{align*}

For $v$, the univariate \ces\ for the constraints on $x$ and $y$ are,
\begin{align*}
    \p{1}{v}(x,y,g^v(x,y)) &= g^v(x,y) + 5 - \cos(\pi y) - g^v(0,y)\\
    \p{2}{v}(x,y,g^v(x,y),g^u(x,y)) &= g^v(x,y) + 5 - g^v(x,0) - u(x,0,g^u(x,y)),
\end{align*}
where monomials have been used as the support functions. The full multivariate \ce\ for $v$ is,
\begin{align*}
    v(x,y,&g^v(x,y),g^u(x,y)) = g^v(x,y)-g^v(x,0)-g^v(0,y)+g^v(0,0)-u(x,0,g^u(x,y))\\
    &+u(0,0,g^u(x,y))-\cos (\pi  y)+5.
\end{align*}

The constrained expressions satisfy the constraints for any valid $g^u(x,y)$ and $g^v(x,y)$. Figure \ref{fig:MultiExLinear} shows the constrained expressions for $u$ and $v$ evaluated using $g^u(x,y) = xy+\sin(x)+y^2$ and $g^v(x,y) = x^2y\cos(y)e^x$. The value constraints that are easily visualized, the constraints on $u(x,0)$ and $v(x,0)$, are shown as black lines. The remaining constraints are not easily visualized and are therefore not shown, but they are satisfied nonetheless.
\begin{figure}[H]
    \centering
    \begin{subfigure}{0.49\linewidth}
        \centering
        \includegraphics[width=\linewidth]{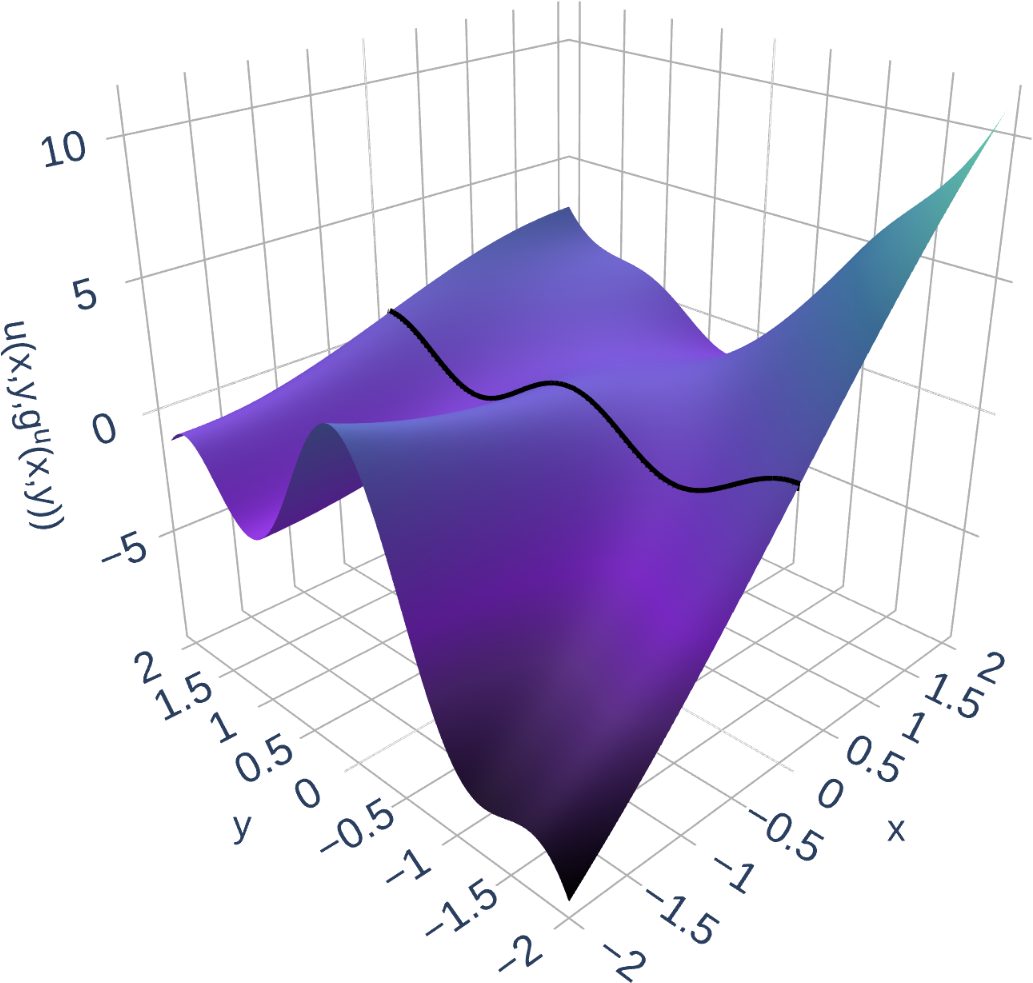}
        \caption{$u(x,y,g^u(x,y))$}
    \end{subfigure}
    \begin{subfigure}{0.49\linewidth}
        \centering
        \includegraphics[width=\linewidth]{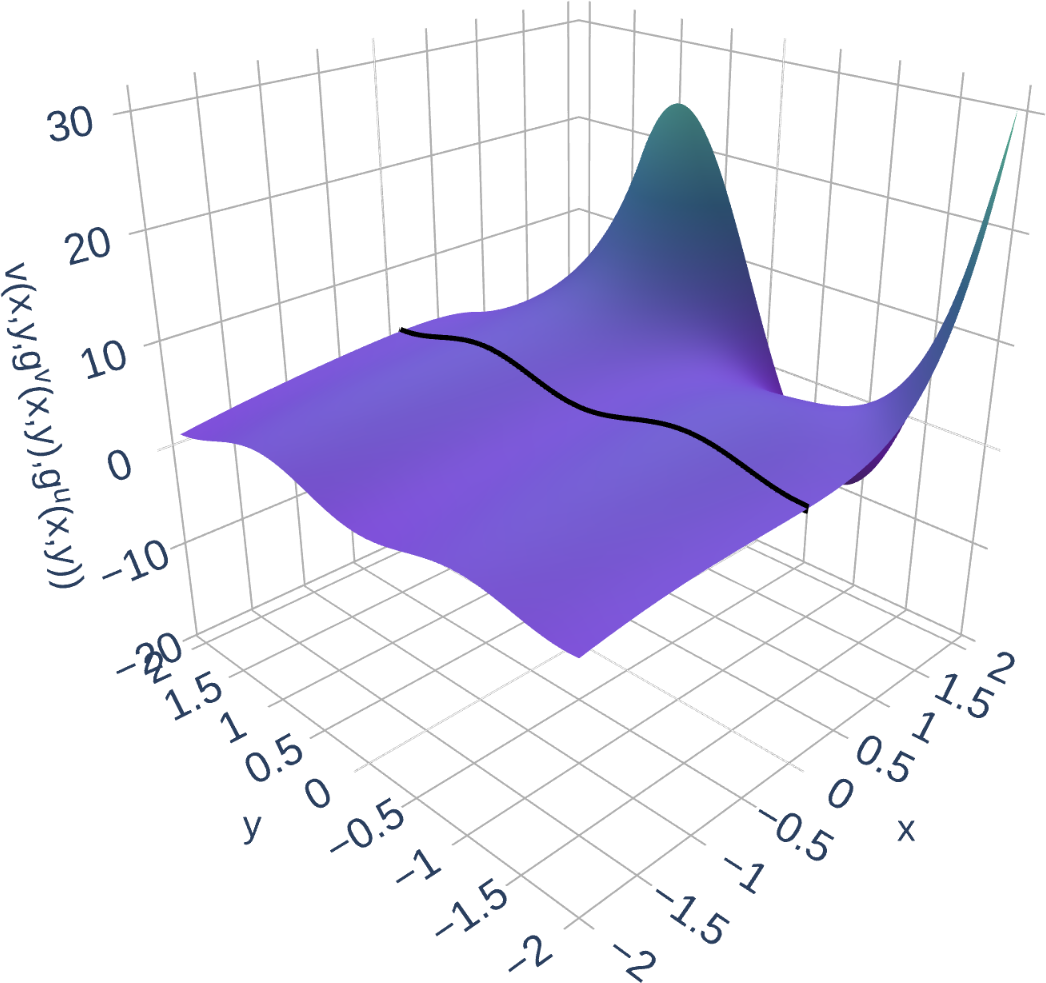}
        \caption{$v(x,y,g^v(x,y),g^u(x,y))$}
    \end{subfigure}
    \caption{Constrained expressions evaluated using $g^u(x,y) = xy+\sin(x)+y^2$ and $g^v(x,y) = x^2y\cos(y)e^x$. The value constraints on $u(x,0)$ and $v(x,0)$ are shown as black lines. The remaining constraints are not easily visualized and are therefore not shown but are satisfied nonetheless.}
    \label{fig:MultiExLinear}
\end{figure}
\end{example}

\subsection{Multivariate Constrained Expression Theorems}
This section introduces theorems for multivariate \ces\ that parallel the theorems for univariate \ces\ given in Section \ref{subsec:UniProofs}. Note that Theorems \ref{thrm:ConsructingNdCes} and \ref{thrm:IntegralConstraintModification} for multivariate \ces\ given earlier parallel Theorem \ref{thrm:UniCe} for univariate \ces, and will not be repeated here.

Theorem \ref{thrm:MultiGExists} shows that the \ce\ functional represents the family of all possible functions that satisfy the constraints.

\begin{theorem}\label{thrm:MultiGExists}%
For any function satisfying the constraints, $f (\B{x})\colon\mathbb{R}^n\mapsto\mathbb{R}$, there exists at least one free function, $g (\B{x})$, such that the \ce\ $u(\B{x},g(\B{x})) = f(\B{x})$. In other words, \ces\ are surjective functionals whose domain is all free functions and whose codomain is all functions that satisfy the constraints.

\proof
Note that the processing order used in this proof is chosen arbitrarily, and could be modified to use the processing order required for any set of constraints and still produce the same final result. Based on Theorem \ref{thrm:UniGExists}, the univariate constrained expression will return the free function if the free function satisfies the constraints. Let $ \p{1}{u}(\B{x},g(\B{x}))$ represent the univariate constrained expression for the independent variable $x_1$ that uses the free function $g(\B{x})$, $\p{2}{u}(\B{x},\p{1}{u}(\B{x},g(\B{x})))$ represent the univariate constrained expression for the independent variable $x_2$ that uses the free function $\p{1}{u}(\B{x},g(\B{x}))$, and so on up to $ \p{n}{u}(\B{x},\p{n-1}{u}(\B{x},g(\B{x})))$, which is simply the constrained expression $u(\B{x},g(\B{x}))$. If one chooses $g(\B{x}) = f(\B{x})$, then based on Theorem \ref{thrm:UniGExists} $\p{1}{u}(\B{x},f(\B{x})) = f(\B{x})$. Applying Theorem \ref{thrm:UniGExists} recursively leads to $\p{2}{u}(\B{x},\p{1}{u}(\B{x},g(\B{x}))) = f(\B{x})$ and so on until $u(\B{x},f(\B{x})) = f(\B{x})$. Hence, for any function satisfying the constraints, $f(\B{x})$, there exists a free function, $g(\B{x})=f(\B{x})$, such that the multivariate constrained expression is equal to the function satisfying the constraints, i.e., $u(\B{x},f(\B{x})) = f(\B{x})$. 
\end{theorem}

Based on the univariate \ce\ theorems, one is inclined to guess that the free function obtained in the previous theorem is not unique. As Theorem \ref{thrm:NonMultiG} shows, this inclination is correct.

\begin{theorem}\label{thrm:NonMultiG}%
For a given function satisfying the constraints, $f (\B{x})\colon\mathbb{R}^n\mapsto\mathbb{R}$, the free function, $g (\B{x})$, in the \ce\ $u(\B{x},g(\B{x})) = f(\B{x})$ is not unique. In other words, \ces\ are not injective functionals over the domain of all free functions and codomain of all functions that satisfy the constraints.

\proof
Since each expression  $\p{i}{u}(\B{x},g(\B{x}))$ used in deriving the multivariate \ce\ is derived through the univariate formulation, the results of the proof of Theorem \ref{thrm:NonUniG} apply for each each $ \p{i}{u}(\B{x},g(\B{x}))$, and therefore, the free function $g(\B{x})$ is not unique.
\end{theorem}

Like in the univariate case, this proof has immediate implications when using the constrained expression for optimization. Through the recursive application of the univariate TFC approach, for cases with no integral constraints, any terms in $g(\B{x})$ that are linearly dependent to the the support functions, $s_i(x_1)$, $s_j(x_2)$, ... , $s_k(x_n)$, will not contribute to the solution. In the multivariate case, this also includes products of the support functions that include one and exactly one support function from each independent variable, e.g., $s_i(x_1)s_j(x_2)...s_k(x_n)$. For example, suppose the support functions $s_i(x) = \{1,x,x^2\}$ and $s_j(y) = \{y,y^2\}$ were used when deriving a multivariate \ce. Then, any terms in the free function linearly dependent to any of the following functions $\{1,x,x^2,y,y^2,xy,x^2y,xy^2,x^2y^2\}$ can be removed, as they do not affect the output of the \ce.

Integral constraints add a slight complication, as they change the number of support functions used. In general, when using integral constraints in the multivariate case, the expression,
\begin{equation*}
    \alpha_{ki}S_{ij} = \delta_{jk},
\end{equation*}
is no longer true. Thus, some functions that are linearly dependent to the support functions may still be significant when included in the free function. Example \ref{ex:MultiIntDependentTerms} shows how to determine which functions linearly dependent to the support functions still have a significant contribution when included in the free function using the constraints from Example \ref{ex:MultiInt}. 

\begin{example}{Integral constraint linear dependence}\label{ex:MultiIntDependentTerms}
Consider the constraints from Example \ref{ex:MultiInt}, copied below for the reader's convenience:
\begin{equation*}
    u(x,0) = 2 u_y(x,1), \quad u(x,2) = \sin(x), \andd \int_{-1}^1 u(2,y) \dd y = 5.
\end{equation*}
In Example \ref{ex:MultiInt}, the two \ces\ were derived as,
\begin{align*}
    \p{1}{u}(x,y,g(x,y)) &= g(x,y) +\frac{2+11y-6y^2}{4}\Big(2g_y(x,1)-g(x,0)\Big)\\
    &\quad+\frac{2+13y-6y^2}{4}\Big(\sin(x)-g(x,2)\Big)\\
    \p{2}{u}(x,y,g(x,y)) &= g(x,y) + \frac{1}{2} \bigg(5 - \int_{-1}^1 g(2,\tau) \dd \tau\bigg),
\end{align*}
using the support functions $s_1(x) = 1$ for $x$ and $s_1(y)=1$, $s_2(y)=y$, and $s_3(y)=y^2$ for $y$. Theorem \ref{thrm:NonUniG} applies without modification to $\p{2}{u}$, as $\p{2}{u}$ was created using the regular univariate theory. In contrast, the switching functions of $\p{1}{u}$ were modified to include $\int_{-1}^1 \p{y}{\phi}_k(y) \dd y = 0$ for $k\in\{1,2\}$. As mentioned earlier, this means that Theorem \ref{thrm:NonUniG} must be modified slightly as $\alpha_{ki}S_{ij} \neq \delta_{jk}$ in this case. 

For the reader's convenience, the last few lines of Theorem \ref{thrm:NonUniG} have been copied below:
\begin{align*}
    y(x) &= f(x) + \beta_j \Big(\delta_{jk} - \alpha_{ki} \, \mathbb{S}_{ij}\Big) s_k(x)\\
    y(x) &= f(x) + \beta_j \Big(\delta_{jk} - \delta_{jk} \Big) s_k(x)\\
    y(x) &= f(x).
\end{align*}

For the constraints given in this example, the last few lines are re-derived. Let $B_{jk}$ be defined by
\begin{equation*}
    B_{jk} = \alpha_{ki} \mathbb{S}_{ij} = \begin{bmatrix} 1 & 1 \\ -2 & 2 \\ -4 & 4\end{bmatrix} \begin{bmatrix} \frac{1}{2} & \frac{11}{4} & -\frac{3}{2} \\ \frac{1}{2} & \frac{13}{4} & -\frac{3}{2}\end{bmatrix} = \begin{bmatrix} 1 & 6 & -3 \\ 0 & 1 & 0 \\ 0 & 2 & 0\end{bmatrix}.
\end{equation*}
Suppose that $\beta_j = \begin{Bmatrix} a, & b, & c \end{Bmatrix}$, then,
\begin{align*}
     \beta_j \Big(\delta_{jk} - \alpha_{ki} \mathbb{S}_{ij}\Big) s_k(x) &= \beta_j \Big(\delta_{jk} - B_{jk}\Big) s_k(x)\\
     &= \begin{Bmatrix} a, & b, & c \end{Bmatrix} \Bigg( \begin{bmatrix} 1 & 0 & 0 \\ 0 & 1 & 0 \\ 0 & 0 & 1 \end{bmatrix} - \begin{bmatrix} 1 & 6 & -3 \\ 0 & 1 & 0 \\ 0 & 2 & 0\end{bmatrix} \Bigg) \begin{Bmatrix} 1 \\ y \\ y^2 \end{Bmatrix}\\
      &= \begin{Bmatrix} a, & b, & c \end{Bmatrix} \Bigg(\begin{bmatrix} 0 & -6 & 3 \\ 0 & 0 & 0 \\ 0 & -2 & 1 \\ \end{bmatrix} \Bigg) \begin{Bmatrix} 1 \\ y \\ y^2 \end{Bmatrix}\\
      &= 3 a y^2-6 a y+c y^2-2 c y.
\end{align*}
Hence, only the constants $a$ and $c$ affect the final results. However, notice that the first and third row of $\delta_{jk}-B_{jk}$ are linearly dependent. Consequently, the effect of $a$ and $c$ on the final solution differs only by a constant. Therefore, one concludes that in this case, any functions linearly dependent to $y$ do not affect the final solution, and the effect of any functions linearly dependent to $1$ on the final solution will be linearly dependent with the effect of any functions linearly dependent to $y^2$ on the final solution.  

Consequently, one can remove any functions linearly dependent to two different functions, either $1$ and $y$, or $y$ and $y^2$, from the free function $g(x,y)$ without changing the final result. This can be interpreted intuitively, as the number of functions to be removed from $g(x,y)$ due to constraints on $y$ matches the number of constraints on $y$: two. Moreover, because this is a multivariate case, one can also remove any products of the support functions that include one and exactly one support function from each independent variable. However, in this case, the only support function used for the $x$ constrained expression is $1$. Therefore, this does not contribute to the terms to be removed from $g(x,y)$. 
\end{example}

As Example \ref{ex:MultiIntDependentTerms} shows, one must examine the matrix $\delta_{jk} - \alpha_{ki} \mathbb{S}_{ij}$ to calculate which terms linearly dependent to the support functions contribute to the non-uniqueness of the free function. In this matrix, a row of all zeros indicates that the corresponding support function does not contribute to the final result; for the optimization process, this corresponds to removing terms linearly dependent to that support function from the free function. Linear dependence between rows of the matrix indicates that the effect of the corresponding support functions on the final result differs only by a constant; for the optimization process, this corresponds to removing the terms in the free function linearly dependent to one of the support functions associated with the linearly dependent rows in the matrix.  

As in the univariate case, the multivariate \ces\ can be shown to be projection functionals: this is done in Theorem \ref{thrm:ProjMult}. 

\begin{theorem}\label{thrm:ProjMult}%
The multivariate \ce\ is a projection functional.

\proof
To prove Theorem \ref{thrm:ProjMult}, one must show that $u(\B{x},u(\B{x},g(\B{x}))) = u(\B{x},g(\B{x}))$. Theorems \ref{thrm:ConsructingNdCes} and \ref{thrm:IntegralConstraintModification} show that constrained expression returns a function that satisfies the constraints. In other words, for any $g(\B{x})$, $u(\B{x},g(\B{x}))$ is a function that satisfies the constraints. From Theorem \ref{thrm:MultiGExists}, if the free function used in the \ce\ satisfies the constraints, then the \ce\ returns that free function exactly. Hence, if the \ce\ function is given itself as the free function, it will simply return itself.
\end{theorem}

In addition, just as in the univariate case, Theorems \ref{thrm:MultiGExists}, \ref{thrm:NonMultiG}, and \ref{thrm:ProjMult}  allow for a more rigorous definition of the multivariate \ce. The multivariate \ce\ is a surjective, projection functional whose domain is the set of all free functions and whose codomain is the set of all functions that satisfy the constraints.

\subsection{Tensor Form}
Recursive applications of univariate TFC lead to \ces\ that lend themselves nicely to mathematical proofs, such as those in the previous section. However, at times it may be more convenient to express the constrained expression in a more compact form. Conveniently, multivariate \ces\ that are formed from recursive applications of univariate TFC can be succinctly expressed in the following tensor form,
\begin{equation*}
    u(\B{x}) = g(\B{x})+\mathcal{M}(\rho(\B{x},g(\B{x}))_{i_1i_2\dots i_n}\Phi_{i_1}(x_1)\Phi_{i_2}(x_2)\dots\Phi_{i_n}(x_n)
\end{equation*}
where $i_1, i_2, \dots, i_n$ are $n$ indices associated with the $n$-dimensions that have constraints, $\mathcal{M}$ is an $n$-dimensional tensor whose elements are based on the projection functionals, $\rho(\B{x},g(\B{x}))$, and the $n$ vectors $\Phi_{i_k}$ are vectors whose elements are based on the switching functions for the associated dimension. 

The $\mathcal{M}$ tensor can be constructed using a simple two-step process. Note that the arguments of functions and functionals are dropped in this explanation for clarity. 
\begin{enumerate}
    \item The elements of the first order sub-tensors of $\mathcal{M}$ acquired by setting all but one index equal to one are a zero followed by the projection functionals for the dimension associated with that index. Mathematically,
    \begin{equation*}
        \mathcal{M}_{1\dots i_k \dots 1} = \begin{Bmatrix} 0, & \p{k}{\rho}_1, & \cdots, & \p{k}{\rho}_{\ell_k} \end{Bmatrix},
    \end{equation*}
    where $\p{k}{\rho}_j$ indicates the $j$-th projection functional of the $k$-th independent variable and $\ell_k$ is the number of constraints associated with the $k$-th independent variable. 
    \item The remaining elements of the $\mathcal{M}$ tensor, those that have more than one index not equal to one, are the geometric intersection of the associated projection functionals multiplied by a sign ($-$~or~$+$). Mathematically, this can be written as,
    \begin{equation}\label{eq:GenMElement}
        \mathcal{M}_{i_1i_2 \dots i_n} = \pC{j}{i_j-1}\bigg[\pC{k}{i_k-1}\Big[ \cdots \big[\p{h}\rho_{i_h-1}\big] \cdots \Big]\bigg](-1)^{m+1},
    \end{equation}
    where $i_j$, $i_k$, $\dots$, $i_h$ are the indices of $\mathcal{M}_{i_1i_2 \dots i_n}$ that are not equal to one and $m$ is equal to the number of non-one indices. If no integral constraints are present, i.e., the processing order of the independent variables does not matter, then by multiple applications of Clairaut's Theorem the variables associated with the constraint operators and projection functional in Equation \eqref{eq:GenMElement} can be freely permuted \cite{M-TFC2,M-TFC}. For example, if no integral constraints are present, then Equation \eqref{eq:GenMElement} could be re-written as,
    \begin{equation*}
        \mathcal{M}_{i_1i_2 \dots i_n} = \pC{h}{i_h-1}\bigg[\pC{j}{i_j-1}\Big[ \cdots \big[\p{k}\rho_{i_k-1}\big] \cdots \Big]\bigg](-1)^{m+1}.
    \end{equation*}
    If integral constraints are present, then the processing order of the associated elements of the $\mathcal{M}$ tensor must match the processing order used in the recursive formulation.
\end{enumerate}

The elements of the vectors $\Phi_{i_k}$ are composed of a $1$ followed by the switching functions associated with the $k$-th independent variable. Mathematically,
\begin{equation*}
    \Phi_{i_k} = \begin{Bmatrix} 1, & \p{k}{\phi}_1, & \cdots, & \p{k}{\phi}_{\ell_k}\end{Bmatrix},
\end{equation*}
where $\p{k}{\phi}_j$ denotes the $j$-th switching function of the $k$-th independent variable. 

To solidify the reader's understanding of the tensor form explained above, some of the previous examples' constrained expressions are re-derived below in Examples \ref{ex:MultiExSimpleTensor}, \ref{ex:MultiIntTensor}, and \ref{ex:MultiLinearConstraintsTensor}.

\begin{example}{Non-integral constraints in tensor form}\mbox{}\label{ex:MultiExSimpleTensor}
Consider the constraints from Example \ref{ex:MultiExSimple}
\begin{gather*}
   u(0,y) = y^2\sin(\pi y), \quad u(1,y)+u(2,y) = y\sin(\pi y),\\
   u_y(x,0) = 0, \andd u(x,0) = u(x,1).
\end{gather*}
The first step of the two-step process yields the first order sub-tensors of $\mathcal{M}$.
\begin{equation*}
    \mathcal{M}_{ij}(x,y,g(x,y)) = \begin{bmatrix} 0 & -g_y(x,0) & g(x,1)-g(x,0) \\
    y^2 \sin(\pi y) -g(0,y) & \text{-} & \text{-} \\
    y \sin(\pi y) -g(2,y)-g(1,y) & \text{-} & \text{-}\end{bmatrix}
\end{equation*}
Then, the elements of $\mathcal{M}$ associated with more than one, non-one index can be found using step two. For example,
\begin{align*}
    M_{22} &= (-1)^3\pC{1}{1}[\p{2}{\rho}_1] = -[-g_y(x,0)]\Big|_{x=0} = g_y(0,0)\\
    &= (-1)^3\pC{2}{1}[\p{1}{\rho}_1] = -\frac{\partial [y^2\sin(\pi y)-g(0,y)]}{\partial y}\Big|_{y=0} = g_y(0,0).
\end{align*}
Hence, the full $\mathcal{M}$ tensor can be written as,
\small
\begin{align*}
    &\mathcal{M}_{ij}(x,y,g(x,y)) =\\ &\begin{bmatrix} 0 & -g_y(x,0) & g(x,1)-g(x,0) \\
    y^2 \sin(\pi y) -g(0,y) & g_y(0,0) & g(0,0)-g(0,1) \\
    y \sin(\pi y) -g(2,y)-g(1,y) & g_y(2,0) + g_y(1,0) & g(2,0) + g(1,0) - g(2,1)-g(1,1)\end{bmatrix}.
\end{align*}
\normalsize

The $\Phi$ vectors are built using the switching functions from the univariate cases,
\begin{equation*}
    \Phi_i(x) = \begin{Bmatrix} 1, & \frac{3-2x}{3}, & \frac{x}{3} \end{Bmatrix} \andd  \Phi_j(y) = \begin{Bmatrix} 1, & y-y^2, & -y^2 \end{Bmatrix}.
\end{equation*}
Using the $\mathcal{M}$ tensor and the $\Phi$ vectors, the full \ce\ is,
\begin{equation*}
    u(x,y,g(x,y)) = g(x,y) + \mathcal{M}_{ij}(x,y,g(x,y))\Phi_i(x) \Phi_j(y).
\end{equation*}
Expanding this expression and simplifying yields,
\begin{equation*}
\begin{aligned}
    u(x,y,&g(x,y)) = g(x,y)+\left(y-y^2\right) \Big(\frac{3-2 x}{3} g_y(0,0)-\frac{x}{3} \left(-g_y(1,0)-g_y(2,0)\right)\\
    &-g_y(x,0)\Big)-y^2 \Big(\frac{3-2 x}{3} g(0,0)-\frac{3-2 x}{3} g(0,1)-\frac{x}{3}  (-g(1,0)-g(2,0))\\
    &+\frac{x}{3}  (-g(1,1)-g(2,1)) -g(x,0)+g(x,1)\Big)+\frac{3-2 x}{3} \left(y^2 \sin (\pi  y)-g(0,y)\right)\\
    &+\frac{x}{3}  \Big(-g(1,y) -g(2,y)+y \sin (\pi  y)\Big),
\end{aligned}
\end{equation*}
the same result as in Example \ref{ex:MultiExSimple}.
\end{example}

\begin{example}{Integral constraints in tensor form}\label{ex:MultiIntTensor}
Consider the constraints from Example \ref{ex:MultiInt}
\begin{equation*}
    u(x,0) = 2 u_y(x,1), \quad u(x,2) = \sin(x), \andd \int_{-1}^1 u(2,y) \dd y = 5.
\end{equation*}
Using the same two step process the $\mathcal{M}$ tensor is constructed,
\begin{equation*}
    \mathcal{M}_{ij}(x,y,g(x,y)) = \begin{bmatrix} 0 & 2g_y(x,1)-g(x,0) & \sin(x)-g(x,2) \\ 5-\int_{-1}^1 g(2,\tau) \dd{\tau} & \int_{-1}^1 g(2,\tau) \dd{\tau} -5 & \int_{-1}^1 g(2,\tau) \dd{\tau} -5\end{bmatrix}.
\end{equation*}
Since all the elements of the $\mathcal{M}$ tensor with more than one, non-one index contain intersections including integral constraints, they must be processed in a specific order. For example,
\begin{align*}
    M_{22} = (-1)^3\pC{2}{1}[\p{1}{\rho}_1] &= 2 \frac{\partial \Big(5-\int_{-1}^1 g(2,\tau) \dd{\tau}\Big)}{\partial y}\Big|_{y=1}-\Big( 5-\int_{-1}^1 g(2,\tau) \dd{\tau}\Big)\Big|_{y=0} \\
    &= \int_{-1}^1 g(2,\tau) \dd{\tau}-5,
\end{align*}
produces the correct result that leads to a valid \ce, whereas,
\begin{equation*}
    (-1)^3\pC{1}{1}[\p{2}{\rho}_1] = \int_{-1}^1 \Big( g(2,0) - 2g_y(x,1) \Big)\dd{\tau},
\end{equation*}
does not.

The $\Phi$ vectors are built using the switching functions from the univariate cases,
\begin{equation*}
    \Phi_i(x) = \begin{Bmatrix} 1, & \frac{1}{2} \end{Bmatrix} \andd  \Phi_j(y) = \begin{Bmatrix} 1, & \frac{2+11y-6y^2}{4}, & \frac{2+13y-6y^2}{4} \end{Bmatrix}.
\end{equation*}
Using the $\mathcal{M}$ tensor and the $\Phi$ vectors, the full \ce\ is,
\begin{equation*}
    u(x,y,g(x,y)) = g(x,y) + \mathcal{M}_{ij}(x,y,g(x,y))\Phi_i(x) \Phi_j(y).
\end{equation*}
Expanding this expression and simplifying yields,
\begin{align*}
    u(x,y,&g(x,y)) = g(x,y)+\frac{1}{4} \Big(2(2-y)\left((6 y+1) g_y(x,1)+3 y \left(\int_{-1}^1 g(2,\tau ) \, d\tau -5\right)\right)\\
    &+(y-2) (6 y+1) g(x,0)+(y (6 y-13)-2) g(x,2)+(y (13-6 y)+2) \sin (x)\Big),
\end{align*}
the same result as in Example \ref{ex:MultiInt}.
\end{example}

\begin{example}{Linear constraints in tensor form}\label{ex:MultiLinearConstraintsTensor}
Consider the constraints from Example \ref{ex:MultiLinearConstraints},
\begin{gather*}
    u(0,y) = \cos(\pi y), \quad \int_{-1}^2 u(1,y) \dd{y} = e, \quad u(x,1)-u(x,2) = -2,\\
    u(x,0)+v(x,0) = 5, \andd v(0,y) = 5-\cos(\pi y).
\end{gather*}
Using the two step process the $\mathcal{M}$ tensors for $u$ and $v$ are constructed,
\begin{align*}
    \mathcal{M}^u_{ij}(x,y,&g^u(x,y),g^v(x,y)) = \begin{bmatrix} 0 & -2 -g^u(x,1)+g^u(x,2) \\ \cos(\pi y)-g^u(0,y) & 2+g^u(0,1)-g^u(0,2) \\ e-\int_{-1}^2 g^u(1,\tau)\dd{\tau} & 0\end{bmatrix},\\
    \mathcal{M}^v_{ij}(x,y,&g^u(x,y),g^v(x,y)) =\\
    &\begin{bmatrix} 0 & 5-g^v(x,0)-u(x,0,g^u(x,y)) \\ 5-\cos(\pi y)-g^v(0,y) & -5+g^v(0,0)+u(0,0,g^u(x,y))\end{bmatrix}.
\end{align*}
The $\Phi$ vectors are built using the switching functions from the univariate cases,
\begin{align*}
    \Phi^u_i(x) = \begin{Bmatrix} 1, & 1-x, & \frac{x}{3} \end{Bmatrix}, \quad   \Phi^u_j(y) = \begin{Bmatrix} 1, & \frac{1-2y}{2}\end{Bmatrix},\\
    \Phi^v_i(x) = \begin{Bmatrix} 1, & 1\end{Bmatrix}, \andd \Phi^v_j(y) =  \begin{Bmatrix} 1, & 1 \end{Bmatrix}.
\end{align*}
Using the $\mathcal{M}$ tensors and the $\Phi$ vectors, the full \ces\ are,
\begin{align*}
    u(x,y,g^u(x,y)) &= g^u(x,y) + \mathcal{M}_{ij}(x,y,g^u(x,y))\Phi_i^u(x) \Phi_j^u(y),\\
    v(x,y,g^v(x,y),g^u(x,y)) &= g^v(x,y) + \mathcal{M}_{ij}(x,y,g^v(x,y),g^u(x,y))\Phi_i^v(x) \Phi_j^v(y).
\end{align*}
Expanding these expressions and simplifying yields,
\begin{align*}
    u(x,y,&g^u(x,y)) = g^u(x,y)+\frac{1}{3} x \left(e-\int_{-1}^2 g^u(1,\tau ) \dd{\tau} \right)+\frac{1-2 y}{2}\big((1-x) (1+g^u(0,1))\\
    &+(1-x) (1-g^u(0,2))-g^u(x,1)+g^u(x,2)-2\Big)+(1-x) (\cos (\pi  y)-g^u(0,y)),\\
    v(x,y,&g^v(x,y),g^u(x,y)) = g^v(x,y)-g^v(x,0)-g^v(0,y)+g^v(0,0)-u(x,0,g^u(x,y))\\
    &+u(0,0,g^u(x,y))-\cos (\pi  y)+5,
\end{align*}
the same result as in Example \ref{ex:MultiLinearConstraints}.
\end{example}

%% file: Data/DeApplications.tex

\chapter{APPLICATIONS IN DIFFERENTIAL EQUATIONS}\label{chap:deApplications}

The \ces\ introduced in the previous chapter provide a way to analytically embed linear constraints in $n$-dimensions, which has a wide variety of applications, such as Computer-Aided Design (CAD) \cite{CoonsCadBook,CoonsCadSiggraph}, image warping \cite{CoonsImageWarping}, and security pattern design \cite{CoonsSecurityPattern}. However, this dissertation focuses on the application of TFC to differential equations. 

Differential equations are used to model and simulate physics as well as design and refine ideas, objects, systems of objects, and systems of systems. Consequently, due to their general scope, differential equations are used across a diverse range of fields, such as engineering, finance, medicine, biology, and chemistry. Although ordinary differential equations (ODEs) will be discussed, this section's primary focus will be on partial differential equations (PDEs). Due to their wide applicability, a variety of methods exist to approximate the solutions of PDEs: chief among them is the finite element method (FEM) \cite{FEM,FEA1,FEA2,FEA3}. Although FEM has been incredibly successful in solving PDEs, it does have some drawbacks.

FEM discretizes the domain into elements. This works well for low-dimensional cases, but the number of elements grows exponentially with the number of dimensions. Therefore, the discretization becomes prohibitive as the number of dimensions increases. Another issue is that FEM solves the PDE at discrete nodes, but if the solution is needed at locations other than these nodes, an interpolation scheme must be used. Moreover, extra numerical techniques are needed to perform further manipulation of the FEM solution.

Spectral methods \cite{Spectral0, Spectral1, Spectral2}, pseudo-spectral methods \cite{Pseudospectral0}, and many of their variants avoid some of these issues by proposing an analytical solution form with unknown terms that can ultimately be used to reduce the residual of the PDE at a discrete set of training points\footnote{Here, and throughout the dissertation, ``training points'' refers to the points used by the algorithm to estimate the solution of the differential equation.} and simultaneously to reduce the error in the constraints. Since these techniques use an analytical solution form, they do not require an interpolation scheme for obtaining the solution at points not seen during training, and the PDE solution can be further manipulated afterward without any special techniques. However, these methods still rely on a set of basis functions, and as will be shown later, the number of basis functions required to obtain an accurate solution can become prohibitive, e.g., see the Navier-Stokes example in Section \ref{sec:PdeExamples}.

Using neural networks (NNs) to estimate the solution of PDEs can help circumvent this issue \cite{NnDeBook,OrigOdePde,ModernPDE}, as NNs have been proven to be extremely useful when approximating high-dimensional, nonlinear functions: for example, consider computer vision problems that contain thousands of dimensions or the 200 dimensional PDE estimated via NNs in Reference \cite{ModernPDE}. Similar to the spectral and pseudo-spectral methods, the NN techniques represent a closed-form, analytical estimation of the PDE, and therefore, do not require an interpolation scheme or other special techniques to further manipulate the estimated PDE solution. In most of these techniques, the constraints and minimization of the PDE residual are handled via the loss function that is minimized when training the NN. Although some of them do use a functional solution form that automatically satisfies the constraints, these functionals do not have the same mathematical guarantees as the TFC \ces, see the theorems in Chapter \ref{chap:tfcTheory}, and cannot satisfy certain sets of constraints, such as the constraints in the Navier-Stokes problem in Section \ref{sec:PdeExamples}. In addition, the functionals of the other algorithms are multiplicative in nature, whereas TFC \ces\ are additive in nature. 

One property that all of these techniques lack is a mechanized way to satisfy the PDE constraints analytically. Such a technique is particularly important in physics-informed problems and/or when constraint information is known with a high degree of confidence \cite{Deep-TFC,piRaissi}. Fortunately, as shown in previous chapters, TFC provides a mechanized method to analytically satisfy constraints while maintaining a free function. Furthermore, with a little imagination, many of the concepts from the previously introduced PDE solution methods can be adopted and combined with the \ce\ to form useful PDE estimation methods. The benefits of these methods are:
\begin{itemize}
    \item The constraints are analytically satisfied by the \ce, and therefore, do not need to be handled in a separate way, such as augmenting the loss function, using an optimizer that handles constraints, or appending the constraints to the system of equations to be solved. 
    \item A consequence of the previous benefit is that the TFC methods are typically faster than the competing algorithms.
    \item Improved convergence, especially when the initial guess is poor \cite{JohnstonDissertation}.
\end{itemize}

\begin{figure}[!ht]
    \centering
    \includegraphics[width=\linewidth]{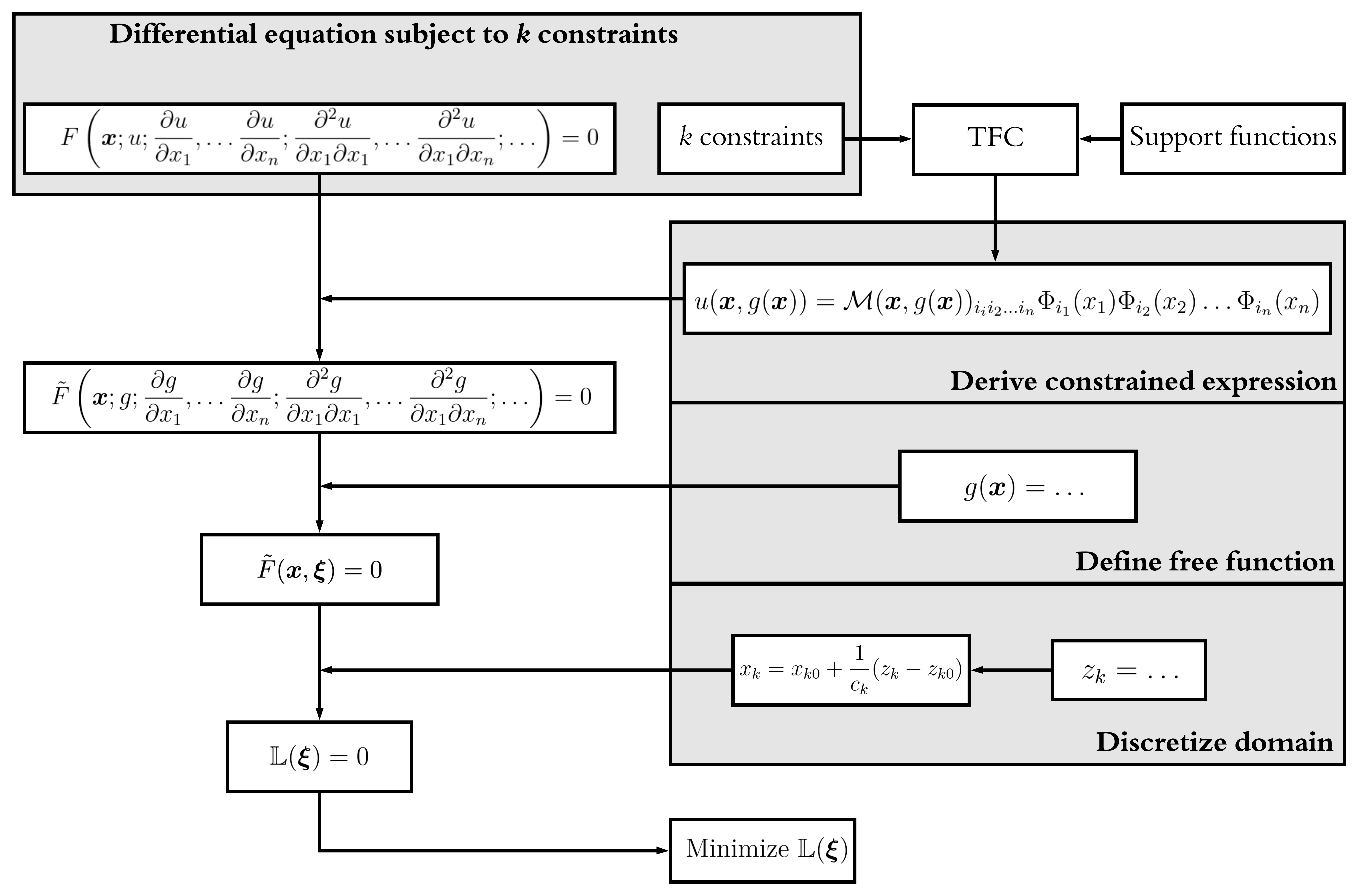}
    \caption{Differential equation solution estimation using TFC.}
    \label{fig:DeFlowchart}
\end{figure}

The methods that utilize the TFC \ce\ differ only in their choice of free function. Hence, a general methodology can be introduced that can be used for any differential equation with linear constraints, and one can switch between the PDE estimation methods simply by varying the free function. The general methodology can be summarized in five steps:
\begin{enumerate}
    \item Derive the \ce\ associated with the differential equation's\\ constraints.
    \item Define the free function, $g(\B{x})$.
    \item Discretize the domain.
    \item Formulate the loss function based on the residual of the differential equation. 
    \item Minimize the loss function in the previous step using the free function defined in step two.
\end{enumerate}
Figure \ref{fig:DeFlowchart} visually depicts these steps via a flowchart.

In general, a differential equation can be represented by some function $F$ of the independent variables, $\B{x}$, the dependent variable $u$, and its derivatives, i.e.,
\begin{equation*}
    F\left(\B{x}; u; \frac{\partial u}{\partial x_{1}},\ldots {\frac {\partial u}{\partial x_{n}}};{\frac {\partial ^{2}u}{\partial x_{1}\partial x_{1}}},\ldots {\frac {\partial ^{2}u}{\partial x_{1}\partial x_{n}}};\ldots\right) = 0.
\end{equation*}
The constraints of this differential equation can be used to create a TFC \ce. Then, this constrained expression can be substituted into the differential equation to form a differential equation with no constraints, $\tilde{F}$, that is a function of the free function rather than the dependent variable. 

Next, the free function is defined and substituted into the differential equation. In Figure \ref{fig:DeFlowchart}, the unknown parameters in $g(\B{x})$ are represented by the symbol $\B{\xi}$, e.g., $\B{\xi}$ represents $\theta$ if $g(\B{x})$ is selected as a neural network, $\B{\xi}$ represents $\B{w}$ if $g(\B{x})$ is selected as a LS-SVM, etc., see below for more details on each of these free function choices. Once the free function is substituted, the differential equation, $\tilde{F}$, becomes an algebraic equation that is a function of the independent variables $\B{x}$ and the unknown parameters $\B{\xi}$ only. 

The dependence on the independent variables is removed by discretizing the domain. In general, the domain of the free function may not coincide with the domain of the problem. For example, suppose the free function is selected as a linear combination of Chebyshev orthogonal polynomials which are defined on $[-1, 1]$. Let the free function be defined on $z \in [z_0, z_f]$ and the problem be defined on $x_k \in [x_{k_0}, x_{k_f}]$ where $k$ corresponds to the dimension. In order to use the free function, a map between the basis function domain and problem domain must be created. The simplest map is a linear one,
\begin{equation}\label{eq:linearMapping}
z = z_0 + \frac{z_f-z_0}{x_{k_f}-x_{k_0}}(x - x_{k_0}) \quad \longleftrightarrow \quad x_k = x_{k_0} + \frac{x_{k_f}-x_{k_0}}{z_f-z_0}(z - z_0).
\end{equation}

After discretizing the domain, the resultant set of algebraic equations is now only a function of the unknown parameters $\B{\xi}$: this algebraic set of equations, also known as the loss function, is denoted by $\mathbb{\B{L}}(\B{\xi})$. Thus, $\B{\xi}$ are used to minimize the difference between $\mathbb{\B{L}}(\B{\xi})$ and $\B{0}$. Once the parameters that minimize the difference are found, they can be substituted back into the \ce\ to estimate the solution of the differential equation. Note that because the \ce\ is an analytical expression, it can be easily manipulated afterward, e.g., differentiated, integrated, etc.

\section{Useful Free Function Choices}\label{sec:freeFunctions}
This section explains in detail some useful free function choices that are used in examples in later sections to solidify the reader's understanding of how TFC is applied to differential equations. Note that while it is included in this section for completeness and historical significance, as it paved the way for other machine learning algorithms, the Constrained Support Vector Machine (CSVM) methodology is no longer actively used as a free function choice, because it requires a complex analytical analysis for each new differential equation, and the resultant payoff in terms of solution error is overshadowed by the other free function choices.

\subsection{Linear Combination of Basis Functions}
A natural choice for the free function is a linear combination of basis functions, as this choice is capable of spanning the entire function space that the basis spans as the number of basis functions approaches infinity. For readers unfamiliar with univariate and multivariate basis functions, Appendix \ref{app:BasisFunctions} provides a cursory overview. 

Mathematically, a linear combination of $m$ basis functions can be expressed as,
\begin{equation*}
    g(\B{x}) =  \B{h}\T \B{\xi},
\end{equation*}
where $\B{h}\in\mathbb{R}^{m}$ is a vector of the $m$ basis functions evaluated at $\B{x}$, and $\B{\xi}\in\mathbb{R}^{m}$ is a vector of the unknown coefficients. 
The subsequent derivatives of the free function can be computed,
\begin{equation*}
    \frac{\partial^{n} g}{\partial x_k^{n}} = \left(\frac{\dd z}{\dd x_k}\right)^{n}  \frac{\partial^{n} \B{h}\T}{\partial z^n} \B{\xi}.
\end{equation*}
By defining,
\begin{equation*}
c_k := \frac{\dd z}{\dd x_k} = \frac{z_f - z_0}{x_{k_f} - x_{k_0}},
\end{equation*}
the derivative computations can be written more succinctly as, 
\begin{equation*}
    \frac{\partial^{n} g}{\partial x_k^{n}} = c_k^{n} \frac{\partial^n \B{h}\T}{\partial z^n} \B{\xi}.
\end{equation*}
It follows that a partial derivative with respect to multiple independent variables, e.g., $x_1$ and $x_2$, can be written as,
\vspace{-6pt}
\begin{equation*}
    \frac{\partial^{2} g}{\partial x_1 \partial x_2} = c_1 c_2 \frac{\partial \B{h}\T}{\partial x_1 \partial x_2} \B{\xi}.
\end{equation*}
This process applies to any derivative of the free function.

Throughout this dissertation, whenever the free function is taken to be a linear combination of basis functions, either the Chebyshev or Legendre orthogonal polynomials are used. Thus, it is useful to mention that their optimal\footnote{Optimal here refers to minimizing the condition number of the matrix to invert when minimizing the residual of the differential equation via least-squares.} discretization scheme is the Chebyshev-Gauss-Lobatto nodes \cite{Collocation1,Collocation2}. For $N$ points, the Chebyshev-Gauss-Lobatto nodes are calculated using,
\begin{equation*}
    z_j = -\cos\left(\frac{j \pi}{N-1}\right) \quad \text{for} \quad j = 0, 1, 2, \cdots, N-1.
\end{equation*}
If least-squares is used as the optimization scheme, then the collocation point distribution results in a much slower increase, relative to the uniform distribution, of the condition number of the matrix to be inverted as the number of basis functions increases. The collocation points can be realized in the problem domain through the relationship provided in Equation \eqref{eq:linearMapping}.

A linear expansion of basis functions was the first free function used for solving differential equations using TFC and has been used extensively to solve ODEs \cite{Integral-TFC,TFC-Selected,LDE,NDE,TFC-8thOrder,TFC-WeightedLS} and PDEs \cite{M-TFC2,TFC-VanillaPDE}.
However, one drawback of this free function choice is that it will become computationally prohibitive as
the dimension increases. Compelling alternatives can be found in the machine learning community.

\subsection{Support Vector Machines}
Support vector machines (SVMs) were originally introduced to solve classification problems \cite{SvmVapnik} like determining which class a given input $x$ belongs to, where there are two possible classes $x$ may belong to. The proposed solution was to find a decision boundary surface that separates the two classes. The equation of the separating boundary depended only on a few input vectors called the support~vectors.

The training data is assumed to be separable by a linear decision boundary. Hence, a separating hyperplane, $H$, with equation $\B{w}\T \B{\varphi} (\B{x}) + b = 0$, is sought. The parameters are rescaled such that the closest training point to the hyperplane $H$, $(\B{x}_k,u_k)$, is on a parallel hyperplane $H_1$ with equation $\B{w}\T \B{\varphi} (x) + b = 1$. By using the formula for orthogonal projection, if $\B{x}$ satisfies the equation of one of the hyperplanes, then the signed distance from the origin of the space to the corresponding hyperplane is given by $\B{w}\T\B{\varphi} (\B{x}) /\B{w}\T\B{w}$. Since $\B{w}\T\B{\varphi} (\B{x}) $ equals $-b$ for $H$, and $1-b$ for $H_1$, it follows that the distance between the two hyperplanes, called the ``separating margin,'' is $1/\B{w}\T\B{w}$. Thus, to find the largest separating margin, one needs to minimize $ \B{w}\T\B{w}$. The optimization problem becomes,
\begin{equation*}
    \min \dfrac{1}{2} \left(\B{w}\T \B{w} \right) \quad \text{subject to } \, u_i(\B{w}\T \B{\varphi} (\B{x}_i) + b) \geq 1, \quad i =1,\dots,n.
\end{equation*}

If a separable hyperplane does not exist, the problem is reformulated by taking into account the classification errors, or slack variables, $\Gamma_i$, and a linear or quadratic expression is added to the cost function.  The optimization problem in the non-separable case is, 
\begin{equation*}
    \min \dfrac{1}{2} \left(\B{w}\T \B{w} \right) +C \left(\sum \Gamma_i\right) \quad \text{subject to } \, u_i(\B{w}\T \B{\varphi} (\Gamma_i) + b) \geq 1-\Gamma_i.
\end{equation*}

When solving the optimization problem by using Lagrange multipliers, the function $\B{\varphi} (x)$ always shows up as a dot product with itself; thus, the kernel trick \cite{KernelTrick} can be applied. In this dissertation, the kernel function chosen is the radial basis function (RBF) kernel proposed in \cite{LS_SVP}. Hence, the function $\B{\varphi} (\B{x})$ can be written using the kernel \cite{KernelTrick},
\begin{equation*}
    K(\B{x}_i,\B{x}_j)=\B{\varphi}(\B{x}_i)\T\B{\varphi}(\B{x}_j)=\exp\left(-\dfrac{\left(\B{x}_j-\B{x}_i\right)^2}{\sigma^2}\right),
\end{equation*}
and its partial derivatives \cite{LS_SVP,SvmReviewer1Article},
\begin{align*}
    K (\B{x}_i, \B{x}_j) =  \B{\varphi} (\B{x}_i)\T \B{\varphi} (\B{x}_j) &= \exp\left(-\dfrac{(\B{x}_i - \B{x}_j)^2}{\sigma^2}\right) \nonumber\\
    K_1 (\B{x}_i,  \B{x}_j) = \B{\varphi}' (\B{x}_i)\T \B{\varphi} (\B{x}_j) &= -\dfrac{2(\B{x}_i - \B{x}_j)}{\sigma^2} \exp\left(-\dfrac{(\B{x}_i - \B{x}_j)^2}{\sigma^2}\right) \nonumber\\
    K_1\T (\B{x}_i, \B{x}_j) = \B{\varphi} (\B{x}_i)\T \B{\varphi'} (\B{x}_j) &= \dfrac{2(\B{x}_i - \B{x}_j)}{\sigma^2} \exp\left(-\dfrac{(\B{x}_i - \B{x}_j)^2}{\sigma^2}\right) \\
    K_{11} (\B{x}_i, \B{x}_j) = \B{\varphi}' (\B{x}_i)\T \B{\varphi'} (\B{x}_j) &= \dfrac{2}{\sigma^2} - \dfrac{4 (\B{x}_i - \B{x}_j)^2}{\sigma^4} \exp\left(-\dfrac{(\B{x}_i - \B{x}_j)^2}{\sigma^2}\right)\nonumber.
\end{align*}

The SVM free function choice was inspired by least-squares SVMs (LS-SVMs) and their success in solving differential equations \cite{LS_SVP}. They can be written mathematically as,
\begin{equation*}
    g(\B{x}) = \B{w}\T\B{\varphi}(\B{x}),
\end{equation*}
where $\B{w}$ is a vector of weights used in the optimization process to reduce the residual of the differential equation and $\B{\varphi}(\B{x})$ is defined in terms of the kernel function. This free function choice was used in Reference \cite{SVM-TFC} to solve ODEs and PDEs.

\subsection{Neural Networks}
The architecture of neural networks is a rich topic, and one could spend a lot of time analyzing how different architecture choices ultimately affect the accuracy of the differential equation solution. The author has chosen to leave this research effort to future work and instead focus on just one of the simpler architectures, a fully connected neural network. Each layer of a fully connected neural network consists of a nonlinear activation function composed with a linear transformation of the form $\mathcal{A} = W\cdot \B{x} + \B{b}$, where $W$ is a matrix of the neuron weights, $\B{b}$ is a vector of the neuron biases, and $\B{x}$ is a vector of inputs from the previous layer (or the inputs to the neural network if it is the first layer). Then, each layer is composed to form the entire network. For the fully connected neural networks used in this dissertation, the last layer is simply a linear output layer. For example, a neural network with three hidden layers that each use the nonlinear activation function $\psi$ and a linear output layer can be written mathematically as,
\begin{equation*}
    \mathcal{N}(\B{x};\theta) = W_4\cdot\psi\bigg(W_3\cdot\psi\Big(W_2\cdot\psi\big(W_1\cdot \B{x}+\B{b}_1\big)+\B{b}_2\Big)+\B{b}_3\bigg)+\B{b}_4,
\end{equation*}
where $\mathcal{N}$ is the neural network function, $\B{x}$ is the vector of inputs, $W_k$ are the weight matrices, $\B{b}_k$ are the bias vectors, and $\theta$ is a symbol that represents all trainable parameters of the neural network: the weights and biases of each layer constitute the trainable parameters. Note that the notation $\mathcal{N}(x,y,\dots;\theta)$ is also used in this dissertation for independent variables $x,y,\dots$ and trainable parameters $\theta$. In this dissertation, all neural networks' weights are initialized using the Glorot uniform initialization \cite{GlorotInitialization}, and the biases are initialized as zeros. Whenever a neural network is used as the free function in a \ce\ to solve a differential equation, the overall technique is referred to as Deep-TFC. This technique was used in Reference \cite{Deep-TFC} to solve a variety of PDEs.

\subsection{Extreme Learning Machines}
Extreme learning machines (ELMs) are a learning algorithm for single-hidden layer neural networks that randomly selects the hidden layer's input weights and biases and computes the output weights via least-squares \cite{ELM}. Since the weights and biases of the hidden layer are not tuned during the training, i.e., they are not trainable parameters, the neural network is linear with respect to the trainable parameters; thus, they can be computed via least-squares. In terms of the neural network description given in the previous section, an ELM can be expressed mathematically as,
\begin{equation}\label{eq:ElmForm}
    \mathcal{N}(\B{x};\theta) = W_2\cdot\psi\big(W_1\cdot \B{x}+\B{b}_1\big)
\end{equation}
where $\theta$ consists of $W_2$ only. In this dissertation, the hidden layer's weights and biases, $W_1$ and $b_1$, respectively, are initialized using the uniform distribution $U(-10,10)$ when solving ODEs and the uniform distribution on $U(-1,1)$ when solving PDEs. Whenever an ELM is used as the free function in a \ce, the overall technique is referred to as X-TFC. X-TFC was used in Reference \cite{X-TFC} to solve various ODEs and PDEs. 

\section{Useful Optimization Options}
The optimization/minimization methods introduced in this section do not constitute an exhaustive list of optimizers that can be used with TFC. Rather, they form a short list of the optimization/minimization methods used for the problems and examples given in this dissertation. Many other optimization/minimization schemes could be used in conjunction with TFC to estimate the solutions of differential equations, and exploring them is a topic of future work.

\subsection{Least-Squares}
When using basis functions, SVMs, or ELMs as the free function, the resultant minimization problem, $\mathbb{L}(\B{\xi}) = 0$, can be solved via least-squares. For linear differential equations, the loss function can be written as,
\begin{equation*}
    \mathbb{\B{L}}(\B{\xi}) = \mathbb{A}\B{\xi}-\B{b} = 0,
\end{equation*}
and a linear least-squares technique can be employed to solve,
\begin{equation*}
    \mathbb{A}\B{\xi} = \B{b}.
\end{equation*}
Appendix \ref{app:LinearLeastSquares} describes some common methods to solve the linear least-squares problem. Of the approaches presented in Appendix \ref{app:LinearLeastSquares}, the scaled QR method tends to have the lowest condition number and is thus the most numerically stable. In the specific case of ELMs, the $\mathbb{A}$ matrix tends to be ill-conditioned as the number of basis functions increases: as the number of basis functions increases, the probability of selecting nearly identical values for the weights and biases of two different neurons increases, which results in linearly dependent or nearly linearly dependent columns in $\mathbb{A}$. In this case, taking the pseudo-inverse using one of the previous techniques typically results in an inaccurate solution. To remedy this, a least-squares technique designed for ill-conditioned matrices is employed, such as the \verb"lstsq" function in the NumPy package available for Python or the \verb"lsqminnorm" function available in MATLAB.

If the differential equation is nonlinear, then a nonlinear least-squares, also known as iterative least-squares, method can be used. In this case, the loss function for the $j$-th iteration is approximated using the first two terms of the Taylor series, 
\begin{equation*}
    \mathbb{\B{L}}(\B{\xi}) \approx \mathbb{\B{L}}(\B{\xi}_j) + \mathcal{J}\Delta\B{\xi},
\end{equation*}
where $\Delta\B{\xi} = \B{\xi}-\B{\xi}_j$ and 
\begin{equation*}
    \mathcal{J}(\B{\xi}_j) = \frac{\partial \mathbb{\B{L}}(\B{\xi})}{\partial \B{\xi}}\Big|_{\B{\xi}_j},
\end{equation*}
is the Jacobian matrix of $\mathbb{\B{L}}(\B{\xi})$. The goal is to drive the loss function to zero, i.e., set $\B{\xi} = \B{\xi}_{j+1}$ in the truncated Taylor series and set it equal to zero. The result is,
\begin{equation*}
    \B{\xi}_{j+1} = \B{\xi}_j + \Delta\B{\xi},
\end{equation*}
where 
\begin{equation*}
    \mathbb{\B{L}}(\B{\xi}_j) + \mathcal{J}(\B{\xi}_j)\Delta\B{\xi} = 0
\end{equation*}
is used to determine the value of $\Delta\B{\xi}$. Notice that the solution for $\Delta\B{\xi}$ can be re-written as
\begin{equation*}
    \mathcal{J}(\B{\xi}_j)\Delta\B{\xi} = -\mathbb{\B{L}}(\B{\xi}_j),
\end{equation*}
which can be solved using the previously described linear least-squares techniques. This method is repeated until the user-specified termination condition(s) are met. For a description of the termination conditions used in this dissertation, see Section \ref{appSec:NLLS} of Appendix \ref{app:JaxCode}.

\subsection{Quasi-Newton Methods}
Quasi-Newton methods approximate the Jacobian matrix---some also approximate the Hessian---of a given function and use them to minimize the given function. To use them with TFC, the loss function is first converted to a scalar by taking its norm: in this dissertation, that norm is either the $L_1$ or $L_2$ norm. Then, the new loss function is minimized using the Quasi-Newton method. In this dissertation, the only Quasi-Newton method used is the limited-memory Broyden-Fletcher-Goldfarb-Shanno \cite{BfgsBook} (L-BFGS) algorithm.

\subsection{Gradient Descent}
Gradient descent algorithms are an iterative method used to minimize a scalar loss function---as with Quasi-Newton methods, the TFC loss functions are converted to scalar functions by taking either their $L_1$ or $L_2$ norm---via its Jacobian matrix. In its simplest form,
\begin{equation*}
    \B{\xi}_{j+1} = \B{\xi}_j + \Delta\B{\xi}
\end{equation*}
where
\begin{equation*}
    \Delta\B{\xi} = -\lambda\mathcal{J}(\B{\xi}_j)
\end{equation*}
and $\lambda$ is some positive constant frequently referred to as the learning rate. However, the gradient descent can become more complex by, for example, randomly sampling a subset of the training points for each iteration as is done in stochastic gradient descent \cite{Sgd1,Sgd2}, adapting the learning rate based on data from previous iterations as is done in AdaGrad \cite{AdaGrad}, or using a concept analogous to linear momentum to inform the update along with gradient descent \cite{MomentumGradientDescent}. Variants of these techniques are utilized by the Adaptive Momentum Estimation (Adam) algorithm \cite{AdamOptimizer}, which is the gradient descent algorithm used in this dissertation.  

\subsection{Constrained Least-Squares Support Vector Machines}
The inspiration to use Least-Squares Support Vector Machines (LS-SVMs) stems from References \cite{LS_SVP} and \cite{LsSvmPde}, which used LS-SVMs to solve ODEs and PDEs, respectively. Essentially, this technique uses an LS-SVM to reduce the residual of the differential equation via least-squares, where Lagrange multipliers are used to enforce the differential equation constraints. When TFC is used with an SVM as the free function to solve the differential equation using the LS-SVM optimization technique, the overall methodology is referred to as a Constrained Support Vector Machine (CSVM) \cite{SVM-TFC}. This methodology is best understood via an example. 

\begin{example}{CSVM applied to a linear, first-order ODE}
Consider a first-order, linear ODE,
\begin{equation*}
    \dot{y}- p (t) y = r (t), \quad \text{subject to} \quad y (t_0) = y_0,
\end{equation*}
where $\dot{y} = \frac{\text{d} y}{\text{d} t}$. The TFC \ce\ for the constraint is,
\begin{equation*}
    y (t, g(t)) = g (t) + y_0 - g(0).
\end{equation*}
Further, let the free function be defined as an SVM,
\begin{equation*}
    g (t) = \B{w}\T \B{\varphi} (t),
\end{equation*}
so the constrained expression becomes,
\begin{equation}\label{y_TFC}
    y (t, \B{w}) = \B{w}\T \big(\B{\varphi} (t) - \B{\varphi} (t_0)\big) + y_0.
\end{equation}
Notice that a least-squares technique cannot be formed directly from the residual of the differential equation because $\B{\varphi} (t)$ is only defined via the kernel trick, i.e., only dot products with itself are defined. Hence, a loss function is constructed,
\begin{equation*}
    \min\dfrac{1}{2}\left(\B{w}\T\B{w} + \gamma\B{e}\T\B{e}\right)
\end{equation*}
where $\gamma$ is a positive, user-specified constant and 
\begin{equation*}
    \B{e} = e_i = \dot{y}(t_i,\B{w})- p (t_i) y_i(t_i,\B{w}) - r(t_i).
\end{equation*}
This loss function contains the original term used to find the largest separating margin, $\B{w}\T\B{w}$, as well as a term to reduce the residual of the ODE, $\gamma\B{e}\T\B{e}$. Since the \ce\ embeds the differential equation constraint, the only constraint that remains is the one associated with the error term, $\B{e}$. Hence, the optimization problem is,
\begin{align*}
    &\min\dfrac{1}{2}\left(\B{w}\T\B{w} + \gamma\B{e}\T\B{e}\right)\\
    &\text{subject to} \quad \B{w}\T \dot{\B{\varphi}} (t_i) - p(t_i) \big(\B{w}\T \B{\varphi} (t_i) - \B{w}\T \B{\varphi} (t_0) + y_0\big) - r (t_i) - e_i = 0.
\end{align*}

Using the method of Lagrange multipliers a term is introduced for the constraint on the residuals, which leads to the expression,
\begin{align*}
    \mathbb{\B{L}}(\B{w}, \B{e}, \B{\alpha}) =& \dfrac{1}{2} \left(\B{w}\T \B{w} + \gamma \B{e}\T \B{e}\right) \\
    &-\sum_{i=1}^N \alpha_i \big[\B{w}\T \dot{\B{\varphi}} (t_i) - p (t_i) \left(\B{w}\T \B{\varphi} (t_i) - \B{w}\T \B{\varphi} (t_0) + y_0\right) - r (t_i) - e_i\big],
\end{align*}
where $\alpha_i$ are the Lagrange multipliers. The values that force the gradients of $\mathbb{\B{L}}$ to be equal to zero give candidates for the minimum,
\begin{equation*}
\begin{aligned}
    \dfrac{\partial\mathbb{\B{L}}}{\partial\B{w}} = 0 &\qquad\to\qquad \B{w} = \sum_{i=1}^N \alpha_i \left[\dot{\B{\varphi}} (t_i) - p(t_i) \left(\B{\varphi} (t_i) - \varphi(t_0)\right)\right] \\ 
    \dfrac{\partial\mathbb{\B{L}}}{\partial e_i} = 0 &\qquad\to\qquad e_i = -\frac{\alpha_i}{\gamma} \\
    \dfrac{\partial\mathbb{\B{L}}}{\partial\alpha_i} = 0 &\qquad\to\qquad 0 = \B{w}\T\dot{\B{\varphi}} (t_i) -  p(t_i)\left(\B{w}\T\left(\B{\varphi} (t_i) - \B{\varphi} (t_0)\right) + y_0\right) - r (t_i) - e_i.
\end{aligned}
\end{equation*}
Using,
\begin{equation*}
    \B{w} = \ds \sum_{j=1}^N \alpha_j \left[\dot{\B{\varphi}} (t_j) - p(t_j) \left(\B{\varphi} (t_j) - \varphi(t_0)\right)\right],
\end{equation*}
one obtains a new formulation of the approximate solution given by Equation \eqref{y_TFC} that can be expressed in terms of the kernel and its derivatives. One can combine the three equations obtained by setting the gradients of $\mathbb{\B{L}}$ equal to zero together to create a linear system with unknowns $\alpha_j$,
\begin{equation*}
    M_{ij} \alpha_j = r (t_i) + p (t_i) y_0.
\end{equation*}
The coefficient matrix, $M_{ij}$, is given by,
\begin{equation*}
\begin{aligned}
    M_{ij} = K_{11} (t_i, t_j) &- p (t_j)\left[K_1 (t_i, t_j) - K_1 (t_i, t_0)\right]  - p (t_i) K_y (t_i, t_j) + \delta_{ij}/\gamma,
\end{aligned}
\end{equation*}
where,
\begin{equation*}
\begin{aligned}
    K_4 (t_i,t_j) &= K (t_i, t_j) - K (t_j, t_0) - K (t_i, t_0) + 1, \\
    K_y (t_i, t_j) & = K_1 (t_j, t_i) - K_1 (t_j, t_0) - p (t_j) K_4 (t_i, t_j).
\end{aligned}
\end{equation*}
Finally, in terms of the kernel matrix, the approximate solution at the training points, $t_i$, is given by,
\begin{equation*}
    y(t_i) = \ds\sum_{j=1}^N \alpha_j K_y (t_i, t_j) + y_0,
\end{equation*}
and a formula for the approximate solution at an arbitrary point $t$ is given by,
\begin{equation*}
    y(t) = \ds\sum_{j=1}^N \alpha_j K_y (t, t_j) + y_0.
\end{equation*}
\end{example}

The CSVM technique creates a loss function based on the residual of the differential equation that can ultimately be solved via least-squares. Moreover, the least-squares system and the constrained expression can be rewritten in terms of the kernel function and its derivatives. Although not utilized in this dissertation, a similar derivation for first-order, nonlinear ODEs is included in Appendix \ref{app:nonLinearSvm} for completeness. 

\section{Numerical Implementation}
Even for simple PDEs, taking the derivatives necessary to implement TFC analytically is tedious and error-prone. The errors in taking the derivatives can be reduced by using a symbolic programming paradigm, but the results must still be copied into another framework, Python, MATLAB, etc., which is time-consuming and error-prone as well. Of course, TFC could be implemented directly in a symbolic program, but the computation speed would suffer: one of TFC's main benefits. If these were the only options, applying TFC to differential equations would be cumbersome, frustrating, and slow, and users would most likely choose other differential equation solution methods due to this pitfall. Fortunately, automatic differentiation alleviates the issues that plague the other implementation options. 

Automatic differentiation utilizes the chain rule of differential calculus and modifies the variable types to calculate derivatives \cite{AdInMachineLearning}. As such, automatic differentiation can be applied to code with minimal changes and can evaluate ``derivatives at machine-level precision with only a small constant factor of overhead and ideal asymptotic efficiency'' \cite{AdInMachineLearning}. For example, Reference \cite{AdNilponentAlgebraSchutte} utilizes a nilpotent algebra to calculate arbitrary order derivatives; this technique is an example of forward-mode automatic differentiation, where the derivative is calculated alongside the primary value. A second type is reverse-mode automatic differentiation \cite{ReverseMode}, where the chain rule is traversed from the output backwards towards the input to calculate the derivative: also known as back-propagation. In general, a good rule of thumb is that reverse mode differentiation should be used for functions $f\colon\mathbb{R}^n\mapsto\mathbb{R}^m$ where $n \gg m$; otherwise, forward mode differentiation should be used. Hence, forward mode automatic differentiation will typically be used when implementing TFC. An in-depth understanding of automatic differentiation is not required to understand TFC's numerical implementation, so it will not be discussed here. However, if the reader is interested in learning more, they should consult Reference \cite{AdInMachineLearning}.

JAX \cite{JaxOriginalPaper,JaxGithub} is a framework for Python that combines the automatic differentiation power of Autograd \cite{autograd} with XLA (Accelerated Linear Algebra) \cite{Tensorflow} to produce fast, composable transformations of NumPy/Python code. Moreover, a just-in-time compiler (JIT) allows one to easily convert their code into XLA-optimized kernels. Ultimately, this allows the user to easily compute the derivatives necessary to apply TFC to differential equations, and JIT-ing the resultant code makes the run time fast: many of the differential equations in this dissertation were estimated via TFC in less than a second. Numerically implementing TFC in a JIT-able way via JAX was not a simple plug-and-play. Rather, the author wrote over 8,000 lines of C++ and Python code that interface with JAX to produce the final product. This code forms some general-use classes that can be used to apply TFC to a large variety of differential equations. The intricacies of this code are not germane to the topics covered in the body of this dissertation, and therefore are not included here; a more detailed description of the codebase can be found in Appendix \ref{app:JaxCode}, and the \href{https://tfc-documentation.readthedocs.io/en/latest/}{\textcolor{blue}{\underline{reference documentation}}} contains a complete description. Furthermore, this general-use code is publicly available on the \href{https://github.com/leakec/tfc}{\textcolor{blue}{\underline{TFC GitHub}}} and includes the scripts used to generate many of the examples and results found in this dissertation \cite{TfcGithub}.

The computations for all examples and results in this dissertation were performed in Python on a desktop computer running Ubuntu 20.04 with an Intel\textsuperscript{\textregistered} Core\texttrademark \, i5-2400 and 16 GB of RAM. All run times were calculated using the \verb"process_timer" function from the Python \verb"time" package.

\section{Simple PDE Example}\label{sec:PdeExamples}
\begin{figure}[!ht]
    \centering
    \includegraphics[width=0.75\linewidth]{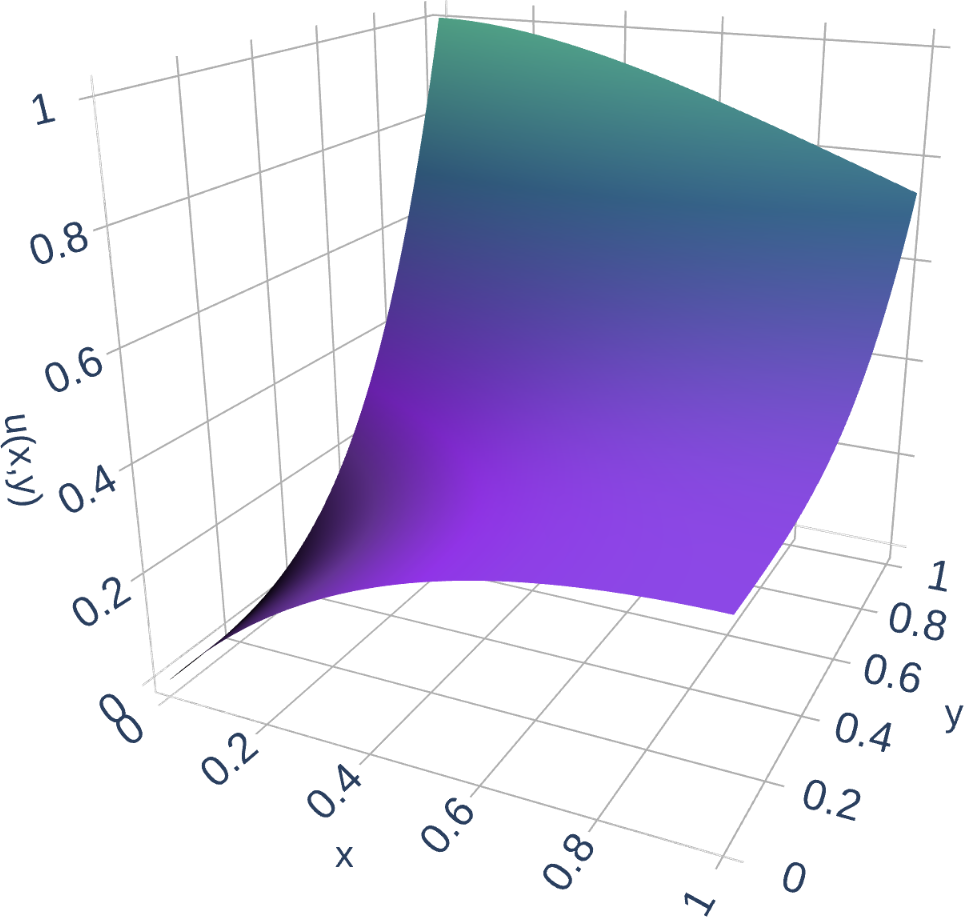}
    \caption{Analytical solution of the simple PDE.}
    \label{fig:LaragisPdeSoln}
\end{figure}

To better understand each of the previously introduced free functions, this section applies each one to the same linear PDE:
\begin{equation*}
    u_{xx} (x,y)+u_{yy} (x,y) = e^{-x}(x - 2 + y^3 + 6y)
\end{equation*}
where $x,y \in [0,1]$ and subject to,
\begin{eqnarray*}
    u(0,y) &=& y^3\\
    u(1,y) &=& (1+y^3)e^{-1}\\
    u(x,0) &=& xe^{-x}\\
    u(x,1) &=& e^{-x}(x+1),
\end{eqnarray*}
which has the true solution $u(x,y) = e^{-x}(x + y^3)$. The true solution is shown in Figure \ref{fig:LaragisPdeSoln}.
Following the step-by-step method given earlier, the PDE can be re-written as,
\begin{equation*}
    F(x,y,u,u_{xx},u_{yy}) = u_{xx} (x,y)+u_{yy} (x,y) - e^{-x}(x - 2 + y^3 + 6y) = 0.
\end{equation*}
The constraints can be embedded into a \ce,
\begin{align*}
    \p{1}{u}(x,y,g(x,y)) &= g(x,y) + (1-x)\Big(y^3-g(0,y)\Big) + x\Big((1+y^3)e^{-1}-g(1,y)\Big) \\
    \p{2}{u}(x,y,g(x,y)) &= g(x,y) + (1-y)\Big(xe^{-x}-g(x,0)\Big) + y\Big(e^{-x}(x+1)-g(x,1)\Big)
\end{align*}
where $\p{1}{u}$ can be used as the free function in $\p{2}{u}$ or $\p{2}{u}$ can be used as the free function in $\p{1}{u}$ to create the full \ce. The \ce\ written in tensor form is,
\begin{equation*}
    u(x,y) = g(x,y,g(x,y))+\mathcal{M}(x,y,g(x,y))_{ij}\Phi_i(x)\Phi_j(y)
\end{equation*}
where
\begin{equation*}
    \mathcal{M}(x,y,g(x,y))_{ij} = \begin{bmatrix} 0 & xe^{-x}-g(x,0) & e^{-x}(x+1)-g(x,1) \\ y^3-g(0,y) & g(0,0) & g(0,1)-1 \\ (1+y^3)e^{-1}-g(1,y) & g(1,0)-e^{-1} & g(1,1)-2e^{-1}\end{bmatrix},
\end{equation*}
\begin{equation*}
    \Phi_i(x) = \begin{Bmatrix} 1, & 1-x, & x\end{Bmatrix}, \andd
    \Phi_j(y) = \begin{Bmatrix} 1, & 1-y, & y\end{Bmatrix}.
\end{equation*}
Substituting the \ce\ into $F$ yields $\tilde{F}(x,y,g(x,y))$, which does not have any constraints. Now, the various free function choices introduced earlier will be used to minimize $\tilde{F}$. 

\begin{example}{Simple PDE solved using basis functions \cite{M-TFC2}}\label{ex:SimplePdeBasis}
Let $g(x,y)$ be a linear expansion of Chebyshev orthogonal polynomials, and let $m$ be the maximum degree of said polynomials. Remember, as shown in Appendix \ref{app:BasisFunctions}, that the two-dimensional basis set is just a tensor product of the univariate Chebyshev orthogonal polynomials. Further, recall from the result of Theorem \ref{thrm:NonMultiG} that the basis functions linearly dependent to the support functions must be removed from the expansion: in the multivariate case, this also includes products of the support functions that include exactly one support function from each independent variable, e.g., $s_i(x_1)s_j(x_2)...s_k(x_n)$. To expound, suppose instead that the linear expansion for $g(x,y)$ was simply the set of monomials.
\begin{equation*}
    g(x,y) = \xi_1 + \xi_2 x + \xi_3 y + \xi_4 x^2 + \xi_5 xy + \xi_6 y^2 + \dots
\end{equation*}
In this case, the terms $1$, $x$, $y$, and $xy$ need to be removed from the expansion, as $1$, $x$, and $y$ are used as support functions in the \ce. The same needs to be done for the Chebyshev orthogonal expansion used in this example.

Since the linear expansion is a tensor product and the terms linearly dependent to the support functions have to be removed, the degree of the expansion, $m$, and the number of basis functions in the expansion do not have a simple relationship. Therefore, the degree of the expansion, $m$, and the number of basis functions in the expansion are tabulated for this example in Table \ref{tab:numBasis}.
\begin{table}[H]
\centering
\caption{Tabulated values for the degree of basis expansion and equivalent number of basis functions.}
\label{tab:numBasis}
\begin{tabular}{|c|c|} 
\hline
 \makecell{\textbf{m}} & \makecell{\textbf{Number of Functions}}\\ \hline
5  & 17 \\\hline
10 & 62 \\\hline
15 & 132 \\\hline
20 & 227 \\\hline
25 & 347 \\
\hline
\end{tabular}
\end{table}

Once the free function is substituted into the differential equation, $\tilde{F}=\tilde{F}(x,y,\B{\xi})=0$. Next, the domain is discretized. Since Chebyshev orthogonal polynomials are used, the domain is discretized using Chebyshev-Gauss-Lobatto nodes. Let the number of points per independent variable be given by $n$. For example, a value of $n = 5$ would imply a $5\times 5$ grid or 25 total training points. After the domain is discretized $\tilde{F}$ becomes $\mathbb{\B{L}}(\B{\xi}) = 0$. The PDE in this example is linear, so $\mathbb{\B{L}}$ is linear in $\B{\xi}$, and therefore, linear least-squares can be used to minimize $\mathbb{\B{L}}$. 

Let the test set be a $100\times100$ grid of uniformly spaced points. Table \ref{tab:LagarisCp} shows the maximum test set solution error,
\begin{equation*}
    e = \max_{(x,y) \in \text{test set}} |u(x,y)-u_{\text{true}}(x,y)|,
\end{equation*}
where $u_{\text{true}}(x,y)$ is the true solution given earlier, for different values of $n$ and $m$. Table \ref{tab:LagarisCp} shows that in general as the number of basis functions and training points increases, the maximum test set solution error decreases.
\begin{table}[H]
\caption{Maximum test set solution error using TFC with Chebyshev orthogonal polynomials.}\label{tab:LagarisCp}
\begin{tabular}{c|ccccc}
\noalign{\hrule height 1.0pt}
\hline
\diagbox[height=1.25cm]{\textbf{n}}{\textbf{m}} & \textbf{5} & \textbf{10} & \textbf{15} & \textbf{20} & \textbf{25}\\
\hline
5  & $6.26\times 10^{-4}$ & $\text{-}$   & $\text{-}$   & $\text{-}$   & $\text{-}$ \\
10 & $5.53\times 10^{-4}$ & $1.20\times 10^{-10}$ & $\text{-}$   & $\text{-}$   & $\text{-}$ \\
15 & $5.30\times 10^{-4}$ & $1.17\times 10^{-10}$ & $4.44\times 10^{-16}$ & $\text{-}$   & $\text{-}$ \\
20 & $5.20\times 10^{-4}$ & $1.16\times 10^{-10}$ & $5.00\times 10^{-16}$ & $4.44\times 10^{-16}$ & $\text{-}$ \\
25 & $5.13\times 10^{-4}$ & $1.15\times 10^{-10}$ & $7.22\times 10^{-16}$ & $2.61\times 10^{-15}$ & $5.55\times 10^{-16}$\\
30 & $5.09\times 10^{-4}$ & $1.14\times 10^{-10}$ & $6.66\times 10^{-16}$ & $8.88\times 10^{-16}$ & $3.22\times 10^{-15}$\\  \noalign{\hrule height 1.0pt}
\end{tabular}
\end{table}

Since the TFC method of solving differential equations is closely related to the spectral method---the only real difference is how the constraints are handled---it is worth comparing the two methods. To this end, Table \ref{tab:LagarisSpectral} shows the same results as Table \ref{tab:LagarisCp} but using the spectral method.


\begin{table}[H]
\caption{Maximum test set solution error using spectral method with Chebyshev orthogonal polynomials.}\label{tab:LagarisSpectral}
\begin{tabular}{c|ccccc}
\noalign{\hrule height 1.0pt}
\hline
\diagbox[height=1.25cm]{\textbf{n}}{\textbf{m}} & \textbf{5} & \textbf{10} & \textbf{15} & \textbf{20} & \textbf{25}\\
\hline
5  & $4.25\times 10^{-4}$ & \text{-}   & \text{-}   & \text{-}   & \text{-} \\
10 & $3.40\times 10^{-4}$ & $7.11\times 10^{-11}$ & \text{-}   & \text{-}   & \text{-} \\
15 & $3.16\times 10^{-4}$ & $7.95\times 10^{-11}$ & $1.41\times 10^{-12}$ & \text{-}   & \text{-} \\
20 & $3.04\times 10^{-4}$ & $7.77\times 10^{-11}$ & $4.85\times 10^{-12}$ & $5.75\times 10^{-12}$ & \text{-} \\
25 & $2.97\times 10^{-4}$ & $7.69\times 10^{-11}$ & $3.45\times 10^{-12}$ & $9.91\times 10^{-12}$ & $2.71\times 10^{-11}$\\
30 & $2.92\times 10^{-4}$ & $7.59\times 10^{-11}$ & $3.12\times 10^{-12}$ & $1.19\times 10^{-11}$ & $1.79\times 10^{-11}$\\ 
\noalign{\hrule height 1.0pt}
\end{tabular}
\end{table}

Comparing Tables \ref{tab:LagarisCp} and \ref{tab:LagarisSpectral} reveals that the spectral method is slightly more accurate---less than an order of magnitude---than TFC when the number of basis functions is low. However, as the number of basis functions increases, TFC becomes as many as five orders of magnitude more accurate than the spectral method. The accuracy difference between the two methods for a low number of basis functions stems from the fact that the spectral method can relax the error on the constraints in order to reduce the average error over the domain, whereas TFC is constrained to satisfy the constraints exactly, and so does not have the same freedom. The accuracy difference between the two methods for a high number of basis functions stems from the fact that TFC effectively has more information than spectral method, as it has the exact constraint information over the entire boundary, whereas spectral method only has information about the constraints at discrete points on the associated boundaries. Moreover, TFC is faster than the spectral method, as the matrix that is inverted during the least-squares process is smaller; it is smaller because it does not contain the extra equations that the spectral method needs to satisfy the constraints. 
\end{example}

\begin{example}{Simple PDE solved using CSVM \cite{SVM-TFC}}\label{ex:SimplePdeCsvm}
Let $g(x,y)$ be a SVM, so $\tilde{F}$ becomes $\tilde{F}(x,y,\B{w})$. Then, the CSVM technique must be applied to rewrite the \ce\ and the optimization process in the dual form. For this example only, let superscripts denote a derivative with respect to the superscript variable and a subscript be a normal tensor index: this is done for clarity and compactness. For example, the symbol $A^{xx}_{ij}$ would denote a second-order derivative of the second-order tensor $A_{ij}$ with respect to the variable $x$, i.e., $\frac{\partial^2 A_{ij}}{\partial x^2}$. In the same spirit, for this example only, the arguments of most functions and functionals will be dropped.

The \ce\ shown earlier can be re-written as,
\begin{align*}
    &u = A_{ij}\Phi_i\Phi_j + w_j\varphi_j(x,y) - w_kB_{ijk}\Phi_i\Phi_j, \quad \text{where}\\
    &A_{ij} = \begin{bmatrix} 0 & xe^{-x} & e^{-x}(x+1) \\ y^3 & 0 & -1 \\ (1+y^3)e^{-1} & g(1,0)-e^{-1} & -2e^{-1}\end{bmatrix}\\
    &B_{ijk} = \begin{bmatrix} 0 & \varphi_k(x,0) & \varphi_k(x,1) \\ \varphi_k(0,y) & -\varphi_k(0,0) & -\varphi_k(0,1) \\ \varphi_k(1,y) & -\varphi_k(1,0) & -\varphi_k(1,1) \end{bmatrix}\\
    &\Phi_i = \begin{Bmatrix} 1, & 1-x, & x \end{Bmatrix}, \andd \Phi_j = \begin{Bmatrix} 1, & 1-y, & y \end{Bmatrix}.
\end{align*}
Now, discretize the domain and use Lagrange multiplies to form $\mathbb{\B{L}}$,
\begin{equation*}
   \mathbb{\B{L}} (\B{w}, \B{\alpha}, \B{e}) = \frac{1}{2} w_i w_i + \frac{\gamma}{2} e_I e_I - \alpha_I (u^{xx}_I+u^{yy}_I-f_I-e_I),
\end{equation*}
where $u_I$ is a vector whose elements are $u(x_I,y_I,\B{w})$ where $(x_I,y_I)$ is the $I$-th training point. The gradients of $\mathbb{\B{L}}$ give candidates for the minimum,
\begin{align*}
    &\frac{\partial \mathbb{\B{L}}}{\partial w_k} = w_k - \alpha_I(\varphi^{xx}_{Ik}-B^{xx}_{Iijk}\Phi_i\Phi_j+\varphi^{yy}_{Ik}-B^{yy}_{Iijk}\Phi_i\Phi_j) = 0\\
    & \frac{\partial\mathbb{\B{L}}}{\alpha_I} = \hat{z}^{xx}_I+\hat{z}^{yy}_I-f_I-e_I = 0 \\
    & \frac{\partial\mathbb{\B{L}}}{e_I} = \frac{\gamma}{2}e_I-\alpha_I = 0,
\end{align*}
where $\varphi_{Ik}$ is the second-order tensor composed of the vectors $\varphi_k(x_I,y_I)$ and $B_{Iijk}$ is the fourth-order tensor composed of the third-order tensors $B(x_I,y_I)_{ijk}$. The gradients of $\mathbb{\B{L}}$ can be used to form a system of simultaneous linear equations to solve for the unknowns and write $u$ in the dual form. The system of simultaneous linear equations is,
\begin{equation*}
    {\cal A}_{IJ}\alpha_J = {\cal B}_I,
\end{equation*}
where
\begin{align*}
    {\cal A}_{IJ} =\ &\varphi^{xx}_{Ik}\varphi^{xx}_{Jk}-\varphi^{xx}_{Ik} B^{xx}_{Jijk}\Phi_i\Phi_j+\varphi^{xx}_{Ik}\varphi^{yy}_{Jk}-\varphi^{xx}_{Ik} B^{yy}_{Jijk}\Phi_i\Phi_j-B^{xx}_{Iijk}\Phi_i\Phi_j\varphi^{xx}_{Jk}\\
    &+B^{xx}_{Iijk}\Phi_i\Phi_jB^{xx}_{Jmnk}\Phi_m\Phi_n-B^{xx}_{Iijk}\Phi_i\Phi_j\varphi^{yy}_{Jk}+B^{xx}_{Iijk}\Phi_i\Phi_jB^{yy}_{Jmnk}\Phi_m\Phi_n+\varphi^{yy}_{Ik}\varphi^{xx}_{Jk}\\
    &-\varphi^{yy}_{Ik} B^{xx}_{Jijk}\Phi_i\Phi_j+\varphi^{yy}_{Ik}\varphi^{yy}_{Jk}-\varphi^{yy}_{Ik} B^{yy}_{Jijk}\Phi_i\Phi_j-B^{yy}_{Iijk}\Phi_i\Phi_j\varphi^{xx}_{Jk}\\
    &+B^{yy}_{Iijk}\Phi_i\Phi_jB^{xx}_{Jmnk}\Phi_m\Phi_n-B^{yy}_{Iijk}\Phi_i\Phi_j\varphi^{yy}_{Jk}+B^{yy}_{Iijk}\Phi_i\Phi_jB^{yy}_{Jmnk}\Phi_m\Phi_n+\frac{1}{\gamma}\delta_{IJ}\\
    {\cal B}_I =\ & f_I-A^{xx}_{Iij}\Phi_i\Phi_j-A^{yy}_{Iij}\Phi_i\Phi_j
\end{align*}
where $\Phi_m=\Phi_i$, $\Phi_n=\Phi_j$, and $A_{Iijk}$ is the fourth order tensor composed of the third order tensors $A(x_I,y_I)_{ijk}$. The dual-form of the solution is,
\begin{align*}
    u(&x,y,\B{\alpha}) = A_{ij}\Phi_i\Phi_j\\
    &+\alpha_I\bigg[\varphi^{xx}_{Ik}\varphi(x,y)_k-B^{xx}_{Iijk}\Phi_i\Phi_j\varphi_k(x,y)+\varphi^{yy}_{Ik}\varphi_k(x,y)-B^{yy}_{Iijk}\Phi_i\Phi_j\varphi_k(x,y)\bigg]\\
    &-\alpha_I\bigg[\varphi^{xx}_{Ik}B_{ijk}\Phi_i\Phi_j-B^{xx}_{Iijk}\Phi_i\Phi_jB_{mnk}\Phi_m\Phi_n+\varphi^{yy}_{Ik}B_{ijk}\Phi_i\Phi_j-B^{yy}_{Iijk}\Phi_i\Phi_jB_{mnk}\Phi_m\Phi_n\bigg].
\end{align*}
The system of simultaneous linear equations and the dual form of the solution can be written and were solved using the kernel matrix and its partial derivatives.

The attentive reader will notice two user-specified hyperparameters remain to be selected: $\sigma$, the user-specified parameter that modifies the kernel matrix, and $\gamma$, the parameter that appears in the loss function. In Reference \cite{SVM-TFC}, a grid search was used to find the pair of hyperparameters that minimized the residual of the differential equation on a set of validation points. Using that set of hyperparameters and $100$ training points, the test set's maximum error was $5.561 \times 10^{-8}$.
\end{example}

\begin{example}{Simple PDE solved using Deep-TFC \cite{Deep-TFC}}
Let $g(x,y)$ be a neural network as defined earlier with nonlinear activation function $\psi = \tanh$. Further, let this neural network have six hidden layers with 15 neurons per layer and a linear output layer. Substituting the neural network as the free function into the constrained expression, then substituting the \ce\ into the differential equation, and finally discretizing the domain leads to a loss function, $\mathbb{L}(\theta)$, where $\theta$ are the trainable parameters of the neural network as defined earlier; the domain was discretized using a $10\times10$ grid of uniformly spaced points. 

Neural networks are typically trained using gradient descent algorithms, but the author has found that for solving PDEs using Deep-TFC, the Quasi-Newton algorithms typically perform better. In this example, the L-BFGS algorithm was used to train the network. 

Due to the inherent stochasticity of the Deep-TFC method, the problem was solved ten times, and the best solution was retained; that solution had a maximum error of $2.780\times10^{-7}$ on the test set, a uniformly distributed $100\times100$ grid. The aforementioned stochasticity is easily visualized as a histogram: Figure \ref{fig:nnHist} shows such a histogram for 100 Monte Carlo trials.

\begin{figure}[H]
    \centering
    \includegraphics[width=\linewidth]{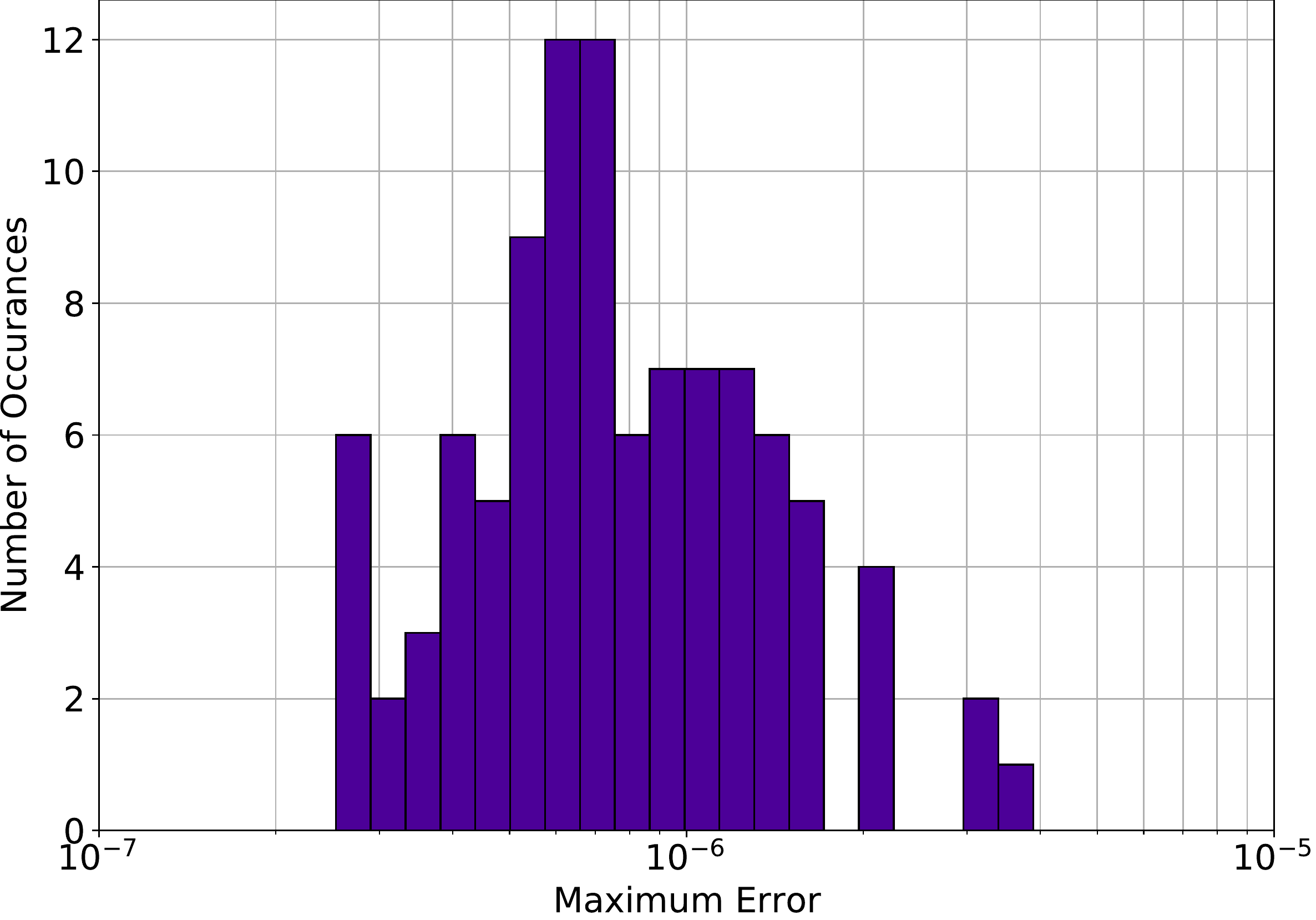}
    \caption{Histogram of the Deep-TFC maximum solution error on the test set for 100 Monte Carlo trials.}
    \label{fig:nnHist}
\end{figure}

Figure \ref{fig:nnHist} shows that Deep-TFC produces a solution at least as accurate as the solution reported earlier approximately $10\%$ of the time; this aligns well with one's intuition, as the reported solution was the best of ten trials. The remaining $90\%$ of the time the solution error is larger, but Figure \ref{fig:nnHist} shows that the Deep-TFC method is consistent: the maximum solution error in the 100 Monte Carlo trials was $3.891\times10^{-6}$, only an order of magnitude larger than the maximum solution error reported earlier \cite{Deep-TFC}.
\end{example}

\begin{example}{Simple PDE solved using X-TFC}\label{ex:SimplePdeXTfc}
Let $g(x,y)$ be an ELM as defined earlier with nonlinear activation function $\psi = \tanh$. Similar to the solution that used basis functions, Example \ref{ex:SimplePdeBasis}, selecting $g(x,y)$ in this way ultimately results in a loss function, $\mathbb{\B{L}}(W_2)$, that can be solved via linear least-squares. As in the previous examples, let the test set of points be a $100\times100$ uniform grid. Table \ref{tab:LagarisXtfc} shows the maximum test set solution error using X-TFC, where $m$ corresponds to the number of neurons in the hidden layer of the ELM: note that the number of neurons in each column of Table \ref{tab:LagarisXtfc} coincides with the number of basis functions in the corresponding columns of Table \ref{tab:LagarisCp} from Example \ref{ex:SimplePdeBasis}.
\begin{table}[H]
\caption{Maximum test set solution error using X-TFC with the $\tanh$ activation function.}\label{tab:LagarisXtfc}
\begin{tabular}{c|ccccc}
\noalign{\hrule height 1.0pt}
\hline
\diagbox[height=1.25cm]{\textbf{n}}{\textbf{m}}
 & \textbf{17} & \textbf{62} & \textbf{132} & \textbf{227} & \textbf{347}\\
\hline
5  & $1.74\times 10^{-5}$ & $\text{-}$   & $\text{-}$   & $\text{-}$   & $\text{-}$ \\
10 & $4.44\times 10^{-6}$ & $1.49\times 10^{-10}$ & $\text{-}$   & $\text{-}$   & $\text{-}$ \\
15 & $4.12\times 10^{-6}$ & $1.11\times 10^{-10}$ & $1.21\times 10^{-12}$ & $\text{-}$   & $\text{-}$ \\
20 & $3.95\times 10^{-6}$ & $5.80\times 10^{-11}$ & $4.40\times 10^{-13}$ & $2.37\times 10^{-13}$ &  $\text{-}$ \\
25 & $3.84\times 10^{-6}$ & $5.96\times 10^{-11}$ & $5.46\times 10^{-13}$ & $2.10\times 10^{-13}$ & $2.14\times 10^{-13}$ \\
30 & $3.77\times 10^{-6}$ & $5.50\times 10^{-11}$ & $5.55\times 10^{-13}$ & $1.66\times 10^{-13}$ & $1.83\times 10^{-13}$\\
\noalign{\hrule height 1.0pt}
\end{tabular}
\end{table}

\noindent Comparing Tables \ref{tab:LagarisCp} and \ref{tab:LagarisXtfc} reveals that when the number of basis functions is low, X-TFC outperforms TFC by as much as two orders of magnitude; however, as the number of basis functions increases, TFC outperforms X-TFC by as much as three orders of magnitude.

\begin{figure}[H]
    \centering
    \includegraphics[width=\linewidth]{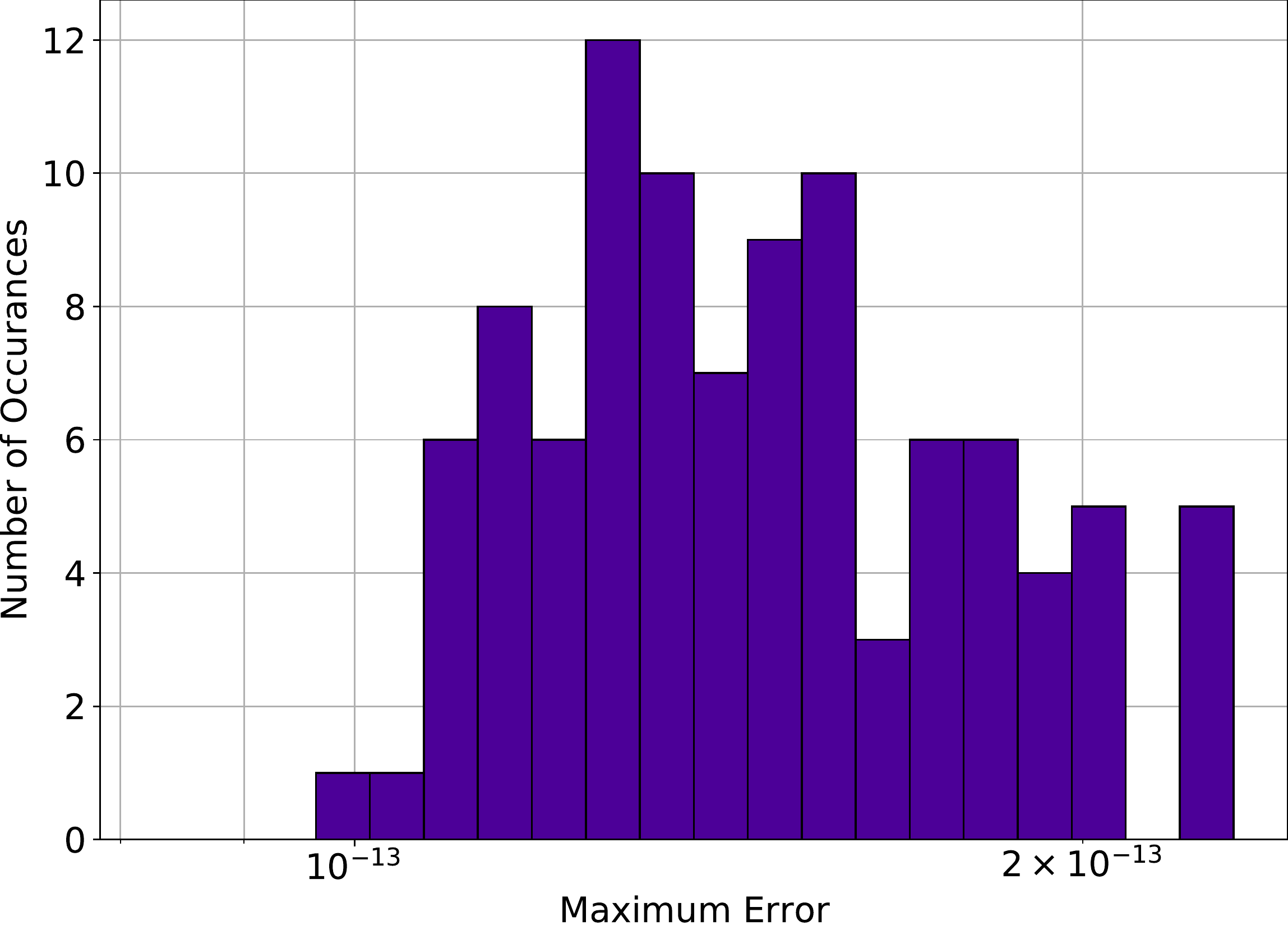}
    \caption{Histogram of X-TFC maximum solution error where $n$ = 30 and $m$ = 347 on the test set for 100 Monte Carlo trials.}
    \label{fig:elmHist}
\end{figure}

Similar to Deep-TFC, X-TFC is inherently stochastic. Hence, Figure \ref{fig:elmHist} shows a histogram of 100 Monte Carlo trials of the final case: $n$ = 30, $m$ = 347. Figure \ref{fig:elmHist} shows that the corresponding value given in Table \ref{tab:LagarisXtfc} is actually at the higher end of the distribution; however, it should be noted here that the deviation between Monte Carlo trials when using X-TFC is only as large as approximately $2\times10^{-13}$. This relative deviation is much smaller than when using Deep-TFC, where cases varied by as much as an order of magnitude.

The small relative difference between X-TFC cases can be attributed to the large number of neurons relative to the initial distribution. In each of the Monte Carlo trials, the weights and biases are chosen using $U(-1,1)$. Since there are $347$ neurons, the sample space is well represented each time. If fewer neurons were used, the relative difference between test cases would be larger. To illustrate, Figure \ref{fig:elmHistSmallM} shows a histogram of 100 Monte Carlo trials when $n$ = 30 and $m$ = 17. 

\begin{figure}[H]
    \centering
    \includegraphics[width=\linewidth]{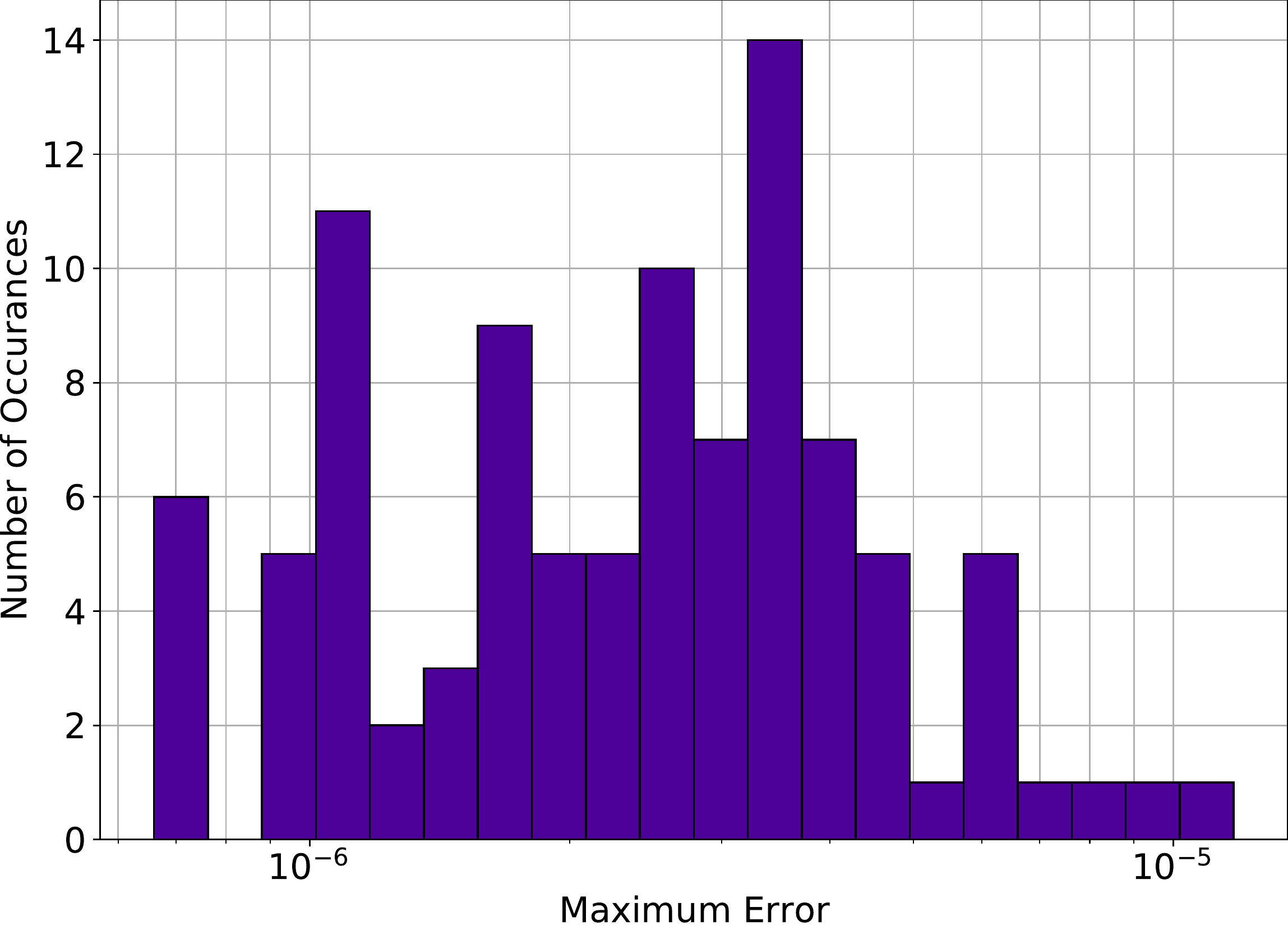}
    \caption{Histogram of X-TFC maximum solution error where $n$ = 30 and $m$ = 17 on the test set for 100 Monte Carlo trials.}
    \label{fig:elmHistSmallM}
\end{figure}

In Figure \ref{fig:elmHist}, the maximum test set error in the worst case was approximately twice as large as the maximum test set error in the best case, whereas in Figure \ref{fig:elmHistSmallM}, the worst case is approximately ten times as large as the best case. Of course, other factors---such as the nonlinear relationship between the trainable parameters of Deep-TFC versus the linear relationship between the trainable parameters of X-TFC---play a role in the differences between the histograms of Figures \ref{fig:nnHist} and \ref{fig:elmHist} as well.
\end{example}

\section{A Juxtaposition of TFC, CSVM, X-TFC, and Deep-TFC}

The previous simple PDE examples---Examples \ref{ex:SimplePdeBasis} through \ref{ex:SimplePdeXTfc}---highlight some of the strengths and weaknesses of the four free function options discussed earlier; this section analyzes those strengths and weaknesses further. To aid in that analysis, Table \ref{tab:SimplePdeComp} shows the maximum training and test set errors when using TFC and each of the four free function choices to solve the simple PDE. In addition, the solution errors of other state-of-the-art algorithms are included for reference: these algorithms include the well-known FEM, a neural-network-based method that analytically satisfies the constraints via a functional\footnote{This functional is multiplicative in nature whereas the TFC \ce\ is additive in nature. Moreover, the multiplicative functional cannot satisfy certain sets of constraints, nor does it have the mathematical guarantees that TFC \ces\ do.} \cite{OrigOdePde}, a Bernstein neural-network-based approach \cite{BNN}, and a Chebyshev neural-network-based approach \cite{CNN}. Table \ref{tab:SimplePdeComp} shows that the TFC methodology outperforms all the others in terms of accuracy on the training and test sets, followed by X-TFC. In terms of the test set error, these two methods are followed by CSVM, Deep-TFC, and then the other state-of-the-art algorithms. As for the training set error, FEM outperforms both Deep-TFC and CSVM, but the other state-of-the-art algorithms do not.

\begin{table}[!ht]
\centering
\caption{Comparison of maximum training set and test set errors between TFC methods and current state-of-the-art techniques.}
\label{tab:SimplePdeComp}
\begin{tabular}{|c|c|c|} 
\hline
\makecell{\textbf{Method}} & \makecell{\textbf{Training Set}\\\textbf{Maximum Error}} & \makecell{\textbf{Test Set}\\\textbf{Maximum Error}}\\\hline
{TFC \cite{M-TFC2}} & {$2.22 \times 10^{-16}$} & {$4.44 \times 10^{-16}$} \\\hline
{X-TFC \cite{X-TFC}} & {$3.8 \times 10^{-13}$} & {$5.1 \times 10^{-13}$} \\\hline
{CSVM \cite{SVM-TFC}} & {$4.4 \times 10^{-8}$} & {$5.6 \times 10^{-8}$} \\\hline
{Deep-TFC \cite{Deep-TFC}} & {$2.7 \times 10^{-7}$} & {$2.8 \times 10^{-7}$} \\\hline
{FEM \cite{OrigOdePde}} & {$2 \times 10^{-8}$} & {$1.5 \times 10^{-5}$} \\\hline
{NN \cite{OrigOdePde}} & {$5 \times 10^{-7}$} & {$5 \times 10^{-7}$} \\\hline
{Bernstein NN \cite{BNN}} & {\text{-}} & {$2.4 \times 10^{-4}$} \\\hline
{Chebyshev NN \cite{CNN}} & {\text{-}} & {$3.2 \times 10^{-2}$} \\\hline
\end{tabular}
\end{table}

As mentioned earlier, the CSVM technique is no longer actively being used as a free function choice because it requires a complex analytical analysis for each new differential equation, and the resultant payoff in terms of solution error is overshadowed by the other free function choices. Example \ref{ex:SimplePdeCsvm} demonstrates well the aforementioned complex analytical analysis: the long expressions containing multiple $4$-th and $5$-th order tensors in the example are daunting enough, but the author reminds readers that these expressions are compact and simple compared to the five or so pages of work it took to derive them. Furthermore, Table \ref{tab:SimplePdeComp} clearly shows the accuracy gained when using TFC or X-TFC rather than CSVM. In addition, because CSVM requires a grid search to find the two hyperparameters, the training time is longer than when using TFC or X-TFC. Although this free function choice has become antiquated, it laid the foundation for the synergy between TFC and machine learning algorithms and is therefore historically significant. 

Looking at the error values alone in Table \ref{tab:SimplePdeComp}, the reader may wonder why Deep-TFC is useful. After all, one of the arguments against CSVM was the error relative to TFC and X-TFC; however, on more complex problems, Deep-TFC actually does better than X-TFC and TFC. Moreover, Deep-TFC does not require the same complex analytical analysis that CSVM did. To highlight Deep-TFC's performance on complex problems, consider low-speed, two-dimensional, developing channel flow governed by the Navier-Stokes equations and the following boundary conditions:
\begin{equation*}
\begin{aligned}
    &\Part{u}{x}+\Part{v}{y} = 0\\
    &\rho \bigg(\Part{u}{t}+u\Part{u}{x}+v\Part{u}{y}\bigg) = -\Part{P}{x}+\mu \bigg(\PartS{u}{x}+\PartS{u}{y}\bigg)\\
    &\rho \bigg(\Part{v}{t}+u\Part{v}{x}+v\Part{v}{y}\bigg) = \mu \bigg(\PartS{v}{x}+\PartS{v}{y}\bigg)\\
    &\text{subject to}\quad \begin{cases}
    &u(0,y,t) = \Part{u}{x}(L,y,t) = u(x,y,0) = 0\\
    &u(x,\frac{H}{2},t) = u(x,-\frac{H}{2},t) = 0\\
    &v(0,y,t)=\Part{v}{x}(L,y,t)=v(x,y,0)=0\\
    &v(x,\frac{H}{2},t)=v(x,-\frac{H}{2},t)=0,
    \end{cases}
\end{aligned}
\end{equation*}
where $u$ and $v$ are velocities in the $x$ and $y$ directions respectively, $H$ is the height of the channel, $P$ is the pressure, $\rho$ is the density, and $\mu$ is the viscosity. For this problem, the values $H=1$ m, $\rho=1$ kg/m$^3$, $\mu=1$ Pa$\cdot$s, and $\Part{P}{x}=-5$ N/m$^3$ were chosen. 

The $u$ and $v$ dependent variables each have the same constraints; therefore, their constrained expressions are the same. Hence, just the constrained expression for $u$ will be shown. In recursive form, the \ce\ for $u$ is,
\begin{align*}
    \p{1}{u}(x,&y,t,g^u(x,y,t)) = g^u(x,y,t)-g^u(0,y,t)-xg_x^u(L,y,t) \\
    \p{2}{u}(x,&y,t,g^u(x,y,t)) = g^u(x,y,t)-\frac{H-2y}{2H}g^u\Big(x,-\frac{H}{2},t\Big)-\frac{H+2y}{2H}g^u\Big(x,\frac{H}{2},t\Big) \\
    \p{3}{u}(x,&y,t,g^u(x,y,t)) = g^u(x,y,t)-g^u(x,y,0),
\end{align*}
where $\p{1}{u}$, $\p{2}{u}$, and $\p{3}{u}$ can be processed in any order to produce the full \ce. In tensor form, the \ce\ is,
\begin{equation*}
    u(x,y,t,g^u(x,y,t)) = g^u(x,y,t)+\mathcal{M}(x,y,t,g^u(x,y,t))_{ijk}\Phi_i(x)\Phi_j(y)\Phi_k(t)
\end{equation*}
where,
\begin{align*}
    \mathcal{M}_{ij1}(x,y,t,g^u(x,y,t)) &= \begin{bmatrix} 0 & -g^u(x,-\frac{H}{2},t) & -g^u(x,\frac{H}{2},t) \\
    -g^u(0,y,t) & g^u(0,-\frac{H}{2},t) & g^u(0,\frac{H}{2},t) \\ 
    -g^u_x(L,y,t) & g^u_x(L,-\frac{H}{2},t) & g^u_x(L,\frac{H}{2},t) \end{bmatrix}\\
    \mathcal{M}_{ij2}(x,y,t,g^u(x,y,t)) &= \begin{bmatrix} -g^u(x,y,0) & g^u(x,-\frac{H}{2},0) & g^u(x,\frac{H}{2},0) \\
    g^u(0,y,0) & -g^u(0,-\frac{H}{2},0) & -g^u(0,\frac{H}{2},0) \\ 
    g^u_x(L,y,0) & -g^u_x(L,-\frac{H}{2},0) & -g^u_x(L,\frac{H}{2},0) \end{bmatrix}
\end{align*}
and
\begin{equation*}
    \Phi_i(x) = \begin{Bmatrix} 1, & 1, & x\end{Bmatrix}, \quad \Phi_j(y) = \begin{Bmatrix} 1, & \frac{H-2y}{2H}, & \frac{H+2y}{2H}\end{Bmatrix}, \quad \Phi_k(t) = \begin{Bmatrix} 1, & 1 \end{Bmatrix}.
\end{equation*}

For Deep-TFC, the training set used was $2,000$ independently and identically distributed (i.i.d.) points sampled from $x\in U(0,15)$, $y\in U(-H/2,H/2)$, and $t\in U(0,3)$. For X-TFC and TFC, the training set was a grid of $10\times10\times10$ uniformly spaced points. For each method, the test set consisted of a grid of $100\times100$ evenly spaced points in $x$ and $y$ at three different times: $t=0.01$, $t=0.1$ and $t=3$. The test set can be analyzed in two different ways: 
\begin{enumerate}
    \item Qualitatively - The solution should be symmetric about the line $y=0$, and the solution should develop spatially and temporally such that after a sufficient amount of time has passed and sufficiently far from the inlet, $x=0$, the $u$-velocity will be equal, or very nearly equal, to the steady-state Poiseuille flow solution.
    \item Quantitatively - The solution at $x=15$ and $t=3$ can be compared to the steady-state Poiseuille flow solution.
\end{enumerate}

The neural network used for the Deep-TFC solution had four hidden layers and 30 neurons per layer, and the nonlinear activation function used was the hyperbolic tangent \cite{Deep-TFC}. The X-TFC solution used 200 neurons and the hyperbolic tangent as the nonlinear activation function; adding additional neurons, up to 300, did not improve the solution over the case with 200 neurons. The TFC solution used Chebyshev orthogonal polynomials up to degree ten; adding additional polynomials, up to degree 15, did not improve the solution over the case with polynomials up to degree ten.

\begin{figure}[!hbtp]
    \centering
    \begin{minipage}[t]{0.5\linewidth}
       \centering\includegraphics[width=\linewidth]{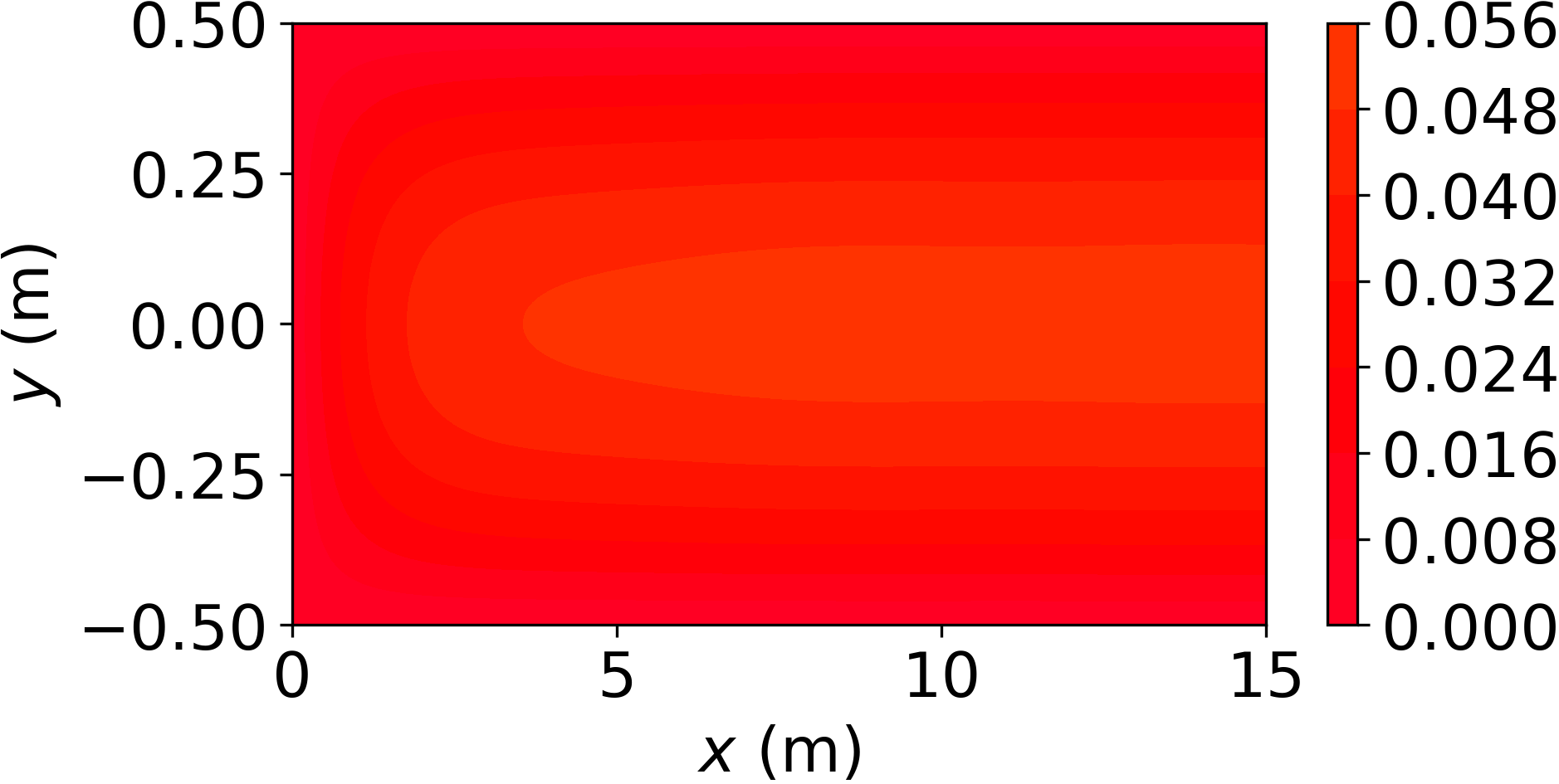}
        \caption{TFC solution at $t=0.01$.}
       \label{fig:NsTFC0}
    \end{minipage}%
    \begin{minipage}[t]{0.5\linewidth}
       \centering\includegraphics[width=\linewidth]{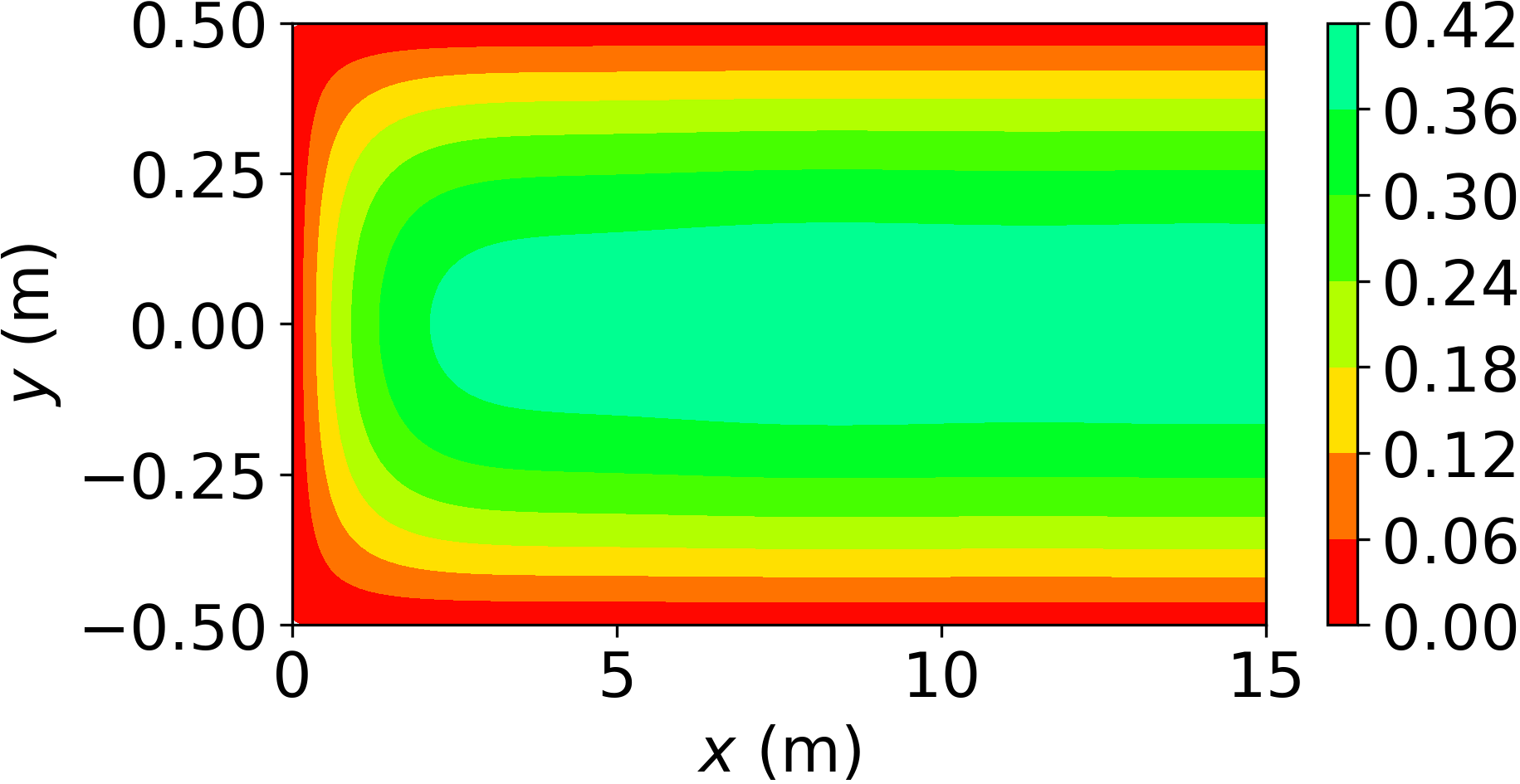}
        \caption{TFC solution at $t=0.1$.}
    \end{minipage}%
    \allowdisplaybreaks\newline%
    \begin{minipage}[t]{0.5\linewidth}
       \centering\includegraphics[width=\linewidth]{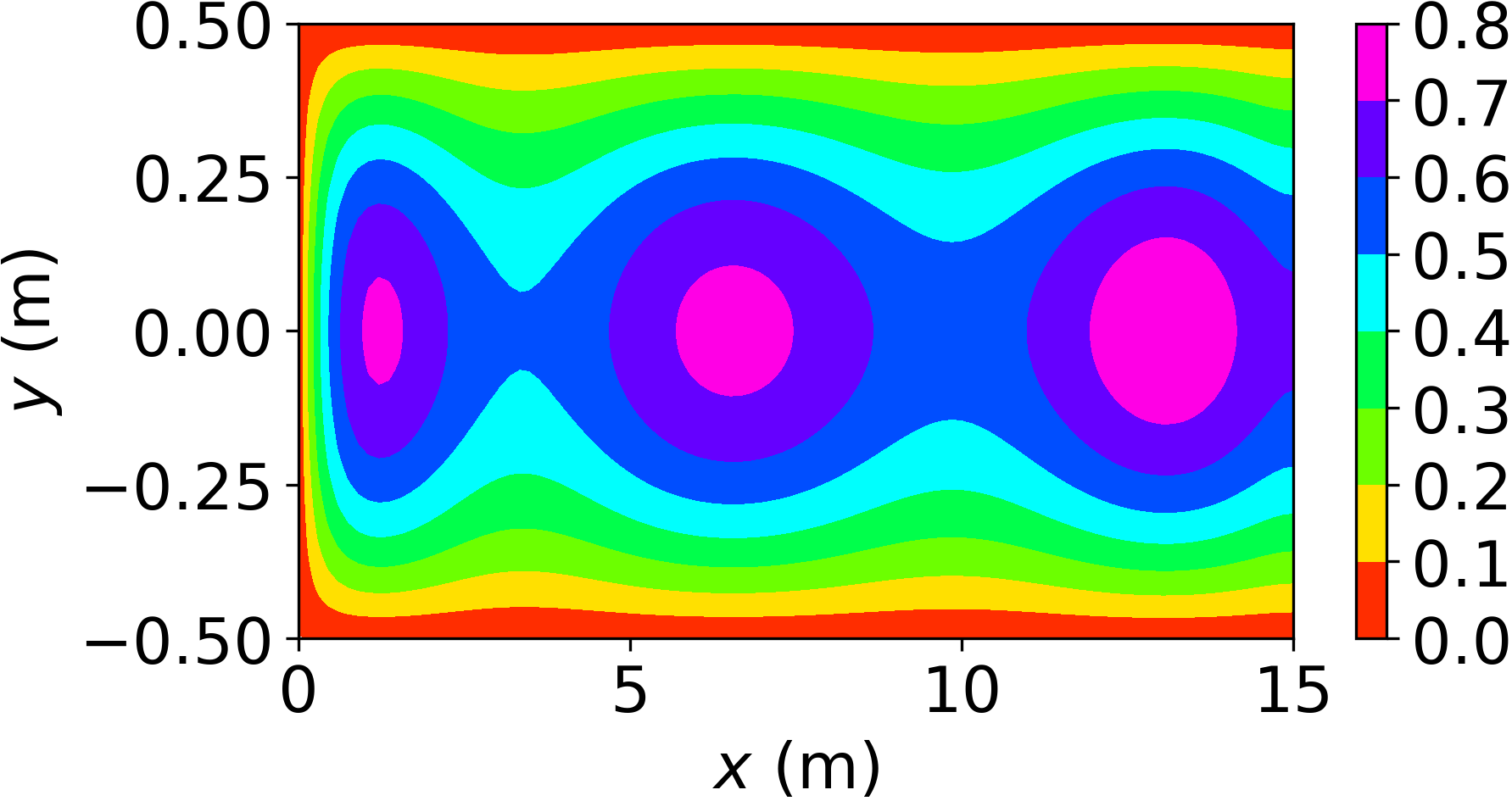}
        \caption{TFC solution at $t=3.0$.}
        \label{fig:NsTFC2}
    \end{minipage}%
\end{figure}

\begin{figure}[!hbtp]
    \centering
    \begin{minipage}[t]{0.5\linewidth}
       \centering\includegraphics[width=\linewidth]{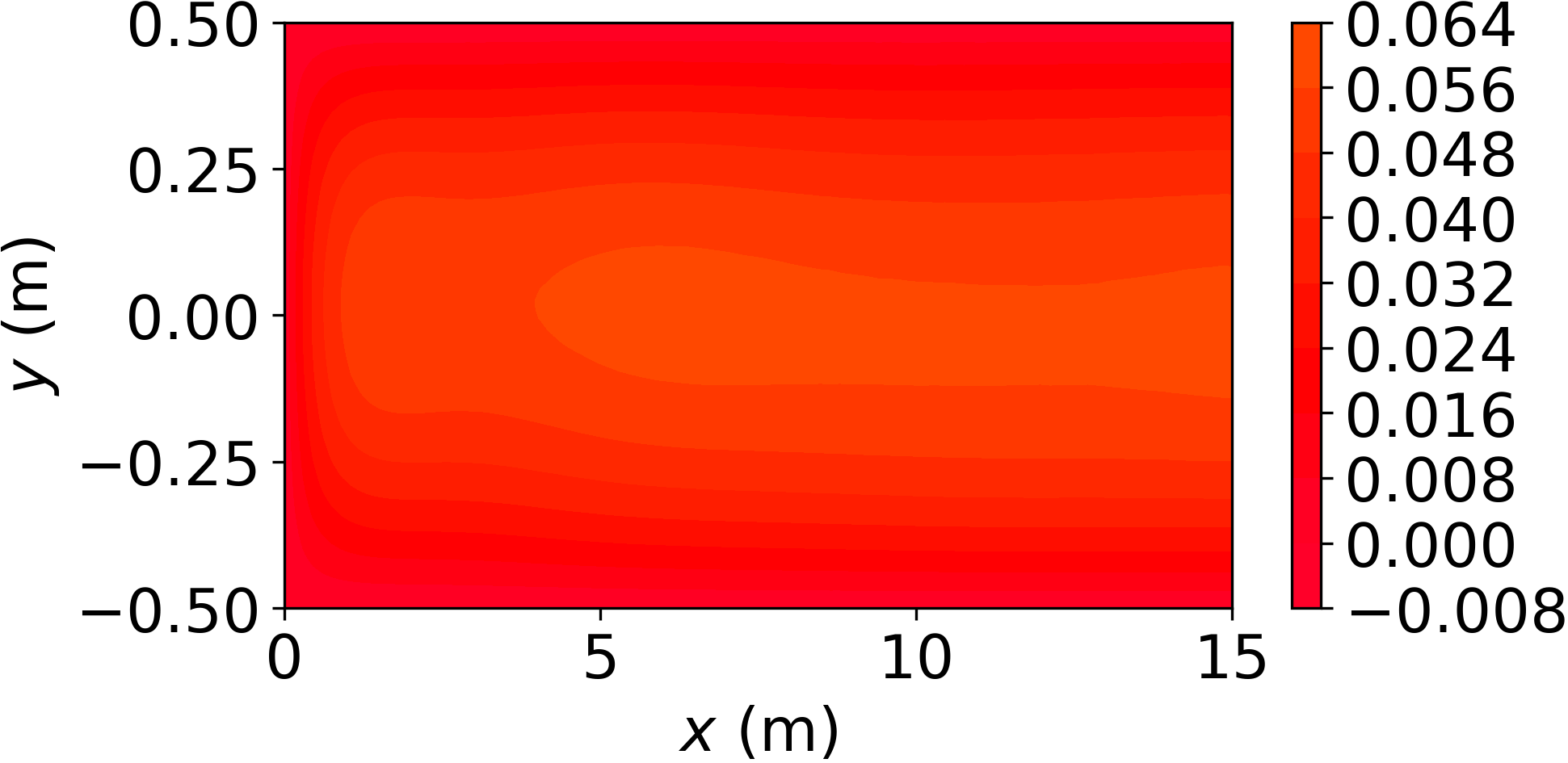}
        \caption{X-TFC solution at $t=0.01$.}
       \label{fig:NsXTFC0}
    \end{minipage}%
    \begin{minipage}[t]{0.5\linewidth}
       \centering\includegraphics[width=\linewidth]{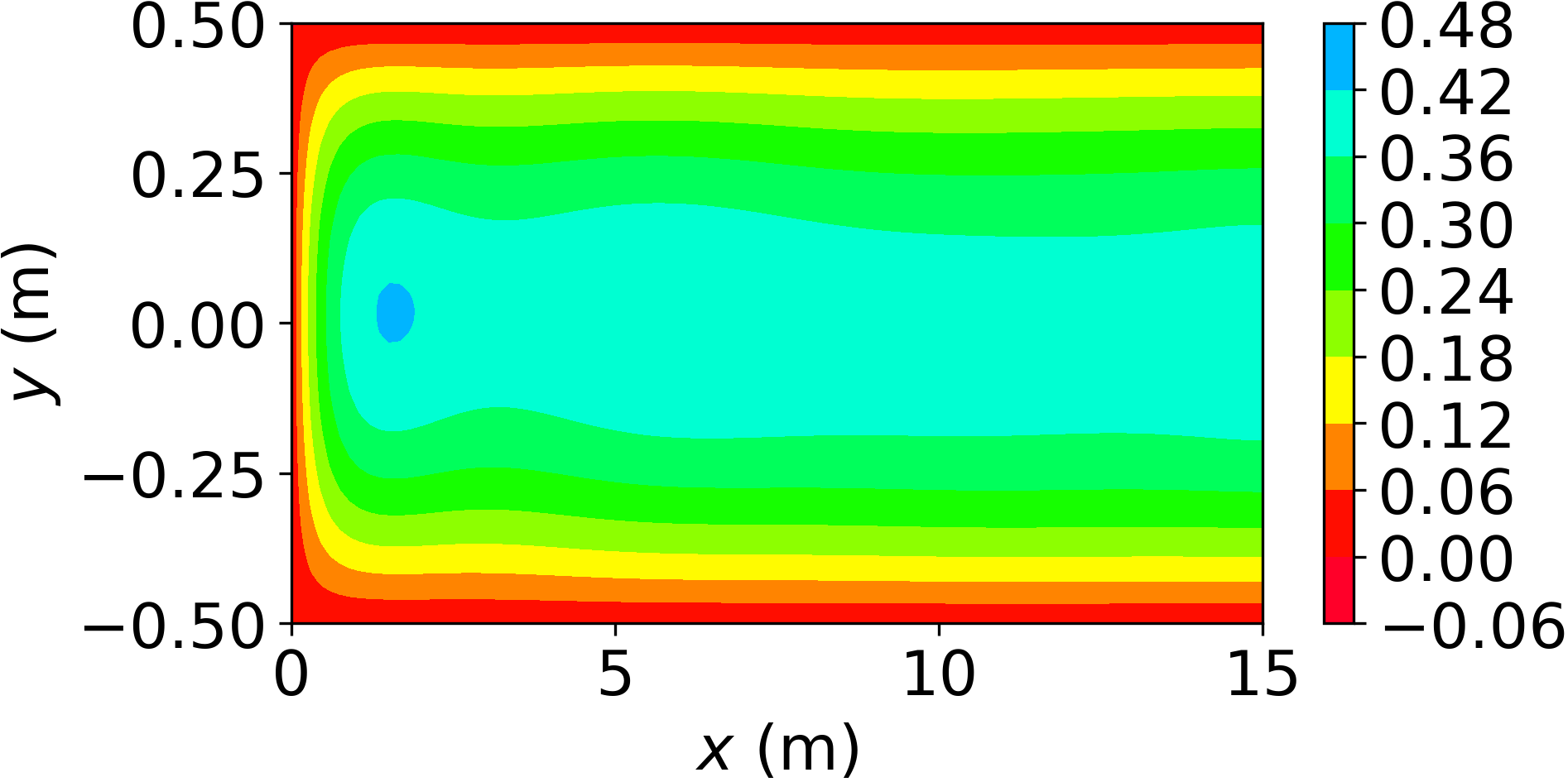}
        \caption{X-TFC solution at $t=0.1$.}
    \end{minipage}%
    \allowdisplaybreaks\newline%
    \begin{minipage}[t]{0.5\linewidth}
       \centering\includegraphics[width=\linewidth]{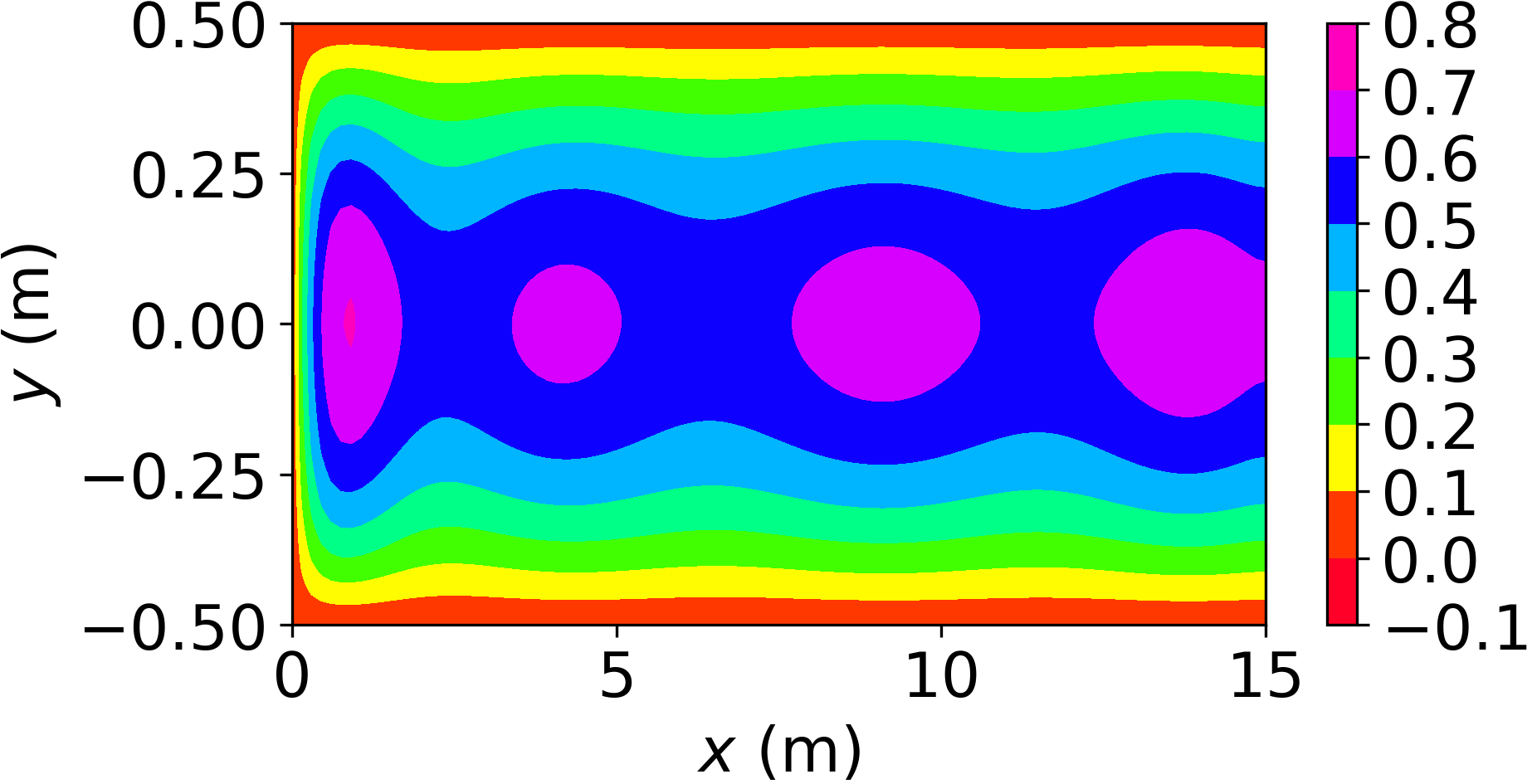}
        \caption{X-TFC solution at $t=3.0$.}
        \label{fig:NsXTFC2}
    \end{minipage}%
\end{figure}

\begin{figure}[!hbtp]
    \centering
    \begin{minipage}[t]{0.5\linewidth}
       \centering\includegraphics[width=\linewidth]{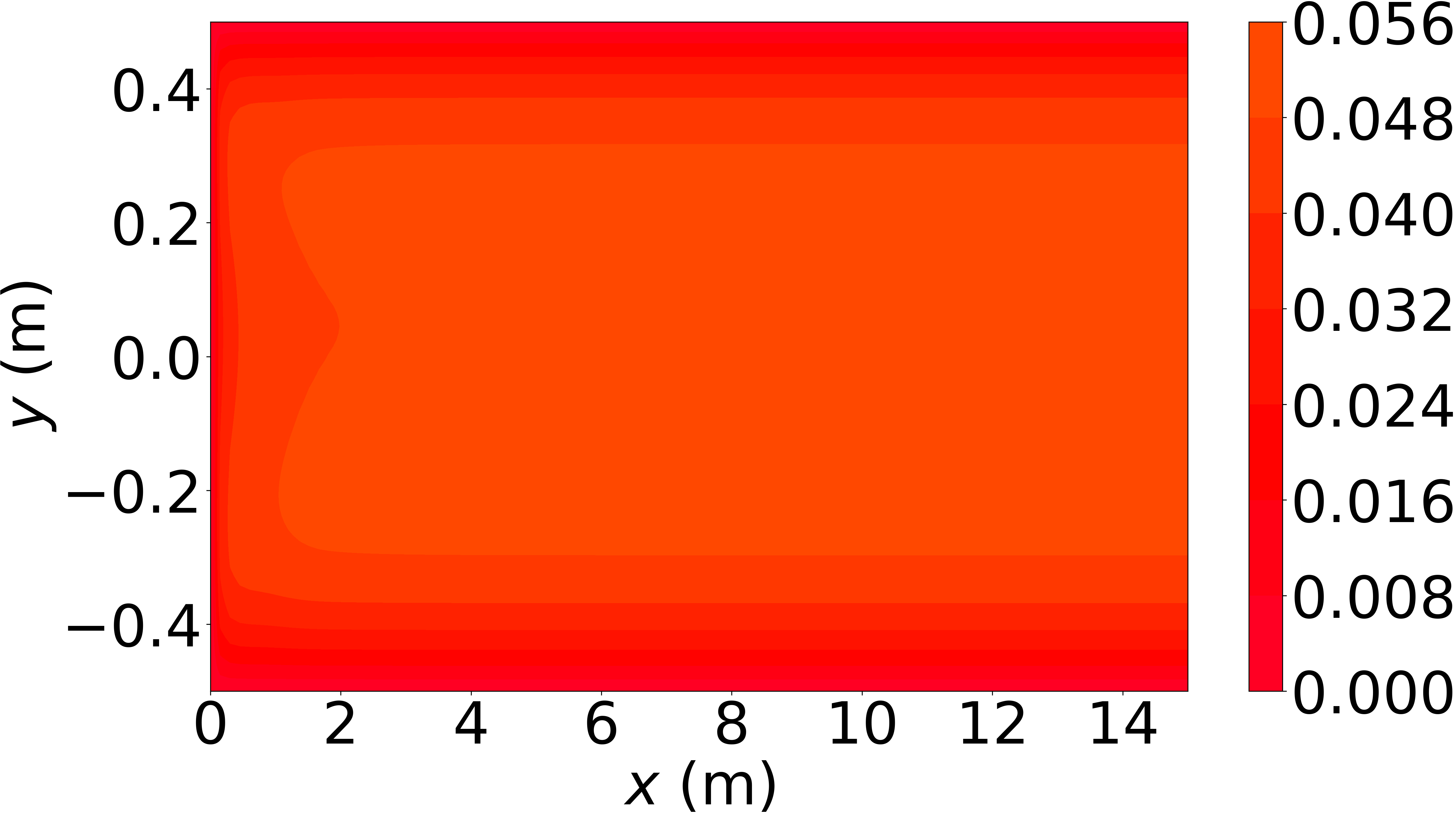}
        \caption{Deep-TFC solution at $t=0.01$.}
       \label{fig:NsDTFC0}
    \end{minipage}%
    \begin{minipage}[t]{0.5\linewidth}
       \centering\includegraphics[width=\linewidth]{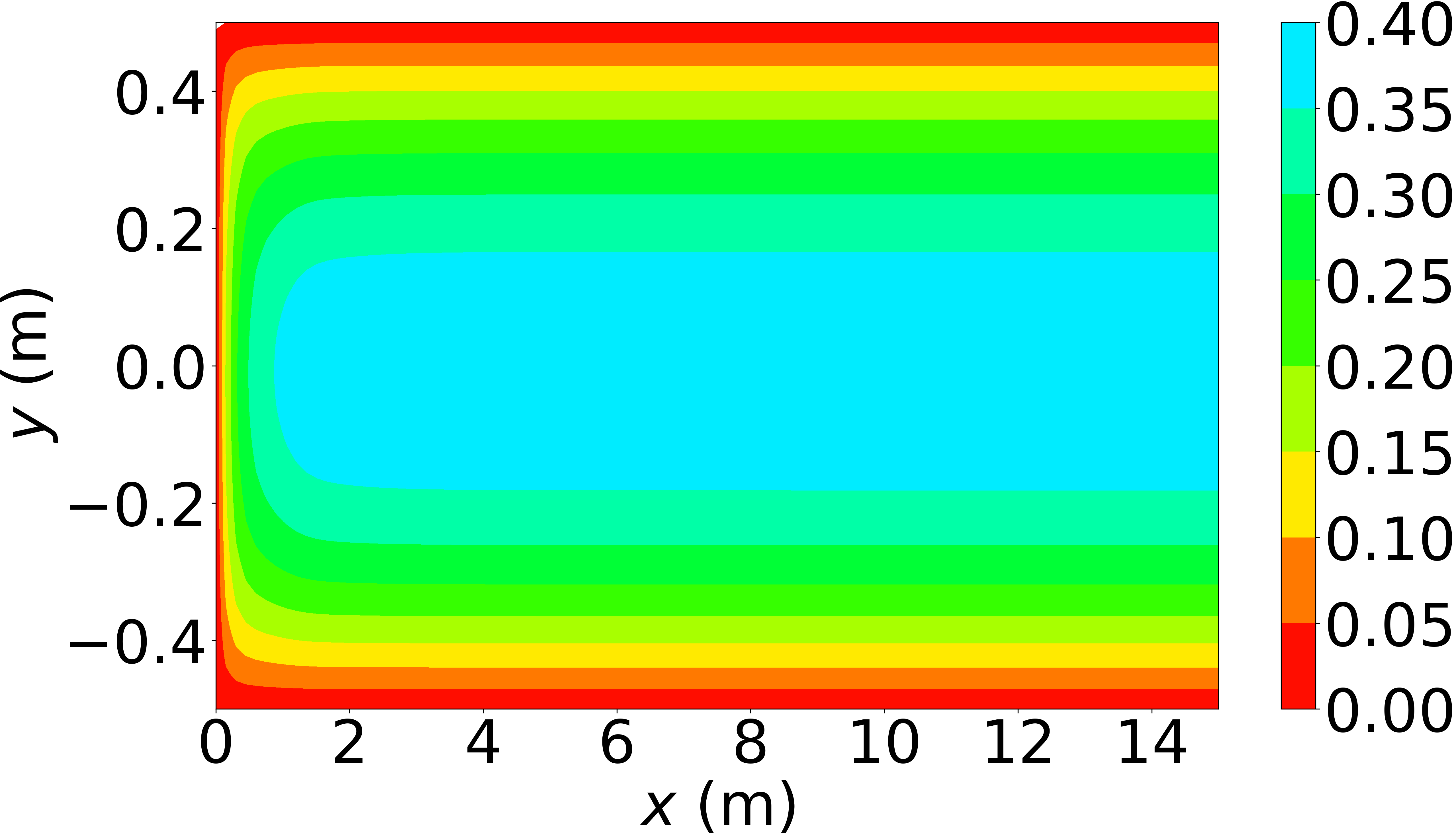}
        \caption{Deep-TFC solution at $t=0.1$.}
    \end{minipage}%
    \allowdisplaybreaks\newline%
    \begin{minipage}[t]{0.5\linewidth}
       \centering\includegraphics[width=\linewidth]{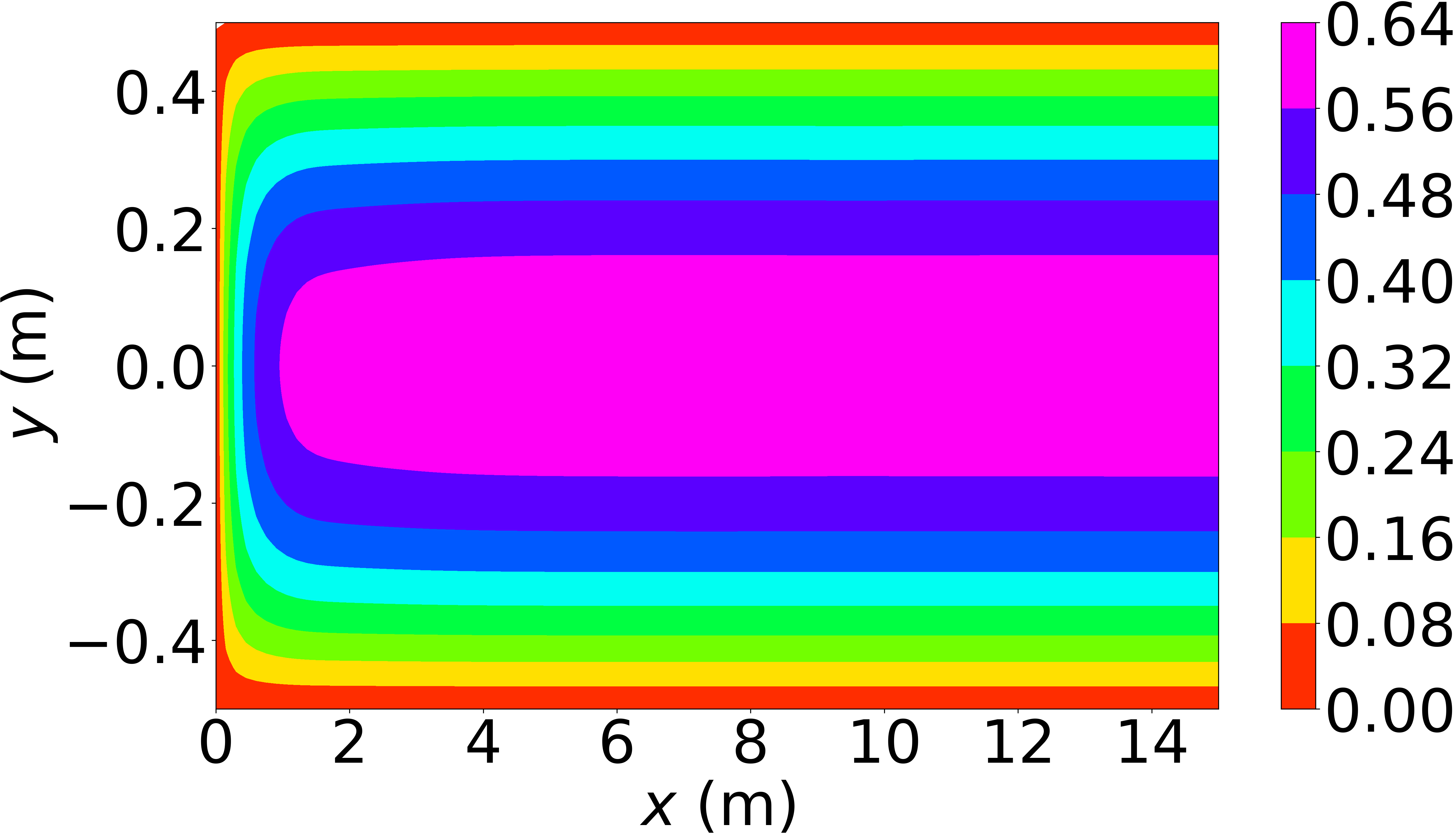}
        \caption{Deep-TFC solution at $t=3.0$.}
        \label{fig:NsDTFC2}
    \end{minipage}%
\end{figure}

A quantitative comparison of the three methods' errors with respect to the steady-state Poiseuille flow solution at $x=15$ and $t=3$ is shown in Table \ref{tab:NsComp}. In addition, a qualitative comparison is illustrated via Figures \ref{fig:NsTFC0} through \ref{fig:NsDTFC2}: Figures \ref{fig:NsTFC0} through \ref{fig:NsTFC2} correspond to the TFC solution, Figures \ref{fig:NsXTFC0} through \ref{fig:NsXTFC2} correspond to the X-TFC solution, and Figures \ref{fig:NsDTFC0} through \ref{fig:NsDTFC2} correspond to the Deep-TFC solution.

\begin{table}[!hbtp]
\centering
\caption{Comparison of maximum and mean test set errors between TFC, X-TFC, and Deep-TFC.}
\label{tab:NsComp}
\begin{tabular}{|c|c|c|} 
\hline
\makecell{\textbf{Method}} & \makecell{\textbf{Test Set}\\\textbf{Maximum Error}} & \makecell{\textbf{Test Set}\\\textbf{Average Error}}\\\hline
{TFC} & {$5.59 \times 10^{-3}$} & {$3.68 \times 10^{-3}$} \\\hline
{X-TFC} & {$4.02 \times 10^{-3}$} & {$1.89 \times 10^{-3}$}\\\hline
{Deep-TFC \cite{Deep-TFC}} & {$5.38 \times 10^{-4}$} & {$3.12 \times 10^{-4}$} \\\hline
\end{tabular}
\end{table}

Table \ref{tab:NsComp} shows that the X-TFC solution does marginally better than the TFC solution in terms of error, and the Deep-TFC solution is approximately an order of magnitude better than X-TFC and TFC. This is reflected qualitatively in the figures as well. 

The Deep-TFC figures, Figures \ref{fig:NsDTFC0} through \ref{fig:NsDTFC2}, match the qualitative expectation given earlier. In contrast, the TFC and X-TFC solutions' figures do not: this difference is highlighted most in figures for $t=3.0$. 

In summation, this Navier-Stokes example demonstrates the utility of Deep-TFC as problems become sufficiently complex. These results coupled with those of Table \ref{tab:SimplePdeComp} make good arguments for using Deep-TFC for complex problems and TFC for simpler problems. What about X-TFC?

For certain problems, X-TFC outperforms TFC: the two-dimensional wave equation in Section \ref{sec:2dWaveEqn} is one example. Moreover, comparing the solution errors of the two methods---see Table \ref{tab:waveEquation2DComparison} or compare the first two columns of Tables \ref{tab:LagarisCp} and \ref{tab:LagarisXtfc}---shows that the difference between the two methods is the most significant for a low number of trainable parameters. This is another benefit of the X-TFC framework: lower solution error than other methods when a lower number of parameters is used. 

Hence, for multidimensional problems in general, if a problem has a complex residual, use Deep-TFC; if the problem is simple, use X-TFC or TFC, and if the application is memory restrictive, i.e., a lower number of trainable parameters is required, use X-TFC. Oftentimes X-TFC or TFC is sufficient to estimate the solution, and the user is encouraged to try out both frameworks for their specific problem. In fact, the numerical implementation makes it extremely simple to switch between these two methods.

%% file: Data/FlexBodyApplications.tex

\chapter{APPLICATIONS IN FLEXIBLE BODY PROBLEMS}\label{chap:flexBodyApplications}

The previous chapter explained how to apply TFC to differential equations. Consequently, TFC can be applied to a wide variety of problems spanning multiple fields and multiple disciplines within those fields. However, the author is particularly interested in some of the differential equations appearing in flexible body problems, and therefore, this chapter is dedicated to them exclusively. The flexible body problems solved in this chapter include:
\begin{itemize}
    \item Natural tandem balloon shape - A set of four, coupled, first-order, nonlinear ODEs wherein both ends of the domain are themselves unknowns that must be solved simultaneously alongside the ODEs. 
    \item One-dimensional wave equation - A two-dimensional (one space, one time), second-order, linear PDE.
    \item Two-dimensional wave equation - A three-dimensional (two space, one time), second-order, linear PDE. 
    \item Biharmonic equation, Cartesian coordinates - A two-dimensional, fourth-order, linear PDE.
    \item Biharmonic equation, polar coordinates - A two-dimensional, fourth-order, linear PDE.
\end{itemize}

\section{Natural Tandem Balloon Shape}
\begin{figure}[!htbp]
    \centering
    \includegraphics[width=0.5\linewidth]{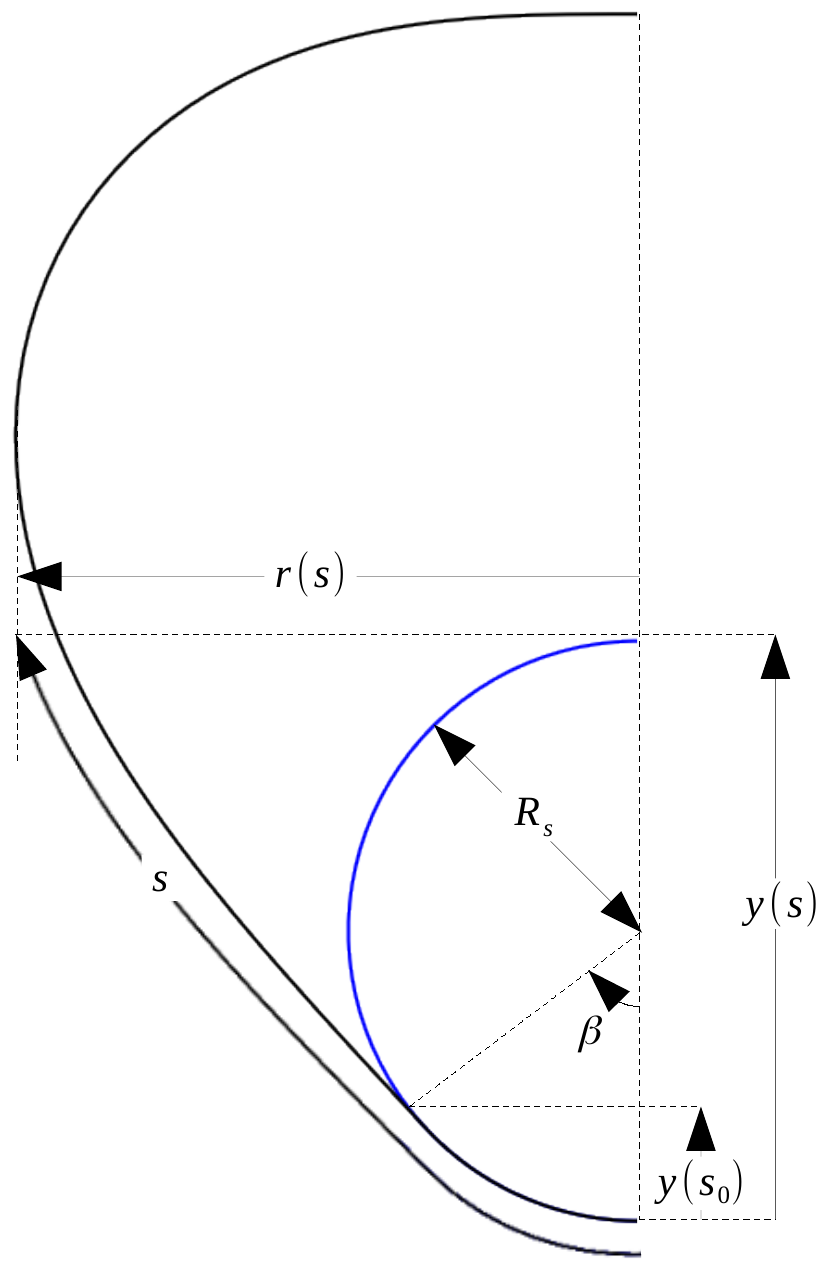}
    \caption{Tandem balloon diagram.}
    \label{fig:TandemBalloonDiagram}
\end{figure}

Tandem balloons are useful scientific vessels for collecting terrestrial atmospheric data and are being considered for the same task on other planets and moons \cite{BalloonVenusInit}. Naturally, the shape of these balloons is a critical component involved in simulating their trajectories. This section will provide a general overview of the problem \cite{JonathanCameronNaturalBalloonShapes}; a more detailed description can be found in References \cite{NaturalBalloonShapeBaginski} and \cite{BalloonSmalley}. Figure \ref{fig:TandemBalloonDiagram} is a diagram of the tandem balloon and the coordinate system used to describe its shape. Table \ref{tab:BalloonNomenclature} provides the nomenclature used to describe the tandem balloon.
\begin{table}[!ht]
    \centering
    \caption{Tandem balloon nomenclature.}
    \begin{tabularx}{\linewidth}{|c|X|}
        \hline
        Symbol & Description \\\hline
        $A_s$ & Surface area of the super-pressure balloon\\\hline
        $M_{atm}$ & Molecular weight of the atmosphere \\\hline
        $M_g$ & Molecular weight of the lifting gas \\\hline
        $R_s$ & Radius of the super-pressure balloon\\\hline
        $T_0$ & Total vertical load felt at $s_0$ \\\hline
        $V_s$ & Volume of the super-pressure balloon\\\hline
        $b$ & Specific density of the lifting gas \\\hline
        $g$ & Acceleration due to gravity \\\hline
        $m_{sg}$ & Total mass of the gas in the super pressure balloon\\\hline
        $r$ & Coordinate that describes position perpendicular to axisymmetric line \\\hline
        $s$ & Coordinate that describes position along the balloon film\\\hline
        $s_0$ & Point where the zero-pressure balloon and super-pressure balloon come into contact \\\hline
        $w$ & Zero-pressure balloon film mass per unit area\\\hline
        $w_s$ & Super-pressure balloon film mass per unit area\\\hline
        $y$ & Coordinate that describes position parallel to axisymmetric line\\\hline
        $\beta$ & Angle measured from the center of the super-pressure balloon between the vertical and the point where the zero-pressure balloon and super-pressure balloon surfaces diverge\\\hline    
        $\ell_d$ & Length of the balloon film \\\hline
        $\rho$ & Atmospheric density \\\hline
        $\sigma_c$ & Circumferential stress\\\hline
        $\sigma_m$ & Meridional stress\\\hline
        $\theta$ & Angle in $[-\frac{\pi}{2},\frac{\pi}{2}]$ between the vertical and a line tangent to the balloon's surface \\\hline
    \end{tabularx}
    \label{tab:BalloonNomenclature}
\end{table}

The differential equations that govern the balloon are summarized in Equation \eqref{eq:BalloonDesOrig} \cite{NaturalBalloonShapeBaginski,BalloonSmalley}.
\begin{equation}\label{eq:BalloonDesOrig}
\begin{aligned}
    \frac{\dd \theta}{\dd s} &= \frac{1}{\sigma_m} \Big(\frac{\sigma_c}{r}\cos(\theta)-w \sin(\theta)-b\big(y-y(s_0)\big)\Big)\\
    \frac{\dd \sigma_m}{\dd s} &= \frac{\sigma_c}{r}\sin(\theta)+ w\cos(\theta)-\frac{\sigma_m}{r}\sin(\theta)\\
    \frac{\dd r}{\dd s} &= \sin(\theta)\\
    \frac{\dd y}{\dd s} &= \cos(\theta),
\end{aligned}
\end{equation}
subject to the boundary constraints,
\begin{align*}
    s_0 &= R_s\beta\\
    \theta(s_0) &= \frac{\pi}{2}-\beta\\
    \sigma_m(s_0) &= \frac{T_0}{2\pi r \sin(\theta)}\\
    r(s_0) &= R_s\sin(\beta)\\
    y(s_0) &= R_s(1-\cos(\beta))\\
    \theta(\ell_d) &= -\frac{\pi}{2} \\
    r(\ell_d) &= 0,
\end{align*}
where 
\begin{align*}
    &T_0 = L+g(w+w_s)A_{s0}+g\Big(\frac{V_{s0}}{V_s}m_{sg}-\rho V_{s0}\Big)\\
    &\text{if }\beta < \frac{\pi}{2}\rightarrow\begin{cases} A_{s0} = 2\pi R_s y(s_0)\\ V_{s0} = \frac{\pi}{3}z^2(s_0)\Big(3R_s-y(s_0)\Big)\end{cases}\\
    &\text{if }\beta \geq \frac{\pi}{2}\rightarrow\begin{cases} h_0 = 2R_s-y(s_0)\\ A_{s0} = A_s-2\pi R_s h_0\\ V_{s0} = V_s - \frac{\pi}{3}h_0\Big(3R_s-h_0\Big)\end{cases}\\
    A_s &= 4\pi R_s^2 \\
    V_s &= \frac{4}{3}\pi R_s^3
\end{align*}
and
\begin{align*}
   b &= g\rho\Big(1-\frac{M_g}{M_{atm}}\Big).
\end{align*}

The simplest version of these differential equations is the natural balloon shape, which has zero circumferential stress, i.e., $\sigma_c = 0$. Even with this simplification, the set of coupled, nonlinear differential equations is challenging to solve as the problem domain, $s\in[s_0,\ell_d]$, is variable on both ends: $\beta$ and  $\ell_d$ are variables to be solved alongside the differential equations, i.e., $s_0$ and $\ell_d$ are unknown. However, since TFC must map the free function domain to the problem domain anyway, the mapping parameter can be used in the least-squares when reducing the residual. 

The meridional stress generates another complication, as $r(\ell_d) = 0$ and $\pm\frac{\dd \sigma_m}{\dd s}\rightarrow\infty$ as $r\rightarrow0$ for non-zero $\sigma_m$, where the sign, $\pm$, depends on the sign of $\sigma_m$. Of course, this singularity does not exist in real life, as an infinite stress would rip the balloon apart, rather, it is due to the coordinate system chosen and assumptions made when deriving the differential equations. Fortunately, there is a change of variables that prevents a singularity in the dependent variables \cite{NaturalBalloonShapeBaginski}. Let $q = \frac{1}{\sigma_m r}$, then, Equation \eqref{eq:BalloonDesOrig} can be rewritten as,
\begin{align*}
    \frac{\dd \theta}{\dd s} &=  q\sigma_c\cos(\theta)-qrw \sin(\theta)-qrb(y-y(s_0))\\
    \frac{\dd q}{\dd s} &= -q^2\Big(\sigma_c\sin(\theta)+ wr\cos(\theta)\Big)\\
    \frac{\dd r}{\dd s} &= \sin(\theta)\\
    \frac{\dd y}{\dd s} &= \cos(\theta).
\end{align*}

The \ces\ that embed the boundary constraints given previously are shown in Equation \eqref{eq:BalloonCes}. Note that theses \ces\ are written for the domain of the free function, because as mentioned earlier, the differential equations will ultimately be written in terms of the free function domain, $z$, so the mapping parameter can be used to solve for $\beta$ and $\ell_d$. For this problem, Chebyshev orthogonal polynomials are used, which have a domain of $z\in[-1,1]$. Further, note that the domain is the only part of the problem being modified, so the right-hand side of the boundary conditions remains unchanged, e.g., $\theta(z=-1) = \theta(s_0)$.
\begin{equation}\label{eq:BalloonCes}
\begin{aligned}
    \theta(z,g^\theta(z)) &= g^\theta(z) + \frac{1-z}{2}\Big(\theta(s_0)-g^\theta(-1)\Big) + \frac{z+1}{2}\Big(\theta(\ell_d)-g^\theta(1)\Big)\\
    r(z,g^r(z)) &= g^r(z) + \frac{1-z}{2}\Big(r(s_0)-g^r(-1)\Big) + \frac{z+1}{2}\Big(r(\ell_d)-g^r(1)\Big)\\
    q(z,g^q(z)) &= g^q(z) + q(s_0)-g^q(-1) \\
    y(z,g^y(z)) &= g^y(z) + y(s_0)-g^y(-1)
\end{aligned}
\end{equation}

The differential equations written on this domain are,
\begin{align*}
    c\frac{\dd \theta}{\dd z} &=  q\sigma_c\cos(\theta)-qrw \sin(\theta)-qrb(y-y(s_0))\\
    c\frac{\dd q}{\dd z} &= -q^2\Big(\sigma_c\sin(\theta)+ wr\cos(\theta)\Big)\\
    c\frac{\dd r}{\dd z} &= \sin(\theta)\\
    c\frac{\dd y}{\dd z} &= \cos(\theta),
\end{align*}
where $c(\beta,\ell_d) = 2/\big(\ell_d-s_0(\beta)\big)$ is the mapping parameter. In this form, the residuals of the differential equations include $\beta$ and $\ell_d$ as well as the $\B{\xi}$ vectors for each of the independent variables. Therefore, the equations are in a form such that an optimization technique can be used to minimize the residuals: in this case, nonlinear least-squares was used. 

\begin{figure}[!htbp]
    \centering
    \includegraphics[width=\linewidth]{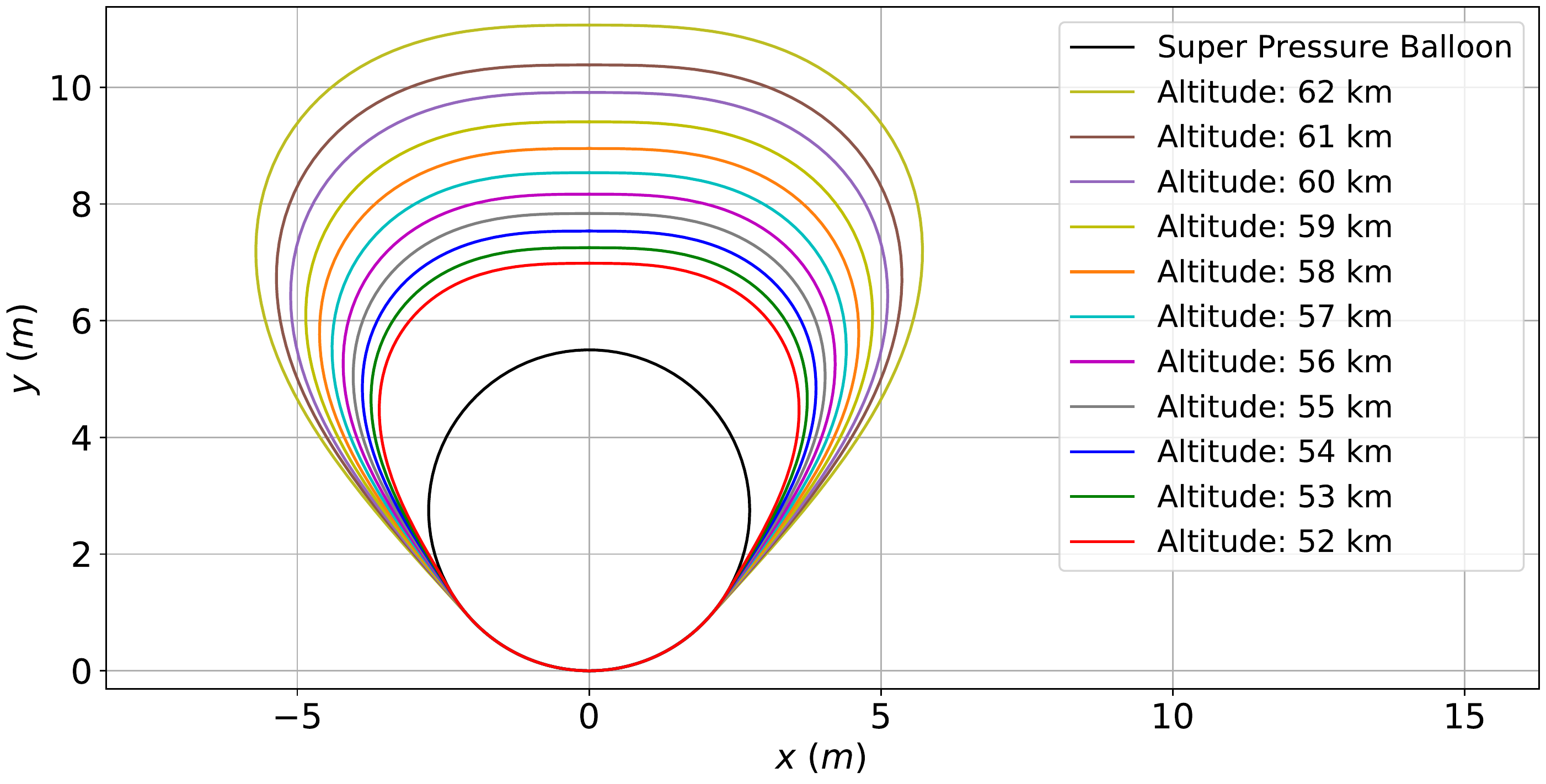}
    \caption{Natural balloon shapes on Venus for a range of altitudes from 52 km to 62 km.}
    \label{fig:NaturalBalloonShapes}
\end{figure}

Consider the Venus atmospheric data in Table \ref{tab:BalloonAtmData}, which was collected as part of the Venus Variable Altitude Aerobot project at JPL, and let the balloon constants be those given in Table \ref{tab:BalloonData}.
\begin{table}[!ht]
    \centering
    \caption{Tandem balloon atmospheric data.}
    \label{tab:BalloonAtmData}
    \begin{tabular}{|c|c|c|c|}
        \hline
        Altitude (km) & \makecell{Atmospheric\\Density ($\frac{kg}{m^3}$)} & \makecell{Super Pressure\\Balloon Gas Mass (kg)} & Gravity ($\frac{m}{s^2}$)\\\hline
        52 & 1.28 & 11.62 & 8.719 \\\hline
        53 & 1.15 & 10.74 & 8.716 \\\hline
        54 & 1.03 & 9.97 & 8.713 \\\hline
        55 & 0.921 & 9.29 & 8.71 \\\hline
        56 & 0.818 & 8.67 & 8.707 \\\hline
        57 & 0.721 & 8.12 & 8.704 \\\hline
        58 & 0.629 & 7.58 & 8.702 \\\hline
        59 & 0.545 & 7.14 & 8.699 \\\hline
        60 & 0.469 & 6.812 & 8.696 \\\hline
        61 & 0.41 & 6.675 & 8.693 \\\hline
        62 & 0.341 & 6.2675 & 8.69 \\\hline
    \end{tabular}
\end{table}
\begin{table}[!ht]
    \centering
    \caption{Tandem balloon constants.}
    \label{tab:BalloonData}
    \begin{tabular}{|c|c|}
        \hline
        Parameter & Value and Units\\\hline
        $w\ \Big(\frac{kg}{m^2}\Big)$ & $0.095$\\\hline
        $w_s\ \Big(\frac{kg}{m^2}\Big)$ & $0.215$\\\hline
        $M_g\ \Big(\frac{kg}{mol}\Big)$ & $4\times10^{-3}$\\\hline
        $M_{atm}\ \Big(\frac{kg}{mol}\Big)$ & $4.34\times10^{-2}$\\\hline
        $L\ (N)$ & $208g$\\\hline
    \end{tabular}
\end{table}
Using this data and the differential equations above, TFC was used to solve the natural balloon shapes: the results are shown in Figure \ref{fig:NaturalBalloonShapes}. The average solution time for the natural balloon shapes in Figure \ref{fig:NaturalBalloonShapes} was $0.65$ seconds, and the residual of the differential equation at all points was $\mathcal{O}(10^{-15})$. To compare, solving the same problem in Octave via a shooting method that uses \verb"fzero" and \verb"ode45" takes an average of 9.85 seconds per case. Of course, this comparison is not exactly one-to-one, because as mentioned earlier, TFC is implemented in JAX, and many of the functions have been JIT-ed.

\begin{figure}[!htbp]
    \centering
    \includegraphics[width=\linewidth]{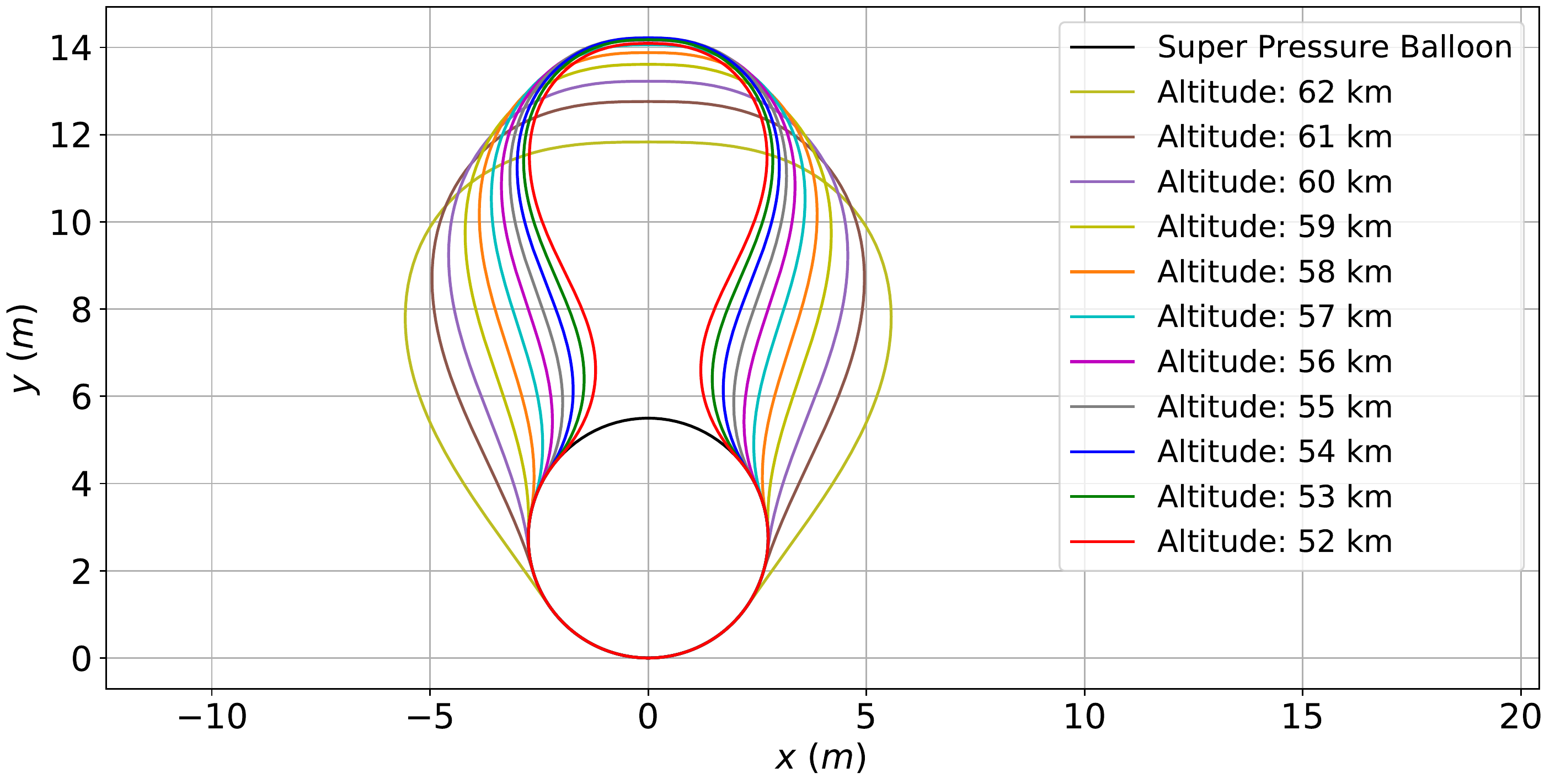}
    \caption{Constant circumferential stress balloon shapes with $\ell_d = 18$ meters on Venus for a range of altitudes from 52 km to 62 km.}
    \label{fig:ConstantStressBalloonShapes}
\end{figure}

Another classic balloon shape scenario is one wherein the circumferential stress is constant. In this case, one can trade the unknown $\ell_d$ for $\sigma_c$ if desired, i.e., the balloon has a fixed material length, and the circumferential stress is some unknown constant. Let the fixed material length $\ell_d = 18$ meters; Figure \ref{fig:ConstantStressBalloonShapes} shows the balloon shapes for this case. The average solution time for each case in Figure \ref{fig:ConstantStressBalloonShapes} was 0.75 seconds, and the residual of the differential equation at all points was $\mathcal{O}(10^{-15})$. Solving the same problem in Octave via a shooting method that uses \verb"fsolve" and \verb"ode45" takes an average of 45.7 seconds per case.

\section{Wave Equation}
The wave equation is a well-known PDE that describes the propagation of waves, such as those found in a vibrating string or $n$-dimensional membrane. This section applies TFC to the one-dimensional (one spatial dimension and one time dimension) wave equation and the two-dimensional (two spatial dimensions and one time dimension) wave equation.

\subsection{One-Dimensional Wave Equation}
\begin{figure}[!ht]
    \centering
    \ifDarkTheme
        \includegraphics[width=0.8\linewidth]{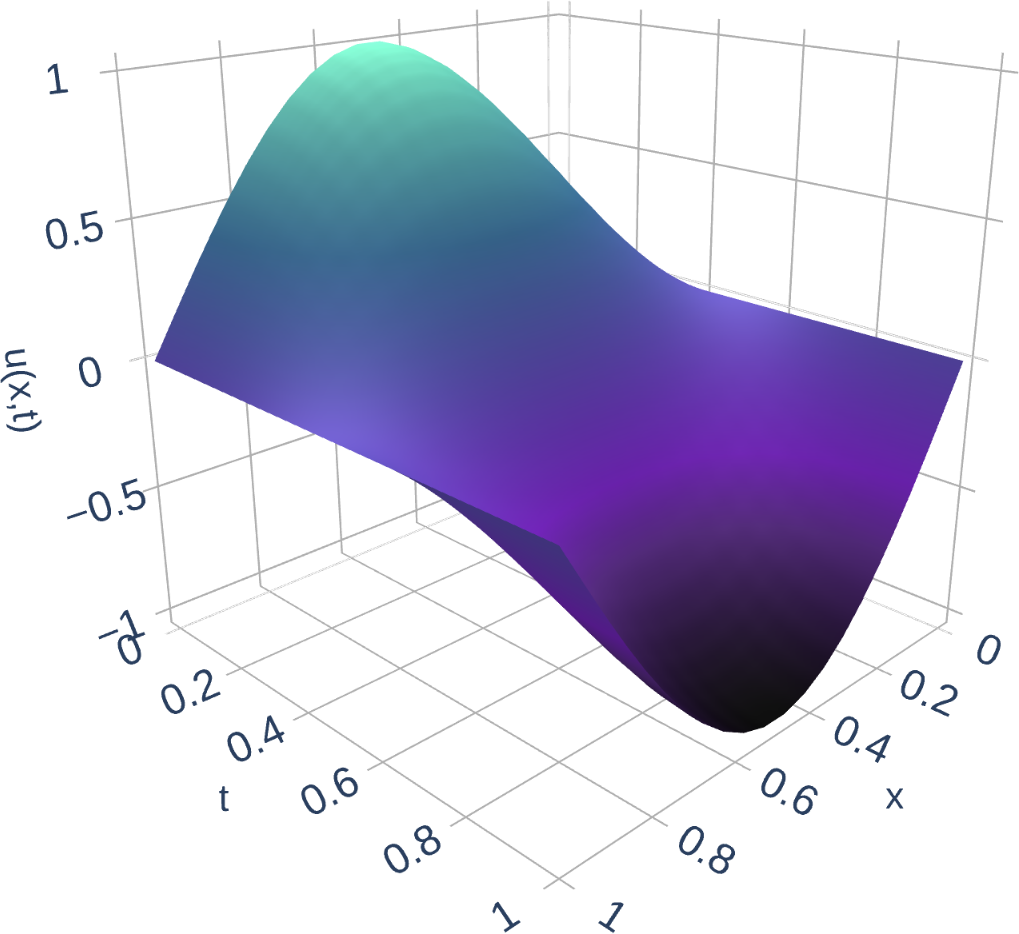}
    \else
        \includegraphics[width=0.8\linewidth]{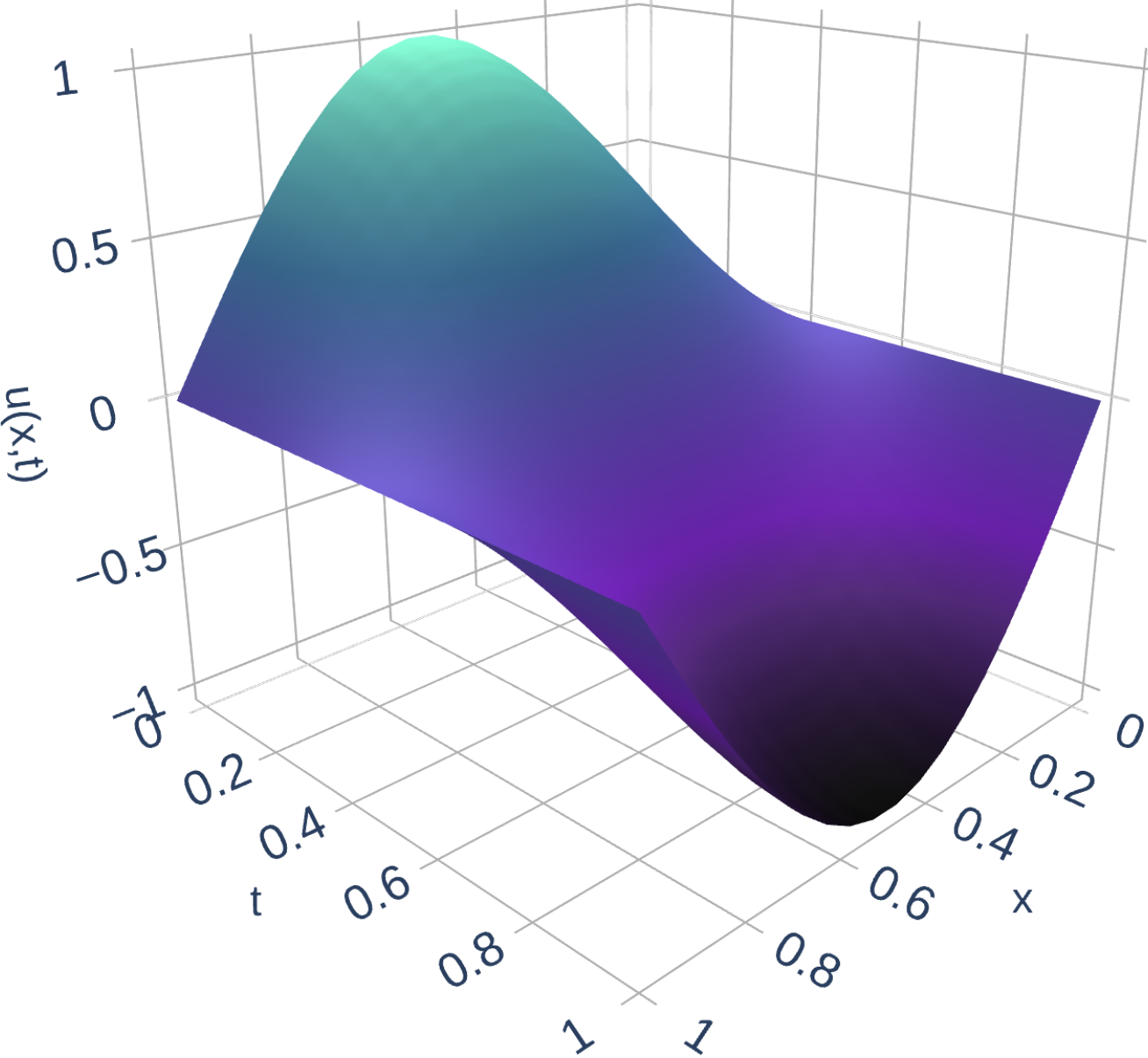}
    \fi
    \caption{Analytical solution for the one-dimensional wave equation.}
    \label{fig:WaveEqnAnalyticalSoln}
\end{figure}

Consider the wave equation for a one-dimensional object,
\begin{equation*}
    u_{xx} = k^2 u_{tt},
\end{equation*}
for some constant $k$ on the domain $(x,t)\in[0,1]\times[0,1]$ with the following boundary conditions,
\begin{equation*}
    u(0,t) = 0, \quad u(1,t) = 0, \quad u(x,0) = \sin(\pi x), \andd u_t(x,0) = 0.
\end{equation*}
One can physically imagine these boundary conditions as describing the string on a musical instrument that is fixed at both ends and free to vibrate with initial displacement $u(x,0) = \sin(\pi x)$. Let the constant $k=1$; then, the analytical solution is,
\begin{align*}
    u(x,t) = \sin(\pi x)\cos(\pi t).
\end{align*}
The analytical solution is shown in Figure \ref{fig:WaveEqnAnalyticalSoln}.

The \ce\ written in recursive form is,
\begin{align*}
    \p{1}u(x,t,g(x,t)) &= g(x,t) -(1-x)g(0,t) -x g(1,t) \\
    \p{2}u(x,t,g(x,t)) &= g(x,t) + \sin(\pi x)-g(x,0) - t g_t(x,0)
\end{align*}
where $\p{1}{u}$ can be used as the free function in $\p{2}{u}$ or $\p{2}{u}$ can be used as the free function in $\p{1}{u}$ to create the full \ce. The \ce\ written in tensor form is,
\begin{equation*}
    u(x,t,g(x,t)) = g(x,t)+\mathcal{M}_{ij}(x,t,g(x,t))\Phi_i(x)\Phi_j(t),
\end{equation*}
where
\begin{equation*}
    \mathcal{M}_{ij}(x,t,g(x,t)) = \begin{bmatrix} 0 & \sin(\pi x)-g(x,0) & -g_t(x,0) \\
    -g(0,t) & g(0,0) & g_t(0,0) \\ 
    -g(1,t) & g(1,0) & g_t(1,0) \end{bmatrix},
\end{equation*}
\begin{equation*}
    \Phi_i(x) = \begin{Bmatrix} 1, & 1-x, & x\end{Bmatrix} \andd
    \Phi_j(t) = \begin{Bmatrix} 1, & 1, & t\end{Bmatrix}.
\end{equation*}

Using Legendre orthogonal polynomials up to degree 20 as the free function and a grid of $30\times30$ training points, the PDE solution was estimated using the TFC method. The solution was obtained in 0.49 seconds, and the average error on a test set of $100\times100$ evenly spaced training points was $1.044\times10^{-15}$.

\subsection{Two-Dimensional Wave Equation}\label{sec:2dWaveEqn}
\begin{figure}[!ht]
    \centering
    \ifDarkTheme
        \includegraphics[width=0.8\linewidth]{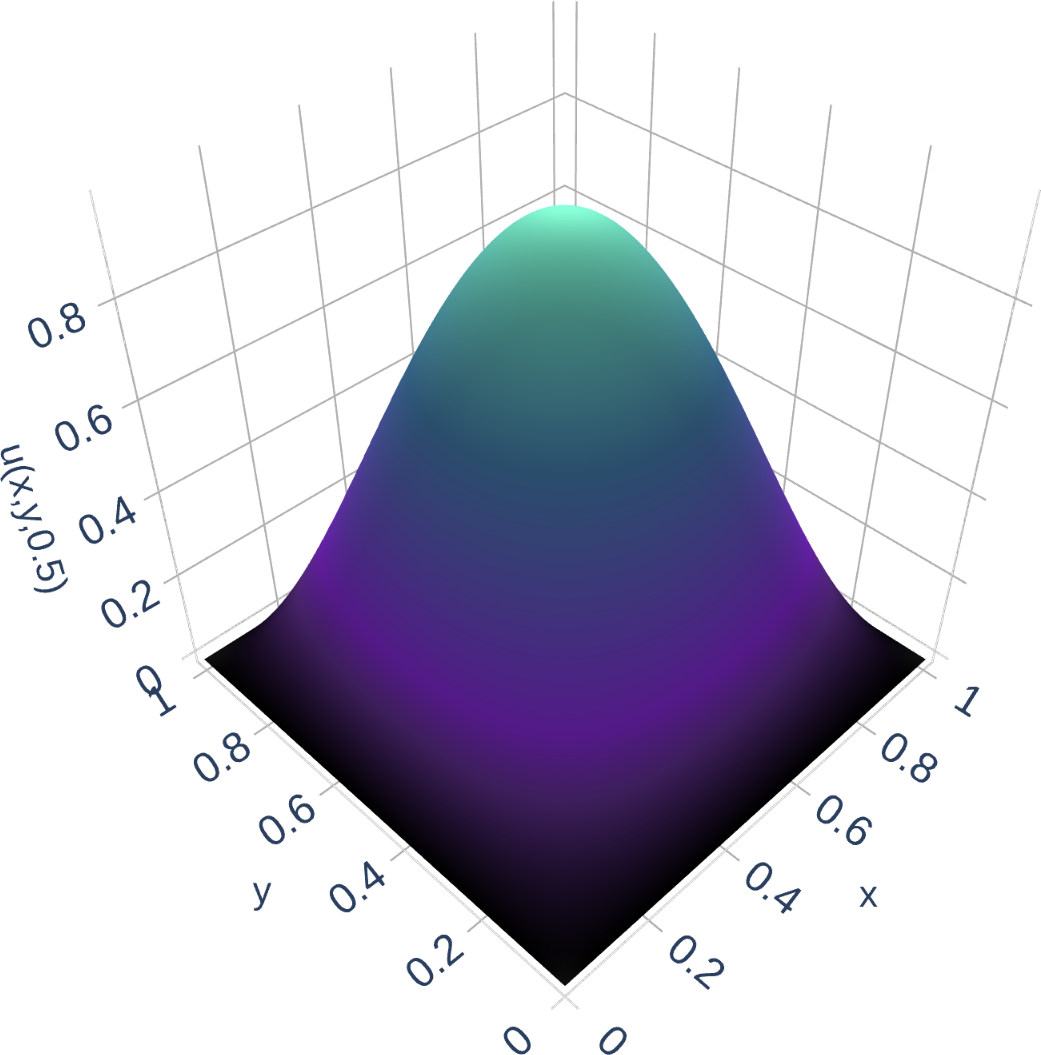}
    \else
        \includegraphics[width=0.8\linewidth]{Figures/Wave2DAnalyticalSoln.png}
    \fi
    \caption{Two-dimensional wave equation analytical solution at $t=0.5$.}
    \label{fig:WaveEqnAnalyticalSoln2D}
\end{figure}

The two-dimensional wave equation can be used to describe objects such as a flexible two-dimensional membrane. Consider such a membrane clamped at all sides with an initial deformation $u(x,y,0) = \sin(\pi x)\sin(\pi y)$. Then, the governing PDE can be written as,
\begin{equation*}
    u_{xx}+u_{yy} = k^2 u_{tt},
\end{equation*}
on the domain $(x,t)\in[0,1]\times[0,1]\times[0,1]$ with the following boundary conditions,
\begin{equation*}
\begin{gathered}
    u(0,y,t) = 0, \quad u(1,y,t) = 0, \quad u(x,0,t) = 0, \quad u(x,1,t) = 0,\\
    \quad u(x,y,0) = \sin(\pi x)\sin(\pi y), \andd u_t(x,y,0) = 0.
\end{gathered}
\end{equation*}
Let $k=8$, then the analytical solution is,
\begin{equation*}
    u(x,y,t) = \sin(\pi x)\sin(\pi y)\cos\left(\frac{\pi\sqrt{2}}{8} t\right).
\end{equation*}
The analytical solution at $t=0.5$ is shown in Figure \ref{fig:WaveEqnAnalyticalSoln2D}.

The \ce\ written in recursive form is,
\begin{align*}
    \p{1}{u}(x,y,t,g(x,y,t)) &= g(x,y,t)-(1-x)g(0,y,t)-xg(1,y,t)\\
    \p{2}{u}(x,y,t,g(x,y,t)) &= g(x,y,t)-(1-y)g(x,0,t)-yg(x,1,t)\\
    \p{3}{u}(x,y,t,g(x,y,t)) &= g(x,y,t)+\sin(\pi x)\sin(\pi y)-g(x,y,0)-tg_t(x,y,0)\\
\end{align*}
where $\p{1}{u}$, $\p{2}{u}$, and $\p{3}{u}$ can be processed in any order to produce the full \ce. The tensor form of the \ce\ is,
\begin{equation*}
    u(x,y,t,g(x,y,t)) = g(x,y,t)+\mathcal{M}_{ijk}(x,y,t,g(x,y,t))\Phi_i(x)\Phi_j(y)\Phi_k(t),
\end{equation*}
where
\begin{align*}
    \mathcal{M}_{ij1}(x,y,t,g(x,y,t)) &= \begin{bmatrix} 0 & -g(x,0,t) & -g(x,1,t) \\
    -g(0,y,t) & g(0,0,t) & g(0,1,t) \\ 
    -g(1,y,t) & g(1,0,t) & g(1,1,t) \end{bmatrix}\\
    \mathcal{M}_{ij2}(x,y,t,g(x,y,t)) &= \begin{bmatrix} \sin(\pi x)\sin(\pi y)-g(x,y,0) & g(x,0,0) & g(x,1,0) \\
    g(0,y,0) & -g(0,0,0) & -g(0,1,0) \\ 
    g(1,y,0) & -g(1,0,0) & -g(1,1,0) \end{bmatrix}\\
    \mathcal{M}_{ij3}(x,y,t,g(x,y,t)) &= \begin{bmatrix} -g_t(x,y,0) & g_t(x,0,0) & g_t(x,1,0) \\
    g_t(0,y,0) & -g_t(0,0,0) & -g_t(0,1,0) \\ 
    g_t(1,y,0) & -g_t(1,0,0) & -g_t(1,1,0) \end{bmatrix}\\
\end{align*}
and
\begin{align*}
    \Phi_i(x) &= \begin{Bmatrix} 1, & 1-x, & x\end{Bmatrix}, \\
    \Phi_j(y) &= \begin{Bmatrix} 1, & 1-y, & y\end{Bmatrix}, \\
    \Phi_k(t) &= \begin{Bmatrix} 1, & 1, & t\end{Bmatrix}.
\end{align*}

For the two-dimensional wave equation, choosing ELMs as the free function led to a better estimate of the solution than Chebyshev or Legendre orthogonal polynomials on average: although the two were similar. To illustrate, Table \ref{tab:waveEquation2DComparison} shows the maximum and mean errors when using the two methods for different numbers of basis functions: the number of basis functions corresponds to the number of Chebyshev polynomials there are on this problem when keeping all polynomials up to degree 3, 6, 9, 12, 15, and 18. The TFC method used Chebyshev orthogonal polynomials, and the X-TFC method used the hyperbolic tangent as the activation function. Each method used a grid of $11\times11\times11$ training points and a test set of $15\times15\times15$ uniformly spaced points.
\begin{table}[!ht]
\centering
\caption{TFC and X-TFC solution errors for various numbers of basis functions when solving the two-dimensional wave equation.}
\label{tab:waveEquation2DComparison}
\begin{tabular}{|c|c|c|c|c|} 
\hline
 \makecell{$m$} & \multicolumn{2}{c|}{\textbf{TFC}} & \multicolumn{2}{c|}{\textbf{X-TFC}}\\ \hline
  & \makecell{Maximum Error} & \makecell{Mean Error} & \makecell{Maximum Error} & \makecell{Mean Error} \\\hline
12 & $5.32$ & $5.16\times10^{-1}$ & $6.24\times10^{-3}$& $6.98\times10^{-4}$\\ \hline
76 & $8.07\times10^{-3}$ & $1.04\times10^{-3}$ & $4.89\times10^{-3}$ & $4.56\times10^{-4}$\\ \hline
212 & $1.64\times10^{-1}$ & $1.55\times10^{-2}$ & $2.42\times10^{-3}$ & $2.29\times10^{-4}$\\ \hline
447 & $2.22\times10^{-2}$ & $2.09\times10^{-3}$ & $9.34\times10^{-3}$ & $8.79\times10^{-4}$\\ \hline
808 & $3.91\times10^{-3}$ & $3.67\times10^{-4}$ & $3.32\times10^{-3}$ & $3.12\times10^{-4}$\\ \hline
1322 & $3.90\times10^{-3}$ & $3.67\times10^{-4}$ & $3.34\times10^{-3}$ & $3.02\times10^{-4}$\\ \hline
\end{tabular}
\end{table}

Table \ref{tab:waveEquation2DComparison} shows that while the two methods are similar, the X-TFC method performs slightly better. The differences between the two methods are the most pronounced when a lower number of basis functions is used. Consequently, ELMs were used to estimate the solution of this differential equation. For one particular run using 650 neurons, the solution was obtained in 18.4 seconds, and the average solution error on the test set was $2.124\times10^{-5}$.

\section{Biharmonic Equation}
The biharmonic equation is a fourth-order linear PDE that appears in linear elasticity theory \cite{continuumBook}. The PDE is given by,
\begin{equation*}
    \nabla^4 u(\B{x}) = \nabla^2(\nabla^2 u(\B{x}) ) = f(\B{x}) ,
\end{equation*}
where $u$ is the dependent variable of interest, $f$ is a forcing term, and $\nabla^2$ is the Laplacian operator. In two-dimensional plate problems, the variable $u$ is related to the stress experienced by the plate, and $f(\B{x})$ is related to the body forces acting on the plate. 

\subsection{Cartesian Coordinates}
\begin{figure}[!hb]
    \centering
    \ifDarkTheme
        \includegraphics[width=0.8\linewidth]{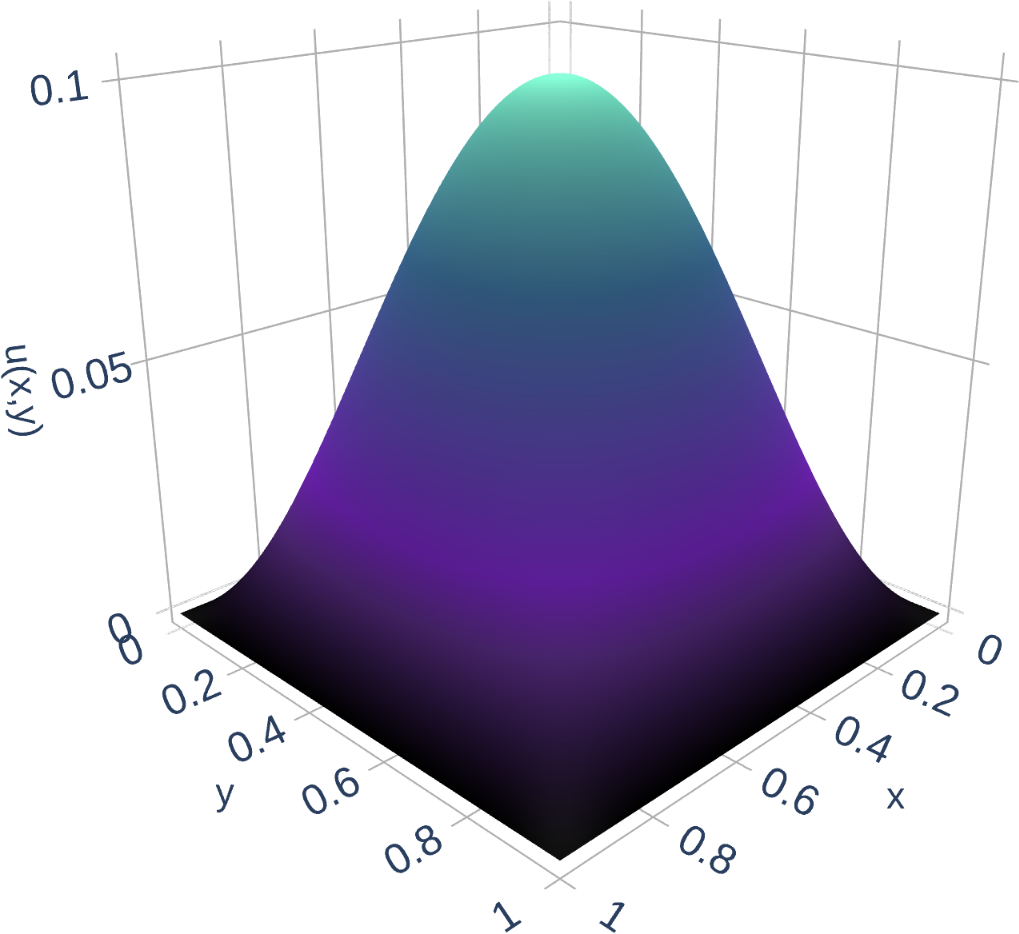}
    \else
        \includegraphics[width=0.8\linewidth]{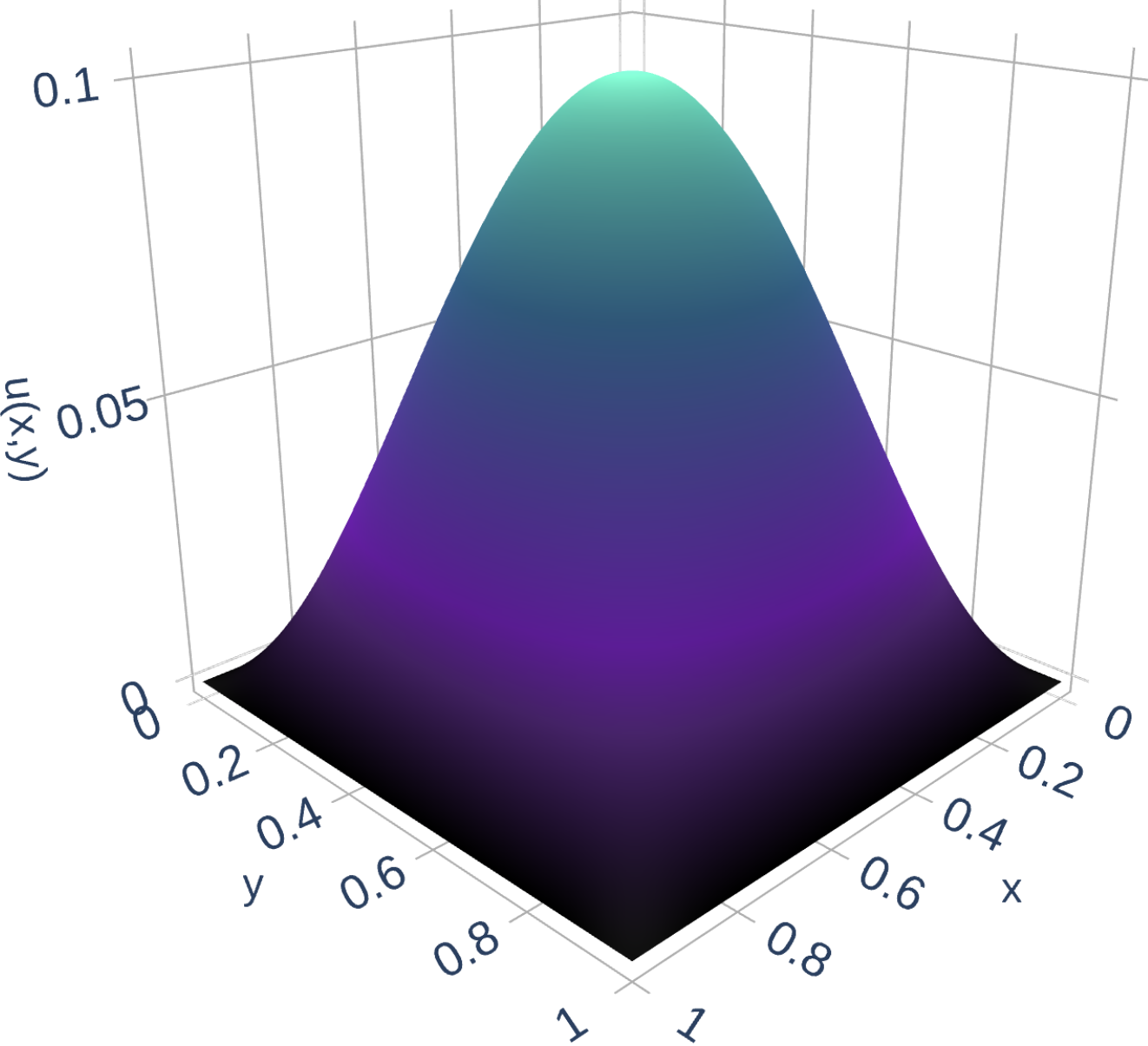}
    \fi
    \caption{Biharmonic equation analytical solution in Cartesian coordinates.}
    \label{fig:biharmonic}
\end{figure}

Consider the following forcing function,
\begin{equation*}
    \nabla^4 u(x,y) = 4 \pi^2 \sin(\pi x)\sin(\pi y),
\end{equation*}
on the domain $(x,y)\in[0,1]\times[0,1]$ with the following boundary conditions,
\begin{gather*}
    u(0,y) = u(1,y) = u(x,0) = u(y,0) = 0\\
    u_{xx}(0,y) = u_{xx}(1,y) = u_{yy}(x,0) = u_{yy}(y,0) = 0
\end{gather*}
The analytical solution to this problem is \cite{FEniCSBook,DolfinCode},
\begin{equation*}
    u(x,y) = \frac{1}{\pi^2} \sin(\pi x)\sin(\pi y),
\end{equation*}
and is shown in Figure \ref{fig:biharmonic}.

The \ce\ for this problem in recursive form is,
\begin{align*}
    \p{1}{u}(x,y,g(x,y)) =\ &g(x,y)-(1-x) g(0,y)-x g(1,y)\\
    &-\frac{-x^3+3 x^2-2 x}{6} g_{xx}(0,y)-\frac{x^3-x}{6} g_{xx}(1,y),\\
    \p{2}{u}(x,y,g(x,y)) =\ &g(x,y)-(1-y) g(x,0)-y g(x,1)\\
    &-\frac{y^3-y}{6} g_{yy}(x,1)-\frac{-y^3+3 y^2-2 y}{6} g_{yy}(x,0),
\end{align*}
where $\p{1}{u}$ can be used as the free function in $\p{2}{u}$ or $\p{2}{u}$ can be used as the free function in $\p{1}{u}$ to create the full \ce. In tensor form the \ce\ is,
\begin{equation*}
    u(x,y,g(x,t)) = g(x,y)+\mathcal{M}_{ij}(x,y,g(x,y))\Phi_i(x)\Phi_j(y),
\end{equation*}
where
\begin{equation*}
    \mathcal{M}_{ij}(x,y,g(x,y)) = \begin{bmatrix}
    0 & -g(x,0) & -g(x,1) & -g_{yy}(x,0) & -g_{yy}(x,1) \\
    -g(0,y) & g(0,0) & g(0,1) & g_{yy}(0,0) & g_{yy}(0,1) \\
    -g(1,y) & g(1,0) & g(1,1) & g_{yy}(1,0) & g_{yy}(1,1) \\
    -g_{xx}(0,y) & g_{xx}(0,0) & g_{xx}(0,1) & g_{xxyy}(0,0) & g_{xxyy}(0,1) \\
    -g_{xx}(1,y) & g_{xx}(1,0) & g_{xx}(1,1) & g_{xxyy}(1,0) & g_{xxyy}(1,1) \end{bmatrix}
\end{equation*}
and
\begin{align*}
    \Phi_i(x) &= \begin{Bmatrix} 1, & 1-x, & x, & \frac{-x^3+3 x^2-2 x}{6}, & \frac{x^3-x}{6} \end{Bmatrix}, \\
    \Phi_i(y) &= \begin{Bmatrix} 1, & 1-y, & y, & \frac{-y^3+3 y^2-2 y}{6}, & \frac{y^3-y}{6} \end{Bmatrix}.
\end{align*}


Using Chebyshev orthogonal polynomials up to degree 26 as the free function and a grid of $20\times20$ training points, the PDE solution was estimated using the TFC method. The solution was obtained in 0.94 seconds, and the average error on a test set of $100\times100$ uniformly spaced points was $1.661\times10^{-16}$.

\subsection{Polar Coordinates}
\begin{figure}[!ht]
    \centering
    \ifDarkTheme
        \textattachfile{Figures/BiharmonicPolarAnalyticalSolution.html}{\includegraphics[width=0.7\linewidth]{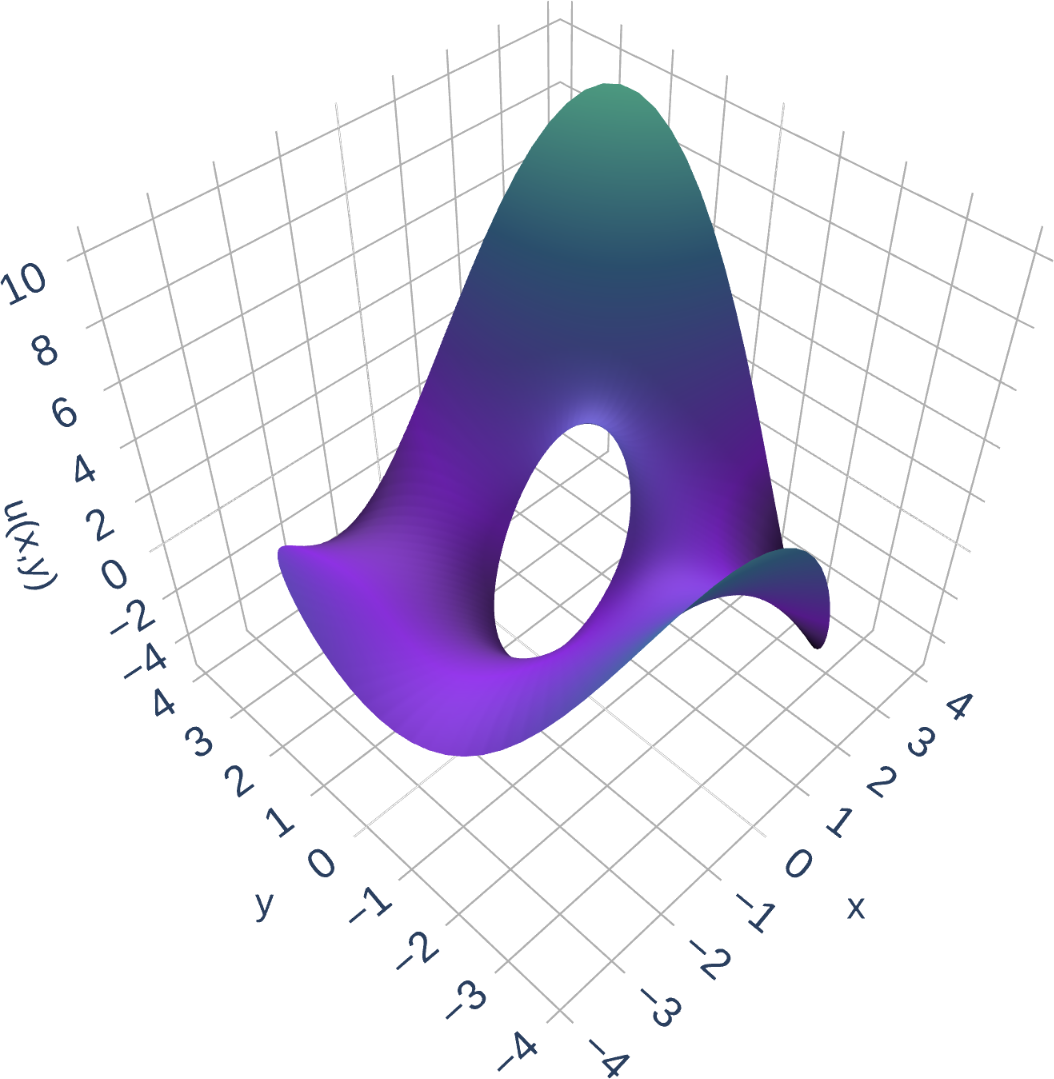}}
    \else
        \textattachfile{Figures/BiharmonicPolarAnalyticalSolution.html}{\includegraphics[width=0.7\linewidth]{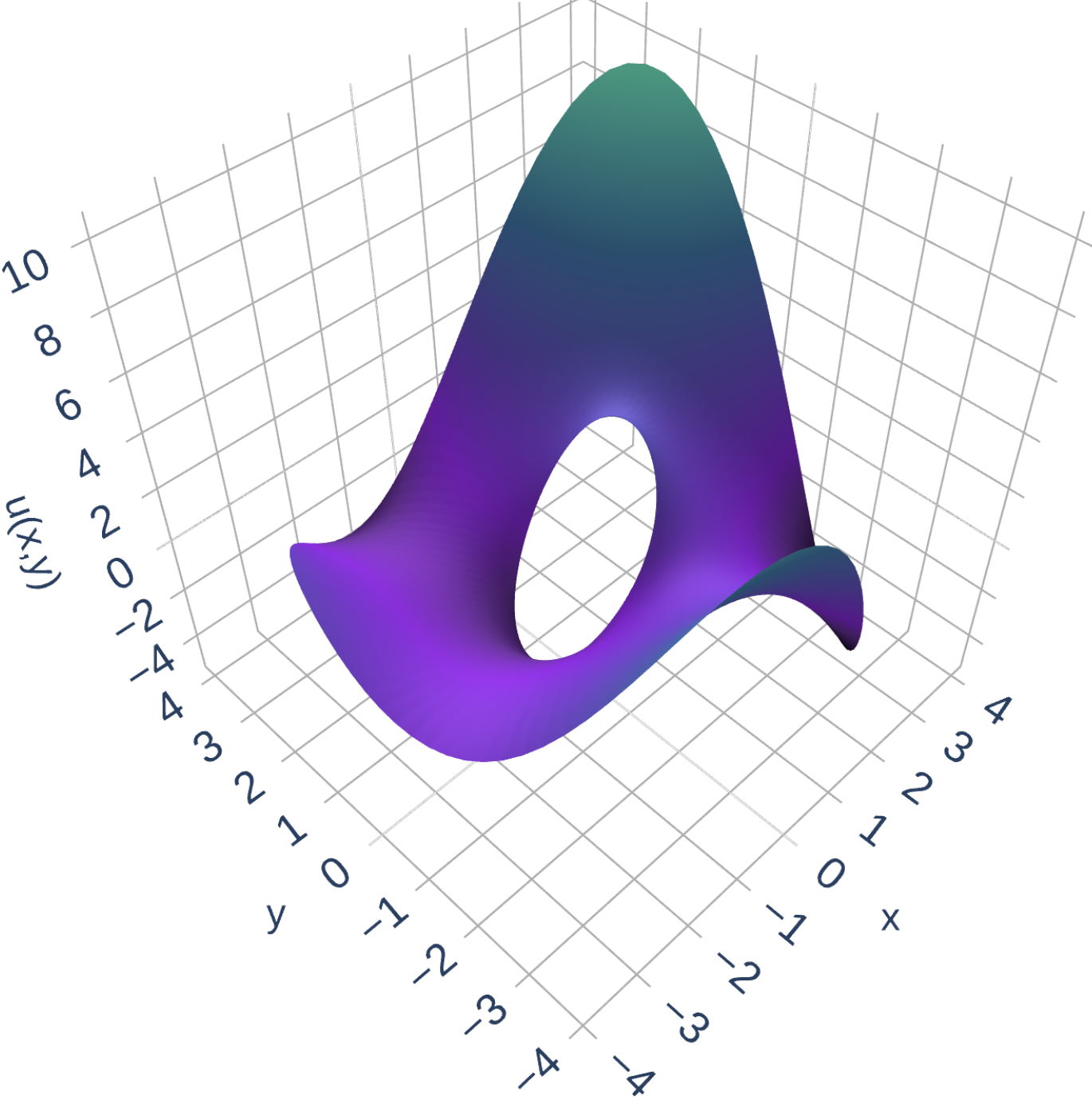}}
    \fi
    \caption{Biharmonic equation analytical solution in polar coordinates. Note, this figure contains an embedded, standalone HMTL version of the plot that can be viewed/downloaded by clicking on it. Doing so may require a dedicated PDF viewer such as Adobe Acrobat or Okular.}
    \label{fig:biharmonicPolar}
\end{figure}

Consider the following forcing function,
\begin{equation*}
    \nabla^4 u(r,\theta) = 0,
\end{equation*}
on the domain $(r,\theta)\in[1,4]\times[0,2\pi]$ with the following boundary conditions,
\begin{gather*}
    u(1,\theta) = \frac{1}{4} \sin (2 \theta )+\frac{1}{16} \sin (3 \theta )+\pi  \cos (\theta )+\frac{1}{8} \\
    u(4,\theta) = 4 \sin (2 \theta )+4 \sin (3 \theta )+\frac{1}{4} \pi  \cos (\theta )+2 \\
    u_{rr}(1,\theta) = \frac{1}{2} \sin (2 \theta )+\frac{3}{8} \sin (3 \theta )+2 \pi  \cos (\theta )+\frac{1}{4} \\
    u_{rr}(4,\theta) = \frac{1}{2} \sin (2 \theta )+\frac{3}{2} \sin (3 \theta )+\frac{1}{32} \pi  \cos (\theta )+\frac{1}{4} \\
    u(r,0) = u(r,2\pi) \\
    u_{\theta}(r,0) = u_{\theta}(r,2\pi) \\
    u_{\theta\theta}(r,0) = u_{\theta\theta}(r,2\pi) \\
    u_{\theta\theta\theta}(r,0) = u_{\theta\theta\theta}(r,2\pi).
\end{gather*}
The analytical solution to this problem is \cite{PolarBiharmonicSolution},
\begin{equation*}
    u(r,\theta) = \frac{r^3}{16} \sin (3 \theta )+\frac{r^2}{4} \sin (2 \theta )+\frac{r^2}{8}+\frac{\pi  \cos (\theta )}{r},
\end{equation*}
and is shown in Figure \ref{fig:biharmonicPolar}.
Note that the Laplacian in polar coordinates is,
\begin{equation*}
    \nabla^2 u(r,\theta) = u_{rr}+\frac{1}{r}u_r + \frac{1}{r^2} u_{\theta\theta},
\end{equation*}
so the biharmonic operator in polar coordinates is \cite{PolarBiharmonicOperator},
\begin{equation*}
    \nabla^4 u(r,\theta) = u_{rrrr}+\frac{2}{r^2}u_{rr\theta\theta}+\frac{1}{r^4}u_{\theta\theta\theta\theta}+\frac{2}{r}u_{rrr}-\frac{2}{r^3}u_{r\theta\theta}-\frac{1}{r^2}u_{rr}+\frac{4}{r^4}u_{\theta\theta}+\frac{1}{r^3}u_r.
\end{equation*}

The \ce\ for this problem given in recursive form is,
\begin{align*}
    \p{1}{u}(r,&\theta,g(r,\theta)) = g(r,\theta)\\
    &+\frac{1}{3} (r-1) \left(4 \sin (2 \theta )+4 \sin (3 \theta )+\frac{1}{4} \pi  \cos (\theta )+2 -g(4,\theta )\right)\\
    &+\frac{1}{3} (4-r) \left(\frac{1}{4} \sin (2 \theta )+\frac{1}{16} \sin (3 \theta )+\pi  \cos (\theta )+\frac{1}{8} -g(1,\theta )\right) \\
    &+\frac{1}{18} \left(-r^3+12 r^2-39 r+28\right) \left(-g_{rr}(1,\theta )+\frac{1}{2} \sin (2 \theta )+\frac{3}{8} \sin (3 \theta )+2 \pi  \cos (\theta )+\frac{1}{4}\right)\\
    &+\frac{1}{18} \left(r^3-3 r^2-6 r+8\right) \left(\frac{1}{2} \sin (2 \theta )+\frac{3}{2} \sin (3 \theta )+\frac{1}{32} \pi  \cos (\theta )+\frac{1}{4}-g_{rr}(4,\theta )\right),\\
    \p{2}{u}(r,&\theta,g(r,\theta)) = g(r,\theta)-\frac{\theta}{2 \pi }  \Big(g(r,2\pi)-g(r,0)\Big)+\frac{2 \pi  \theta -\theta ^2 }{4 \pi }\Big(g_r(r,2\pi)-g_r(r,0)\Big)\\
    &+\frac{-\theta ^3+3 \pi  \theta ^2-2 \pi ^2 \theta }{12 \pi } \Big(g_{rr}(r,2\pi)-g_{rr}(r,0))\Big)\\
    &+\frac{-\theta ^4+4 \pi  \theta ^3-4 \pi ^2 \theta ^2 }{48 \pi } \Big(g_{rrr}(r,2\pi)-g_{rrr}(r,0)\Big) ,
\end{align*}
where $\p{1}{u}$ can be used as the free function in $\p{2}{u}$ or $\p{2}{u}$ can be used as the free function in $\p{1}{u}$ to create the full \ce. In tensor form the \ce\ is,
\begin{equation*}
    u(r,\theta,g(r,\theta)) = g(r,\theta)+\mathcal{M}_{ij}(r,\theta,g(r,\theta))\Phi_i(r)\Phi_j(\theta),
\end{equation*}
where
\begin{align*}
    \mathcal{M}_{ij}(r,\theta,g(r,\theta)) &= \begin{bmatrix}
     0\\
     -g(1,\theta )+\frac{1}{4} \sin (2 \theta )+\frac{1}{16} \sin (3 \theta )+\pi  \cos (\theta )+\frac{1}{8}\\
     -g(4,\theta )+4 \sin (2 \theta )+4 \sin (3 \theta )+\frac{1}{4} \pi  \cos (\theta )+2 \\
     -g_{rr}(1,\theta )+\frac{1}{2} \sin (2 \theta )+\frac{3}{8} \sin (3 \theta )+2 \pi  \cos (\theta )+\frac{1}{4} \\
     -g_{rr}(4,\theta )+\frac{1}{2} \sin (2 \theta )+\frac{3}{2} \sin (3 \theta )+\frac{1}{32} \pi  \cos (\theta )+\frac{1}{4}
    \end{bmatrix}, \quad i=1\\
    \mathcal{M}_{ij}(r,\theta,g(r,\theta)) &= \begin{bmatrix}
    g(r,2 \pi )-g(r,0) & g_{\theta}(r,2 \pi )-g_{\theta}(r,0) \\
    g(1,0)-g(1,2 \pi ) & g_{\theta}(1,0)-g_{\theta}(1,2 \pi )\\
    g(4,0)-g(4,2 \pi ) & g_{\theta}(4,0)-g_{\theta}(4,2 \pi ) \\
    g_{rr}(1,0)-g_{rr}(1,2 \pi ) & g_{rr\theta}(1,0)-g_{rr\theta}(1,2 \pi ) \\
    g_{rr}(4,0)-g_{rr}(4,2 \pi ) & g_{rr\theta}(4,0)-g_{rr\theta}(4,2 \pi )
    \end{bmatrix}, \quad i\in\{2,3\}\\
    \mathcal{M}_{ij}(r,\theta,g(r,\theta)) &= \begin{bmatrix}
     g_{\theta\theta}(r,2 \pi )-g_{\theta\theta}(r,0) & g_{\theta\theta\theta}(r,2 \pi )-g_{\theta\theta\theta}(r,0) \\
     g_{\theta\theta}(1,0)-g_{\theta\theta}(1,2 \pi ) & g_{\theta\theta\theta}(1,0)-g_{\theta\theta\theta}(1,2 \pi ) \\
     g_{\theta\theta}(4,0)-g_{\theta\theta}(4,2 \pi ) & g_{\theta\theta\theta}(4,0)-g_{\theta\theta\theta}(4,2 \pi ) \\
     g_{rr\theta\theta}(1,0)-g_{rr\theta\theta}(1,2 \pi ) & g_{rr\theta\theta\theta}(1,0)-g_{rr\theta\theta\theta}(1,2 \pi ) \\
     g_{rr\theta\theta}(4,0)-g_{rr\theta\theta}(4,2 \pi ) & g_{rr\theta\theta\theta}(4,0)-g_{rr\theta\theta\theta}(4,2 \pi )
    \end{bmatrix}, \quad i\in\{4,5\}
\end{align*}
and
\begin{align*}
    \Phi_i(r) &= \begin{Bmatrix} 1,\frac{4-r}{3},\frac{r-1}{3},\frac{-r^3+12 r^2-39 r+28}{18},\frac{r^3-3 r^2-6 r+8}{18} \end{Bmatrix}, \\
    \Phi_i(\theta) &= \begin{Bmatrix} 1,-\frac{\theta }{2 \pi },\frac{2 \pi  \theta -\theta ^2}{4 \pi },\frac{-\theta ^3+3 \pi  \theta ^2-2 \pi ^2 \theta }{12 \pi },\frac{-\theta ^4+4 \pi  \theta ^3-4 \pi ^2 \theta ^2}{48 \pi } \end{Bmatrix}.
\end{align*}

Using Chebyshev orthogonal polynomials up to degree 30 as the free function and a grid of $30\times30$ training points, the PDE solution was estimated using the TFC method. The solution was obtained in 10.67 seconds, and the average error on a test set of $100\times100$ uniformly spaced points was $1.535\times10^{-8}$.

%% file: Data/conclusions.tex

\chapter{SUMMARY AND CONCLUSIONS \label{chap:Conclusions}}

This dissertation is titled ``The Multivariate Theory of Functional Connections: An $n$-Dimensional Constraint Embedding Technique Applied to Partial Differential Equations'' because it presents two main ideas: (1) the derivation and analysis of multivariate TFC \ces\ (2) the estimation of PDE solutions using TFC. 

The first of these main ideas is conveyed in Chapter \ref{chap:tfcTheory}, which is a self-contained presentation of multivariate TFC. It includes the derivation of multivariate \ces\ starting from the univariate theory and all the associated mathematical theorems that have been proven to date. Although most readers will likely only be interested in using \ces\ for value constraints, derivative constraints, and linear combinations thereof, as these are the most common types of constraints found in differential equations, integral and component constraints are included as well. Because the chapter is comprehensive, it serves as a convenient and useful reference for any reader interested in TFC, regardless of their familiarity with the subject. In addition, extensions of the theory to inequality constraints, nonlinear constraints, parallelotope domains, lower-dimensional constraints in $n$-dimensions, and to any field, i.e., beyond real numbers to other fields such as complex numbers, are covered in Appendices \ref{app:SimpleNLConstraints}, \ref{app:Inequality},  and \ref{app:Miscellaneous}.

The second idea is conveyed in Chapter \ref{chap:deApplications}, which describes how to apply TFC to DEs, in particular, PDEs. In addition to describing the general methodology for solving DEs via TFC, this chapter discusses the common free function and optimization choices as well as their strengths and weaknesses. Similar to Chapter \ref{chap:tfcTheory}, Chapter \ref{chap:deApplications} is also self-contained;  consequently, it is a convenient and useful reference for any reader interested in solving DEs via TFC. 

Although Chapters \ref{chap:tfcTheory} and \ref{chap:deApplications} contain all the information necessary to apply TFC to DEs, they show few complicated examples: such examples are useful, as they highlight and clarify some of the nuances of TFC discussed in these chapters. To that end, Chapter \ref{chap:flexBodyApplications} contains some complex problems in a field of particular interest to the author: flexible body dynamics. Furthermore, these complex problems showcase the power and convenience of the numerical implementation. As examples:
\begin{enumerate}
    \item The natural tandem balloon shape problem is a complex system of four ODEs wherein both ends of the domain are themselves unknowns. The code for this problem is simple to read and write due to automatic differentiation---none of the derivatives have to be written out explicitly, which would clutter the code significantly---and the optimization of the free functions is simple too, despite the number of unknowns---a $\xi$ vector for each of the four dependent variables and the two unknowns associated with the ends of the domain.
    \item The polar biharmonic equation is a linear, fourth-order PDE with relative constraints up to the third derivative in one of the two independent variables. The numerical implementation allows the \ce\ to be written using the recursive format, which is short and easy to read. Furthermore, automatic differentiation makes creating the residual easy and straightforward even though it contains ninth-order partial derivatives (fifth-order partial derivatives in the \ce\ plus the fourth-order partial derivatives appearing in the residual itself).
\end{enumerate}
The code for both of these problems and most of the problems and examples in this dissertation can be found for free on the \href{https://github.com/leakec/tfc}{\textcolor{blue}{\underline{TFC GitHub}}} \cite{TfcGithub}. 

This dissertation showed that TFC is useful for solving differential equations. In many of the examples shown, the solution error found using TFC is multiple orders of magnitude lower than competing state-of-the-art methods. Moreover, most of the TFC solutions are found in seconds or fractions of a second, except for those found using Deep-TFC, which typically solves problems on the order of minutes. However, in its current state, TFC can only be applied to rectangular domains\footnote{Rectangular domains here means with respect to the coordinates being used, e.g., problems using polar coordinates such as the polar bi-harmonic problem appear to be on a cylindrical domain when viewed on a Cartesian grid, but are rectangular from the perspective of the polar coordinates.} and a limited number of non-rectangular domains \cite{TFC-Selected, BijectiveMapping-TFC}: this restriction is not present in many of the competing state-of-the-art methods.

\section{Future Work}
In terms of constraints, TFC can currently embed value, derivative, integral, component, and linear constraints, and it can embed any number of these constraints on any number of dimensions. However, it cannot embed sets of integral constraints whose integration variables refer to one another, such as, 
\begin{equation*}
    \int_0^1u(x,0) \dd{x} = 1 \andd \int_0^1u(0,y) \dd{y} = 1.
\end{equation*}
Finding a way to embed these constraints is a topic of future work. In addition, inequality constraints can currently be used in conjunction with value constraints only, see Appendix \ref{app:Inequality}. Integrating inequality constraints fully into the theory, i.e., finding a way to combine inequality constraints and the remaining types of linear constraints, remains a topic of further study. Also, as noted in the comparison with other methods, TFC is restricted to rectangular domains and a handful of irregular domains: extending TFC to all irregular domains is a topic of future research.

This dissertation focused on applying TFC to differential equations. However, there are a plethora of other applications that have yet to be explored, such as Computer-Aided Design (CAD) \cite{CoonsCadBook,CoonsCadSiggraph}, image warping \cite{CoonsImageWarping}, and security pattern design \cite{CoonsSecurityPattern}. Yet, even in the application of differential equations, there are numerous directions future research can pursue:
\begin{itemize}
    \item Hybrid basis functions - Combining two or more sets of basis functions and using the result as the free function, e.g., Fourier basis functions and Legendre orthogonal polynomials. Note that one could even combine X-TFC and basis functions in this way because both are simply a linear combination of functions.
    \item Optimizers - Only four optimizers are utilized in this dissertation: least-squares, L-BFGS, Adam, and CSVM. However, many optimizers could have been used instead, some of which may outperform those used here.
    \item NN architecture - Deep-TFC has only used fully connected NNs up to this point. Like the optimizers, there are a variety of NN architectures that could be used, some of which may outperform the fully connected NNs used here. 
\end{itemize}
This list is by no means exhaustive, and there are almost certainly research directions that have been excluded. Hence, the reader is encouraged to consider and pursue the research ideas that appear here as well as those that the author has not considered.

%% file: Data/appendices.tex
%

\begin{appendices}
\titleformat{\chapter}{\centering\normalsize}{APPENDIX \thechapter}{0em}{\vskip .5\baselineskip\centering}
\renewcommand{\appendixname}{APPENDIX}

\input{Data/appendixGraphTheory}
\input{Data/appendixSimpleNL}
\input{Data/appendixInequality}
\input{Data/appendixDividingDomain}
\input{Data/appendixOthBasisFunc}
\input{Data/appendixLS}
\input{Data/appendixJaxCode}

\input{Data/appendixNonLinearSvm}
\input{Data/appendixMiscellaneous}

\end{appendices}

%% file: Data/appendixGraphTheory.tex

\chapter{GRAPH THEORY}\label{app:GraphTheory}
This appendix provides a cursory overview of the graph theory concepts germane to determining the processing order of \ces; readers who would like a more in-depth presentation of these topics should consult Reference \cite{GraphTheory}. First, a few different types of graphs are defined.

\begin{figure}[!ht]
    \begin{subfigure}{0.49\linewidth}
        \centering
        \includegraphics[width=1.0\linewidth]{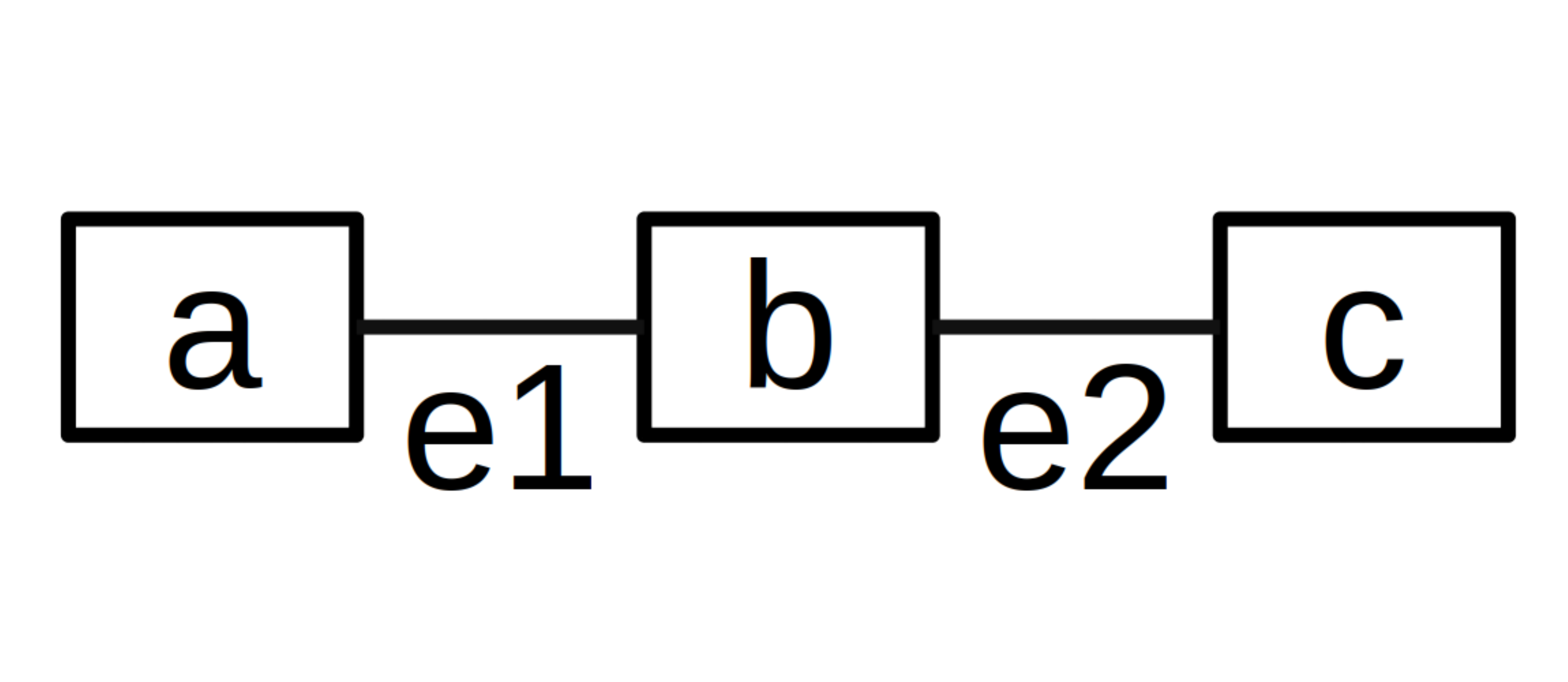}
        \caption{Graph}
    \end{subfigure}\hspace*{\fill}
    \begin{subfigure}{0.49\linewidth}
        \centering
        \includegraphics[width=1.0\linewidth]{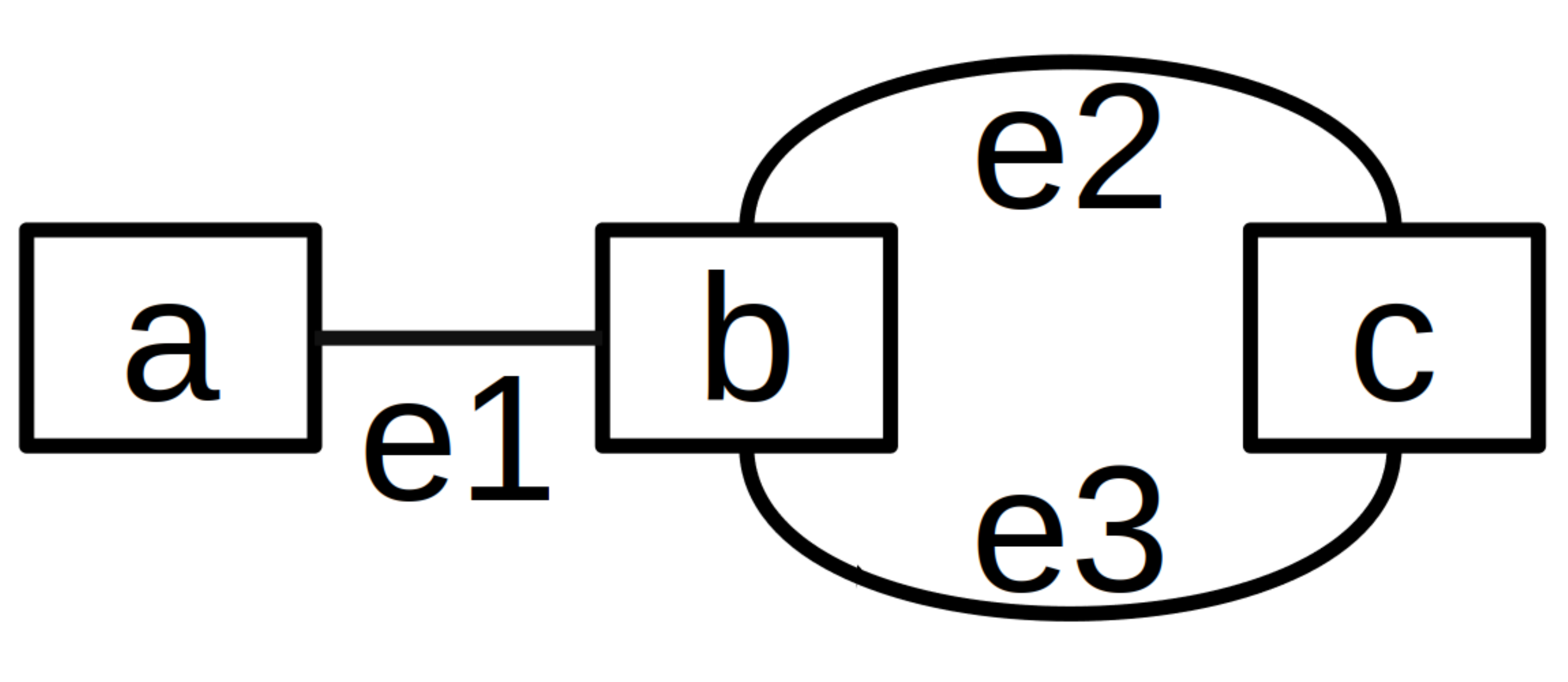}
        \caption{Multigraph}
    \end{subfigure}
    \begin{subfigure}{0.49\linewidth}
        \centering
        \includegraphics[width=1.0\linewidth]{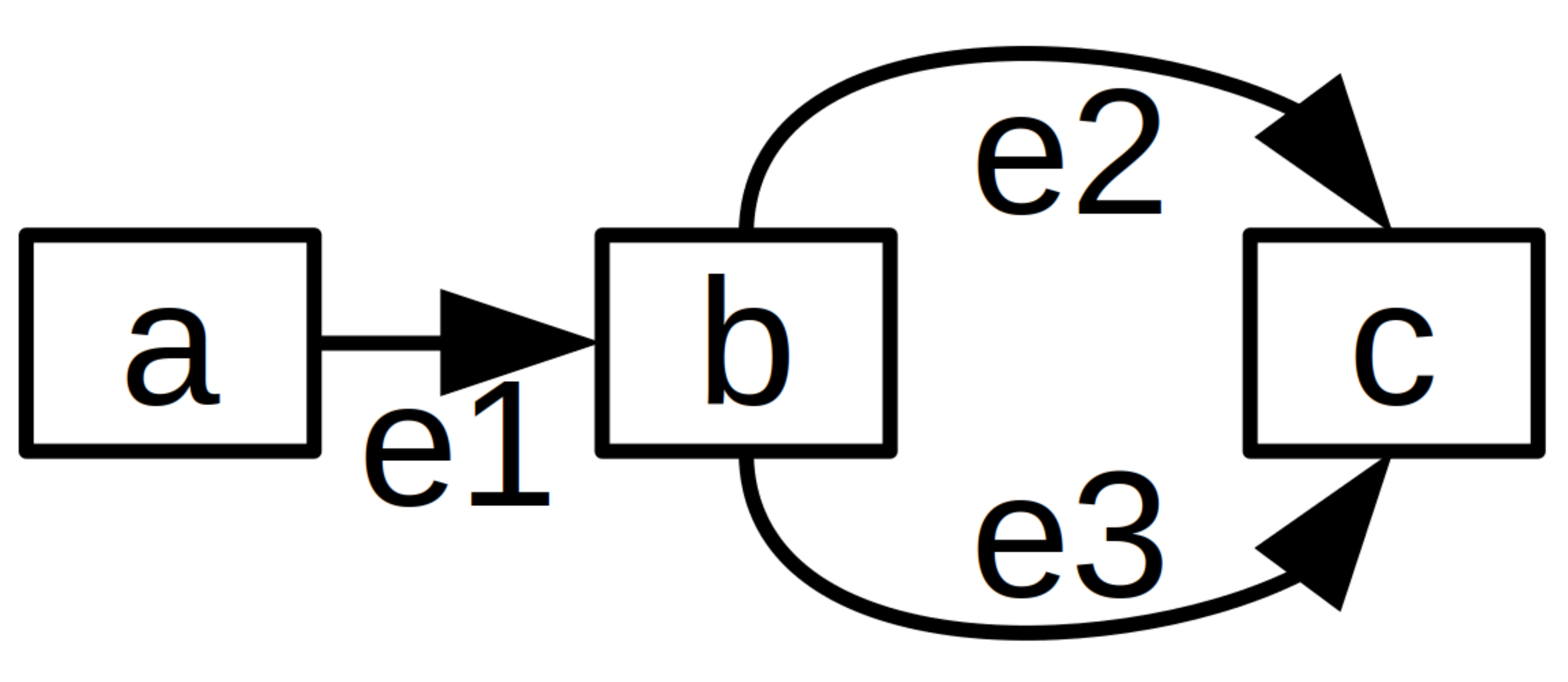}
        \caption{Directed acyclic multigraph}
    \end{subfigure}\hspace*{\fill}
   \begin{subfigure}{0.49\linewidth}
        \centering
        \includegraphics[width=1.0\linewidth]{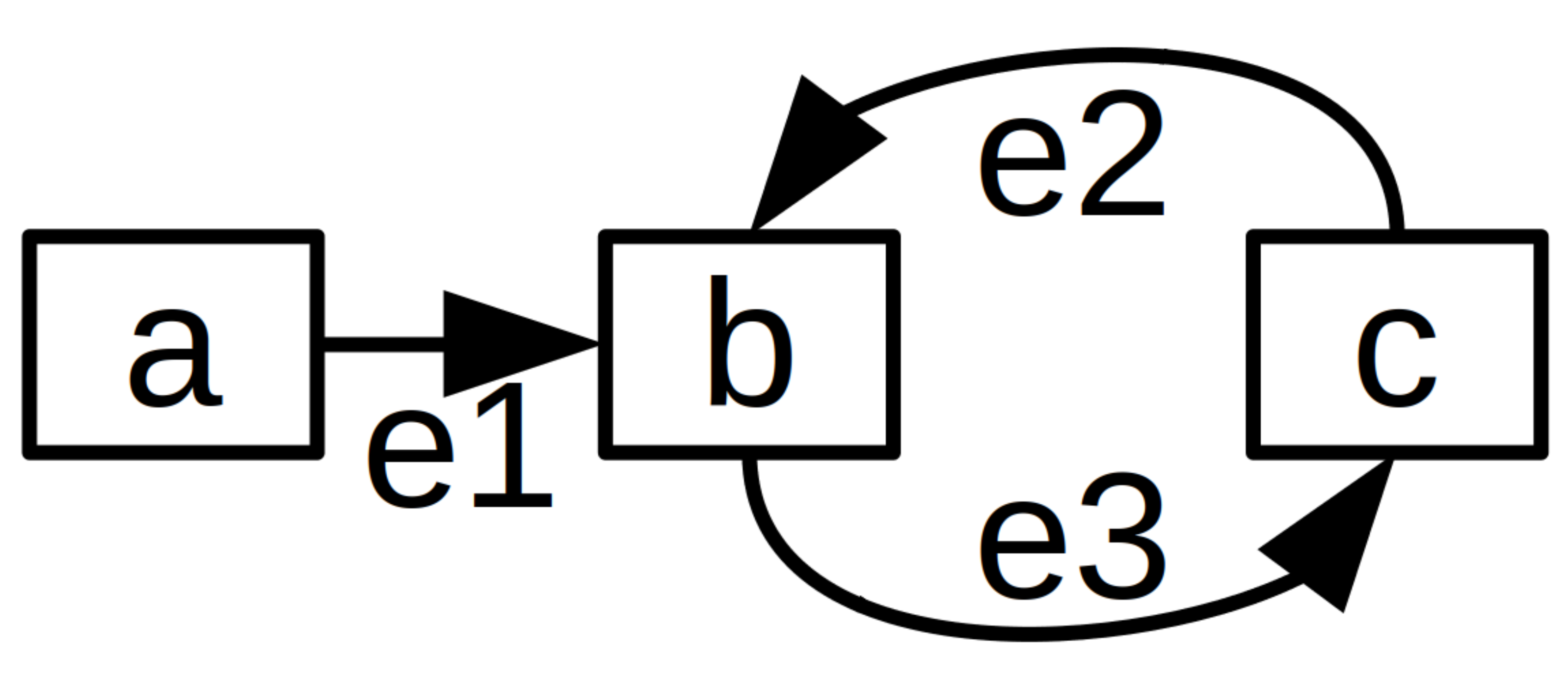}
        \caption{Directed cyclic multigraph}
    \end{subfigure} 
    \caption{Example graphs.}
    \label{fig:GraphExamples}
\end{figure}

\begin{definition}
A graph is a set of nodes connected by edges.
\end{definition} 
\noindent Figure \ref{fig:GraphExamples}(a) shows an example of a graph.

\begin{definition}
A multigraph is a graph wherein at least one pair of nodes is connected by more than one edge.
\end{definition}
\noindent Figure \ref{fig:GraphExamples}(b) shows an example of a multigraph.

\begin{definition}
A directed graph is a graph wherein the edges have direction.
\end{definition}
\noindent Figure \ref{fig:GraphExamples}(c) shows an example of a directed graph: more specifically, a directed multigraph. Arrows denote the directions of the edges. For a directed edge, the target node is the node at the head of the arrow, and the source node is the node at the tail of the arrow. 

\begin{definition}
A cyclic graph is a graph that contains at least one cycle.
\end{definition}
\begin{definition}
An acyclic graph is a graph that contains no cycles.
\end{definition}
\begin{definition}
A cycle is a path on the graph wherein the only repeated nodes are the first and last nodes.
\end{definition}
\noindent Figure \ref{fig:GraphExamples}(d) shows an example of a cyclic graph: more specifically, a directed cyclic multigraph.

It is also convenient to define some node types.
\begin{definition}
A root node is a node in a directed graph that is not the target of any edges.
\end{definition}
\begin{definition}
A leaf node is a node in a directed graph that is not the source of any edges.
\end{definition}
\noindent In Fig. \ref{fig:GraphExamples}, node $a$ is a root node in (c) and (d), and node $c$ is a leaf node in (c) but not in (d).
\begin{definition}
A parent node of node $i$ is any node $j$ in a directed graph such that an edge exists where $i$ is the target and $j$ is the source.
\end{definition}
\begin{definition}
A child node of node $i$ is any node $j$ in a directed graph such that an edge exists where $i$ is the source and $j$ is the target.
\end{definition}
\noindent In Fig. \ref{fig:GraphExamples}, node $a$ is a parent of node $b$ and node $b$ is a child of $a$ in (c) and (d).

Lastly, the concept of the adjacency matrix, $A_{ij}$, for a directed graph is introduced. The adjacency matrix can be constructed using,
\begin{equation*}
    A_{ij} = \begin{cases} 1, &\text{if node $i$ is a parent of node $j$}\\
                           0, &\text{otherwise}. \end{cases}
\end{equation*}
The adjacency matrix has many uses, but in the context of TFC, the adjacency matrix is used to determine if a graph is acyclic or not. If $A_{ij}$ is nilpotent, then the graph is acyclic \cite{GraphTheory}.

%% file: Data/appendixSimpleNL.tex
\chapter{EXTENSION TO NONLINEAR CONSTRAINTS}\label{app:SimpleNLConstraints}

This appendix extends TFC to simple nonlinear constraints and parameterized nonlinear constraints. The extension is accomplished by transforming the nonlinear constraints into linear constraints by introducing extra variables into the constrained expression, which are found in the $\kappa$ terms. These extra variables are constants from the perspective of the constraint operators. Hence, the rest of the theory---derivation of switching functions and projection functionals, extension to $n$-dimensions, and associated mathematical theorems---remains unchanged and still applies as presented in Chapter \ref{chap:tfcTheory}. 

This is one of the critical aspects of abstracting the \ce\ into the symbols associated with the switching-projection form: if one can rewrite constraints as $\C{}[y] = \kappa$ and maintain the relevant mathematical properties, the rest of the theory still applies. This core idea is used in many of the extensions to the TFC theory, e.g., extending beyond the field of real numbers to all mathematical fields, as shown in Appendix \ref{app:Miscellaneous}.

\section{Simple Nonlinear Constraints}
\begin{definition}
    Simple nonlinear constraints are those that can be written as,
    \begin{equation*}
        \psi\Big[\C{}[y]\Big] = \hat{\kappa}
    \end{equation*}
    for some nonlinear function or operator $\psi$ that has a well-defined inverse, i.e., $\psi^{-1}$ exists and can be calculated. 
\end{definition}

Notice that these types of constraints can be rewritten as linear constraints by applying $\psi^{-1}$ to each side:
\begin{equation*}
    \C{}[y] = \psi^{-1}[\hat{\kappa}] = \kappa,
\end{equation*}
where $\kappa = \psi^{-1}[\hat{\kappa}]$. However, doing so may result in multiple solutions. That is, there may be multiple $\kappa$ values that satisfy $\kappa = \psi^{-1}[\hat{\kappa}]$ or even an infinite number of $\kappa$ values.

One option is to write a constrained expression for each $\kappa$ solution. However, this is impossible for an infinite number of $\kappa$ solutions and quickly becomes burdensome when there are multiple simple nonlinear constraints; the total number of \ces\ required for multiple nonlinear constraints is $\prod_i \text{num}(\kappa_i)$ where $\text{num}(\kappa_i)$ is the number of $\kappa$ solutions for the $i$-th nonlinear constraint. 

Fortunately, by introducing new variables into the constrained expression, one can condense the set of \ces\ into one \ce. The following examples highlight some common cases: solutions with $\pm$, a finite number of solutions, and a countably infinite number of solutions.

\begin{example}{Constraints with solutions that contain $\pm$}
Consider the constraint $y^2(0) = 3$. For this constraint, the nonlinear operator raises the function to the power of two, e.g., $\psi[f(x)] = f^2(x)$. Inverting the nonlinear operator results in the solutions $y(0) = \pm\sqrt{3}$. One could write the set of \ces\ that satisfy these constraints,
\begin{align}
    y(x,g(x)) &= g(x) + \sqrt{3} - g(0) \label{eq:Optim_positive} \\
    y(x,g(x)) &= g(x) - \sqrt{3} - g(0). \label{eq:Optim_negative}
\end{align}
However, by including a variable $n$ in the \ce, these two \ces\ can be combined into one,
\begin{equation*}
    y(x,n,g(x)) = g(x) + (-1)^{\mathds{1}_0(n)}\sqrt{3} - g(0),
\end{equation*}
where $\mathds{1}_0(x)$ is the unit step function where the step occurs at $x=0$. If $n < 0$, then Equation \eqref{eq:Optim_negative} is recovered, and if $n > 0$, then Equation \eqref{eq:Optim_positive} is recovered. Hence, $n\in\mathbb{R}$ is now just a variable, and its value dictates which \ce\ branch is used.

Figure \ref{fig:RootExNonlinear} shows the constrained expression plotted with randomly chosen values of $n$ and free functions chosen as polynomials with random coefficients. The solutions of $y(0)$ that satisfy the constraint are shown as black points.
\begin{figure}[H]
    \centering
    \includegraphics[width=0.75\linewidth]{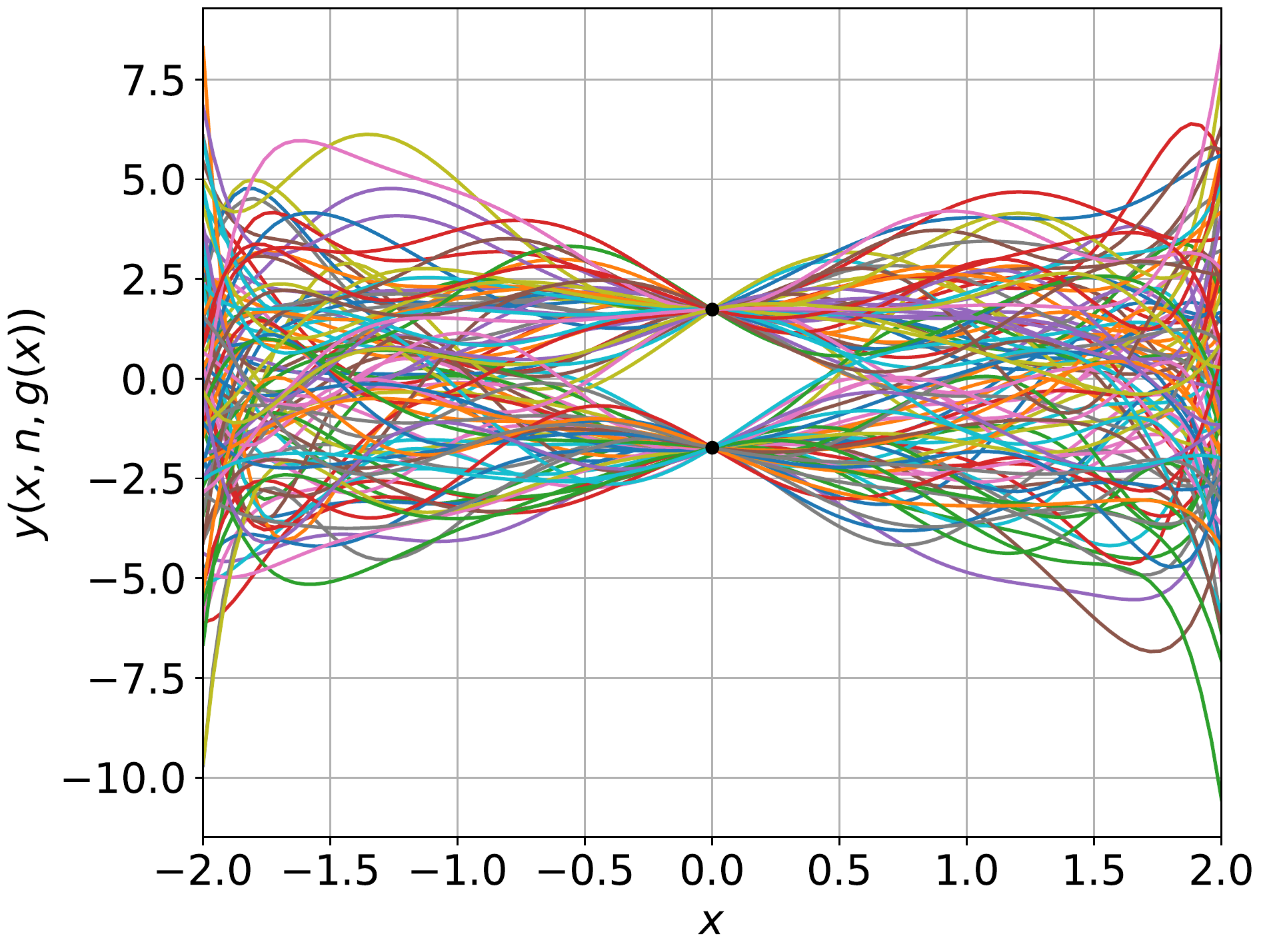}
    \caption{Squared constraint example for randomly chosen $g(x)$ and $n$.}
    \label{fig:RootExNonlinear}
\end{figure}
\end{example}

\begin{example}{Constraints with a finite number of solutions}
Consider the constraint 
\begin{equation*}
    y^3(0)-6y^2(0)+11y(0)=6.
\end{equation*}
Applying $\psi_1^{-1}[\hat{\kappa}]$ yields $\kappa = \{1,2,3\}$, a finite set. The set of \ces\ that satisfy these solutions can be written compactly as,
\begin{equation*}
    y(x,n,g(x)) = g(x) - \kappa[n] - g(0),
\end{equation*}
where $n\in\mathbb{Z}/3\mathbb{Z}$ and $\kappa[n]$ is the $n$-th solution in the set of $\kappa$ that satisfies the constraints---one may think of $\kappa[n]$ as the indexing operation of a zero-indexed array $\kappa$ that contains the solutions of $\psi_1^{-1}[\hat{\kappa}]$. Alternatively, if one prefers to keep $n\in\mathbb{R}$, then this \ce\ can be rewritten as,
\begin{equation*}
    y(x,n,g(x)) = g(x) - \kappa\Big[\nint{\bmod(n,2)}\Big] - g(0).
\end{equation*}

\begin{figure}[H]
    \centering
    \includegraphics[width=0.75\linewidth]{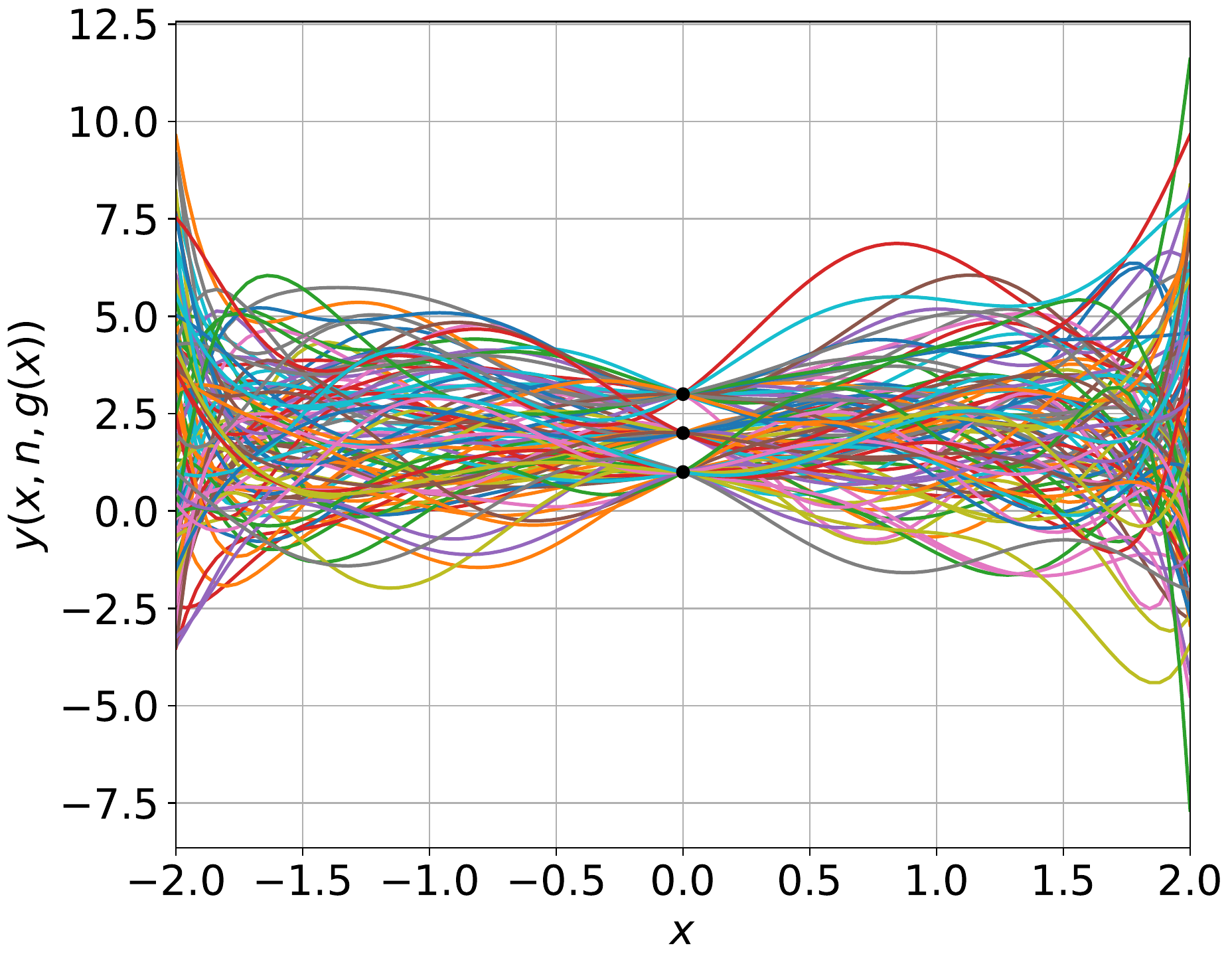}
    \caption{Polynomial constraint example for randomly chosen $g(x)$ and $n$.}
    \label{fig:PolyExNonlinear}
\end{figure}

\noindent The term $\nint{(\bmod(n,2)}$, where $\nint{x}$ rounds $x$ to the nearest integer, forces $n\in\mathbb{R}$ to lie on $\mathbb{Z}/3\mathbb{Z}$.

Figure \ref{fig:PolyExNonlinear} shows the constrained expression plotted with randomly chosen values of $n$ and free functions chosen as polynomials with random coefficients. The solutions of $y(0)$ that satisfy the constraint are shown as black points.
\end{example}

\begin{example}{Constraints with a countably infinite number of solutions}
Consider the constraint,
\begin{equation*}
    \sin\big(y(1)\big) = 0.
\end{equation*}
Applying $\psi_1^{-1}[\hat{\kappa}]$ yields $\kappa = n\pi$ where $n\in\mathbb{Z}$: a countably infinite number of solutions. Similar to the previous example, one can write the \ce\ as,
\begin{figure}[H]
    \centering
    \includegraphics[width=0.75\linewidth]{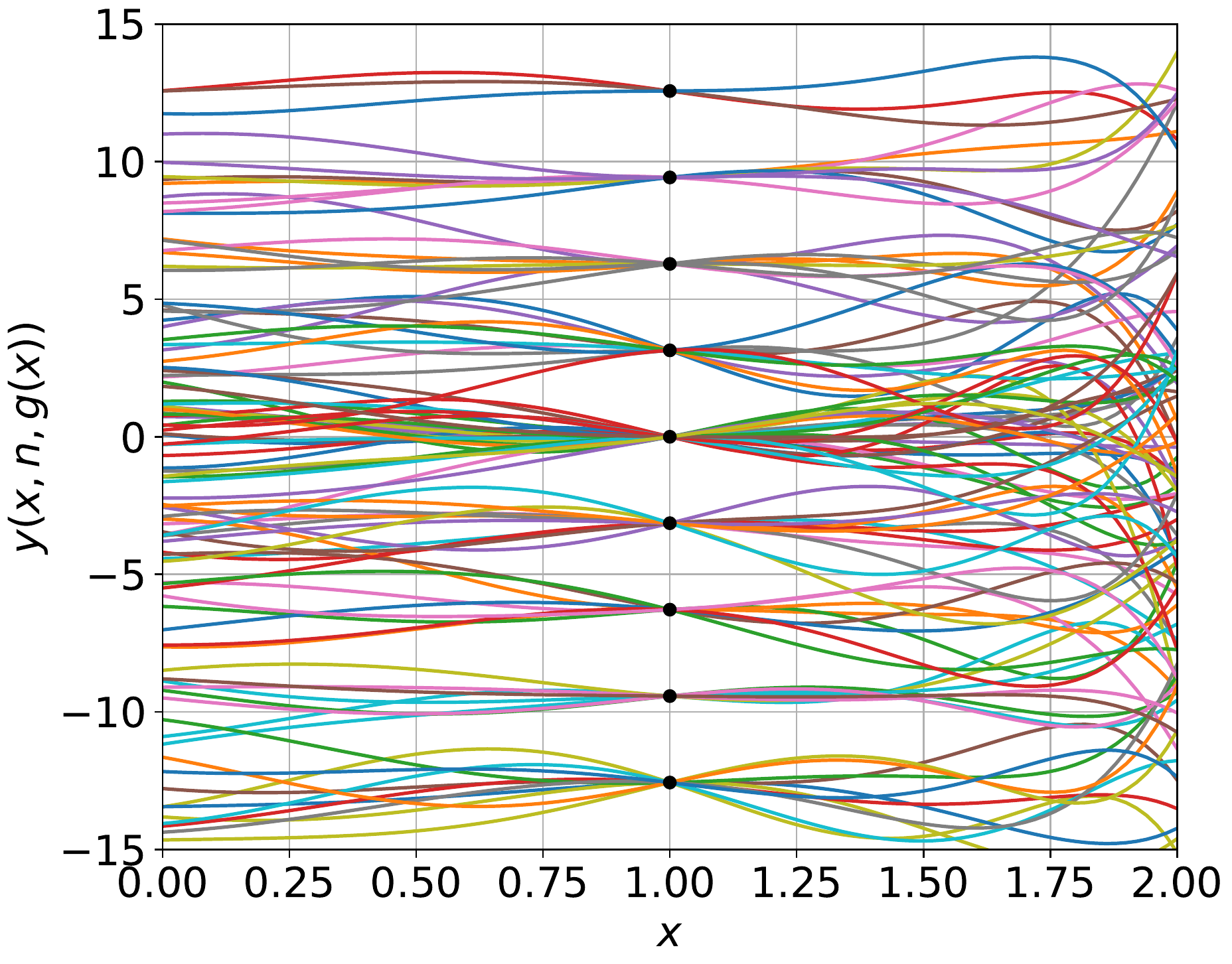}
    \caption{Sine constraint example for randomly chosen $g(x)$ and $n$.}
    \label{fig:SineExNonlinear}
\end{figure}
\begin{equation*}
    y(x,n,g(x)) = g(x) - n\pi - g(1),
\end{equation*}
\noindent where $n\in\mathbb{Z}$ or as,
\begin{equation*}
    y(x,n,g(x)) = g(x) - \nint{n}\pi - g(1),
\end{equation*}
where $n\in\mathbb{R}$.

Figure \ref{fig:SineExNonlinear} shows the constrained expression plotted with randomly chosen values of $n$ and free functions chosen as polynomials with random coefficients. A subset of the solutions of $y(1)$ that satisfy the constraint is shown via black points.

\end{example}

\section{Parameterized Nonlinear Constraints}
Oftentimes, one can parameterize the solution space of the nonlinear constraints and rewrite them as a set of linear constraints. Examples \ref{ex:nonLinParamConstraint} and \ref{ex:nonLinParamConics} demonstrate the idea.

\begin{example}{Simple parameterized nonlinear constraint}\label{ex:nonLinParamConstraint}
Consider the nonlinear constraint,
\begin{equation*}
    y^2(0) + (y_x(0) - 4)^2 = 9.
\end{equation*}
The solutions to this nonlinear constraint lie on a circle that can be parameterized in terms of an unknown, $\theta$:
\begin{equation*}
    y(0) = 3\sin(\theta) \andd y_x(0) = 3\cos(\theta) + 4.
\end{equation*}
Then, these linear constraints can be embedded into a \ce\ using the usual method,
\begin{equation*}
    y(x,\theta,g(x)) = g(x) + 3\sin(\theta) - g(0) + x\Big(3\cos(\theta) + 4 - g_x(0)\Big),
\end{equation*}
where $\theta\in\mathbb{R}$.
\end{example}

\begin{example}{Parameterized nonlinear constraints on conics}\label{ex:nonLinParamConics}
Consider the following nonlinear constraints,
\begin{equation*}
     x^2(0) + \frac{16}{9} y^2(0) + \frac{1}{4} z^2(0) = 1 \andd 100 \big(x(3)-3\big)^2 + 100 y^2(3) - \frac{100}{9}z^2(3) = -1.
\end{equation*}
The solutions of these nonlinear constraints form an ellipsoid at $t=0$ and a hyperboloid of two sheets at $t=3$. Hence, theses nonlinear constraints can be parameterized as,
\begin{align*}
    x(0) &= \sin(\phi)\cos(\theta) &&x(3) = \frac{1}{10}\sinh(|v|)\cos(\psi) + 3 \\
    y(0) &= \frac{3}{4}\sin(\phi)\sin(\theta) &&y(3) = \frac{1}{10}\sinh(|v|)\sin(\psi) \\
    z(0) &= 2\cos(\phi) &&z(3) = (-1)^{\mathds{1}_0(n)}\frac{3}{10}\cosh(|v|),
\end{align*}
where $\phi,\theta,v,\psi,n \in \mathbb{R}$. These parameterized constraints can be embedded into \ces:
\begin{align*}
    x(t,\phi,\theta,v,\psi,g(t)) &= g(t) + \frac{3-t}{3}\Big(\sin(\phi)\cos(\theta) - g(0)\Big) \\
    &\quad + \frac{t}{3}\Big(\frac{1}{10}\sinh(|v|)\cos(\psi) + 3 - g(3)\Big) \\
    y(t,\phi,\theta,v,\psi,g(t)) &= g(t) + \frac{3-t}{3}\Big(\frac{3}{4}\sin(\phi)\sin(\theta) - g(0)\Big) \\
    &\quad+ \frac{t}{3}\Big(\frac{1}{10}\sinh(|v|)\sin(\psi) - g(3)\Big) \\
    z(t,\phi,\theta,v,\psi,g(t)) &= g(t) + \frac{3-t}{3}\Big(2\cos(\phi) - g(0)\Big) \\
    &\quad + \frac{t}{3}\Big((-1)^{\mathds{1}_0(n)}\frac{3}{10}\cosh(|v|) - g(3)\Big).   
\end{align*}

Figure \ref{fig:nonLineConics} shows the \ces\ plotted with randomly chosen values of $\phi$, $\theta$, $v$, $\psi$, and $n$ and free functions that were chosen as polynomials with random coefficients. The nonlinear constraint surfaces are shown as a black ellipsoid and purple hyperboloid.
\begin{figure}[H]
    \centering
    \textattachfile{Figures/Conics.html}{\includegraphics[width=0.75\linewidth]{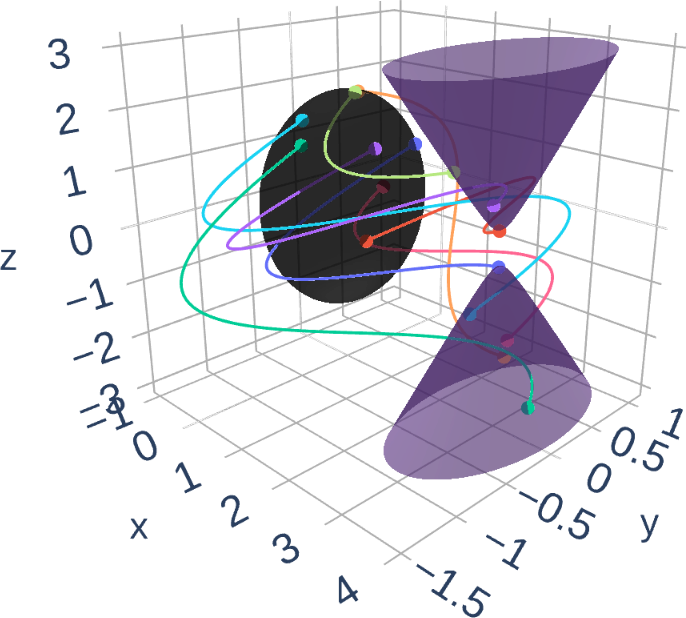}}
    \caption{Parameterized nonlinear constraints on conics. Note, this figure contains an embedded, standalone HMTL version of the plot that can be viewed/downloaded by clicking on it. Doing so may require a dedicated PDF viewer such as Adobe Acrobat or Okular.}
    \label{fig:nonLineConics}
\end{figure}
\end{example}

%% file: Data/appendixInequality.tex

\chapter{ADDING INEQUALITY CONSTRAINTS TO CONSTRAINED EXPRESSIONS}\label{app:Inequality}

Inequality constraints such as $y(x) < f_u(x) \ \forall x\in\Omega$ where $y$ is the dependent variable, $f_u(x)$ is some function that specifies an upper bound, and $\Omega$ is the domain of interest, can also be added to TFC \ces. However, at the time this dissertation is written, the method for incorporating inequality constraints cannot be used in combination with all other constraint types introduced previously, nor have all the mathematical theorems related to the \ce\ been extended to include inequality constraints; thus, they are included here as an appendix, rather than in the main body of the text. 

Early attempts to incorporate inequality constraints utilized the sigmoid function to satisfy the inequality constraints approximately \cite{Inequality-TFC}, but the introduction of an automatic differentiation framework allows them to be incorporated exactly. The enabling component of the automatic differentiation framework is the concept of primitives, which allow a user to specify both a function and its derivative as black boxes: meaning that the derivative specified does not have to be the actual mathematical derivative of the original function. The Heaviside function is a prime example \cite{JaxGithub}; the function value is encoded as,
\begin{equation*}
    \mathds{1}(x,x_1) = \begin{cases} 0,& x<x_1 \\ x_1,& x=0 \\ 1,& x > x_1 \end{cases}
\end{equation*}
but the derivative is encoded as,
\begin{equation*}
    \frac{\dd \mathds{1}}{\dd x} = 0.
\end{equation*}
Mathematically this is not correct as the derivative is really the Dirac delta function \cite{HeavisideDefinition}, but encoding it in this way allows one to incorporate inequality constraints into \ces.

Let $\mathds{1}_0(x) = \mathds{1}(x,0)$, which is equivalent to the unit step function where the step occurs at $x=0$. One can think of this function as the mathematical equivalent of a switch or gate when it is composed with multiplication. For example,
\begin{equation*}
    f(x) = g(x)\mathds{1}(x,0) = g(x)
    \mathds{1}_0,
\end{equation*}
will return $g(x)$ when $x > 0$ and $0$ otherwise. This switching behavior is exactly the desired behavior needed to implement inequality constraints. 

Consider the following two inequality constraints,
\begin{equation*}
    y(x) > f_\ell(x) \andd y(x) < f_u(x),
\end{equation*}
where $f_\ell(x)$ and $f_u(x)$ are the lower and upper bound functions respectively. Using the switch-like behavior of $\mathds{1}_0$, it is straightforward to write a functional that maintains a free function, i.e., a \ce, and obeys these two inequality constraints,
\begin{equation}\label{eq:InequalityForm}
    y(x,g(x)) = g(x) + \Big(f_u(x)-g(x)\Big)\mathds{1}_0\Big(g(x)-f_u(x)\Big) + \Big(f_\ell(x)-g(x)\Big)\mathds{1}_0\Big(f_\ell(x)-g(x)\Big).
\end{equation}
One can prove that this form satisfies the two inequality constraints via brute-force by checking the three possible cases:
\begin{enumerate}
    \item $g(x) < f_\ell(x) \to y(x,g(x)) = g(x) + \Big(f_u(x)-g(x)\Big) (0) + \Big(f_\ell(x)-g(x)\Big) (1) = f_\ell(x)$
    \item $f_\ell(x) \leq g(x) \leq f_u(x) \to y(x,g(x)) = g(x) + \Big(f_u(x)-g(x)\Big) (0) + \Big(f_\ell(x)-g(x)\Big) (0) \\ \text{\phantom{$f_\ell(x) \leq g(x) \leq f_u(x) \to y(x,g(x)) $}}= g(x)$
    \item $f_u(x) < g(x) \to y(x,g(x)) = g(x) + \Big(f_u(x)-g(x)\Big) (1) + \Big(f_\ell(x)-g(x)\Big) (0) = f_u(x)$
\end{enumerate}
Furthermore, the Heaviside derivative overridden by the automatic differentiation program produces the desired behavior in the derivative of the bounded \ce,
\begin{equation*}
    y_x(x,g(x)) = \begin{cases} \frac{\dd f_\ell}{\dd x}(x),& g(x) < f_\ell(x) \\ 
                                g_x(x),& f_\ell(x) \leq g(x) \leq f_u(x)\\
                                \frac{\dd f_u}{\dd x}(x),& f_u(x) < g(x);
                  \end{cases}
\end{equation*} 
that is, the derivative of the \ce\ is equal to the derivatives of the lower and upper bound functions when they are active and equal to the derivative of $g(x)$ otherwise. 

Although inequality constraints cannot yet be combined with all the constraint types introduced earlier, they can be combined with the most commonly occurring constraint: point constraints. Let $\hat{y}(x,g(x))$ be a \ce\ satisfying some set of point constraints that are consistent with the inequality constraints,
\begin{equation*}
    y(x) > f_\ell(x) \andd y(x) < f_u(x).
\end{equation*}
Then, $\hat{y}(x,g(x))$ can be substituted as the free function into the inequality \ce\ given earlier,
\begin{equation}\label{eq:IneqAndEq}
\begin{aligned}
    y(x,g(x)) =\ &\hat{y}(x,g(x)) + \Big(f_u(x)-\hat{y}(x,g(x))\Big)\mathds{1}_0\Big(\hat{y}(x,g(x))-f_u(x)\Big)\\
    &+ \Big(f_\ell(x)-\hat{y}(x,g(x))\Big)\mathds{1}_0\Big(f_\ell(x)-\hat{y}(x,g(x))\Big),
\end{aligned}
\end{equation}
and the result satisfies both the equality and inequality constraints. The inequality constraints are satisfied because Equation \eqref{eq:InequalityForm} satisfies them for any free function, including $\hat{y}(x,g(x))$, and the equality constraints can be shown to be satisfied by simply applying the constraint operator to the constrained expression (similar to the proof of Theorem \ref{thrm:UniCe}),
\begin{align*}
    \C{i}[y(x,g(x))] &= \C{i}[\hat{y}(x,g(x))] + \C{i}\Big[\Big(f_u(x)-g(x)\Big) (0)\Big] + \C{i}\Big[\Big(f_\ell(x)-g(x)\Big)(0)\Big] \\
    &= \kappa_i + 0 + 0 = \kappa_i.
\end{align*}
Note that for any consistent set of constraints, the $\mathds{1}_0$ functions must yield zero at the equality constraint locations. The theorem on the existence of $g(x)$ (Theorem \ref{thrm:UniGExists}) can easily be extended as well. Let $f(x)$ be some function that satisfies the constraints, then,
\begin{align*}
    y(x,f(x)) &= \hat{y}(x,f(x)) + \Big(f_u(x)-f(x)\Big) (0) + \Big(f_\ell(x)-f(x)\Big)(0) \\
    &= f(x) + 0 + 0 \\
    &= f(x).
\end{align*}
Thus, for any function satisfying the constraints, $f(x)$, there exists at least one free function $g (x) = f (x)$, such that the constrained expression is equal to the function satisfying the constraints, i.e., $y (x, f (x)) = f (x)$.

\begin{figure}[!ht]
    \centering
    \begin{minipage}{0.49\linewidth}
    \includegraphics[width=\linewidth]{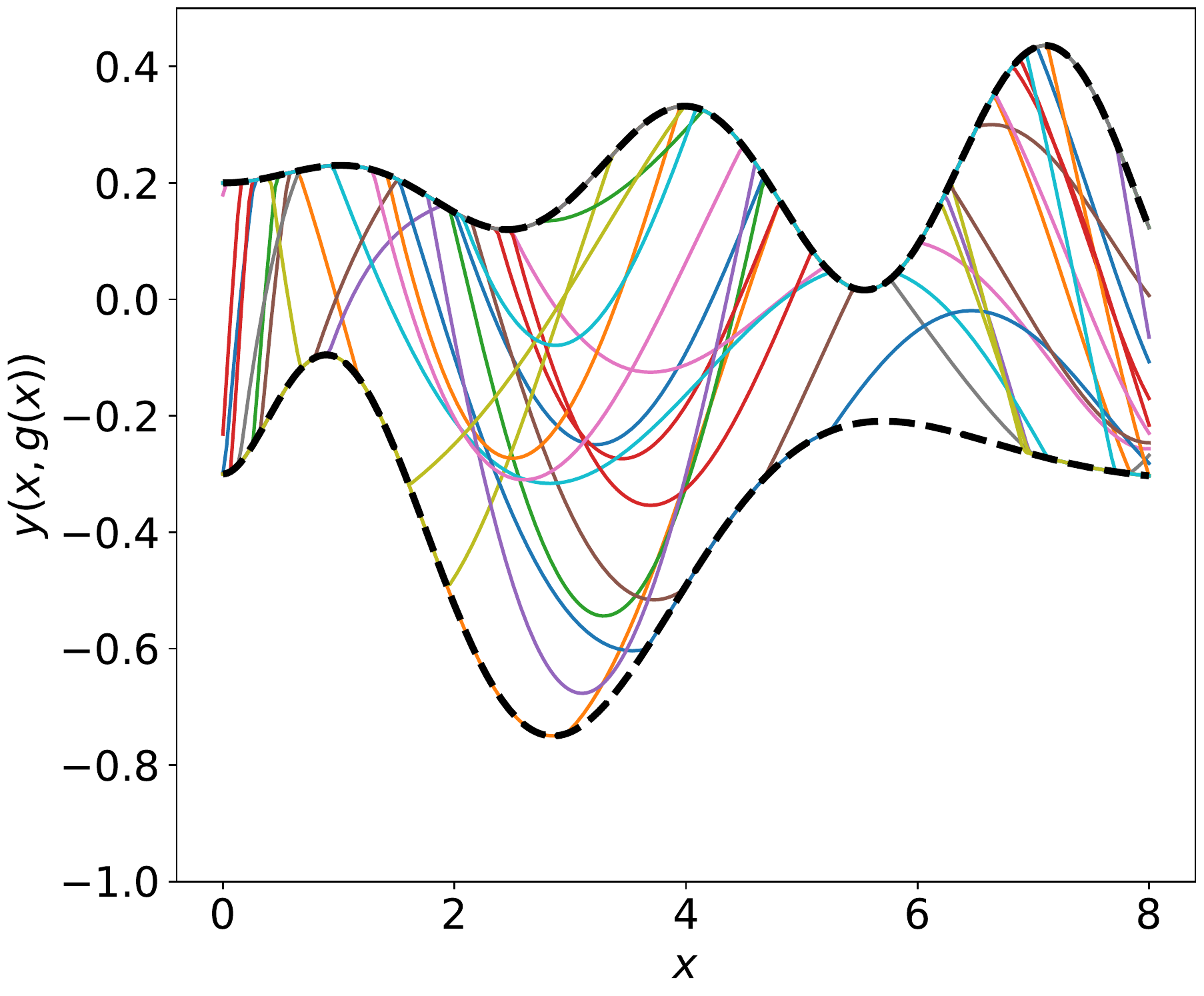}
    \caption{Inequality constraints only.}
    \label{fig:ineqaulityOnly}
    \end{minipage}%
    \begin{minipage}{0.49\linewidth}
    \includegraphics[width=\linewidth]{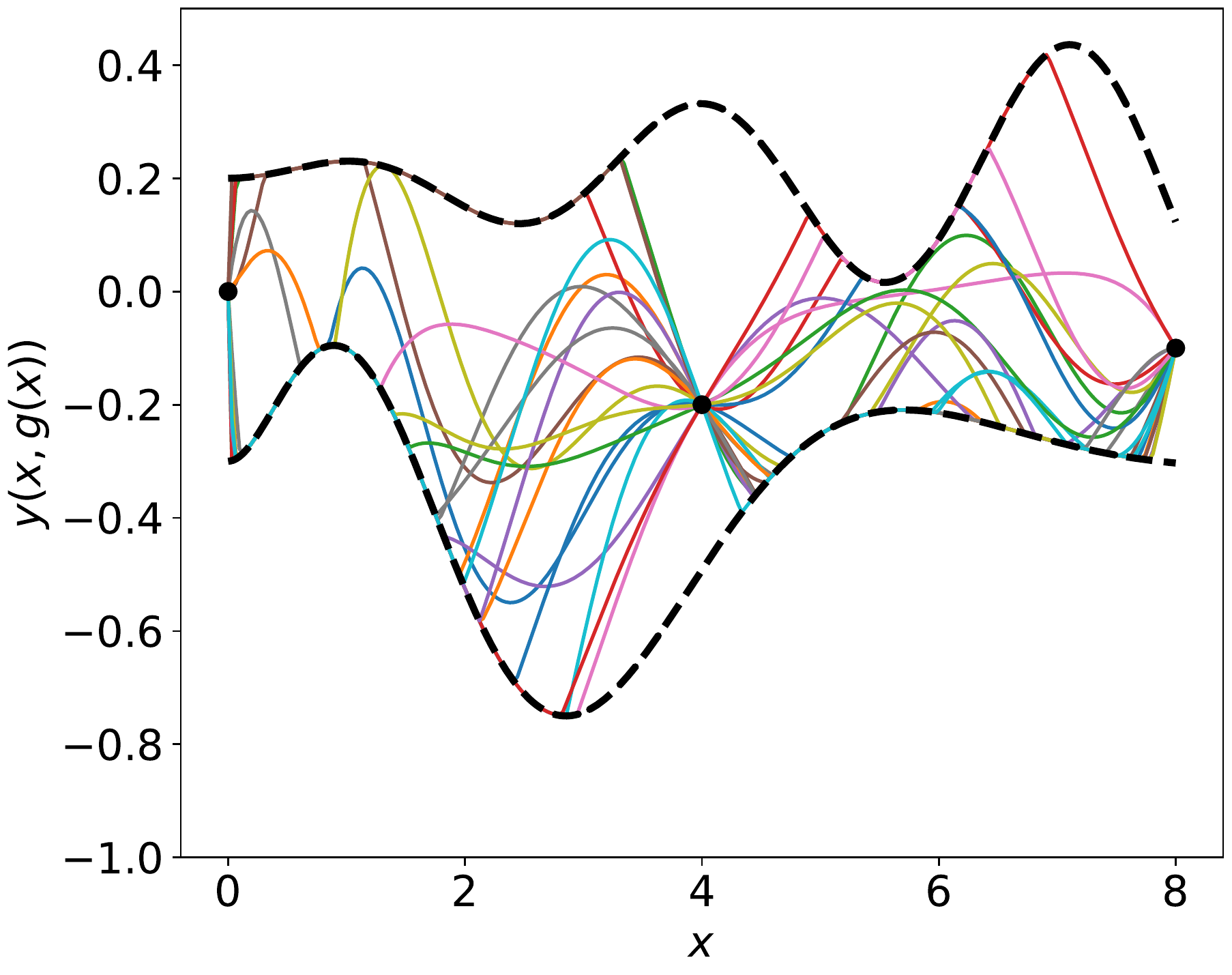}
    \caption{Inequality and value-level\\equality constraints.}
    \label{fig:ineqaulityAndEquality}
    \end{minipage}
\end{figure}

Figures \ref{fig:ineqaulityOnly} and \ref{fig:ineqaulityAndEquality} visually show that Equations \ref{eq:InequalityForm} and \ref{eq:IneqAndEq} work, respectively. Each figure contains randomly generated upper and lower bounds, shown as black, dashed lines, and randomly generated free functions, shown as solid, colored lines. In addition, Figure \ref{fig:ineqaulityAndEquality} contains randomly generated point constraints shown via black dots.

%% file: Data/appendixDividingDomain.tex

\chapter{SPLITTING THE DOMAIN}\label{app:DivideDomain}

When the solution of a differential equation has steep gradients, it is oftentimes difficult to describe the solution over the whole domain using one \ce. In these cases, it is convenient to split the domain into sections and enforce continuity at the intersections via the \ces. For example, consider the following differential equation that describes convection and diffusion processes:
\begin{equation*}
    y_{xx} - P_e y_x = 0
\end{equation*}
subject to,
\begin{equation*}
    y(0) = 1 \andd y(1) = 0,
\end{equation*}
where $x\in[0,1]$, $P_e$ is the Peclet number, and the analytical solution is,
\begin{equation*}
    y = \frac{1 - e^{P_e(x-1)}}{1 - e^{-P_e}}.
\end{equation*}

\begin{figure}[!ht]
    \centering
    \includegraphics[width=0.7\linewidth]{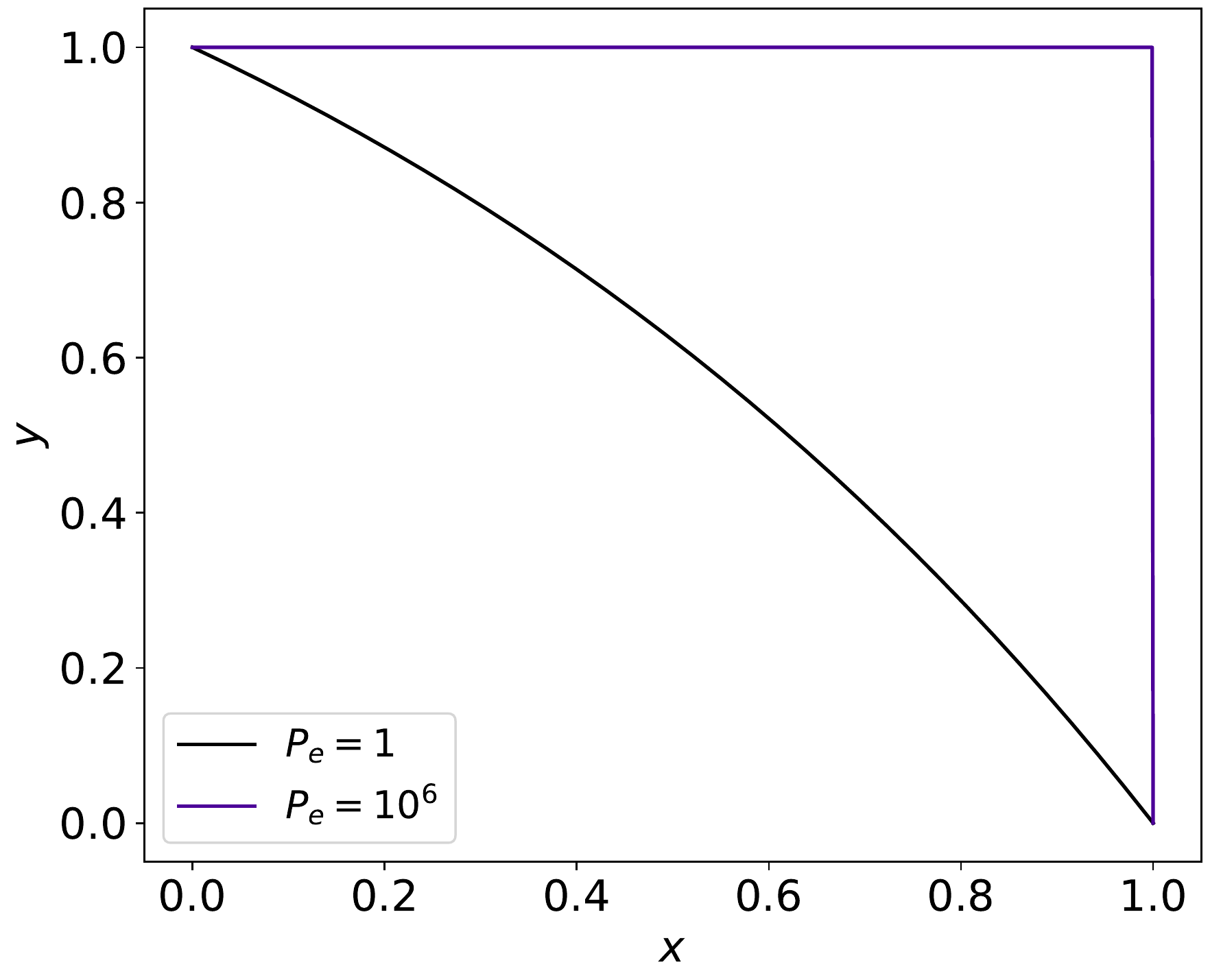}
    \caption{Analytical solutions of the convection-diffusion equation with different Peclet numbers.}
    \label{fig:convectionDiffusionAnalytical}
\end{figure}

As the Peclet number increases, the solution's gradients become larger. To illustrate, Figure \ref{fig:convectionDiffusionAnalytical} shows the analytical solutions to the convection-diffusion equation with two different Peclet numbers, $P_e=1$ and $P_e=10^6$.
When the Peclet number is low enough, one can estimate the solution well using only one \ce\ for the entire domain,
\begin{equation*}
    y(x,g(x)) = g(x) + (1-x)(1-g(0)) - x g(1).
\end{equation*}
However, as the Peclet number increases, a better estimation is obtained with two \ces. These two \ces\ are for the subdomains $x_1\in[0,x_p]$ and $x_2\in[x_p,1]$, where $x_p$ is the point of intersection between the two. For this differential equation, choosing the point $x_p$ arbitrarily does not increase the accuracy noticeably; hence, the point $x_p$ is an unknown that will be solved as part of the solution process. 

Since the domains of the two \ces\ are dependent on $x_p$, it is simplest to write them on the basis function domain,
\begin{align*}
    \p{1}{y}(z,g(z)) &= \p{1}{g}(z) + \frac{1-2z+z^2}{4}\Big(1-\p{1}{g}(z_0)\Big) + \frac{3+2z-z^3}{4}\Big(y_p-\p{1}{g}(z_f)\Big)\\
    &\quad +\frac{z^2-1}{2}\Big(\dd y_p/c_1-\frac{\dd \p{1}{g}}{\dd z}(z_f)\Big)\\
    \p{2}{y}(z,g(z)) &= \p{2}{g}(z) + \frac{3-2z-z^2}{4}\Big(y_p-\p{2}{g}(z_0)\Big) - \frac{1+2z+z^2}{4}\p{2}{g}(z_f)\\
    &\quad +\frac{1-z^2}{2}\Big(\dd y_p/c_2-\frac{\dd \p{2}{g}}{\dd z}(z_0)\Big)
\end{align*}
where $\p{1}{g}(z)$ is the free function for the first \ce, $\p{2}{g}(z)$ is the free function for the second constrained expression, $y_p$ and $\dd y_p$ are the value and derivative of the intersection point, $z\in[z_0,z_f]$ is the free function domain, and $c_1$ and $c_2$ are the constants in the linear maps from the problem domain to the basis function domain; for this problem, Legendre orthogonal polynomials are used, so $z\in[-1,1]$. In addition, the constants in the mapping function can be expressed as,
\begin{equation*}
    c_1 = \frac{2}{x_p} \andd c_2 = \frac{2}{1-x_p}.
\end{equation*}
Since the two \ces\ are written on the basis function domain, the differential equation must be modified,
\begin{equation*}
    c_k^2y_{zz}-c_kP_ey_z = 0,
\end{equation*}
where $k=1$ if $x < x_p$ and $k=b$ if $x > x_p$.

Now, the unknown coefficients in the two free functions, $\p{1}{\B{\xi}}$ and $\p{2}{\B{\xi}}$,  and the intersection point and derivative values, $x_p$, $y_p$, and $\dd y_p$, can be used to reduce the residual of the differential equation at each point in the discretized domain: the domain here is broken up into $200$ points per \ce, for a total of $400$ points. However, there is one more nuance to this problem: solving the problem as-is with nonlinear least-squares tends to diverge unless a good initial guess is provided. This divergent behavior is related to trying to solve for $x_p$, which frequently takes on values outside the domain if left unchecked. The author has identified two options that fix this divergent behavior: 
\begin{enumerate}
    \item Remove $x_p$ from the nonlinear least-squares optimizer and estimate it using a separate, exterior optimization scheme such as a genetic algorithm. 
    \item Modify the nonlinear least squares to bound the values that $x_p$ can take.
\end{enumerate}

In this section, the second option is used, but the author has verified that the first method also works. When using the second option, it is tempting to simply perform an update after each iteration of the nonlinear least-squares that bounds $x_p$ to values within the domain. However, this simple change still results in divergent behavior much of the time. In terms of convergence, a better option is to use concepts from inequality constraint embedding (see appendix \ref{app:Inequality} for more details) to bound $x_p$. For this differential equation, $x_p$ was chosen to be,
\begin{equation*}
    x_p = \hat{x}_p + (f_u-\hat{x}_p)\mathds{1}_0(\hat{x}_p-f_u) + (f_\ell-\hat{x}_p)\mathds{1}_0(f_\ell-\hat{x}_p),
\end{equation*}
where $x_p$ is the value used in the \ce\ and $\hat{x}_p$ is the unknown used in the nonlinear least squares; $f_\ell=1\times{10}^{-3}$ and $f_u = 1-1\times{10}^{-3}$ are the lower and upper bounds, respectively, on $x_p$. This simple change results in a nonlinear least-squares that converges.

To demonstrate the benefits of the domain splitting technique, Table \ref{tab:convectionDiffusionError} shows the maximum and mean error on a test set of $1,000$ evenly spaced data points per constrained expression when using a single constrained expression for the whole domain and when using two constrained expression and a split domain as described above. The results are shown for $P_e = 1$ and $P_e = 10^6$. Each case used $200$ training points per constrained expression, and Legendre polynomials up to degree $190$ as the free function.
\begin{table}[!ht]
\centering
\caption{Convection-diffusion equation error: whole vs. split domain.}
\label{tab:convectionDiffusionError}
\begin{tabular}{|c|c|c|c|c|} 
\hline
 \makecell{$P_e$} & \multicolumn{2}{c|}{Whole Domain} & \multicolumn{2}{c|}{Split Domain}\\ \hline
  & \makecell{Maximum Error} & \makecell{Mean Error} & \makecell{Maximum Error} & \makecell{Mean Error} \\\hline
 $1$ & $2.22\times10^{-16}$ & $5.62\times10^{-17}$ & $4.44\times10^{-16}$ & $8.33\times10^{-17}$\\ \hline
 $10^6$ & $1.00$ & $4.99\times10^{-1}$ & $8.61\times10^{-12}$ & $1.10\times10^{-14}$\\ \hline
\end{tabular}
\end{table}
Table \ref{tab:convectionDiffusionError} shows that the split domain and whole domain approaches produce similar error values for the $P_e=1$ case when the solution's gradients are relatively small. However, when $P_e=10^6$ and the gradients are larger, the split domain's errors are orders of magnitude lower than when using the whole domain: the maximum error is $12$ orders of magnitude lower, and the average error is $13$ orders of magnitude lower. 

%% file: Data/appendixOthBasisFunc.tex

\chapter{ORTHONORMAL BASIS FUNCTIONS}\label{app:BasisFunctions}

This appendix provides the reader with an elementary understanding of orthogonal basis functions. Any reader interested in this subject may refer to Reference \cite{MultiVarOrthPolyBook} for a more in-depth understanding. In essence, basis functions are for a function space what vectors are for a vector space. In other words, a linear combination of basis functions spans the function space, just as a linear combination of basis vectors spans the vector space. Thus, a linear combination of basis functions is a useful free function choice for optimization problems.

\section{Mathematical Preliminaries}

This section introduces some mathematical preliminaries needed to understand the properties of basis functions, and in particular, the properties of orthogonal basis functions. The content introduced here is designed to give the reader a basic understanding and will only scratch the surface of this field of mathematics. As such, when appropriate, references will be provided so that the reader can delve deeper into these topics if desired. Moreover, this section assumes the reader is familiar with the properties of vector spaces. If the reader is unfamiliar with these topics, then they may consider reading Reference \cite{TriebelFunctionSpacesBook} for function spaces and reviewing the portion of Reference \cite{GilbertStrangLinearAlgebra} dedicated to vector spaces. 

This dissertation is primarily concerned with function spaces that can be used to describe continuous, non-infinite functions, as these will be particularly useful for describing the solutions of differential equations: the extended Lebesgue spaces, also known as $L^{pe}$ spaces, are the function spaces that contain these functions. The extended Lebesgue spaces are defined based on a generalization of the $p$-norm used to describe vector spaces. Recall that the $p$-norm for a vector is
\begin{equation*}
    ||\B{x}||_p = \bigg(\sum_{k=1}^n |x_k|^p\bigg)^{1/p},
\end{equation*}
where $\B{x}\in\mathbb{R}^n$ is an arbitrary vector, $x_k$ are the components of $\B{x}$, and $p \geq 1$. The $p$-norm of functions is,
\begin{equation*}
    ||f(z)||_p = \bigg( \int_\Omega |f(z)|^p \dd z \bigg)^{1/p},
\end{equation*}
for some arbitrary function $f (z)$ defined on the domain $\Omega$. Note that this $p$-norm can also be defined with a measure $\dd\mu (z)$, in which case the $p$-norm is written as
\begin{equation*}
    ||f(z)||_p = \bigg(\int_\Omega |f(z)|^p \dd \mu (z) \bigg)^{1/p}.
\end{equation*}
The rigorous mathematical definition of a measure will not be discussed here; the interested reader can refer to Reference \cite{MeasureTheoryBook} for more information. For the material in this dissertation, it is sufficient to note that the measure $\dd\mu(z) = W(z) \dd z$ where $W(z) \geq 0 \, \forall\ z\in\Omega$. The measure for a function is analogous to the weights in a weighted vector norm. An arbitrary function, $f(z)$, defined over the domain $\Omega$ is part of the $L^{pe}(\Omega,\mu)$ space if
\begin{equation*}
    ||f(z)||_p = \bigg(\int_\Omega |f(z)|^p \dd \mu (z) \bigg)^{1/p} < \infty.
\end{equation*}
This appendix will focus on basis functions in the $L^{2e}$ space, i.e., for $p=2$.

The generalization of the $p$-norm is sufficient for describing which functions are in the $L^{pe}$ space. However, the $p$-norm gives no information about the orthogonality of two functions. For this, an inner product is needed. Fortunately, the $L^{2e}(\Omega,\mu)$ space already comes equipped with an inner product,
\begin{equation*}
    \langle f, g \rangle = \int_\Omega f(z) \, g(z) \dd \mu (z),
\end{equation*}
where $f (z)$ and $g (z)$ are arbitrary functions in the $L^{2e}(\Omega,\mu)$ space, and $\langle f,g \rangle$ is used to denote an inner product between these functions. The functions $f$ and $g$ are considered orthogonal if $\langle f, g \rangle = 0$. Just as orthogonal basis vectors can be convenient for describing an arbitrary vector in a vector space, so too are orthogonal basis functions for describing an arbitrary function in a function space\footnote{In addition, choosing orthogonal basis functions can also guarantee certain solution properties.}.

Of course, spanning the entirety of $L^{2e}$ space would require an infinite number of basis functions, as the dimension of the $L^{2e}$ space is infinite. Thus, to make problems computationally tractable, a finite number, $m$, of basis functions is used. In general, as the number $m$ increases, the error between the estimated and actual solution will decrease. Finally, note that the basis set domain need not coincide with the domain of the problem. If a bijective map exists that transforms the basis function domain into the problem domain, then that basis may be used to describe the problem's solution. This notion is used frequently throughout this dissertation. 

Based on the description of orthogonal basis function sets thus far, one has two parameters that can be used to describe a basis set for $L^{2e}$:
\begin{enumerate}
\item The domain on which the basis is defined, $\Omega$.
\item The measure used for the basis, $\mu$.
\end{enumerate}
In the following sections, some frequently used orthogonal basis sets will be presented. The presentation will include the domain and measure for each set and recursive generating functions for the set if they exist. The section that follows explains how to extend these basis sets to the multivariate case and concludes with a table that summarizes all the basis functions presented.

\section{Chebyshev Orthogonal Polynomials}

Chebyshev orthogonal polynomials are two sets of basis functions, the first and the second kind. They are usually indicated as $T_k (z)$ and $U_k (z)$, respectively. This section summarizes the main properties of the first kind, $T_k (z)$, only, which are defined on the domain $z\in[-1,+1]$ and with the measure $\dd \mu (z) = \dfrac{1}{\sqrt{1 - z^2}} \dd z$. These polynomials can be generated using the following useful recursive function,\footnote{Note that in this recursive formulation and those that follow, the $z$ argument is dropped for clarity, i.e., $T_k(z)\to T_k$.}
\begin{equation*}
    T_{k + 1} = 2 \, z \, T_k - T_{k - 1} \qquad \text{starting from:} \; \begin{cases} T_0 = 1, \\ T_1 = z.\end{cases}
\end{equation*}
Also, all the derivatives of Chebyshev orthogonal polynomials can be computed recursively, starting from
\begin{equation*}
    \dfrac{\dd T_0}{\dd z} = 0, \quad \dfrac{\dd T_1}{\dd z} = 1 \qquad \text{or} \qquad \dfrac{\dd^d T_0}{\dd z^d} = \dfrac{\dd^d T_1}{\dd z^d} = 0 \quad \forall \; d > 1,
\end{equation*}
and then using,
\begin{align*}
    \dfrac{\dd T_{k+1}}{\dd z} &= 2 \, \left(T_k + z \, \dfrac{\dd T_k}{\dd z}\right) - \dfrac{\dd T_{k-1}}{\dd z} \\ 
    \dfrac{\dd^2 T_{k+1}}{\dd z^2} &= 2 \left(2 \, \dfrac{\dd T_k}{\dd z} + z \, \dfrac{\dd^2 T_k}{\dd z^2}\right) - \dfrac{\dd^2 T_{k-1}}{\dd z^2} \\
    &\ \ \vdots \\ 
    \dfrac{\dd^d T_{k+1}}{\dd z^d} &= 2 \left( d \, \dfrac{\dd^{d-1} T_k}{\dd z^{d-1}} + z \, \dfrac{\dd^d T_k}{\dd z^d}\right) - \dfrac{\dd^d T_{k-1}}{\dd z^d} \quad \forall \; d \ge 1.
\end{align*}
for $k \ge 1$. 
The integral of $ T_k (z)$ has the following useful property,
\begin{equation*}
    \int_{-1}^{+1} T_k (z) \dd z = \left\{\begin{array}{lcl} = 0 & {\rm if} & k = 1 \\ = \dfrac{(-1)^k + 1}{1 - k^2} & {\rm if} & k \ne 1\end{array}\right.
\end{equation*}
while the inner product of two Chebyshev orthogonal polynomials satisfies the orthogonality property,
\begin{equation*}
    \langle T_i (z), T_j (z)\rangle = \int_{-1}^{+1} T_i (z) \, T_j (z) \, \dfrac{1}{\sqrt{1 - z^2}} \dd z = \begin{cases} = 0 & {\rm if} \quad i \ne j \\ = \pi & {\rm if} \quad i = j = 0 \\ = \pi/2 & {\rm if} \quad i = j \ne 0\end{cases}.
\end{equation*}

\begin{figure}[ht]
    \centering\includegraphics[width=0.9\linewidth]{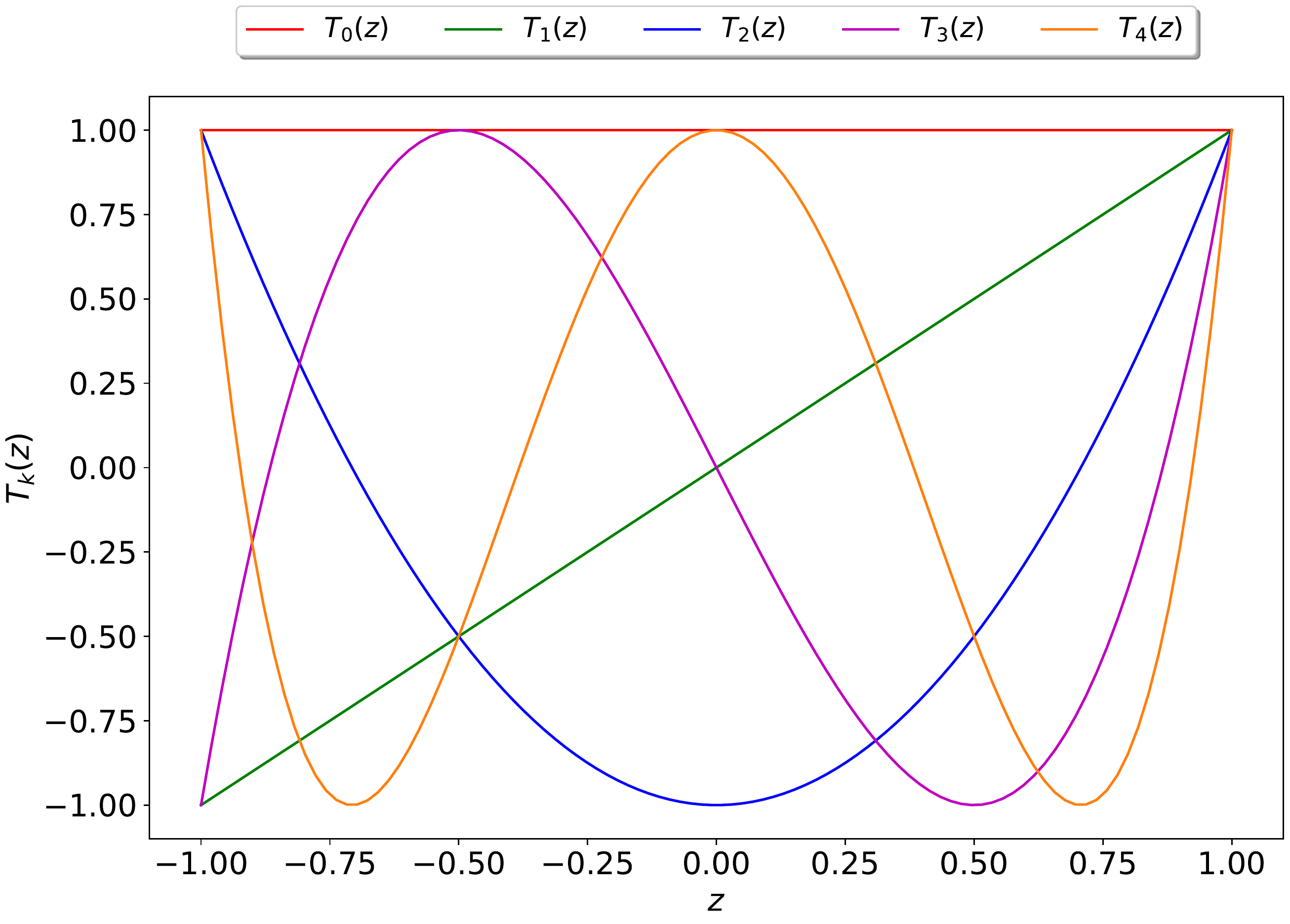}
    \caption{First five Chebyshev orthogonal polynomials.}
    \label{fig:ICOPfig}
\end{figure}

Figure \ref{fig:ICOPfig} shows the first five Chebyshev orthogonal polynomials.

\section{Legendre Orthogonal Polynomials}

\begin{figure}[ht]
    \centering\includegraphics[width=0.9\linewidth]{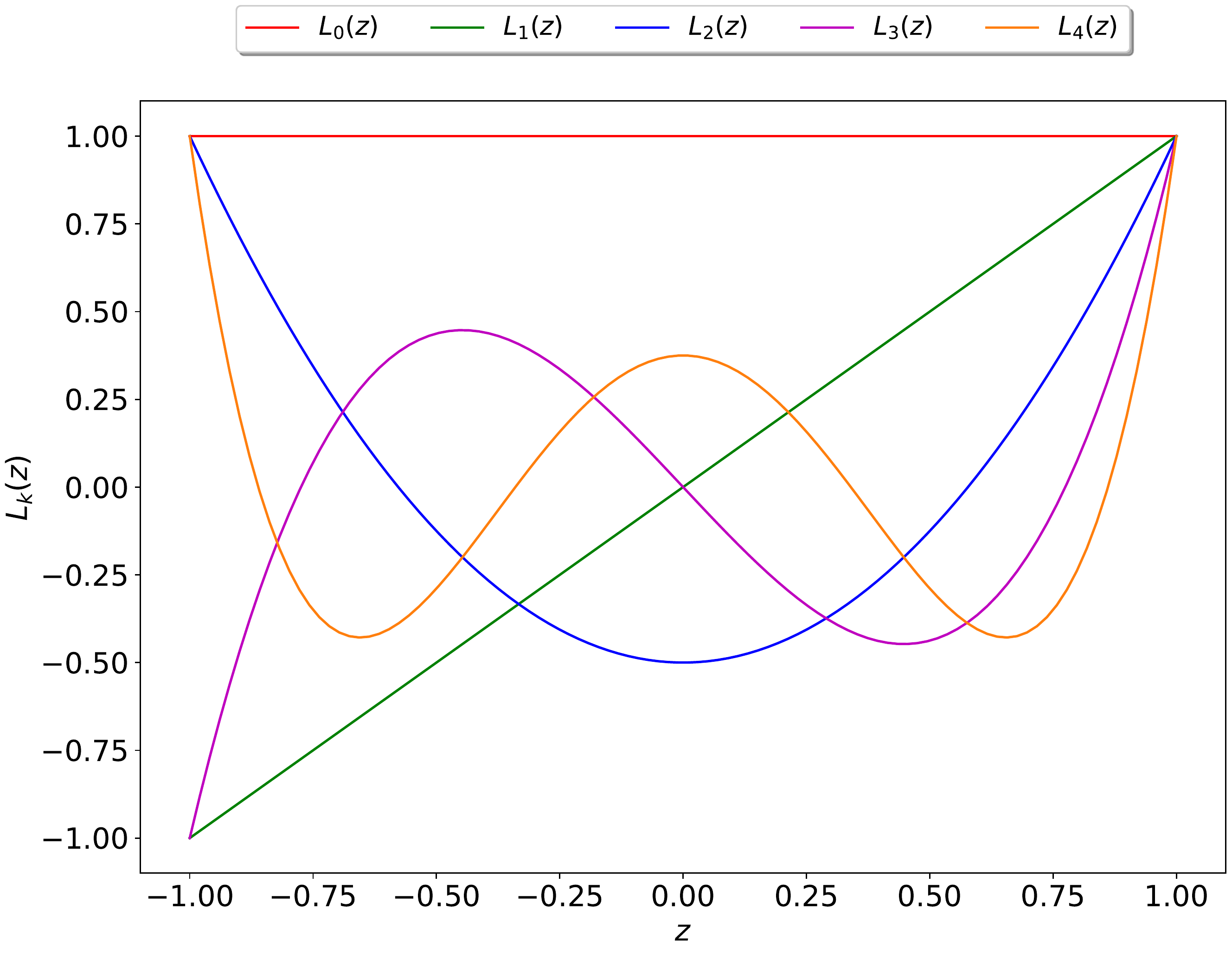}
    \caption{First five Legendre orthogonal polynomials.}
    \label{fig:I_LeP}
\end{figure}

The Legendre orthogonal polynomials, $L_k (z)$, are defined on the domain $z\in[-1, +1]$ with measure $\dd \mu (z) = \dd z$. These polynomials can also be generated recursively by,
\begin{equation*}
    L_{k+1} = \dfrac{2k+1}{k+1} \, z \, L_k - \dfrac{k}{k+1} \, L_{k-1} \qquad \text{starting with:} \; \begin{cases} L_0 = 1 \\ L_1 = z.\end{cases}
\end{equation*}
All derivatives of Legendre orthogonal polynomials can be computed in a recursive way, starting from,
\begin{equation*}
    \dfrac{\dd L_0}{\dd z} = 0, \quad \dfrac{\dd L_1}{\dd z} = 1 \qquad \text{or} \qquad \dfrac{\dd^d L_0}{\dd z^d} = \dfrac{\dd^d L_1}{\dd z^d} = 0 \quad \forall \; d > 1,
\end{equation*}
and continuing with,
\begin{align*}
        \dfrac{\dd L_{k+1}}{\dd z} &= \dfrac{2k+1}{k+1} \left(L_k + z \dfrac{\dd L_k}{\dd z}\right) - \dfrac{k}{k+1} \dfrac{\dd L_{k-1}}{\dd z} \\
        \dfrac{\dd^2 L_{k+1}}{\dd z^2} &= \dfrac{2k+1}{k+1} \left(2\dfrac{\dd L_k}{\dd z} + z \dfrac{\dd^2 L_k}{\dd z^2}\right) - \dfrac{k}{k+1} \dfrac{\dd^2 L_{k-1}}{\dd z^2} \\
        & \ \ \vdots \\
        \dfrac{\dd^d L_{k+1}}{\dd z^d} &= \dfrac{2k+1}{k+1} \left(d\dfrac{\dd^{d-1} L_k}{\dd z^{d-1}} + z \dfrac{\dd^d L_k}{\dd z^d}\right) - \dfrac{k}{k+1} \dfrac{\dd^d L_{k-1}}{\dd z^d} \quad \forall \; d \ge 1,
\end{align*}
for $k \ge 1$. In addition, the inner products of the Legendre polynomials highlight their orthogonality,
\begin{equation*}
    \langle L_i (z), L_j (z)\rangle = \int_{-1}^{+1} L_i (z) \, L_j (z) \, \dd z = \dfrac{2}{2 i + 1} \, \delta_{ij}.
\end{equation*}
Figure \ref{fig:I_LeP} shows the first five Legendre orthogonal Polynomials.

\section{Laguerre Orthogonal Polynomials}

Laguerre orthogonal polynomials, $L_k (z)$, are defined on the domain $[0,\infty)$ and by the measure $\dd\mu (z) = e^{-z} \dd z$. They are generated using the recursive function,
\begin{equation*}
    L_{k+1} = \dfrac{2k + 1 - z}{k + 1} \, L_k - \dfrac{k}{k + 1} \, L_{k-1} \qquad \text{starting with:} \; \begin{cases} L_0 =& 1, \\ L_1 =& 1 - z.\end{cases}
\end{equation*}
All derivatives of Laguerre orthogonal polynomials can be computed recursively, starting from
\begin{equation*}
    \dfrac{\dd L_0}{\dd z} = 0, \quad \dfrac{\dd L_1}{\dd z} =-1 \qquad \text{or} \qquad \dfrac{\dd^d L_0}{\dd z^d} = \dfrac{\dd^d L_1}{\dd z^d} = 0 \quad \forall \; d > 1,
\end{equation*}
then using
\begin{align*}
        \dfrac{\dd L_{k+1}}{\dd z} &= \dfrac{2k + 1 - z}{k + 1} \dfrac{\dd L_k}{\dd z} - \dfrac{1}{k + 1} L_k - \dfrac{k}{k + 1} \dfrac{\dd L_{k-1}}{\dd z} \\
        \dfrac{\dd^2 L_{k+1}}{\dd z^2} &= \dfrac{2k + 1 - z}{k + 1} \dfrac{\dd^2 L_k}{\dd z^2} - \dfrac{2}{k + 1} \dfrac{\dd L_k}{\dd z} - \dfrac{k}{k + 1} \dfrac{\dd^2 L_{k-1}}{\dd z^2} \\
        & \ \ \vdots \\
        \dfrac{\dd^d L_{k+1}}{\dd z^d} &= \dfrac{2k + 1 - z}{k + 1} \dfrac{\dd^d L_k}{\dd z^d} - \dfrac{d}{k + 1} \dfrac{\dd^{d-1} L_k}{\dd z^{d-1}} - \dfrac{k}{k + 1} \dfrac{\dd^d L_{k-1}}{\dd z^d} \quad \forall \; d \ge 1,
\end{align*}
for $k \geq 1$.

\begin{figure}[ht]
    \centering\includegraphics[width=0.9\linewidth]{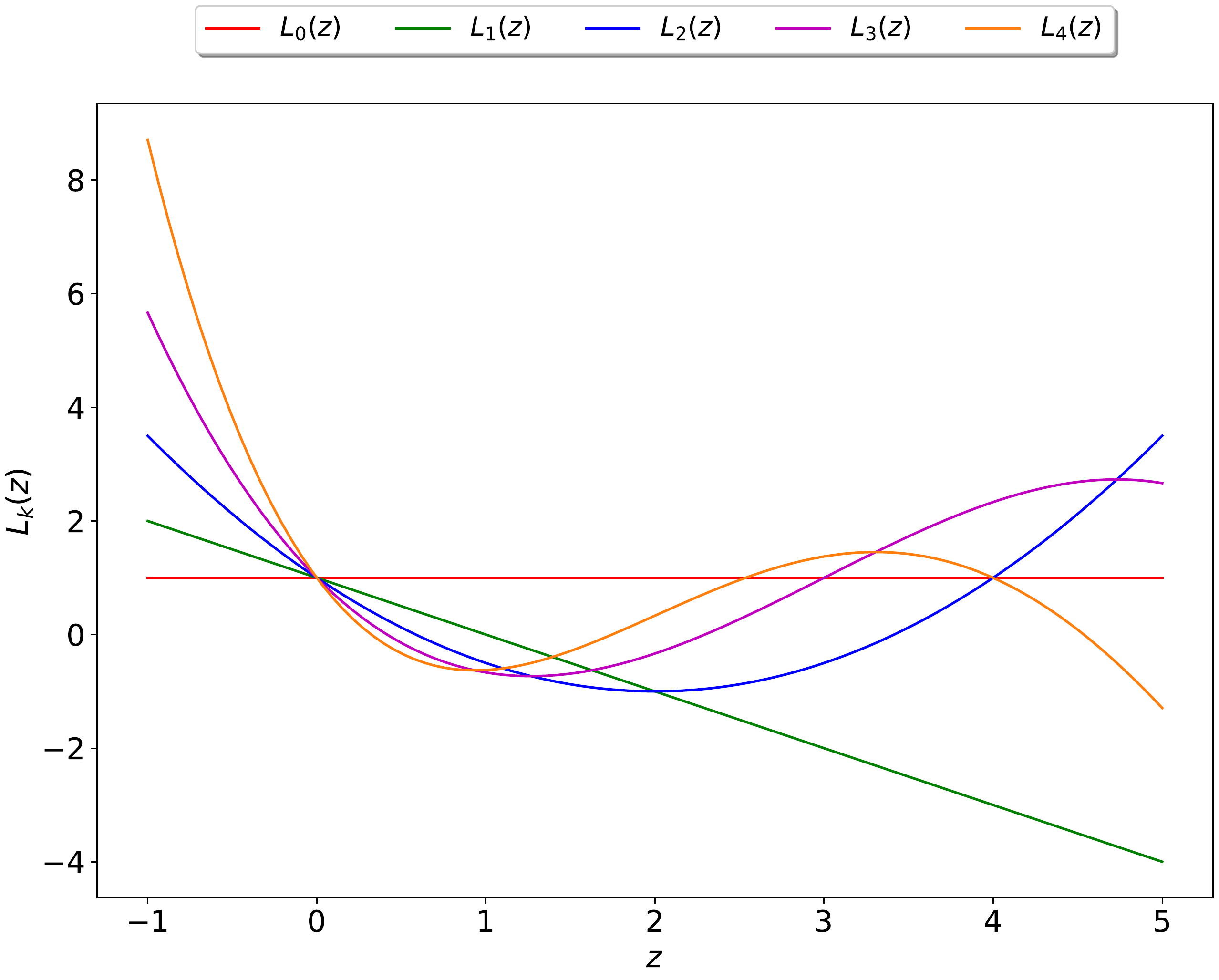}
    \caption{First five Laguerre orthogonal polynomials.}
    \label{fig:I_LaP}
\end{figure}

Figure \ref{fig:I_LaP} shows the first five Laguerre orthogonal Polynomials.

\section{Hermite Orthogonal Polynomials}

There are two Hermite orthogonal polynomials, the probabilists, indicated by $E_k (z)$, and the physicists, indicated by $H_k (z)$. The probabilists are defined on the domain $z\in(-\infty,\infty)$ and with the measure $\dd \mu(z) = e^{-(z^2/2)} \dd z$, and the physicists are defined on the domain $z\in(-\infty,\infty)$ and with the measure $\dd\mu(z) = e^{-z^2}\dd z$. They are both generated using recursive functions.

The probabilists' polynomials can be defined recursively by,
\begin{equation*}
    E_{k + 1} = z \, E_k - k E_{k - 1} \qquad \text{starting with:} \; \begin{cases} E_0 =& 1 \\ E_1 =& z.\end{cases}
\end{equation*}
All derivatives can be computed recursively, starting from
\begin{equation*}
    \dfrac{\dd E_0}{\dd z} = 0, \quad \dfrac{\dd E_1}{\dd z} = 1 \qquad \text{or} \qquad \dfrac{\dd^d E_0}{\dd z^d} = \dfrac{\dd^d E_1}{\dd z^d} = 0 \quad \forall \; d > 1,
\end{equation*}
then using,
\begin{align*}
        \dfrac{\dd E_{k+1}}{\dd z} &= E_k + z \dfrac{\dd E_k}{\dd z} - k \dfrac{\dd E_{k-1}}{\dd z} \\ 
        \dfrac{\dd^2 E_{k+1}}{\dd z^2} &= 2 \dfrac{\dd E_k}{\dd z} + z \dfrac{\dd^2 E_k}{\dd z^2} - k \dfrac{\dd^2 E_{k-1}}{\dd z^2} \\
        & \ \ \vdots \\
        \dfrac{\dd^d E_{k+1}}{\dd z^d} &= d \dfrac{\dd^{d-1} E_k}{\dd z^{d-1}} + z \dfrac{\dd^d E_k}{\dd z^d} - k \dfrac{\dd^d E_{k-1}}{\dd z^d} \quad \forall \; d \ge 1,
\end{align*}
for $k \geq 1$.

The physicists' polynomials can be defined by the recursive relationship,
\begin{equation*}
    H_{k + 1} = 2 z \, H_k - 2 k \, H_{k - 1} \qquad \text{starting with:} \; \begin{cases} H_0 =& 1 \\ H_1 =& 2z.\end{cases}
\end{equation*}
All derivatives can be computed recursively, starting from
\begin{equation*}
    \dfrac{\dd H_0}{\dd z} = 0, \quad \dfrac{\dd H_1}{\dd z} = 2 \qquad \text{or} \qquad \dfrac{\dd^d H_0}{\dd z^d} = \dfrac{\dd^d H_1}{\dd z^d} = 0 \quad \forall \; d > 1,
\end{equation*}
then using,
\begin{align*}
        \dfrac{\dd H_{k+1}}{\dd z} &= 2 H_k + 2 z \dfrac{\dd H_k}{\dd z} - 2k \dfrac{\dd H_{k-1}}{\dd z} \\
        \dfrac{\dd^2 H_{k+1}}{\dd z^2} &= 4 \dfrac{\dd H_k}{\dd z} + 2 z \dfrac{\dd^2 H_k}{\dd z^2} - 2k \dfrac{\dd^2 H_{k-1}}{\dd z^2} \\
        & \ \ \vdots \\
        \dfrac{\dd^d H_{k+1}}{\dd z^d} &= 2 d \dfrac{\dd^{d-1} H_k}{\dd z^{d-1}} + 2 z \dfrac{\dd^d H_k}{\dd z^d} - 2 k \dfrac{\dd^d H_{k-1}}{\dd z^d} \quad \forall \; d \ge 1,
\end{align*}
for $k \geq 1$.

\begin{figure}[ht]
    \centering\includegraphics[width=\linewidth]{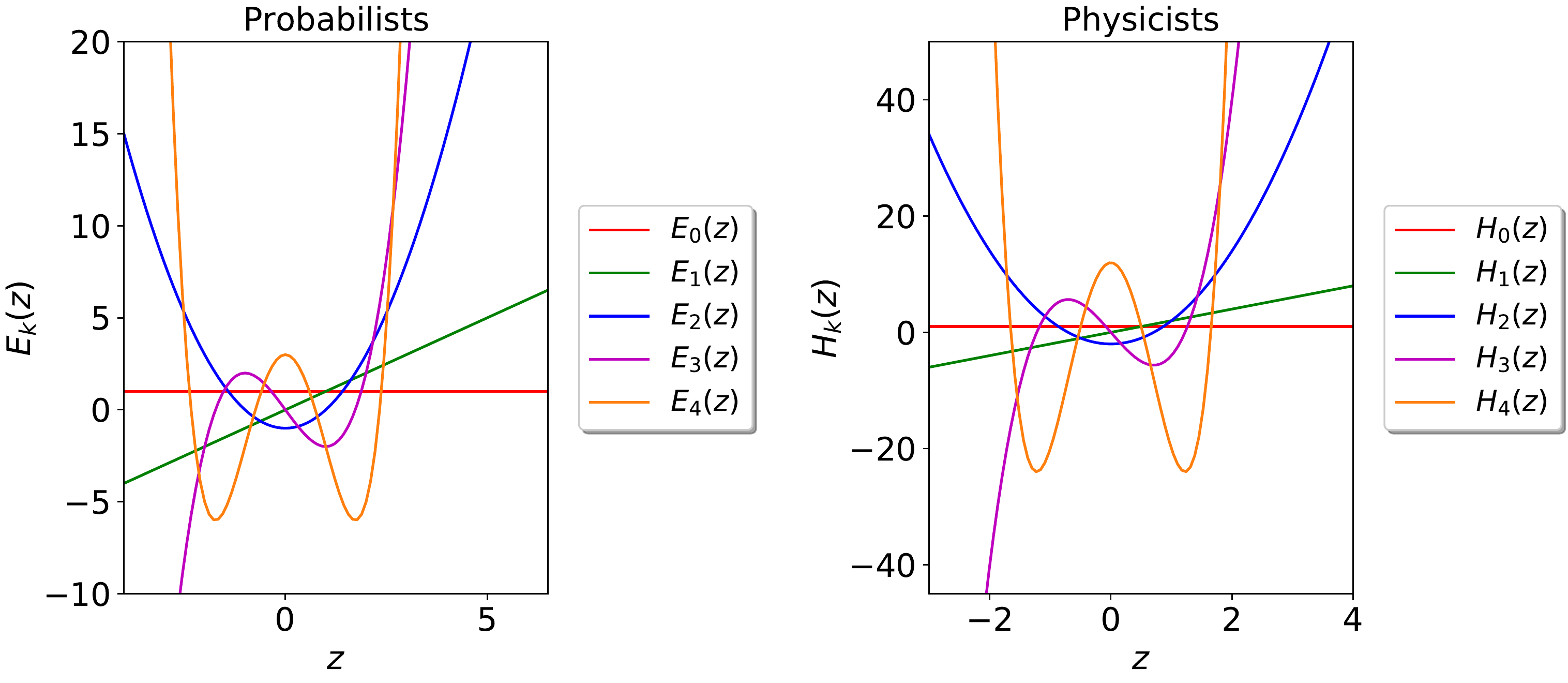}
    \caption{First five Hermite orthogonal polynomials.}
    \label{fig:I_HoP}
\end{figure}

Figure \ref{fig:I_HoP} shows the first five probabilists' and physicists' Hermite orthogonal polynomials.

\section{Fourier Basis}

\begin{figure}[ht]
    \centering\includegraphics[width=0.9\linewidth]{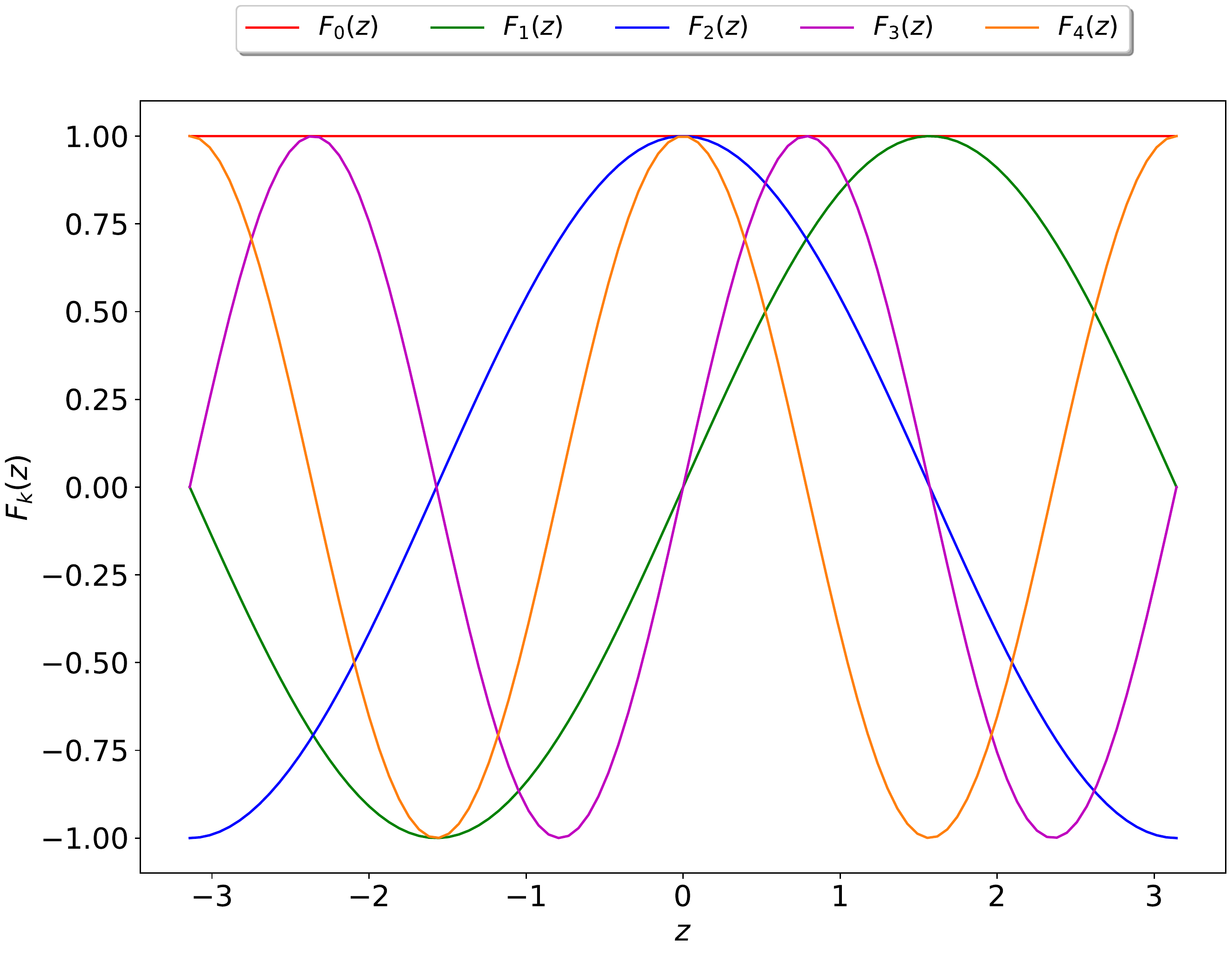}
    \caption{First five Fourier basis functions.}
    \label{fig:I_FS}
\end{figure}

The Fourier basis is defined on the domain $z\in[-\pi,\pi]$ and with the measure $\dd \mu (z) = \dd z$. The basis does not have a recursive generating function. Rather, the basis can be mathematically written as,
\begin{equation*}
    g_k(z) = \begin{cases}1, &k=0\\\cos(\ceil{k/2}z), &\text{$k$ is even}\\\sin(\ceil{k/2}z), &\text{$k$ is odd}\end{cases} 
\end{equation*}
where $\ceil{x}$ rounds $x$ to the next largest integer and $k=0,\dots,m$. There is no recursive relationship to compute the subsequent derivatives of Fourier bases. However, the $n$-th derivative can be computed using,
\begin{equation*}
    \dfrac{\dd^d g (z)}{\dd z^d}  = \begin{cases} \begin{cases}0, &k=0\\\ceil{k/2}^d\cos(\ceil{k/2}z), &\text{$k$ is even}\\\ceil{k/2}^d\sin(\ceil{k/2}z), &\text{$k$ is odd}\end{cases} &\mod(d,4) = 0 \\ \begin{cases}0, &k=0\\-\ceil{k/2}^d\sin(\ceil{k/2}z), &\text{$k$ is even}\\\ceil{k/2}^d\cos(\ceil{k/2}z), &\text{$k$ is odd}\end{cases} &\mod(d,4) = 1 \\ \begin{cases}0, &k=0\\-\ceil{k/2}^d\cos(\ceil{k/2}z), &\text{$k$ is even}\\-\ceil{k/2}^d\sin(\ceil{k/2}z), &\text{$k$ is odd}\end{cases} &\mod(d,4) = 2 \\ \begin{cases}0, &k=0\\\ceil{k/2}^d\sin(\ceil{k/2}z), &\text{$k$ is even}\\-\ceil{k/2}^d\cos(\ceil{k/2}z), &\text{$k$ is odd}\end{cases} &\mod(d,4) = 3\end{cases}
\end{equation*}
whenever $d > 0$. Figure \ref{fig:I_FS} shows the first five Fourier basis functions.

\section{Extension to Multivariate Domains}

In general, multivariate orthogonal basis sets can be created by taking all possible products of functions in the basis sets that make up the individual variables. The measure that makes up this new basis set will be the product of measures of the individual basis sets, and the domain of the multivariate basis set will be the union of the domains that make up the individual basis sets. More details and insights on the 2-dimensional and $n$-dimensional orthogonal basis functions are contained in Reference \cite{Ye} and References \cite{MultiVarOrthPolyBook, Xu}, respectively.

Consider $n$ independent variables in the vector $\B{x} = \{x_1, x_2, \cdots, x_n\}\T$. Moreover, let the orthogonal basis set for each of these independent variables be denoted by $\p{k}{B}_j$, where the subscript $j$ denotes the $j$-th basis function and the pre-superscript $k$ denotes the $k$-th independent variable. For example, the third basis function for $x_2$ would be  $\p{2}{B}_3$. The measure of the multivariate basis set will be denoted by $\mu (\B{x}) = \ds\prod_{k=1}^n \p{k}{\mu}(x_k)$ where $\p{k}{\mu}(x_k)$ is the measure for the $k$-th independent variable. The domain of the multivariate basis will be denoted by $\Omega = \p{1}{\Omega} \times \p{2}{\Omega} \times \cdots \times \p{n}{\Omega}$, where the generic $\p{k}{\Omega}$ denotes the domain of the $k$-th basis set. Then, an arbitrary basis function for the multivariate domain can be written as,
\begin{equation}\label{eq:nDbasisAsTensorProduct}
    \mathcal{B}_{i_1i_2\dots i_n} = \p{1}{B}_{i_1} \p{2}{B}_{i_2} \cdots \p{n}{B}_{i_n},
\end{equation}
where $i_1,\cdots ,i_n \in \mathbb{Z^+}$. In other words, Equation \eqref{eq:nDbasisAsTensorProduct} generates a multivariate basis via a tensor product of univariate basis functions \cite{nDbasisFunctions}. If one were to use all possible products of the functions in the individual basis sets which span $L^{2e}(\p{k}{\Omega},\p{k}{\mu})$, i.e., use all possible combinations of $i_1,\cdots ,i_n \in \mathbb{Z^+}$, an infinite set, then the resulting multivariate basis would span the multivariate function space $L^{2e}(\Omega,\mu)$. Of course, in practice this is not possible, so a finite number of basis functions from the set is used.

Consider the inner product of two different basis functions $\mathcal{B}_{i_1\dots i_n}$ and $\mathcal{B}_{j_1\dots j_n}$ where at least one $i_k\neq j_k$, 
\begin{equation}\label{eq:MultiVarInnerProd}
    \langle \mathcal{B}_{i_1\dots i_n}, \mathcal{B}_{j_1\dots j_n} \rangle = \int_\Omega \mathcal{B}_{i_1\dots i_n} \, \mathcal{B}_{j_1\dots j_n} \dd\mu = \ds\prod_{k=1}^n \int_{\Omega_k} \p{k}{B}_{i_k}  \, \p{k}{B}_{j_k} \dd\mu_k.
\end{equation}
Since these are different basis functions, there must be some $k=\kappa$ such that $i_{\kappa} \neq j_{\kappa}$. For $k=\kappa$, the integral
\begin{equation*}
    \int_{\Omega_\kappa} \p{\kappa}{B}_{i_{\kappa}} \p{k}{B}_{j_{\kappa}} \ \dd \mu_\kappa = 0,
\end{equation*}
and thus, the product of integrals in Equation \eqref{eq:MultiVarInnerProd} is equal to zero. It follows that,
\begin{equation*}
    \langle \mathcal{B}_{i_1\dots i_n}, \mathcal{B}_{j_1\dots j_n} \rangle = 0 \quad \text{if} \quad \exists\ \kappa \mid i_\kappa \neq j_\kappa.
\end{equation*} 
Hence, the resulting multivariate basis set is orthogonal.

Just as in the univariate case, the problem being solved must be made tractable by choosing basis functions up to some finite degree $m$. All the multivariate basis functions of order $m$ are defined by choosing $i_1, \cdots, i_n$ to be on the set,
\begin{equation*}
    \{\B{i} \mid i_k\in\mathbb{Z}^{+}, \sum_{k=1}^n (i_k-1) = m\},
\end{equation*}
where $i_k$ denotes the elements of $\B{i}$. 

Table \ref{tab:BasisSummary} summarizes the orthogonal basis sets described in this section.
\begin{table}[!ht]
    \centering
    \caption{Univariate orthogonal basis functions summary.}
    \begin{tabular}{|l|c|c|}
    \hline
    Basis function name & Domain, $\Omega$ & Measure, $\dd \mu (z)$ \\
    \hline\hline
    Chebyshev polynomials & $[-1,1]$ & $\dfrac{1}{1 - z^2} \dd z$ \\
    Legendre polynomials & $[-1,1]$ & $\dd z$ \\
    Laguerre polynomials & $[0,\infty)$ & $e^{-z} \dd z$ \\
    Hermite probabilists polynomials & $(-\infty,\infty)$ & $e^{-(z^2/2)} \dd z$ \\
    Hermite physicists polynomials & $(-\infty,\infty)$ & $e^{-z^2} \dd z$ \\
    Fourier series & $[-\pi,\pi]$ & $\dd z$ \\
    \hline
    \end{tabular}
    \label{tab:BasisSummary}
\end{table}


%% file: Data/appendixLS.tex

\chapter{LINEAR LEAST-SQUARES METHODS}\label{app:LinearLeastSquares}

There are different numerical techniques to compute the linear least-squares (LS) solution of $\mathbb{A} \, \B{\xi} = \B{b}$. These are:
\begin{itemize}
\item The Moore–Penrose inverse,
    \begin{equation*}
        \B{\xi} = (\mathbb{A}\T \, \mathbb{A})^{-1} \, \mathbb{A}\T \, \B{b}.
    \end{equation*}
\item QR decomposition,
    \begin{equation*}
         \mathbb{A} = Q \, R \qquad \to \qquad \B{\xi} = R^{-1} \, Q\T \, \B{b},
    \end{equation*}
    where $Q$ is an orthogonal matrix and $R$ an upper triangular matrix.
\item SVD decomposition,
    \begin{equation*}
         \mathbb{A} = U \,  \Sigma \, V\T \qquad \to \qquad \B{\xi} = \mathbb{A}^+ \, \B{b} = V \,  \Sigma^+ \, U\T \, \B{b}
    \end{equation*}
     where $U$ and $V$ are two orthogonal matrices and $\Sigma^+$ is the pseudo-inverse of $\Sigma$, which is formed by replacing every non-zero diagonal entry by its reciprocal and transposing the resulting matrix.
\item Cholesky decomposition,
    \begin{equation*}
        \mathbb{A}\T \mathbb{A} \, \B{\xi} = U\T U \B{\xi} = \mathbb{A}\T \, \B{b} \qquad \to \qquad \B{\xi} = U^{-1} \left(U^{-\mbox{\tiny T}} \mathbb{A}\T \, \B{b}\right),
    \end{equation*}
    where $U$ is an upper triangular matrix, and consequently, $U^{-1}$ and $U^{-\mbox{\tiny T}}$ are easy to compute.
\end{itemize}

One can reduce the condition number of the matrix to be inverted by scaling the columns of $\mathbb{A}$,
\begin{equation*}
    \mathbb{A} \left(S S^{-1}\right) \B{\xi} = \left(\mathbb{A} S\right) \left(S^{-1} \B{\xi}\right) = \mathbb{B} \, \B{\eta} = \B{b} \; \to \; \B{\xi} = S \, \B{\eta} = S \, (\mathbb{B}\T \mathbb{B})^{-1} \mathbb{B}\T \B{b},
\end{equation*}
where $S$ is the $m\times m$ scaling diagonal matrix whose diagonal elements are the inverse of the norms of the corresponding columns of $\mathbb{A}$: $s_{kk} = |\B{a}_k|^{-1}$ or the maximum absolute value, $s_{kk} = \max\limits_{i} |a_{ki}|$. 

In this dissertation, the least-squares problem is solved using two methods: (1) the SVD decomposition introduced above (2) a combination of QR decomposition and the previously mentioned scaling, called the scaled QR approach. This approach performs the QR decomposition of the scaled matrix,
\begin{equation*}
    \mathbb{B} = \mathbb{A} \, S = Q \, R \qquad \to \qquad \B{\xi} = S \, R^{-1} \, Q\T \, \B{b}.
\end{equation*}

%% file: Data/appendixJaxCode.tex

\chapter{TFC NUMERICAL IMPLEMENTATION IN JAX}\label{app:JaxCode}

This appendix provides a more detailed description of the code package used to implement TFC in JAX than was given in the main body of the text and some of the major challenges in doing so. In addition, a summary of the main classes available in the package is provided. For a more detailed explanation and tutorials on how to use them see either the \href{https://tfc-documentation.readthedocs.io/en/latest/}{\textcolor{blue}{\underline{code documentation}}} or the \href{https://github.com/leakec/tfc}{\textcolor{blue}{\underline{TFC GitHub}}} \cite{TfcGithub}. Note that while not discussed in detail here, this package also contains some convenience classes and functions that the reader may find useful; for example, the \verb"MakePlot" class assists the user in creating journal-ready plots, and the \verb"Latex" class can be used to export \verb"NumPy" arrays to LaTeX tables.

\section{Basis Function Classes}\label{sec:BasisFuncClasses}
\begin{figure}[ht]
    \centering
    \includegraphics[width=\linewidth]{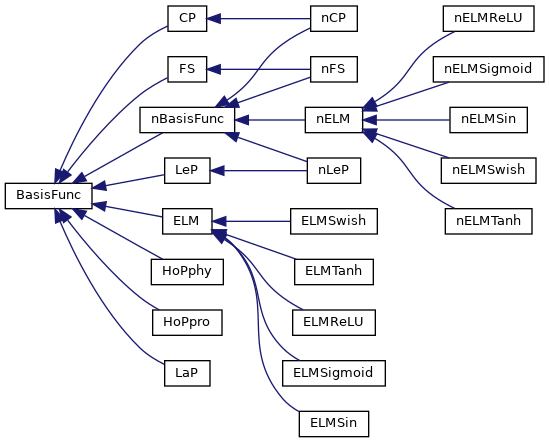}
    \caption{Basis function class hierarchy.}
    \label{fig:TfcBasisClassHierarchy}
\end{figure}

When applying TFC to a differential equation, the free function is used to minimize the differential equation's residual. As described in Section \ref{sec:freeFunctions}, two popular choices for the free function are a linear combination of $n$-dimensional basis functions and ELMs. Notice that these two free function choices can both be written as,
\begin{equation*}
    g(\B{x}) = \B{h}\T \B{\xi},
\end{equation*}
where the calculation for $\B{h}$ evaluated at $\B{x}$ is dictated by the basis function set chosen, if using basis functions, or by the activation function, if using an ELM. Hence, these two free function choices have a lot in common, and in the code are derived from the same abstract, parent class called \verb"BasisFunc". Figure \ref{fig:TfcBasisClassHierarchy} shows the inheritance diagram for the basis function classes.
\noindent From this base class, the univariate basis function classes are derived---\verb"CP" for Chebyshev polynomials, \verb"LeP" for Legendre polynomials, \verb"FS" for Fourier series, \verb"LaP" for Laguerre polynomials, \verb"HoPphy" for the physicists' Hermite polynomials, and \verb"HoPpro" for the probabilists' Hermite polynomials. In addition, an $n$-dimensional abstract class is derived from the base class, from which $n$-dimensional versions of Chebyshev polynomial, Legendre polynomial, and Fourier series classes are derived, called \verb"nCP", \verb"nLeP", and \verb"nFS" respectively. Furthermore, an ELM abstract base class is derived from \verb"BasisFunc", from which five concrete univariate ELM classes are derived; each of these five classes implements a different activation function: \verb"ELMSin" implements the sine activation function, \verb"ELMSwish" implements the swish activation function, \verb"ELMTanh" implements the hyperbolic tangent activation function, \verb"ELMSigmoid" implements the sigmoid activation function, and \verb"ELMReLU" implements the rectified linear activation function. The $n$-dimensional versions of the five univariate ELM classes are derived from an abstract parent class called \verb"nELM", which is derived from the abstract $n$-dimensional basis function class mentioned earlier. The only real difference between the ELM classes and the basis function classes, besides the functions they are implementing, is the addition of randomly generated weights and biases for the hidden layer, i.e., $W_1$ and $b_1$ in Equation \eqref{eq:ElmForm}. 

All of the classes shown in Figure \ref{fig:TfcBasisClassHierarchy} are programmed in C++. The attentive reader may wonder why not just write them in Python, as the scripts that use this core code are ultimately compiled via a JIT anyway? The underlying reasons lie within the restrictions JAX has on JIT. To best understand, the reader must first recall what these basis function classes need to do:
\begin{itemize}
    \item Calculate the basis functions and their derivatives. Recall that this requires using the linear map from $x\to z$, as the problem domain, $x$, and the basis function domain, $z$, do not generally coincide. Furthermore, recall that many of the basis functions are created via a recursion, see Appendix \ref{app:BasisFunctions}. 
    \item While the default behavior is to neglect the terms linearly independent to the support functions when calculating the basis function matrix and its derivatives, there may be times when the user wants to ignore this default behavior. This functionality is accomplished through the argument \verb"full" that shows up as a required input argument in the C++ API and an optional keyword argument in the Python API. 
\end{itemize}

Naturally, the argument \verb"full" requires an if statement to be implemented. Unfortunately, the JIT does not allow tracing through if statements; hence, the basis function classes must be implemented as primitives. Furthermore, if built using JAX directly, the primitive recursions used by many of the basis functions would require using lax operations, as one cannot modify JAX arrays in-place. These lax operations are pure functions; rather than modifying arrays in-place, they return a new array with the updated indices. Hence, the basis function array would be copied to a new array on each iteration of the recursion: not an efficient solution. One may be tempted to implement these using the original NumPy library, but doing so would not provide a function to XLA, which is needed for the JIT. Therefore, to make the final result JIT-able, one must write the basis functions in a lower-level language: the author chose C++.

Writing the basis functions in C++ is a fairly straightforward task, but integrating the result with the JAX JIT compiler, which is XLA, was non-trivial. Since JAX uses XLA, there is a convenient \verb"register_custom_call_target" function that can be used to integrate C++ functions with XLA and ultimately allow them to be JIT-able. However, using this custom call function requires wrapping the C++ function into a PyCapsule object, which means the object being placed in the PyCapsule must be of type \verb"void*"; hence, the C++ function must be cast as a \verb"void*". However, the C++ functions to be wrapped are methods of classes, which means they are not regular functions. Therefore, they do not have standard function pointers that can be easily cast to a \verb"void*". 

Fortunately, this predicament can be overcome through the use of polymorphism. Notice that all the TFC free function classes are derived from \verb"BasisFunc", see Figure \ref{fig:TfcBasisClassHierarchy}. Hence, all basis function classes can be referenced using a \verb"BasisFunc" pointer. Therefore, a static \verb"std::vector" was added to the \verb"BasisFunc" class that contains a pointer to each \verb"BasisFunc" class as it is created. Furthermore, each class contains a unique integer identifier that corresponds to its pointer's position in said \verb"std::vector". That means that given this integer, one could access the correct element of the aforementioned \verb"std::vector" and call the correct member function. Thus, a simple C++ wrapper function was created with the correct function signature that takes in said integer and the required inputs for the member function, calls the member function, and returns the result. Moreover, this simple wrapper function has a standard function pointer that can easily be cast to a \verb"void*", which subsequently can be added to a PyCapsule, and finally added to the set of JIT-able functions via \verb"register_custom_call_target".

\section{Univariate TFC Class}
The univariate TFC class, called \verb"utfc", is used to create univariate TFC expressions and solve ODEs. The class's inputs are:
\begin{itemize}
    \item $N$ - Number of points to use when discretizing the domain.
    \item $nC$ - Number of basis functions to remove from the linear expansion. This variable is used to account for basis functions that are linearly dependent to the support functions used in the construction of the constrained expression. The constraints for each dimension can be expressed in one of two ways. Note that a value of -1 is used to indicate that no constraints exist for a particular dimension.
    \begin{enumerate}
        \item As an integer. When expressed as an integer, the first $nC$ basis functions are removed from the free function.
        \item As a set of integers. The basis functions corresponding to the numbers given in the set are removed from the free function.
    \end{enumerate}
    \item \verb"m" - Degree of the basis function expansion. This number is one less than the number of basis functions used.
    \item \verb"basis" - This optional string argument specifies the basis functions that will be used as the free function. The default is Chebyshev orthogonal polynomials.
    \item $x_0$ - This optional argument specifies the beginning of the DE domain. The default value of $0$ will result in a DE domain that begins at $0$.
    \item $x_f$ - This required keyword argument specifies the end of the DE domain.
\end{itemize}

The class creates a set of $N$ discretized points in the domains $x\in[x_0,x_f]$ and $z\in[z_0,z_f]$, where the values of $z_0$ and $z_f$ are automatically chosen based on the user-specified value of the optional keyword argument \verb"basis". The problem domain values, $x$, and basis domain values, $z$, are both made available to the user as public variables associated with the class. In addition, the TFC class creates an instance of the basis function class specified by \verb"basis" and creates the required JAX primitives to take gradients, Jacobians, and JIT the basis functions. The basis function values and their derivatives are made available to the user via methods associated with the TFC class. These methods contain the optional keyword \verb"full" discussed earlier, see section \ref{sec:BasisFuncClasses}.

\section{Multivariate TFC Class}
The multivariate TFC class is the multidimensional extension of the univariate TFC class. The input arguments to the class are the same, but they must be specified for each dimension, e.g., $N$ specifies the number of discretization points per dimension via a list or array of the proper size. Furthermore, one additional keyword argument, \verb"dim", is used to specify the number of dimensions: the default is two. 

\section{Elementwise Gradients}
The Autograd package contains a function called \verb"egrad", which stands for elementwise gradient. This function does not exist explicitly in JAX, but can easily be recreated using the same methodology as the original \verb"egrad" function available in Autograd \cite{autograd}: extract the diagonal elements of the Jacobian via a Jacobian-vector product or vector-Jacobian product. Using the tree utilities available in JAX, this function can easily be extended to pytrees, i.e., one can take elementwise gradients with respect to nested sets of Python containers. 

Typically, one can use \verb"vmap" to transform calls to the gradient function, \verb"grad", in JAX into elementwise gradients. However, since the basis function class must output a vector for each input, the \verb"grad" JAX transform cannot be used. Hence, the necessity for a separate elementwise gradient function. 

\section{Extending Ordered Dictionaries}
As mentioned earlier, JAX allows one to take gradients and Jacobians with respect to pytrees. This is particularly useful for coupled systems of differential equations or when splitting the domain into sub-domains, as all the unknowns can be combined into one pytree, and the gradient and/or Jacobian of the loss function with respect to all unknowns can be written in one line. 

However, if one is using an iterative least-squares optimization technique, a complication arises. Ideally, one wants to perform the following iteration,
\begin{equation*}
    \B{\xi}_{j+1} = \B{\xi}_j + \Delta\B{\xi},
\end{equation*}
where 
\begin{equation*}
    \mathbb{\B{L}}(\B{\xi}_j) + \mathcal{J}(\B{\xi}_j)\Delta\B{\xi} = 0,
\end{equation*}
but the Jacobian in JAX computed from a dictionary actually returns a dictionary of Jacobians. This can be easily overcome using list comprehension, such as,
\lstset{language=Python}
\begin{lstlisting}
j = jacfwd(L,0)
J = np.hstack([k for k in j.values()]).
\end{lstlisting}
This allows one to calculate $\Delta\B{\xi}$; however, this action will not be repeatable unless an ordered dictionary is used: without an ordered dictionary, the concatenation of the Jacobians might happen in a different order each time. While an ordered dictionary solves the Jacobian creation problem, it still does not allow for $\B{\xi}_j + \Delta\B{\xi}$ to be performed, as $\Delta\B{\xi}$ is a \verb"NumPy" array and $\B{\xi}_j$ is an ordered dictionary. Therefore, the ordered dictionary is extended via operator overloading to include methods that allow for this in the \verb"TFCDict" class. 

The \verb"TFCDict" class is designed for ordered dictionaries that have flat arrays as values. However, when dealing with vector differential equations, it is often convenient to express the free functions' unknown values associated with each component of the vector in one matrix \cite{JohnstonDissertation}. For these types of situations, a second class, \verb"TFCDictRobust", has been created that is similar to \verb"TFCDict", but works for both flat and non-flat arrays. 

\section{Nonlinear Least-Squares}\label{appSec:NLLS}
Nonlinear least squares is used throughout this dissertation to minimize the residuals of differential equations via the unknowns in the free functions. Since this method is used so often, a function called \verb"NLLS" is included that runs the nonlinear least squares. A class called \verb"NllsClass" exists as well for cases where the nonlinear least-squares needs to be called multiple times; the inputs to these two are similar, so only the \verb"NLLS" function will be covered here.

The inputs to the function are:
\begin{itemize}
    \item \verb"xiInit" - Initial guess for the unknown parameters, $\xi$.
    \item \verb"res" - Loss function, $\mathbb{L}$.
    \item \verb"*args" - Any additional arguments taken by $\mathbb{L}$..
\end{itemize}
In addition, the following are optional keyword arguments,
\begin{itemize}
    \item J - User-specified Jacobian. The default value is the Jacobian of $\mathbb{L}$ with respect to $\xi$.
    \item tol - Tolerance for stopping the while loop. Default is $1\times10^{-13}$.
    \item maxIter - Maximum number of nonlinear least-squares iterations. Default is 50.
    \item method - Method used to invert the matrix at each iteration. The default is \verb"pinv". The two options are:
    \begin{enumerate}
        \item \verb"pinv" - Uses \verb"np.linalg.pinv" to perform the inversion.
        \item  \verb"lstsq" - Uses \verb"np.linalg.lstsq" to perform the inversion.
    \end{enumerate} 
     \item User specified condition function. Default is None, which results in a condition that checks the three stopping conditions described below.
    \item body - User specified body function. Default is None, which results in a body function that performs least-squares using the method provided and updates $\xi$, $\Delta \xi$ and \verb"it", the current number of iterations.
    \item \verb"timer" - Setting this to True will time the non-linear least squares using Python's \verb"time.process_time" timer. Note that doing so adds a slight increase in run time, as one iteration of the non-linear least squares is run first to avoid timing the JAX trace. The default is False.
    \item \verb"printOut" - Setting this to true prints out the iteration number and value of $\max(|\mathbb{L}|_{\infty})$ at each iteration.
    \item \verb"printOutEnd" - This string argument is passed to the \verb"end" keyword argument of the print function used in \verb"printOut". The default value is ``\textbackslash n'' (newline).
\end{itemize}
The outputs of the function are
\begin{enumerate}
    \item $\xi$ - The value of $\xi$ at the end of the nonlinear least squares.
    \item \verb"it" - The number of iterations.
    \item \verb"time" - If the keyword argument \verb"timer = True", then the third output is the time required to run the nonlinear least-squares as measured by the timer; otherwise, there is no third output.
\end{enumerate}
          
When using the default condition function \verb"NLLS" checks the following conditions; if any of the conditions are true, then the nonlinear least-squares stops iterating:
\begin{enumerate}
    \item $|\mathbb{L}|_{\infty} < \verb"tol"$
    \item $|\Delta \xi|_{\infty} < \verb"tol"$
    \item Number of iterations > \verb"maxIter"
\end{enumerate}

%% file: Data/appendixNonLinearSvm.tex

\chapter{NONLINEAR SVM DERIVATION}\label{app:nonLinearSvm}
This appendix shows how the CSVM method can be used to solve a first-order, nonlinear ODE. Consider the first-order nonlinear ODE with an initial value boundary condition,
\begin{equation*}
    \dot{y}(t) = f (t, y), \quad y (t_0) = y_0, \quad t \in [t_0, t_f].
\end{equation*}
Similar to the linear case, the \ce\ is,
\begin{equation*}
    y (t) = \B{w}\T \left[\B{\varphi} (t) - \B{\varphi} (t_0)\right] + y_0,
\end{equation*}
and the domain is discretized into $N$ training points $t_0, t_1, \dots, t_N$. Again, let $e_i$ be the residual at $t_i$,
\begin{equation*}
    e_i = \dot{y} (t_i) - f (t_i, y (t_i)).
\end{equation*}

To minimize the error, the sum of the squares of the residuals is minimized. As in the linear case, the regularization term $\B{w}\T \B{w}$ is added to the expression to be minimized. Now, the problem can be formulated as an optimization problem,
\begin{figure*}[!h]
\begin{equation}\label{eq:LagNonLin}
\begin{aligned}
    {\cal L}(\B{w}, b, \B{e}, \B{y}, \B{\alpha}, \beta, \B{\eta}) =& \frac{1}{2} (\B{w}\T \B{w} + \gamma \B{e}\T \B{e}) -\sum_{i = 1}^N \alpha_i \left[\B{w}\T \B{\varphi}' (t_i) - f (t_i, y_i) - e_i\right]\\
    &\quad - \beta [\B{w}\T \B{\varphi} (t_0) + b - y_0]- \sum_{i=1}^N \eta_i \left[\B{w}\T \B{\varphi} (t_i) + b - y_i\right],
\end{aligned}
\end{equation}
\end{figure*}
where the Lagrange multipliers $\B{\alpha}$, $\beta$, and $\B{\eta}$ are used to enforce the constraints: see Equation \eqref{eq:LagNonLin}. The variables $y_i$ are introduced into the optimization problem to keep track of the nonlinear function $f$ at the values corresponding to the training points.

The values where ${\cal L}$ are zero give candidates for the minimum.
\begin{align*}
    & \frac{\partial {\cal L}}{\partial \B{w}} = \B{0} \qquad\to\qquad \B{w} = \sum_{i=1}^N \alpha_i \B{\varphi}' (t_i) + \sum_{i=1}^N \eta_i \B{\varphi} (t_i) + \beta\B{\varphi} (t_0) \\
    &\frac{\partial {\cal L}}{\partial e_i}=0\qquad\to\qquad\gamma e_i=-\alpha_i\\
    &\frac{\partial {\cal L}}{\partial \alpha_i} = 0\qquad\to\qquad\B{w}\T\B{\varphi}'(t_i)=f(t_i,y_i)+e_i\\
    &\frac{\partial {\cal L}}{\partial \eta_i} = 0 \qquad\to\qquad y_i=\B{w}\T\B{\varphi}(t_i)+b\\
    &\frac{\partial {\cal L}}{\partial\beta} = 0 \qquad\to\qquad\B{w}\T\B{\varphi}(t_0)+b=y_0\\
    &\frac{\partial {\cal L}}{\partial b} = 0 \qquad\to\qquad\beta+\sum_{i=1}^N\eta_i=0\\
    &\frac{\partial {\cal L}}{\partial y_i} = 0 \qquad\to\qquad\alpha_i f_y(t_i,y_i)+\eta_i=0
\end{align*}
A system of equations can be constructed by substituting the results found by differentiating ${\cal L}$ with respect to $\B{w}$ and $e_i$ into the remaining five equations. This leads to a set of $3N+2$ equations and $3N+2$ unknowns, which are $\alpha_i$, $\eta_i$, $y_i$, $\beta$, and $b$: this system of equations is shown in Equation (\ref{eq:sys}).
\begin{equation}\label{eq:sys}
\begin{aligned}
    &\sum_{j=1}^N\alpha_j \B{\varphi}'(t_j)\T \B{\varphi}'(t_i) + \sum_{j=1}^N \eta_j \B{\varphi}(t_j)\T \B{\varphi}'(t_i) + \beta \B{\varphi}(t_0)\T \B{\varphi}'(t_i)+\frac{\alpha_i}{\gamma} = f(t_i,y_i) \\
    &\sum_{j=1}^N\alpha_j \B{\varphi}'(t_j)\T \B{\varphi}(t_i)+\sum_{j=1}^N \eta_j \B{\varphi}(t_j)\T \B{\varphi}(t_i) + \beta\B{\varphi}(t_0)\T\B{\varphi}(t_i)+b-y_i=0\\
    &\sum_{j=1}^N\alpha_j\B{\varphi}'(t_j)\T\B{\varphi}(t_0)+\sum_{j=1}^N \eta_j \B{\varphi}(t_j)\T \B{\varphi}(t_0) + \beta \B{\varphi}(t_0)\T\B{\varphi}(t_0)+b=y_0\\
    &\beta+\sum_{i=j}^N\eta_j=0\\
    &\alpha_if_y(t_i,y_i)+\eta_i=0\\
\end{aligned}
\end{equation}
where $i = 1,...,N$. This system of equations can be written in the dual form, in terms of the Kernel matrix and is derivatives, and can be solved using least-squares. Once the set of equations has been solved, the model solution is given in the dual form~by,
\begin{equation*}
    \hat{y}(t) = \sum_{i=1}^N \alpha_i \B{\varphi}' (t_i)\T \B{\varphi} (t) + \sum_{i=1}^N \eta_i \B{\varphi} (t_i)\T \B{\varphi} (t) + \beta \B{\varphi} (t_0)\T \B{\varphi} (t) + b.
\end{equation*}

%% file: Data/appendixMiscellaneous.tex
\chapter{VARIOUS TFC EXTENSIONS}\label{app:Miscellaneous}

This appendix consists of various TFC concepts and extensions that did not fit elsewhere in the dissertation. Naturally, these concepts and extensions are encountered less frequently than those in the main body of the text.

\section{Extension to Parallelotopes via Affine Transformations}
Theorem 12.7 of Reference \cite{AffineTransformation} shows that affine transformations map parallelotopes ($n$-dimensional parallelograms) to parallelotopes. This affine transformation can serve as a bijective map between a general parallelotope and an $n$-dimensional unit hypercube: where multivariate TFC can be applied.

Consider a general $n$-dimensional parallelotope with $n$ sides emanating from each vertex. Pick a vertex arbitrarily and label it $p_0$. Then, for each of the $n$ sides emanating from $p_0$, label the vertices at the opposite ends $p_1$ through $p_n$; again, the order in which they are labeled is arbitrary. Now, an affine transformation from the unit hypercube to the parallelotope can be defined as,
\begin{equation*}
    \begin{Bmatrix} x_1 \\ x_2 \\ \vdots \\ x_n \end{Bmatrix} = \underbrace{\begin{bmatrix} p_1-p_0 & p_2-p_0 & \cdots & p_n-p_0\end{bmatrix}}_{A}\begin{Bmatrix} X_1 \\ X_2 \\ \vdots \\ X_N \end{Bmatrix} + p_0,
\end{equation*}
where $p_0, \dots, p_n$ are written as column vectors, $x_1, \dots, x_n$ are the coordinates in the parallelotope space, and $X_1, \dots, X_n$ are the familiar Cartesian coordinates of the unit hypercube. The matrix that appears in the affine transformation will appear later and has thus been given the symbol, $A$. Let the entire affine transformation be denoted by $\mathcal{A}$, i.e., $\B{x} = \mathcal{A}(\B{X})$. Additionally, let $\mathcal{A}(X_k)$ denote the $x_k$ portion of $\B{x} = \mathcal{A}(\B{X})$. 

The affine transformation can be used to map functions from the parallelotope to the unit hypercube. For example, suppose there is some function $f(\B{x})$ on the parallelotope, then, the function on the unit hypercube, $F(\B{X})$, is
\begin{equation*}
    F(\B{X}) = (f \circ \mathcal{A})(\B{X}) = f(\mathcal{A}(\B{X})).
\end{equation*}
Similarly, the inverse of the affine transformation can be used to map functions from the unit hypercube to the parallelotope,
\begin{equation*}
    f(\B{x}) = (F \circ \mathcal{A}^{-1})(\B{x}) = F(\mathcal{A}^{-1}(\B{x})).
\end{equation*}

Of significance to TFC are the directions in which constraints can be specified in the parallelotope space that correspond to TFC-embeddable constraints in the unit hypercube space. Recall that the constraint operators in multivariate TFC can only operate on one independent variable--- except integral constraints, see Section \ref{sec:RecursiveIntegralConstraints}---else, the structure of the recursive form breaks down. For value-level constraints, a constraint that operates on $X_i$ corresponds to the direction in parallelotope space parallel to the side containing $p_i$ and $p_0$, which corresponds to the coordinate $x_i$; this comes merely from observing the $i$-th column of the affine transformation. Similarly, the allowed constraint derivative directions in the parallelotope space are those parallel to one of the parallelotope's sides:
\begin{equation}\label{eq:gradientMap}
\begin{aligned}
    \frac{\partial F}{\partial X_i} &= \frac{\partial{f}}{\partial x_j}\frac{\partial x_j}{\partial X_i} \\
    &= \frac{\partial{f}}{\partial x_j} A_{ji}.
\end{aligned}
\end{equation}
Notice that the right-hand side is just the gradient of $f$ dotted with one of the columns of $A$; this is precisely a directional derivative in the parallelotope space where the direction corresponds to the side containing $p_i$ and $p_0$. Similarly, integral constraints must integrate in a direction parallel to one of the parallelotope's sides. It follows that general linear constraints must correspond with directions parallel to one of the parallelotope's sides.

In addition, Equation \ref{eq:gradientMap} has important ramifications with regards to scaling the derivatives and integrals. Let $\B{n}_i = p_i-p_0$, $n_i = ||\B{n}_i||$, and $\hat{\B{n}}_i = \frac{\B{n}_i}{n_i}$, then, Equation \ref{eq:gradientMap} can be rewritten as,
\begin{equation*}
    \frac{\partial F}{\partial X_i} = \nabla f \cdot \B{n}_i
\end{equation*}
and after some algebraic simplification,
\begin{equation*}
    \nabla f \cdot \hat{\B{n}}_i = \frac{1}{n_i} \frac{\partial F}{\partial X_i}.
\end{equation*}
Hence, a directional derivative in the $\hat{\B{n}}_i$ direction of the parallelotope space corresponds to a derivative in the unit hypercube space in the direction of $X_i$ and scaled by $\frac{1}{n_i}$. Similarly, integrating in the parallelotope space will require scaling the corresponding integral by $n_i$,
\begin{equation*}
    \int f \dd \hat{\B{n}}_i = n_i \int F \dd X_i.
\end{equation*}

This knowledge allows one to rewrite all the constraints on the dependent variables in the parallelotope space as constraints in the unit hypercube space. The \ce\ can be developed for the unit hypercube and then transformed back into the parallelotope space via the affine transform. 

\begin{example}{Parallelotope example}
Throughout this example, capital letters will be used to denote quantities in the unit hypercube space, and lowercase letters will be used to denote quantities in the parallelotope space, e.g., $U$ is the \ce\ in the unit hypercube space while $u = U\circ\mathcal{A}^{-1}$ is the \ce\ in the parallelotope space. Consider the parallelotope and constraints shown in Figure \ref{fig:parallelogramConstraints}.
\begin{figure}[H]
    \centering
    \includegraphics[width=\linewidth]{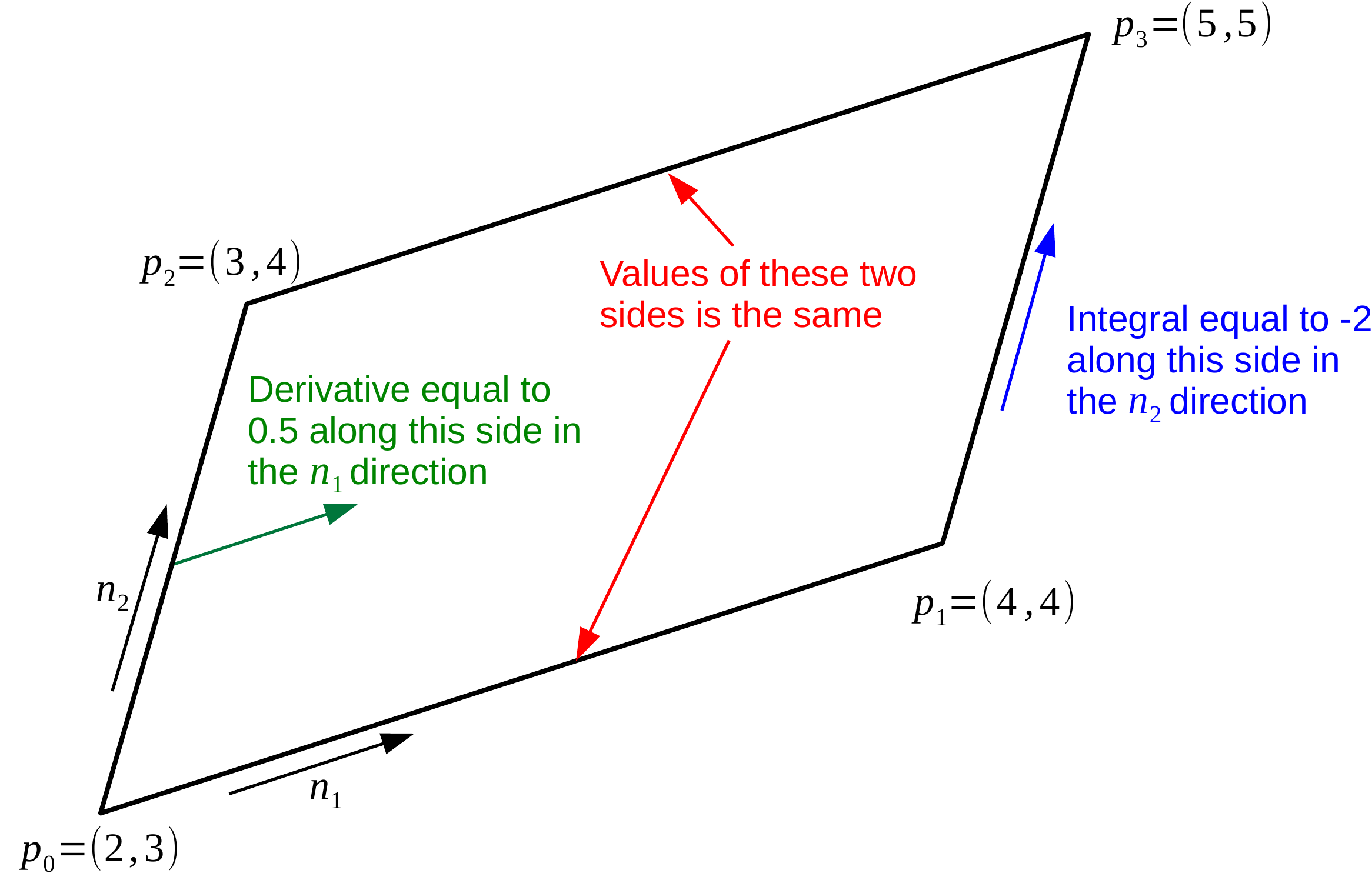}
    \caption{Parallelotope constraints.}
    \label{fig:parallelogramConstraints}
\end{figure}

Using the affine transformation previously described, the parallelotope's constraints in the unit hypercube space can be written as,
\begin{equation*}
    \frac{1}{n_1} U_X(0,Y) = 0.5, \quad n_2\int_0^1 U(1,Y) \dd Y = -2, \andd U(X,0) = U(X,1).
\end{equation*}
Using $S_1(X) = 1$ and $S_2(X)=X$ as the support functions, the univariate \ce\ for the constraints on $X$ is,
\begin{align*}
    \p{1}{U}(X,Y,G(X,Y)) &= G(X,Y) + n_1 (X-1) \Big(0.5-\frac{1}{n_1}G_X(0,Y)\Big) \\
    &\quad + \frac{1}{n_2}\Big(-2 - n_2\int_0^1 G(1,\tau) \dd \tau\Big).
\end{align*}
Utilizing the affine transformation, this \ce\ can be transformed back into the parallelotope space,
\begin{align*}
    \p{1}{u}(x,y,g(x,y)) &=(\p{1}{U}\circ \mathcal{A}^{-1})(x,y,g(x,y))  \\
    &= g(x,y) + n_1(\mathcal{A}^{-1}(x)-1)\Big(0.5-g_{\hat{\B{n}}_1}(\mathcal{A}(0,\mathcal{A}^{-1}(y)))\Big) \\
    &\quad + \frac{1}{n_2}\Big(-2 - \int_0^1 g(\mathcal{A}(1,\mathcal{A}^{-1}(\tau))) \dd \hat{\B{n}}_2 \Big)
\end{align*}
\begin{figure}[H]
    \centering
    \textattachfile{Figures/parallelotopeEx.html}{\includegraphics[width=0.75\linewidth]{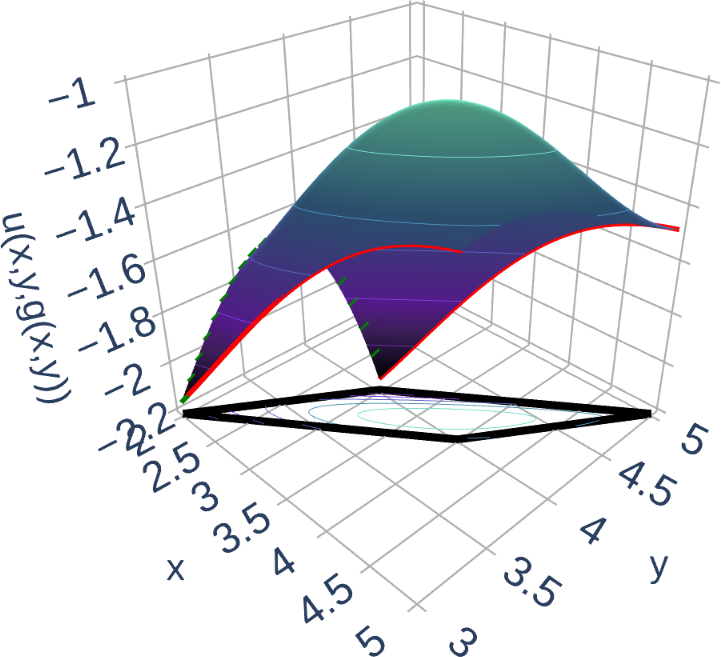}}
    \caption{Parallelotope constrained expression example. Note, this figure contains an embedded, standalone HMTL version of the plot that can be viewed/downloaded by clicking on it. Doing so may require a dedicated PDF viewer such as Adobe Acrobat or Okular.}
    \label{fig:parallelotopeEx}
\end{figure}
\noindent where $g_{\hat{\B{n}}_1}$ represents the derivative of $g$ with respect to $\hat{\B{n}}_1$, $\text{d} \hat{\B{n}}_2$ in the integral constraint is given in terms of the dummy variable $\tau$ (this dummy variable takes the place of $y$), and affine transformations have been simplified where possible, e.g., $G\circ\mathcal{A}^{-1}$ = $g\circ\mathcal{A}\circ\mathcal{A}^{-1} = g$.

Similarly, the univariate \ce\ for the constraints on $Y$ can be transformed into the parallelotope space,
\begin{align*}
    \p{2}{u}(x,y,g(x,y)) &= (\frac{1}{2}-\mathcal{A}^{-1}(y))\Big(g(\mathcal{A}(\mathcal{A}^{-1}(x),1))-g(\mathcal{A}(\mathcal{A}^{-1}(x),0))\Big),
\end{align*}
and the two can be combined using the recursive method to form a multivariate \ce\ that satisfies all of the constraints.

Figure \ref{fig:parallelotopeEx} shows the multivariate \ce\ evaluated using $g(x,y) = \sin(x)\cos(y)$. The constraint on $y$ is shown in red, and the derivative constraint on $x$ is shown via green lines; the integral constraint is not easily visualized but is satisfied nonetheless. The details of the surface can be a bit difficult to distinguish, so the surface's contours have been projected onto the $x$-$y$ plane, and the boundaries of the surface projected onto the $x$-$y$ plane are shown in black. In addition, an interactive HTML version of the plot has been embedded into the PDF and can be opened using a dedicated PDF viewer such as Adobe Acrobat or Okular.
\end{example}

\section{Lower-Dimensional Constraints in \texorpdfstring{$n$}{n}-Dimensions}
There are times in $n$-dimensions when the constraints are not written as $n$-1 dimensional manifolds. For example, consider this point constraint in two-dimensional space: $u(0,0) = 5$. Notice that these types of constraints can be written using a series of constraint operators,
\begin{equation*}
    \ppC{}{k}{i}\Big[\cdots \big[\ppC{}{j}{i}[u]\big] \cdots \Big] = \kappa_i
\end{equation*}
where the pre-subscript in front of the constraint operator represents that variable's contribution to the constraint, e.g., $\ppC{}{k}{i}$ is $x_k$'s contribution to the $i$-th constraint. For example, again consider the constraint $u(0,0) = 5$ and suppose it is the $i$-th constraint, then,
\begin{equation*}
    \C{i}[u(x_1,x_2)] = \ppC{}{1}{i}\Big[\ppC{}{2}{i}[u(x_1,x_2)]\Big] = \ppC{}{1}{i}\Big[u(x_1,0)\Big] = u(0,0).
\end{equation*}
Utilizing lower-dimensional \ces, these lower-dimensional constraints can be embedded into $n$-dimensional \ces\ by modifying the projection functionals. 

Recall that the projection functional for the $i$-th constraint on the $k$-th dimension is written as,
\begin{equation*}
    \p{k}{\rho}_i(\B{x},g(\B{x})) = \p{k}{\kappa}_i - \pC{k}{i}[g(\B{x})].
\end{equation*}
For lower-dimensional constraints, the projection functionals are written as,
\begin{equation*}
    \p{k}{\rho}_i(\B{x},g(\B{x})) = \pp{(k)}{j}{\phi}_i(x_j)\pp{(k)}{j}{\rho}_i(\B{x},\ppC{(k)}{k}{i}[g(\B{x})])
\end{equation*}
where $\pp{(k)}{j}{\phi}_i(x_j)$ and $\pp{(k)}{j}{\rho}_i(\B{x},g(\B{x}))$ are the switching function and projection functional of a lower-dimensional \ce---there is no implied sum over $i$ or $j$ on the right-hand side of this equation as $i$ and $j$ are used here as identifying symbols rather than as indices. This lower-dimensional constrained expression is built using all of the constraint operators associated with the constraint that do not operate on $x_k$, i.e., $\ppC{(k)}{j}{i}$ such that $j \neq k$. Suppose that the constraint is just one dimension lower than the $n$-dimensional \ce. Then, expanding the projection functional gives a form similar to the original,
\begin{equation*}
    \p{k}{\rho}_i(\B{x},g(\B{x})) = \pp{(k)}{j}{\phi}_i(x_j)\kappa_i - \pp{(k)}{j}{\phi}_i(x_j)\ppC{(k)}{j}{i}\Big[\ppC{(k)}{k}{i}[g(\B{x})]\Big] = \pp{(k)}{j}{\phi}_i(x_j) \Big(\kappa_i - \pC{k}{i}[g] \Big),
\end{equation*}
but $\kappa_i$ and $\pC{k}{i}[g(\B{x})]$ are multiplied by $\pp{(k)}{j}{\phi}_i$---again, no sum is implied over $i$ or $j$. 

Notice that if the lower-dimensional constraint is multiple dimensions lower than the $n$-dimensional constrained expression, then the projection functional of the lower-dimensional constrained expression will itself contain an even lower-dimensional constrained expression. Working out the algebra and simplifying results in the following form for the projection functional,
\begin{equation}\label{eq:lowDimRho}
    \p{k}{\rho}_i(\B{x},g(\B{x})) = \p{k}{\rho}_i(\B{x},g(\B{x})) = \Big(\pp{(k)}{j}{\phi}(x_j)\cdots\pp{(k)}{h}{\phi}(x_h)\Big) \Big(\kappa_i - \pC{k}{i}[g] \Big),
\end{equation}
where $j,\dots,h$ are the dimensions associated with the constraint operators $\ppC{(k)}{j}{i},\dots,\ppC{(k)}{h}{i}$ that make up the constraint excluding the $k$-th dimension. Furthermore, the reader should note that if one has multiple lower-dimensional constraints that share the same operator, $\ppC{(k)}{k}{i}$, then these constraints can be written into the same projection functional, i.e., they can be collected into the same lower-dimensional \ce\ that the projection functional projects $g(\B{x})$ to. This will modify Equation \eqref{eq:lowDimRho} by adding an additional two terms for each constraint: a term for the lower-dimensional switching functions and a term for the lower-dimensional projection function.

In the end, these modifications to the projection functional project $g(\B{x})$ to the set of functions that satisfy the constraint value while maintaining two critical properties of the projection functional: (1) the projection functional is constant with respect to $x_k$, i.e., 
\begin{equation*}
    \ppC{(k)}{k}{i}[\p{k}{\phi}_i(x_k)\p{k}{\rho}_i(\B{x},g(\B{x}))] = \ppC{(k)}{k}{i}[\p{k}{\phi}_i(x_k)]\p{k}{\rho}_i(\B{x},g(\B{x}))
\end{equation*} 
and (2) if $g(\B{x})$ satisfies the constraints, then $\p{k}{\rho}_i(\B{x},g(\B{x})) = 0$, i.e., Property \ref{prop:projZero} still holds. These properties are crucial, as they are necessary and sufficient conditions for the proofs of the \ce\ theorems shown in Chapter \ref{chap:tfcTheory} to hold. Hence, all of the \ce\ theorems shown in Chapter \ref{chap:tfcTheory} still apply to these lower-dimensional constraints. 

To help solidify these concepts, the following two examples are provided.
\begin{example}{Single point constraint in two dimensions}
Consider the point constraint proposed earlier: $u(0,0) = 5$. Suppose the constraint is chosen, arbitrarily, to be embedded into the constraints on $x$. Using the support function $s_1(x) = 1$, the switching function for $x$ is $\p{1}\phi_1(x) = 1$, and the resulting \ce\ is,
\begin{equation*}
    u(x,y,g(x,y)) = \p{1}{u}(x,y,g(x,y)) = g(x,y) + \p{1}{\rho}_1(x,y,g(x,y)).
\end{equation*}
Since the constraint is embedded into $x$, $\p{1}{\rho}_1$ needs to project $g(x,y)$ to the univariate \ce\ on $y$. That is,
\begin{equation*}
    \p{1}{\rho}_1(x,y,g(x,y)) = \pp{(1)}{2}{\phi}_1(y) \pp{(1)}{2}{\rho}_1(x,y,\ppC{(1)}{1}{1}[g(x,y)]).
\end{equation*}
Let the support function for this lower-dimensional \ce\ be $s_1(y) = 1$., then, $\pp{(1)}{2}{\phi}_1(y) = 1$. The lower-dimensional projection functional is,
\begin{equation*}
    \pp{(1)}{2}{\rho}_1(x,y,\ppC{(1)}{1}{1}[g(x,y)]) = 5 - \ppC{(1)}{2}{1}\Big[\ppC{(1)}{1}{1}[g(x,y)]\Big] = 5 - \pC{1}{1}[g(x,y)] = 5 - g(0,0).
\end{equation*}
Putting everything together yields the multivariate \ce,
\begin{equation*}
    u(x,y,g(x,y)) = g(x,y) + 5 - g(0,0).
\end{equation*}
For this example, it is simple to verify that the constrained expression satisfies the constraints for any free function. 
\end{example}
\begin{example}{Lower-dimensional constraints in three dimensions}
Consider the following constraints in three-dimensional space,
\begin{equation*}
    u(x,y,1) = \sin(x)\cos(y), \quad u(0,y,0) = e^y, \quad u(1,0,0) = 3, \andd u(1,1,0) = 5.
\end{equation*}
Notice that the last three lower-dimensional constraints all share the same constraint operator $\ppC{(3)}{3}{i}$, where $i = \{2,3,4\}$. Therefore, they can be embedded into the same lower dimensional constraint. Utilizing multivariate TFC, the \ce\ for $u$ can be written as,
\begin{align*}
    u(x,y,z,g(x,y,z)) &= g(x,y,z) + z \Big(\sin(x)\cos(y) - g(x,y,1)\Big) \\
    &\quad+ (1-z) \p{3}{\rho}_2(x,y,z,g(x,y,z))
\end{align*}
where $s_1(z) = 1$ and $s_2(z) = z$ were chosen for the support functions and $\p{3}{\rho}_2(x,y,z,g(x,y,z))$ is the projection functional for the lower-dimensional constraints. 

The lower-dimensional \ce\ is,
\begin{equation*}
    u(x,y,0,g(x,y,z)) = g(x,y,0) + (1-x)\Big(e^y-g(0,y,0)\Big) + x \pp{(3)}{1}{\rho}_2(x,y,0,g(x,y,z))
\end{equation*}
where $s_1(x) = 1$ and $s_2(x) = x$ were chosen for the support functions and $\pp{(3)}{1}{\rho}_2(x,y,0,g(x,y,0))$ is the projection functional for the two point constraints. The univariate \ce\ for the two point constraints is,
\begin{equation*}
    u(1,y,0,g(x,y,z)) = g(1,y,0) + (1-y)\Big(3-g(1,0,0)\Big) + y\Big(5-g(1,1,0)\Big).
\end{equation*}
Utilizing the switching functions and projection functionals from this univariate \ce, the projection functional for the two-dimensional \ce\ can be derived,
\begin{equation*}
    \pp{(3)}{1}{\rho}_2(x,y,0,g(x,y,0)) = (1-y)\Big(3-g(1,0,0)\Big) + y\Big(5-g(1,1,0)\Big).
\end{equation*}
Likewise, the now complete two-dimensional \ce\ can be used to create the projection functional for the full \ce,
\begin{align*}
    \p{3}{\rho}_2(x,y,z,g(x,y,z)) &= (1-x)\Big(e^y-g(0,y,0)\Big) + x\Big( (1-y)\big(3-g(1,0,0)\big) \\
    &\quad + y\big(5-g(1,1,0)\big)\Big),
\end{align*}
so the full \ce\ is,
\begin{align*}
    u(x,y,z,g(x,y,z)) &= g(x,y,z) + z \Big(\sin(x)\cos(y) - g(x,y,1)\Big) \\
    &\quad+ (1-z) \bigg((1-x)\big(e^y-g(0,y,0)\big) + x \Big((1-y)\big(3-g(1,0,0)\big) \\\
    &\quad + y\big(5-g(1,1,0)\big)\Big)\bigg).
\end{align*}
\end{example}

\section{Conjecture: TFC Extends to any Field}
The author of this dissertation believes that TFC \ces\ as presented here are not restricted to the field of real numbers and extend to any mathematical field. However, the author does not feel confident enough in their abstract algebra knowledge to say this for certain, i.e., to write it as a formal theorem. Hence, rather than presenting this as a theorem in the main body of the text, it is presented here as a conjecture with supporting evidence.

The steps to write a univariate \ce\ can be succinctly summarized for a given set of constraints $\kappa_i = \C{i}[y(x)]$ as:
\begin{enumerate}
    \item $\alpha_{ij} = \Big(\C{i}[s_j]\Big)^{-1}$
    \item $\phi_i(x) = s_j(x) \alpha_{ji}$
    \item $\rho_i(x,g(x)) = \kappa_i - \C{i}[g(x)]$
    \item $y(x,g(x)) = g(x) + \phi_i(x) \rho_i(x,g(x))$
\end{enumerate}
These steps and those used to prove the theorems related to the univariate \ce\ shown in Section \ref{subsec:UniProofs} only use operations defined for a field and matrices consisting of elements of that field. Throughout the body of this dissertation, the field of real numbers was used, i.e., $x\in\mathbb{R}$, but the algebra remains the same for $x\in\mathbb{F}$ where $\mathbb{F}$ is any mathematical field. Of course, the types of constraints that can be embedded for a given field are restricted to the types of constraints that can be calculated on that field. For example, one cannot embed integral constraints into a TFC \ce\ for a finite field because one cannot calculate integrals on finite fields. 
Furthermore, multivariate \ces\ are constructed via recursive application of the univariate theory, and the associated proofs also only utilize mathematical operations defined for a field. Hence, multivariate \ces\ and their associated theorems extend to all mathematical fields as well.

To further provide evidence for this conjecture, the following two examples are provided.
\begin{example}{TFC on a finite field}
Consider the finite field containing the four elements $\{0, 1, A, B\}$ with the addition and multiplication tables shown in Tables \ref{tab:FiniteFieldPlus} and \ref{tab:FiniteFieldMult}.\footnote{If the reader is familiar with finite fields, they will recognize this field as $GF(2)[x]/(x^2+x+1)$.} Consider the following constraints,
\begin{equation*}
    u(0,y) = A, \quad u(B,y) = 1, \andd u(x,0) = u(x,B).
\end{equation*}
Utilizing the theory as described in Chapter \ref{chap:tfcTheory}, the multivariate \ce\ that satisfies these constraints can be derived. In this example, the univariate \ce\ that satisfies the constraints on $x$ is derived step by step; notice that the steps have not been modified from those used to derive \ces\ for real numbers. Let the support functions be $s_1(x) = 1$ and $s_2(x) = x$, then,
\begin{table}[H]
    \centering
    \caption{Addition table.}
    \label{tab:FiniteFieldPlus}
    \begin{tabular}{|c|c|c|c|c|}
        \hline
        $+$ & $0$ & $1$ & $A$ & $B$ \\\hline
        $0$ & $0$ & $1$ & $A$ & $B$ \\\hline
        $1$ & $1$ & $0$ & $B$ & $A$ \\\hline
        $A$ & $A$ & $B$ & $0$ & $1$ \\\hline
        $B$ & $B$ & $A$ & $1$ & $0$ \\\hline
    \end{tabular}
\end{table}
\begin{table}[H]
    \centering
    \caption{Multiplication table.}
    \label{tab:FiniteFieldMult}
    \begin{tabular}{|c|c|c|c|c|}
        \hline
         $*$ & $0$ & $1$ & $A$ & $B$ \\\hline
         $0$ & $0$ & $0$ & $0$ & $0$ \\\hline
         $1$ & $0$ & $1$ & $A$ & $B$ \\\hline
         $A$ & $0$ & $A$ & $B$ & $1$ \\\hline
         $B$ & $0$ & $B$ & $1$ & $A$ \\\hline
    \end{tabular}
\end{table}

\begin{align*}
    \alpha_{ij} &= \Big(\pC{1}{i}[s_j(x)]\Big)^{-1} = \begin{bmatrix} 1 & 0 \\ 1 & B \end{bmatrix}^{-1}\\
    \alpha_{ij} &= \begin{bmatrix} 1 & 0 \\ A & A\end{bmatrix}.
\end{align*}
The switching functions are defined as,
\begin{equation*}
    \p{1}{\phi}_i = s_i(x)\alpha_{ij};
\end{equation*}
thus,
\begin{equation*}
    \p{1}{\phi}_1(x) = Ax + 1 \andd \p{1}{\phi}_2(x) = Ax.
\end{equation*}
The projection functionals are defined as
\begin{equation*}
    \p{1}{\rho}_i(\B{x},g(\B{x}) = \kappa_i - \pC{1}{i}[g(\B{x})];
\end{equation*}
thus,
\begin{equation*}
    \p{1}{\rho}_1(x,y,g(x,y)) = A-g(0,y) \andd \p{1}{\rho}_2(x,y,g(x,y)) = 1-g(B,y).
\end{equation*}
Combining the pieces yields the univariate \ce,
\begin{equation*}
    \p{1}{u}(x,y,g(x,y)) = g(x,y)+(Ax+1)(A-g(0,y))+Ax(1-g(B,y)).
\end{equation*}

The univariate \ce\ for the constraints on $y$ is derived in a similar fashion,
\begin{equation*}
     \p{2}{u}(x,y,g(x,y)) = g(x,y)+Ay\big(g(x,B)-g(x,0)\big).
\end{equation*}
Just as with the real numbers, these two univariate \ces\ are combined recursively to yield the multivariate \ce:
\begin{align*}
    u(x,y,g(x,y)) &= \p{2}{u}(x,y,\p{1}{u}(x,y,g(x,y))) \\
    &= g(x,y) + A x (1-g(B,y))+(A x+1) (A-g(0,y))\\
    &\quad +A y \Big((A x+1) (A-g(0,B))+A x (1-g(B,0))+A x (1-g(B,B))\\
    &\quad -(A x+1) (A-g(0,0))+g(x,B)-g(x,0)\Big).
\end{align*}

Table \ref{tab:FiniteFieldExU} shows the output of this \ce\ for $g(x,y) = Ax + xy + y$. The bottom row of the table gives the $x$ value, and the left-most column gives the $y$ value. The remaining table entries give the output, i.e., $u(x,y,g(x,y))$. As expected, the constraints are satisfied.
\begin{table}[H]
    \centering
    \caption{Finite field constrained expression output for $g(x,y) = Ax + xy + y$.}
    \label{tab:FiniteFieldExU}
    \begin{tabular}{|c||c|c|c|c|}
        \hline
        $B$ & $A$ & $B$ & $0$ & $1$ \\\hline
        $A$ & $A$ & $B$ & $0$ & $1$ \\\hline
        $1$ & $A$ & $B$ & $0$ & $1$ \\\hline
        $0$ & $A$ & $B$ & $0$ & $1$ \\\hline\hline
        \diagbox[height=1.25cm,width=1.25cm, dir=NE]{$y$}{$x$} & $0$ & $1$ & $A$ & $B$ \\\hline
    \end{tabular}
\end{table}
\end{example}

\begin{example}{TFC using complex numbers}
Consider the following constraints,
\begin{equation*}
    y\left(\frac{i}{2}\right) = 1+\pi i, \quad y(1) = y(i), \andd y(2+i) + y_x(1) = 2i.
\end{equation*}
Let the support functions be $s_1(x) =1$, $s_2(x) = x$, and $s_3(x) = x^2$. Then,
\begin{align*}
    \alpha_{ij} &= \Big(\pC{1}{i}[s_j(x)]\Big)^{-1} = \begin{bmatrix}  1 & \frac{i}{2} & -\frac{1}{4} \\ 0 & -1+i & -2 \\ 1 & 3+i & 5+4 i \end{bmatrix}^{-1}\\
    \alpha_{ij} &= \begin{bmatrix} \frac{132}{125}-\frac{24 i}{125} & -\frac{82}{125}+\frac{49 i}{125} & -\frac{7}{125}+\frac{24 i}{125} \\ \frac{52}{125}+\frac{36 i}{125} & -\frac{129}{250}-\frac{397 i}{250} & -\frac{52}{125}-\frac{36 i}{125} \\ -\frac{44}{125}+\frac{8 i}{125} & \frac{69}{125}+\frac{67 i}{125} & \frac{44}{125}-\frac{8 i}{125} \end{bmatrix},
\end{align*}
and
\begin{align*}
   \phi_1(x) &= \left(-\frac{44}{125}+\frac{8 i}{125}\right) x^2+\left(\frac{52}{125}+\frac{36 i}{125}\right) x+\left(\frac{132}{125}-\frac{24 i}{125}\right) \\
   \phi_2(x) &= \left(\frac{69}{125}+\frac{67 i}{125}\right) x^2-\left(\frac{129}{250}+\frac{397 i}{250}\right) x+\left(-\frac{82}{125}+\frac{49 i}{125}\right) \\
   \phi_3(x) &= \left(\frac{44}{125}-\frac{8 i}{125}\right) x^2-\left(\frac{52}{125}+\frac{36 i}{125}\right) x+\left(-\frac{7}{125}+\frac{24 i}{125}\right).
\end{align*}
The projection functionals are,
\begin{align*}
    \rho_1(x) &= 1+\pi i-g\left(\frac{i}{2}\right)\\
    \rho_2(x) &= g(1)-g(i) \\
    \rho_3(x) &= 2i - g(2+i)-g_x(1).
\end{align*}
Using the switching functions and projection functionals, the full \ce\ is,
\begin{align*}
    y&(x,g(x)) = g(x) \\
    & + \Bigg(\left(-\frac{44}{125}+\frac{8 i}{125}\right) x^2+\left(\frac{52}{125}+\frac{36 i}{125}\right) x+\left(\frac{132}{125}-\frac{24 i}{125}\right)\Bigg)\Big(1+\pi i-g\left(\frac{i}{2}\right)\Big) \\
    &+\Big(\left(\frac{69}{125}+\frac{67 i}{125}\right) x^2-\left(\frac{129}{250}+\frac{397 i}{250}\right) x+\left(-\frac{82}{125}+\frac{49 i}{125}\right)\Big)\Big(g(1)-g(i)\Big) \\
    &+ \Big(\left(\frac{44}{125}-\frac{8 i}{125}\right) x^2-\left(\frac{52}{125}+\frac{36 i}{125}\right) x+\left(-\frac{7}{125}+\frac{24 i}{125}\right)\Big)\Big(2i - g(2+i)-g_x(1)\Big).
\end{align*}

Figures \ref{fig:complexExRe} and \ref{fig:complexExIm} show the real and imaginary portions of the constrained expression respectively for $g(x) = \frac{1}{4}x+0.3i\cos(x/4)$. In these figures, $\mathrm{Re}[\cdot]$ is used to denote the real portion, and $\mathrm{Im}[\cdot]$ is used to denote the imaginary portion. The first constraint is plotted as a red point, and the second constraint is plotted as green points. The third constraint is harder to visualize but is satisfied nonetheless.
\begin{figure}[H]
    \begin{minipage}{0.5\linewidth}
    \centering
    \includegraphics[width=\linewidth]{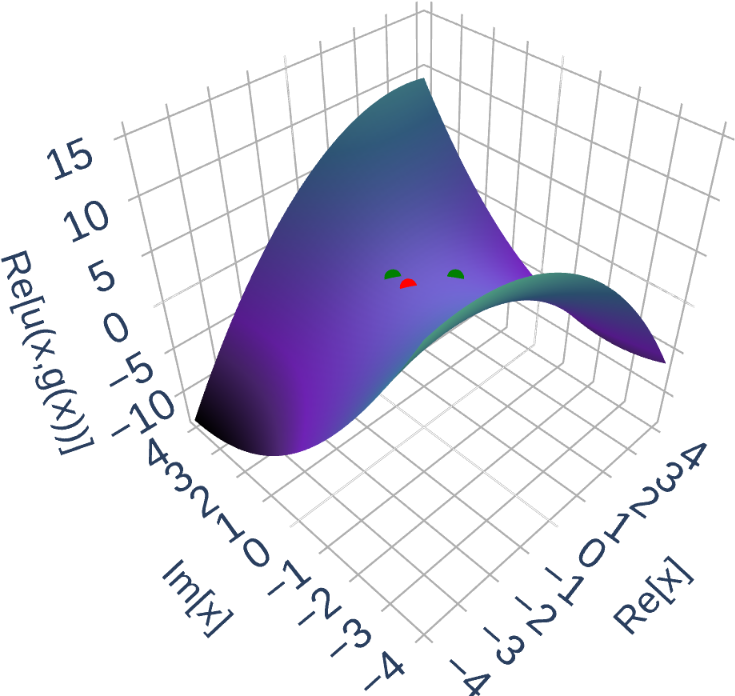}
    \caption{Complex constrained\\expression, real portion.}
    \label{fig:complexExRe}
    \end{minipage}%
    \begin{minipage}{0.5\linewidth}
    \centering
    \includegraphics[width=\linewidth]{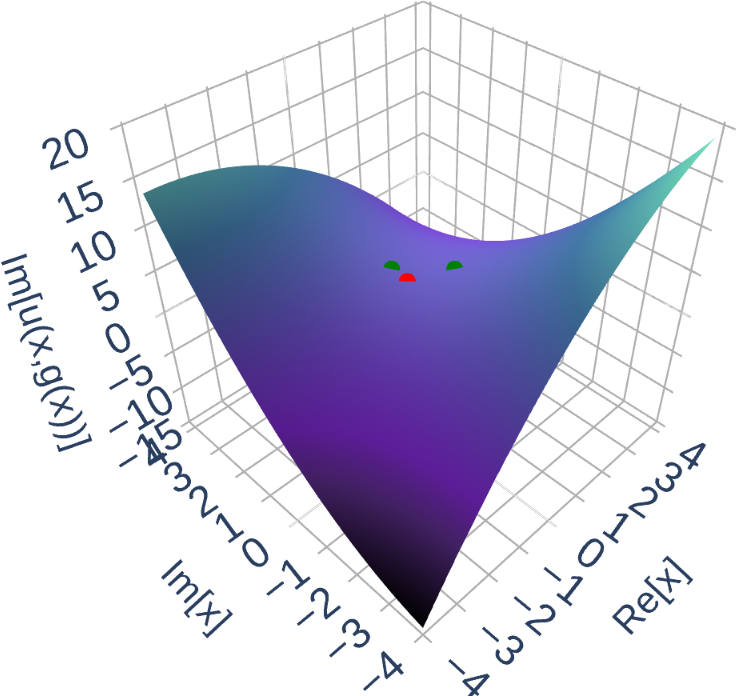}
    \caption{Complex constrained\\expression, imaginary portion.}
    \label{fig:complexExIm}
    \end{minipage}
\end{figure}
\end{example}